\documentclass[12pt]{article}
\usepackage{amsmath}
\usepackage{amssymb}
\usepackage[mathscr]{eucal}
\usepackage{authblk}
\usepackage[usenames]{color}
\usepackage[
bookmarks=true,
backref=true,
colorlinks=true,
linkcolor=red,
urlcolor=green]{hyperref}

\newtheorem{prop}{Proposition}
\newtheorem{thm}{Theorem}

\numberwithin{equation}{section}
\numberwithin{thm}{section}
\numberwithin{lemma}{section}
\numberwithin{prop}{section}
\numberwithin{cor}{section}
\numberwithin{rmk}{section}
\numberwithin{defn}{section}

\newcommand{\gen}[1]{\partial_{#1}}
\newcommand{\lie}{\mathfrak g}

\begin{document}

\title{\Large Symmetry classification of
third-order nonlinear evolution equations.\\ Part II: Solvable algebras.
 }

\author[1]{P.~Basarab-Horwath\thanks{bhrwthp@gmail.com}}
\author[2]{F.~G\"ung\"or\thanks{gungorf@itu.edu.tr}}

\affil[1]{Nygatan 44, 582 27 Link\"oping, Sweden}
\affil[2]{Department of Mathematics, Faculty of Science and Letters, Istanbul Technical University, 34469 Istanbul, Turkey}

\date{}

\maketitle

\abstract{We give a classification of all third-order nonlinear evolution equations
which admit solvable Lie symmetry algebras $\mathsf{A}$ and which are not linearized. We have found that there are 48 types of equations of the form \eqref{eveqn} for $\dim\mathsf{A}=3$, 88 types for $\dim\mathsf{A}=4$ and there are 55 equations for $\dim\mathsf{A}=5$.}

\section{Introduction}
This is a follow-up to our previous paper \cite{BasarabHorwathGuengoerLahno2013} which studied a complete point-symmetry classification of all third-order equations of the form
\begin{equation}\label{eveqn}
u_t= F(t,x,u,u_1, u_2)u_3 + G(t,x,u,u_1, u_2), \quad F\ne 0
\end{equation}
admitting semi-simple symmetry algebras and extensions of these semi-simple Lie algebras by solvable Lie algebras. This paper will be devoted to classification with respect to solvable algebras.  This version  will contain all the computational details. We believe that the researchers working in the symmetry  classification problem of partial differential equations will find them useful.

Unlike semisimple algebras, the classification of solvable Lie algebras is complete only for low dimensions: $\dim\lie\leq 6$. The solution of this  problem was started in the late fifties and sixties in a series of papers \cite{Morozov1958, Mubarakzyanov1963, Mubarakzyanov1963a, Mubarakzyanov1963b, Mubarakzyanov1966} and then extended in \cite{Turkowski1990}. A historical review of the classifications of Lie algebras was given in \cite{BozaFedrianiNunezTenorio2013}.  A summary of the classification and structure results can be found in \cite{BasarabHorwathLahnoZhdanov2001}. For a modern treatment of classification of Lie  algebras the reader is referred to the book \cite{SnobleWinternitz2014}. We mention that a substantial literature exists on classification of certain types of solvable algebras based on the knowledge of existing nilpotent algebras (See references in \cite{SnobleWinternitz2014}).

\section{Some previous results}
We first recall some definitions and  well-known results on solvable algebras \cite{SnobleWinternitz2014}. A Lie algebra $\lie$ is solvable if its derived series $\lie=\lie^{(0)}\supseteq \lie^{(1)}\supseteq\cdots\supseteq\lie^{(k)}\supseteq\cdots$ defined recursively by
$$ \lie^{(k)}=[\lie^{(k-1)},\lie^{(k-1)}]   $$
terminates: $\lie^{(k)}=0$ for some  $k\in \mathbb{N}$.

The lower central series $\lie=\lie^{0}\supseteq \lie^{1}\supseteq\cdots\supseteq\lie^{k}\supseteq\cdots$ is defined recursively by
$$\lie^{k}=[\lie^{k-1},\lie].$$ If the lower central series terminates, namely there exists $k\in \mathbb{N}$ such that $\lie^k=0$, then $\lie$ is called a nilpotent Lie algebra. A nilpotent algebra is also solvable.

Any solvable Lie algebra  $\lie$ contains a unique maximal nilpotent ideal $\mathfrak{n}$ (the nilradical of $\lie$). The dimension of the nilradical satisfies
$$\dim \mathfrak{n}\geq \frac{1}{2}\dim \lie.$$
The derived algebra $\lie'=\lie^{(1)}=\lie^1$ (or ideal) of a solvable algebra is contained in its nilradical $\mathfrak{n}$.

For the basic notions, definitions and relevant discussions of the classification method for the Lie point symmetries of general evolution equations and in particular third order ones we refer to \cite{BasarabHorwathGuengoerLahno2013}. We shall briefly present  some results that will be used throughout this paper.
The classification will be done in terms of vector fields
\begin{equation}\label{gvf}
Q= a(t)\gen t + b(t,x,u)\gen x + c(t,x,u)\gen u
\end{equation}
up to the equivalence transformations
\begin{equation}\label{eqgrp}
t'=T(t),\quad x'=X(t,x,u),\quad u'=U(t,x,u).
\end{equation}
The main strategy used here is to proceed by dimension. Throughout this paper the symbol $\mathsf{A}$ will be used to denote a solvable algebra.

\begin{prop}\label{symmcondition}
The symmetry group of the nonlinear equation \eqref{eveqn} for
arbitrary (fixed) functions $F$ and $G$ is generated by the vector
field \eqref{gvf}
where the functions $a$, $b$ and $c$ satisfy the determining
equations
\begin{subequations}\label{deq}
\begin{equation}\label{deq1}
\begin{array}{ll}
&F\,( -\dot{a}+
     3\,u_1\,b_{u} +
     3\,b_{x} )  +
  [
  u_1(b_{xx}-2c_{xu})+u_1^2(2b_{xu}-c_{uu})+u_1^3b_{uu}+\\
&
  u_2(2b_x-c_u)+
  3u_1u_2b_u-c_{xx}] \,F_{u_2} +\\
&
  [ u_1^2\,b_{u} +
     u_1(b_x-c_u) - c_{x}]\,F_{u_1} -
  c\,F_{u} - a\,F_{t} -
  b\,F_{x}=0,
\end{array}
\end{equation}
and
\begin{equation}\label{deq2}
\begin{array}{ll}
&G\,\left( -\dot{a} -
     u_1\,b_{u} + c_{u}
     \right)  -
  u_1\,b_{t} + c_{t} +\\
&
F\,[u_1(b_{xxx}-3c_{xxu})+3u_1^2(b_{xxu}-c_{xuu})+u_1^3(3b_{xuu}-c_{uuu})
+u_1^4b_{uuu}+\\
&
3u_2(b_{xx}-c_{xu})+3u_1u_2(3b_{xu}-c_{uu})+6u_1^2u_2b_{uu}
+3u_2^2b_u-c_{xxx}]+\\
&
[u_1(b_{xx}-2c_{xu})+u_1^2(2b_{xu}-c_{uu})+u_1^3b_{uu}+u_2(2b_x-c_u)+
3u_1u_2b_u-c_{xx}]\,G_{u_2} +\\
&  [u_1(b_x-c_u)+u_1^2b_u-c_{x}]\,G_{u_1} -
  c\,G_{u} - a\,G_{t} -
  b\,G_{x}=0.
\end{array}
\end{equation}
\end{subequations}
Here the dot over a symbol stands for time derivative.
\end{prop}

\subsection{Admissible Abelian Lie algebras.}

\begin{thm}\label{abeliansymm} (See \cite{BasarabHorwathGuengoerLahno2013} for a proof.) All inequivalent, admissible abelian Lie algebras of vector fields of the form $Q=a(t)\gen t + b(t,x,u)\gen x + c(t,x,u)\gen u$ are given as follows:
\begin{align*}
& \mathsf{A}=\langle \gen t\rangle,\quad \mathsf{A}=\langle \gen u \rangle\;\; (\dim \mathsf{A}=1)\\
& \mathsf{A}=\langle \gen t, \gen u \rangle,\quad \mathsf{A}=\langle \gen x, \gen u \rangle,\quad \mathsf{A}=\langle \gen u, x\gen u\rangle\;\; (\dim \mathsf{A}=2)\\
& \mathsf{A}=\langle \gen t, \gen x, \gen u\rangle,\quad \mathsf{A}=\langle \gen t, \gen u, x\gen u \rangle,\quad \mathsf{A}=\langle \gen u, x\gen u, c(t,x)\gen u \rangle,\; c_{xx}\neq 0\;\; (\dim \mathsf{A}=3)\\
& \mathsf{A}=\langle \gen t, \gen u, x\gen u, c(x)\gen u \rangle,\; c''(x)\neq 0\;\, (\dim \mathsf{A}=4)\\
& \mathsf{A}=\langle \gen u, x\gen u, c(t,x)\gen u, q_1(t,x)\gen u,\dots, q_k(t,x)\gen u \rangle\;\, (\dim \mathsf{A}\geq 4)\\
& c_{xx}\neq 0,\, (q_1)_{xx}\neq 0,\dots, (q_k)_{xx}\neq 0.
\end{align*}
\end{thm}

\noindent Note that
\begin{align*}
\mathsf{A}&=\langle \gen u, x\gen u, c(t,x)\gen u \rangle,\; c_{xx}\neq 0\\
\mathsf{A}&=\langle \gen t, \gen u, x\gen u, c(x)\gen u \rangle,\; c''(x)\neq 0\\
\mathsf{A}&=\langle \gen u, x\gen u, c(t,x)\gen u, q_1(t,x)\gen u,\dots, q_k(t,x)\gen u \rangle\\ &c_{xx}\neq 0,\, (q_1)_{xx}\neq 0,\dots, (q_k)_{xx}\neq 0
\end{align*}
all linearize our third-order evolution equation.

Finally, we mention that we do not list linearizable equations here. When constructing the inequivalent realizations of the solvable Lie algebras we shall make use of the following criteria  for identifying realizations of solvable algebras which lead to linearizable equations.

\begin{thm}[\cite{BasarabHorwathGuengoer2017}]\label{linearizableeqns}
If the evolution equation \eqref{eveqn} admits a rank-one solvable Lie algebra $\mathsf{A}=\langle e_1,e_2,e_3\rangle$ of dimension three or an abelian Lie algebra $\mathsf{A}$ with $\dim\mathsf{A}\geq 4$ as symmetries then the equation is linearizable.
\end{thm}

\section{Two and three-dimensional solvable al\-gebras}
\subsection{Realizations of $\mathsf{A}_{3.2}$}  We have
\[
\mathsf{A}_{3.2}=\mathsf{A}_{2.2}\oplus\mathsf{A}_1.
\]
Here, we note that $\mathsf{A}_{2.2}=\langle e_1, e_2\rangle$ is the solvable Lie algebra with $[e_1, e_2]=e_1$. The following are the canonical admissible realizations of $\mathsf{A}_{2.2}$:

\begin{prop} The admissible two-dimensional solvable Lie algebras $\mathsf{A}_{2.2}$ are given in the following list
\begin{eqnarray*}
\mathsf{A}^1_{2.2} &=& \langle \gen t, t\gen t +x\gen x\rangle\\
\mathsf{A}^2_{2.2} &=& \langle \gen x, t\gen t +x\gen x \rangle\\
\mathsf{A}^3_{2.2} &=& \langle \gen x, x\gen x +u\gen u \rangle\\
\mathsf{A}^4_{2.2} &=& \langle \gen x, x\gen x\rangle.
\end{eqnarray*}
\end{prop}
As can be seen from the list, only the last realization is a rank-one realization, all the others being rank-two realizations.

We shall look at $\mathsf{A}_{3.2}=\langle e_1, e_2, e_3\rangle$ with $[e_1, e_2]=e_1$ and $e_3$ commuting with $e_1$ and $e_2$.

\medskip\noindent $\underline{{\rm rank}\,\mathsf{A}_{3.2}=1}$. Then ${\rm rank}\,\mathsf\mathsf{A}_{2.2}=1$ and so we must have $\displaystyle \mathsf{A}_{2.2}= \langle \gen x, x\gen x\rangle$. Also $e_3=b(t,x,u)\gen x$ because we have ${\rm rank}\,\mathsf{A}_{3.2}=1$. The commutation relations $[e_1, e_3]=[e_2, e_3]=0$ give us $b(t,x,u)=0$ so we have no realization for ${\rm rank}\,\mathsf{A}_{3.2}=1$.

\medskip\noindent $\underline{{\rm rank}\,\mathsf{A}_{3.2}=2}$. Either ${\rm rank}\,\mathsf{A}_{2.2}=1$ or ${\rm rank}\,\mathsf{A}_{2.2}=2$.

\medskip\noindent $\underline{{\rm rank}\,\mathsf{A}_{3.2}=2,\;\;{\rm rank}\,\mathsf{A}_{2.2}=1}$. In this case, $\displaystyle \mathsf{A}_{2.2}= \langle \gen x, x\gen x\rangle$ and $e_3=a(t)\gen t + b(t,x,u)\gen x + c(t,x,u)\gen u$ and $a(t)^2+c(t,x,u)^2\neq 0$ because $e_1\wedge e_3,\;\; e_2\wedge e_3$ are different from zero since ${\rm rank}\,\mathsf{A}_{3.2}=2$. The commutation relations $[e_1, e_3]=[e_2, e_3]=0$ give $b(t,x,u)=0,\; c_x=0$ so that $e_3=a(t)\gen t + c(t,u)\gen u$.

Note that we have the residual equivalence group $\displaystyle \mathscr{E}(\gen x, x\gen x): t'=T(t),\; x'=x,\; u'=U(t,u)$ with $\dot{T}\neq 0,\; U_u\neq 0$. Then under such a transformation we find:
\[
e_3\to e'_3=a(t)\dot{T}(t)\gen {t'} + [a(t)U_t+c(t,u)U_u]\gen {u'}.
\]
If $a(t)\neq 0$ we may always choose $U(t,u)$ so that $a(t)U_t+c(t,u)U_u=0$ and we may choose $T(t)$ so that $a(t)\dot{T}(t)=1$. This gives $e_3=\gen t$ in canonical form.

If $a(t)=0$ then we must have $c(t,u)\neq 0$ and we may choose $U(t,u)$ so that $c(t,u)U_u=1$, giving $e_3=\gen u$ in canonical form.

Hence we have the following realizations of $\mathsf{A}_{3.2}$:
\begin{align*}
\mathsf{A}_{3.2}&=\langle \gen x, x\gen x\rangle\oplus\langle \gen t\rangle\\
\mathsf{A}_{3.2}&=\langle \gen x, x\gen x\rangle\oplus\langle \gen u\rangle.
\end{align*}

\medskip\noindent $\underline{{\rm rank}\,\mathsf{A}_{3.2}=2,\;\;{\rm rank}\,\mathsf{A}_{2.2}=2}$. In this case $e_1\wedge e_2\neq 0$ and $e_1\wedge e_2\wedge e_3=0$. We have three inequivalent realizations of $\mathsf{A}_{2.2}$.

\medskip\noindent$\underline{\mathsf{A}_{2.2}=\langle \gen t, t\gen t + x\gen x\rangle}$. In this case, the residual equivalence group is $\mathscr{E}(e_1, e_2): t'=t,\; x'=p(u)x,\; u'=U(u)$ with $p(u)\neq 0, U'(u)\neq 0$. With $e_3=a(t)\gen t + b\gen x + c\gen u$, we have $c(t,x,u)=0$ because $e_1\wedge e_2\wedge e_3=0$. Then the commutation relation $[e_1, e_3]=0$ gives $\dot{a}(t)=0,\, b_t=0$. Then $[e_2, e_3]=0$ gives $t\dot{a}(t)-a(t)=0,\; xb_x-b=0$ so that $a(t)=0$ and $b=b(u)x$ so that $e_3=b(u)x\gen x$ and $b(u)\neq 0$. If $b'(u)=0$ then we may take $e_3=x\gen x$, but this would then give us $\gen t, t\gen t, x\gen x$ as symmetries and $\gen t, t\gen t$ are incompatible symmetries since they lead to $F=0$, a contradiction. Thus $b'(u)\neq 0$. Under an equivalence transformation, $e_3=b(u)x\gen x$ is transformed to $e'_3=b(u)x'\gen {x'}$ and we choose $U(u)=b(u)$ so that $e'_3=u'x'\gen {x'}$. This gives $e_3=xu\gen x$ in canonical form and

\[
\mathsf{A}_{3.2}=\langle \gen t, t\gen t + x\gen x\rangle\oplus\langle xu\gen x\rangle.
\]

\medskip\noindent$\underline{\mathsf{A}_{2.2}=\langle \gen x, t\gen t + x\gen x\rangle}$. In this case, the residual equivalence group is $\mathscr{E}(e_1, e_2): t'=kt,\; x'=x + p(u)t,\; u'=U(u)$ with $k\neq 0,\; U'(u)\neq 0$. With $e_3=a(t)\gen t + b\gen x + c\gen u$ the condition $e_1\wedge e_2\wedge e_3=0$ gives $c=0$, so $e_3=a(t)\gen t + b\gen x$. Then $[e_1, e_3]=0$ gives $b_x=0$ and $[e_2,  e_3]=0$ gives $tb_t-b=0,\; t\dot{a}(t)-a(t)=0$ so we have $a(t)=mt$ and $b=b(u)t$.

If $m=0$ then $b(u)\neq 0$. If also $b'(u)=0$ then $e_3=t\gen x$ and this is not allowed as we would then have $\gen x,\; t\gen x$ as symmetries, giving a contradiction in the equation for $G$. Thus we must have $b'(u)\neq 0$ in this case. Under an equivalence transformation, $e_3=b(u)t\gen x$ is transformed to $\displaystyle e'_3=\frac{b(u)}{k}t'\gen {x'}$ and we then choose $\displaystyle U(u)=\frac{b(u)}{k}$ so that $e'_3=t'u'\gen {x'}$ giving $e_3=tu\gen x$ in canonical form.

If $m\neq 0$ then we may put $m=1$ (on dividing through by $m$) and then we have  $e_3=t\gen t + b(u)t\gen x$. Under an equivalence transformation $e_3$ is mapped to

\[
e'_3=t'\gen {t'}+ \frac{t'}{k}[b(u)+p(u)]\gen {x'}
\]
and we choose $p(u)$ so that $b(u)+p(u)=0$ giving $e_3=t\gen t$ in canonical form.

So we have the following two realizations of $\mathsf{A}_{3.2}$:
\begin{align*}
\mathsf{A}_{3.2}&=\langle \gen x, t\gen t + x\gen x\rangle\oplus\langle tu\gen x\rangle\\
\mathsf{A}_{3.2}&=\langle \gen x, t\gen t + x\gen x\rangle\oplus\langle t\gen t\rangle.
\end{align*}

\medskip\noindent$\underline{\mathsf{A}_{2.2}=\langle \gen x, x\gen x + u\gen u\rangle}$. In this case, the residual equivalence group is $\mathscr{E}(e_1, e_2): t'=T(t),\; x'=x + p(t)u,\; u'=q(t)u$ with $\dot{T}\neq 0,\; q(t)\neq 0$. The condition $e_1\wedge e_2\wedge e_3=0$ gives $e_3=b\gen x + c\gen u$. The commutation relation $[e_1, e_3]=0$ gives $b_x=c_x=0$, and then $[e_2, e_3]=0$ gives $ub_u-b=0, uc_u-c=0$ so that $b=b(t)u, c=c(t)u$ and then $e_3=b(t)u\gen x + c(t)u\gen u$.

Under an equivalence transformation, $e_3$ is mapped to
\[
e'_3=[b(t)+p(t)c(t)]\frac{u'}{q(t)}\gen {x'} + c(t)u'\gen {u'}.
\]
If $c(t)=0$ then $b(t)\neq 0$ and we may always choose $q(t)=b(t)$ giving $e'_3=u'\gen {x'}$, or $e_3=u\gen x$ in canonical form.

If $c(t)\neq 0$, we choose $p(t)$ so that $c(t)p(t)+b(t)=0$ and we then have $e'_3=c(t)u'\gen {u'}$. if $\dot{c}=0$ then we obtain $e_3=u\gen u$ in canonical form; if $\dot{c}\neq 0$ we choose $T(t)=c(t)$ and we find $e_3=tu\gen u$ in canonical form.

We then have the following realizations of $\mathsf{A}_{3.2}$:
\begin{align*}
\mathsf{A}_{3.2}&=\langle \gen x, x\gen x + u\gen u\rangle\oplus\langle u\gen x\rangle\\
\mathsf{A}_{3.2}&=\langle \gen x, x\gen x + u\gen u\rangle\oplus\langle u\gen u\rangle\\
\mathsf{A}_{3.2}&=\langle \gen x, x\gen x + u\gen u\rangle\oplus\langle tu\gen u\rangle.
\end{align*}

\bigskip\noindent $\underline{{\rm rank}\,\mathsf{A}_{3.2}=3}$. In this case we have $e_1\wedge e_2\wedge e_3\neq 0$. Obviously we must have ${\rm rank}\, \mathsf{A}_{2.2}=2$.

\medskip\noindent$\underline{\mathsf{A}_{2.2}=\langle \gen t, t\gen t + x\gen x\rangle}$. The residual equivalence group is $\mathscr{E}(e_1, e_2): t'=t,\; x'=p(u)x,\; u'=U(u)$ with $p(u)\neq 0, U'(u)\neq 0$. With $e_3=a(t)\gen t + b\gen x + c\gen u$, we have $c(t,x,u)\neq 0$ because $e_1\wedge e_2\wedge e_3\neq 0$. Then the commutation relation $[e_1, e_3]=0$ gives $\dot{a}(t)=0,\, b_t=c_t=0$. Also $[e_2, e_3]=0$ gives $t\dot{a}(t)-a(t)=0,\; xb_x-b=0,\; c_x=0$ so that $a(t)=0$ and $b=b(u)x,\; c=c(u)$ so that $e_3=b(u)x\gen x + c(u)\gen u$ and $c(u)\neq 0$. Under an equivalence transformation, $e_3=b(u)x\gen x + c(u)\gen u$ is transformed to
\[
e'_3=\left[b(u)+\frac{c(u)p'(u)}{p(u)}\right]x'\gen {x'} + c(u)U'(u)\gen {u'}.
\]
Since $c(u)\neq 0$ we may always choose $p(u)$ so that $\displaystyle b(u)+\frac{c(u)p'(u)}{p(u)}=0$ and we choose $U(u)$ so that $c(u)U'(u)=1$. This gives $e_3=\gen u$ in canonical form and we obtain the realization
\[
\mathsf{A}_{3.2}=\langle \gen t, t\gen t + x\gen x\rangle\oplus\langle \gen u\rangle.
\]

\medskip\noindent$\underline{\mathsf{A}_{2.2}=\langle \gen x, t\gen t + x\gen x\rangle}$. The residual equivalence group is $\mathscr{E}(e_1, e_2): t'=kt,\; x'=x + p(u)t,\; u'=U(u)$ with $k\neq 0,\; U'(u)\neq 0$. With $e_3=a(t)\gen t + b\gen x + c\gen u$ the condition $e_1\wedge e_2\wedge e_3\neq 0$ gives $c\neq 0$. The commutation relations $[e_1, e_3]=[e_2, e_3]=0$ give $a(t)=mt,\; b=b(u)t,\; c=c(u)$ and so $e_3=mt\gen t + b(u)t\gen x + c(u)\gen u$. Under an equivalence transformation, $e_3$ is mapped to
\[
e'_3=mt'\gen {t'} + [c(u)p'(u)+b(u)+mp(u)]t\gen {x'} + c(u)U'(u)\gen {u'}.
\]
Since $c(u)\neq 0$ we may always choose $p(u)$ so that $c(u)p'(u)+b(u)+mp(u)=0$, so we may always have $e'_3=mt'\gen {t'} + c(u)U'(u)\gen {u'}$.

If $m=0$ then we choose $U(u)$ so that $c(u)U'(u)=1$ so that $e_3=\gen {u}$ in canonical form.

If $m\neq 0$ then we may put $m=1$ (on dividing through by $m$) and we choose $U(u)$ so that $c(u)U'(u)=U(u)$ then we have  $e_3=t\gen t + u\gen u$ in canonical form.

So we have the following two realizations of $\mathsf{A}_{3.2}$:
\begin{align*}
\mathsf{A}_{3.2}&=\langle \gen x, t\gen t + x\gen x\rangle\oplus\langle \gen u\rangle\\
\mathsf{A}_{3.2}&=\langle \gen x, t\gen t + x\gen x\rangle\oplus\langle t\gen t+u\gen u\rangle.
\end{align*}

\medskip\noindent$\underline{\mathsf{A}_{2.2}=\langle \gen x, x\gen x + u\gen u\rangle}$. The residual equivalence group is $\mathscr{E}(e_1, e_2): t'=T(t),\; x'=x + p(t)u,\; u'=q(t)u$ with $\dot{T}\neq 0,\; q(t)\neq 0$. The condition $e_1\wedge e_2\wedge e_3\neq 0$ gives $a(t)\neq 0$ and so $e_3=a(t)\gen t + b\gen x + c\gen u$. The commutation relation $[e_1, e_3]=0$ gives $b_x=c_x=0$, and then $[e_2, e_3]=0$ gives $ub_u-b=0$, $uc_u-c=0$ so that $b=b(t)u$, $c=c(t)u$ and then $e_3=a(t)\gen t + b(t)u\gen x + c(t)u\gen u$.

Under an equivalence transformation, $e_3$ is mapped to
\[
e'_3=a(t)\dot{T}(t)\gen {t'} + [a(t)\dot{p}(t)+b(t)+p(t)c(t)]\frac{u'}{q(t)}\gen {x'} + \left[\frac{a(t)\dot{q}(t)}{q(t)} + c(t)\right]u'\gen {u'}.
\]
Since $a(t)\neq 0$ we may always choose $p(t)$ and $q(t)$ so that $a(t)\dot{p}(t)+b(t)+p(t)c(t)=0$ and $\displaystyle \frac{a(t)\dot{q}(t)}{q(t)} + c(t)=0$. Also we may choose $T(t)$ so that $a(t)\dot{T}(t)=1$ and this then gives us $e_3=\gen t$ in canonical form.

We then have the following realization of $\mathsf{A}_{3.2}$:
\[
\mathsf{A}_{3.2}=\langle \gen x, x\gen x + u\gen u\rangle\oplus\langle \gen t\rangle.
\]

\subsection{The inequivalent non-linearizing realizations of $\mathsf{A}_{3.2}$:}

\begin{align*}
\mathsf{A}_{3.2}&=\langle \gen x, x\gen x, \gen t\rangle\\
\mathsf{A}_{3.2}&=\langle \gen x, x\gen x, \gen u\rangle\\
\mathsf{A}_{3.2}&=\langle \gen t, t\gen t + x\gen x, xu\gen x\rangle\\
\mathsf{A}_{3.2}&=\langle \gen x, t\gen t + x\gen x, tu\gen x\rangle\\
\mathsf{A}_{3.2}&=\langle \gen x, t\gen t + x\gen x, t\gen t\rangle\\
\mathsf{A}_{3.2}&=\langle \gen x, x\gen x + u\gen u, u\gen x\rangle\\
\mathsf{A}_{3.2}&=\langle \gen x, x\gen x + u\gen u, u\gen u\rangle\\
\mathsf{A}_{3.2}&=\langle \gen x, x\gen x + u\gen u, tu\gen u\rangle\\
\mathsf{A}_{3.2}&=\langle \gen t, t\gen t + x\gen x, \gen u\rangle\\
\mathsf{A}_{3.2}&=\langle \gen x, t\gen t + x\gen x, \gen u\rangle\\
\mathsf{A}_{3.2}&=\langle \gen x, t\gen t + x\gen x, t\gen t+u\gen u\rangle\\
\mathsf{A}_{3.2}&=\langle \gen x, x\gen x + u\gen u, \gen t\rangle.
\end{align*}

\section{Realizations of $\mathsf{A}_{3.3}$:} In this case we have $\mathsf{A}_{3.3}=\langle e_1, e_2, e_3\rangle$ with the commutation relations $[e_1, e_2]=0,\; [e_1, e_3]=0,\; [e_2, e_3]=e_1$. Thus, $\langle e_1, e_2\rangle=A_{2.1}$ is the (two-dimensional) abelian ideal of $\mathsf{A}_{3.3}$.

\medskip\noindent $\underline{{\rm rank}\,\mathsf{A}_{3.3}=1}:$ In this case we must have ${\rm rank}\,\langle e_1, e_2\rangle=1$, so that we have $\langle e_1, e_2\rangle=\langle \gen u, x\gen u\rangle$ in canonical form. Also we must have $e_3=c(t,x,u)\gen u$. First, $[e_1, e_3]=0$ gives $e_3=c(t,x)\gen u$ and then $[e_2, e_3]=e_1$ gives $xc_u=1$ and this is a contradiction. Hence we cannot have ${\rm rank}\,\mathsf{A}_{3.3}=1$.

\medskip\noindent$\underline{{\rm rank}\,\mathsf{A}_{3.3}=2,\;\; {\rm rank}\,\langle e_1, e_2\rangle=1}:$ In this case, we also have $\langle e_1, e_2\rangle=\langle \gen u, x\gen u\rangle$ in canonical form, and $e_3=a(t)\gen t + b\gen x + c\gen u$ with $a(t)^2+b^2\neq 0$. Again, $[e_1, e_3]=0$ gives $b_u=c_u=0$ and then $[e_2, e_3]=e_1$ gives $b=-1$ so that $e_3=a(t)\gen t -\gen x + c(t,x)\gen u$.

The residual equivalence group $\mathscr{E}(e_1, e_2): t'=T(t),\; x'=x,\, u'=u+U(t,x)$. Under such a transformation, $e_3$ is mapped to
\[
e'_3=a(t)\dot{T}(t)\gen {t'} - \gen {x'} + [a(t)U_t(t,x) - U_x(t,x) + c(t,x)]\gen {u'}.
\]
We can always choose $U(t,x)$ so that $a(t)U_t(t,x) - U_x(t,x) + c(t,x)=0$ so that we have $e'_3=a(t)\dot{T}(t)\gen {t'} - \gen {x'}$. If $a(t)=0$ then $e'_3=-\gen {x'}$, giving $e_3=-\gen x$ in canonical form. if $a(t)\neq 0$ then we choose $T(t)$ so that $a(t)\dot{T}(t)=-1$ and then $e_3=-(\gen t + \gen x)$ in canonical form. Thus we have the realizations:
\begin{align*}
\mathsf{A}_{3.3}&=\langle \gen u, x\gen u, -\gen x\rangle\\
\mathsf{A}_{3.3}&=\langle \gen u, x\gen u, -\gen t -\gen x\rangle.
\end{align*}

\medskip\noindent$\underline{{\rm rank}\,\mathsf{A}_{3.3}=2,\;\; {\rm rank}\,\langle e_1, e_2\rangle=2}:$ Here we have three possible canonical forms for $\langle e_1, e_2\rangle$.

\medskip\noindent $\underline{\langle e_1, e_2\rangle=\langle \gen t, \gen x\rangle}:$ Putting $e_3=a(t)\gen t + b\gen x + c\gen u$, the commutation relation $[e_2, e_3]=e_1$ gives $\gen t=0$ which is a contradiction. Thus we have no realization in this case.

\medskip\noindent $\underline{\langle e_1, e_2\rangle=\langle \gen x, \gen t\rangle}:$ Putting $e_3=a(t)\gen t + b\gen x + c\gen u$, then ${\rm rank}\, A_{3.3}=2$ gives $e_1\wedge e_2\wedge e_3=0$ which gives $c=0$ so that $e_3=a(t)\gen t + b\gen x $. The commutation relation $[e_2, e_3]=e_1$ gives $\dot{a}(t)=0,\; b_t=1$ and $[e_1, e_3]=0$ gives $b_x$. So we may assume that $a=0$ without loss of generality and then $e_3=[t+b(u)]\gen x$.

The residual equivalence group $\mathscr{E}(e_1, e_2): t'=t+l,\; x'=x + Y(u),\, u'=U(u)$. Under such a transformation, $e_3$ is mapped to
\[
e'_3=[t'+b(u)-l]\gen {x'}.
\]
If $b'(u)=0$ then we may choose $l$ so that $e'_3=t'\gen {x'}$ so that $e_3=t\gen x$ in canonical form. This is not allowed, since $\gen x$ and $t\gen x$ are incompatible as symmetries: they give a contradiction in the equation for $G$. Thus $b'(u)\neq 0$. In this case, we take $u'=U(u)=b(u)-l$ and we then have $e_3=(t+u)\gen x$ in canonical form. Thus we have the realization
\[
\mathsf{A}_{3.3}=\langle \gen x, \gen t, (t+u)\gen x\rangle.
\]

\medskip\noindent $\underline{\langle e_1, e_2\rangle=\langle \gen x, \gen u\rangle}:$ Putting $e_3=a(t)\gen t + b\gen x + c\gen u$, then ${\rm rank}\, A_{3.3}=2$ gives $e_1\wedge e_2\wedge e_3=0$ which gives $a(t)=0$ and so we have $e_3=b\gen x + c\gen u$. The relation $[e_1, e_3]=0$ gives $b_x=c_x=0$ and then $[e_2, e_3]=e_1$ gives $c_u=0,\; b_u=1$. Hence $c=c(t),\; b=u+b(t)$ and we have $e_3=[u+b(t)]\gen x + c(t)\gen u$.

The residual equivalence group is $\mathscr{E}(e_1, e_2): t'=T(t),\; x'=x+Y(t),\; u'=u+U(t)$. Under such a transformation, $e_3$ is mapped to
\[
e'_3=[u'+b(t)-U(t)]\gen {x'} + c(t)\gen {u'}.
\]
we may always choose $U(t)$ so that $b(t)=U(t)$, giving $e'_3=u'\gen {x'} + c(t)\gen {u'}$. If $c(t)=0$ or if $\dot{c}(t)=0$ then we may take $e_3=u\gen x$ in canonical form. If $\dot{c}(t)\neq 0$ then we may take $T(t)=c(t)$ so that $e'_3=u'\gen {x'} + t'\gen {u'}$ giving $e_3=u\gen x + t\gen u$ in canonical form. Thus, we have the realizations
\begin{align*}
\mathsf{A}_{3.3}&=\langle \gen x, \gen u, u\gen x\rangle\\
\mathsf{A}_{3.3}&=\langle \gen x, \gen u, u\gen x + t\gen u\rangle.
\end{align*}
We note that the realization
\[
\mathsf{A}_{3.3}=\langle \gen x, \gen u, u\gen x\rangle
\]
will give the same evolution equation as
\[
\mathsf{A}_{3.3}=\langle \gen u, x\gen u, -\gen x\rangle.
\]
In fact, the equivalence transformation $t'=t,\; x'=u,\; u'=x$ transforms
\[
\mathsf{A}_{3.3}=\langle \gen u, x\gen u, -\gen x\rangle.
\]
to
\[
\mathsf{A}_{3.3}=\langle \gen x, u\gen x, -\gen u\rangle.
\]
and then the change of basis $e'_1=e_1,\; e'_2=-e_3,\; e'_3=e_2$ gives us
\[
\mathsf{A}_{3.3}=\langle \gen x, \gen u, u\gen x\rangle.
\]
However, they are inequivalent in the sense that there is no equivalence transformation which maps one into the other. Also, as will be seen later, the ordering of the basis elements are of crucial importance when considering extensions of $\mathsf{A}_{3.3}$: in fact the algebras $\mathsf{A}_{4.7}$ and $\mathsf{A}_{4.8}$ are one-dimensional solvable extensions of $\mathsf{A}_{3.3}$.

\medskip\noindent$\underline{{\rm rank}\,\mathsf{A}_{3.3}=3}:$ In this case we must also have ${\rm rank}\,\langle e_1, e_2\rangle=2$.

\medskip\noindent $\underline{\langle e_1, e_2\rangle=\langle \gen t, \gen x\rangle}:$ This case is not admissible because the commutation relation $[e_2, e_3]=e_1$ gives the contradiction $\gen t=0$.

\medskip\noindent $\underline{\langle e_1, e_2\rangle=\langle \gen x, \gen t\rangle}:$ Putting $e_3=a(t)\gen t + b\gen x + c\gen u$ we note that $e_1\wedge e_2\wedge e_3\neq 0$ since ${\rm rank}\,\mathsf{A}_{3.3}=3$, and this gives $c\neq 0$. The commutation relation $[e_2, e_3]=e_1$ gives $\dot{a}(t)=0,\; b_t=1,\; c_t=0$ and $[e_1, e_3]=0$ gives $b_x=c_x=0$. So we may assume that $a=0$ without loss of generality and then $e_3=[t+b(u)]\gen x + c(u)\gen u$ with $c(u)\neq 0$.

The residual equivalence group $\mathscr{E}(e_1, e_2): t'=t+l,\; x'=x + Y(u),\, u'=U(u)$. Under such a transformation, $e_3$ is mapped to
\[
e'_3=[t'+b(u)-l + c(u)Y'(u)]\gen {x'} + c(u)U'(u)\gen {u'}.
\]
Since $c(u)\neq 0$ we may always choose $Y(u)$ so that $b(u)-l + c(u)Y'(u)=0$ and we may choose $U(u)$ so that $c(u)U'(u)=-1$, giving $e'_3=t'\gen {x'}-\gen {u'}$. So $e_3=t\gen x - \gen u$ in canonical form and thus we have the realization
\[
\mathsf{A}_{3.3}=\langle \gen x, \gen t, t\gen x - \gen u\rangle.
\]

\medskip\noindent $\underline{\langle e_1, e_2\rangle=\langle \gen x, \gen u\rangle}:$ Putting $e_3=a(t)\gen t + b\gen x + c\gen u$, then ${\rm rank}\, A_{3.3}=3$ gives $e_1\wedge e_2\wedge e_3\neq 0$ which gives $a(t)\neq 0$ and so we have $e_3=a(t)\gen t + b\gen x + c\gen u$. The relation $[e_1, e_3]=0$ gives $b_x=c_x=0$ and then $[e_2, e_3]=e_1$ gives $c_u=0,\; b_u=1$. Hence $c=c(t),\; b=u+b(t)$ and we have $e_3=a(t)\gen t + [u+b(t)]\gen x + c(t)\gen u$.

The residual equivalence group is $\mathscr{E}(e_1, e_2): t'=T(t),\; x'=x+Y(t),\; u'=u+U(t)$. Under such a transformation, $e_3$ is mapped to
\[
e'_3=a(t)\dot{T}(t)\gen {t'} + [u'+b(t)-U(t) + a(t)\dot{Y}(t)]\gen {x'} + [c(t) + a(t)\dot{U}(t)]\gen {u'}.
\]
Because $a(t)\neq 0$ we choose $U(t)$ so that $c(t)+ a(t)\dot{U}(t)]=0$ and then we may choose $Y(t)$ so that $b(t)-U(t) + a(t)\dot{Y}(t)=0$, and finally we may choose $T(t)$ so that $a(t)\dot{T}(t)=T(t)$, so that $e'_3=t'\gen {t'} + u'\gen {x'}$. Hence $e_3=t\gen t + u\gen x$ in canonical form and we have the realization
\[
\mathsf{A}_{3.3}=\langle \gen x, \gen u, t\gen t + u\gen x\rangle.
\]

\subsection{The inequivalent non-linearizing realizations of $\mathsf{A}_{3.3}$:}

\begin{align*}
\mathsf{A}_{3.3}&=\langle \gen u, x\gen u, -\gen x\rangle\\
\mathsf{A}_{3.3}&=\langle \gen u, x\gen u, -\gen t -\gen x\rangle\\
\mathsf{A}_{3.3}&=\langle \gen x, \gen t, (t+u)\gen x\rangle\\
\mathsf{A}_{3.3}&=\langle \gen x, \gen u, u\gen x\rangle\\
\mathsf{A}_{3.3}&=\langle \gen x, \gen u, u\gen x + t\gen u\rangle\\
\mathsf{A}_{3.3}&=\langle \gen x, \gen t, t\gen x - \gen u\rangle\\
\mathsf{A}_{3.3}&=\langle \gen x, \gen u, t\gen t + u\gen x\rangle.
\end{align*}

\section{Realizations of $\mathsf{A}_{3.4}$:} In this case we have $\mathsf{A}_{3.3}=\langle e_1, e_2, e_3\rangle$ with the commutation relations $[e_1, e_2]=0,\; [e_1, e_3]=e_1,\; [e_2, e_3]=e_1+e_2$. Thus, $\langle e_1, e_2\rangle=A_{2.1}$ is the (two-dimensional) abelian ideal of $\mathsf{A}_{3.4}$.

\medskip\noindent $\underline{{\rm rank}\,\mathsf{A}_{3.4}=1}:$ In this case ${\rm rank}\,\langle e_1, e_2\rangle=1$ and so we have $\langle e_1, e_2\rangle=\langle \gen x, u\gen x\rangle$. Because ${\rm rank}\,\mathsf{A}_{3.4}=1$ we must have $e_3=b(t,x,u)\gen x$. The commutation relation $[e_1, e_3]=e_1$ gives $b_x=1$ and $[e_2, e_3]=e_1+e_2$ gives $ub_x=1+u$. This gives a contradiction: $1=0$, so we have no realization in this case.

\medskip\noindent $\underline{{\rm rank}\,\mathsf{A}_{3.4}=2,\; {\rm rank}\,\langle e_1, e_2\rangle=1}:$ Again, $e_1=\gen x, e_2=u\gen x$ and we take $e_3=a(t)\gen t + b\gen x + c\gen u$ with $a(t)^2+c^2\neq 0$. The commutation relation $[e_1, e_3]=e_1$ gives $b_x=1,\; c_x=0$, and $[e_2, e_3]=e_1+e_2$ gives $ub_x-c=1+u$ so we have $c=-1,\; b=x+b(t,u)$. Hence $e_3=a(t)\gen t + [x+b(t,u)]\gen x - \gen u$.

The residual equivalence group is $\mathscr{E}(e_1, e_2): t'=T(t),\; x'=x+Y(t,u),\; u'=u$. Under such a transformation, $e_3$ is mapped to
\[
e'_3=a(t)\dot{T}(t)\gen {t'} + [x'+b(t,u)-Y(t,u)-Y_u(t,u)+a(t)Y_t(t,u)]\gen {x'} - \gen {u'}.
\]
We may always choose $Y(t,u)$ so that $b(t,u)-Y(t,u)-Y_u(t,u)+a(t)Y_t(t,u)=0$ and so we find that $e'_3=a(t)\dot{T}(t)\gen {t'} + x'\gen {x'} - \gen {u'}$. If $a(t)=0$ then we obtain the canonical form $e_3=x\gen x - \gen u$. If $a(t)\neq 0$ then we choose $T(t)$ so that $a(t)\dot{T}(t)=T(t)$ and we obtain $e_3=t\gen t + x\gen x -\gen u$ in canonical form. We have the realizations
\begin{align*}
\mathsf{A}_{3.4}&=\langle\gen x, u\gen x, x\gen x - \gen u\rangle\\
\mathsf{A}_{3.4}&=\langle\gen x, u\gen x, t\gen t + x\gen x - \gen u\rangle.
\end{align*}

\medskip\noindent $\underline{{\rm rank}\,\mathsf{A}_{3.4}=2,\; {\rm rank}\,\langle e_1, e_2\rangle=2}.$

\medskip\noindent $\underline{\langle e_1, e_2\rangle=\langle \gen t, \gen x\rangle}:$ We have $e_3=a(t)\gen t + b\gen x$ because ${\rm rank}\,\mathsf{A}_{3.4}=2$ so that $e_1\wedge e_2\wedge e_3=0$. Then $[e_2, e_3]=e_1+e_2$ leads to the contradiction $\gen t + \gen x=b_x\gen x$, so we have no realization in this case.

\medskip\noindent $\underline{\langle e_1, e_2\rangle=\langle \gen x, \gen t\rangle}:$ We have $e_3=a(t)\gen t + b\gen x$ because ${\rm rank}\,\mathsf{A}_{3.4}=2$ so that $e_1\wedge e_2\wedge e_3=0$. The commutation relation $[e_1, e_3]=e_1$ gives $b_x=1$ so that $b=x+b(t,u)$. Further, $[e_2, e_3]=e_1+e_2$ gives $\dot{a}(t)\gen t + b_t(t,u)\gen x=\gen t + \gen x$ so that $a(t)=t+k,\; b(t,u)=t+b(u)$ and then $e_3=[t+k]\gen t + [t+x+b(u)]\gen x$.

The residual equivalence group is $\mathscr{E}(e_1, e_2): t'=t+l,\; x'=x+Y(u),\; u'=U(u)$ with $\dot{T}(t)\neq 0,\; U'(u)\neq 0$. Under such a transformation, $e_3$ is mapped to
\[
e'_3=[t'+k-l]\gen {t'} + [x'+t'+b(u)-l-Y(u)]\gen {x'}.
\]
It is clear that we may choose $l=k$ and that we may always choose $Y(u)$ so that $b(u)-l-Y(u)=0$. Thus we obtain $e_3=t\gen t + [t+x]\gen x$ in canonical form and we have the realization
\[
\mathsf{A}_{3.4}=\langle \gen x, \gen t, t\gen t + [t+x]\gen x\rangle.
\]

\medskip\noindent $\underline{\langle e_1, e_2\rangle=\langle \gen x, \gen u\rangle}:$ We have $e_3=b\gen x + c\gen u$ because ${\rm rank}\,\mathsf{A}_{3.4}=2$ so that $e_1\wedge e_2\wedge e_3=0$. The relation $[e_1, e_3]=e_1$ gives $b_x=1,\; c_x=0$. Further, $[e_2, e_3]=e_1+e_2$ gives $b_u=1,\; c_u=1$ and so we find that $e_3=[x+u+b(t)]\gen x + [u+c(t)]\gen u$.

The residual equivalence group is $\mathscr{E}(e_1, e_2): t'=T(t),\; x'=x+Y(t),\; u'=u+U(t)$. Under such a transformation, $e_3$ is mapped to
\[
e'_3=[x'+u'+b(t)-Y(t)-U(t)]\gen {x'} + [u'+c(t)-U(t)]\gen {u'}.
\]
It is clear that we may always choose $Y(t),\; U(t)$ so that $U(t)=c(t),\; Y(t)+U(t)=b(t)$ and so we obtain $e_3=[x+u]\gen x + u\gen u$ in canonical form. We obtain the realization
\[
\mathsf{A}_{3.4}=\langle \gen x, \gen u, [x+u]\gen x, u\gen u\rangle.
\]

\medskip\noindent $\underline{{\rm rank}\,\mathsf{A}_{3.4}=3}.$ In this case, ${\rm rank}\,\langle e_1, e_2\rangle=2$.

\medskip\noindent $\underline{\langle e_1, e_2\rangle=\langle \gen t, \gen x\rangle}.$ With $e_3=a(t)\gen t + b\gen x + c\gen u$ we again have that $[e_2, e_3]=e_1+e_2$ gives us the contradiction $\gen t=0$. So we have no realization in this case.

\medskip\noindent $\underline{\langle e_1, e_2\rangle=\langle \gen x, \gen t\rangle}.$ With $e_3=a(t)\gen t + b\gen x + c\gen u$ we must have $c(t,x,u)\neq 0$ since ${\rm rank}\,\mathsf{A}_{3.4}=3$ gives $e_1\wedge e_2\wedge e_3\neq 0$ and this means that $c\neq 0$. The commutation relations $[e_1, e_3]=e_1$ and $[e_2, e_3]=e_1+e_2$ give $b_x=1,\; c_x=0,\; b_t=1,\; \dot{a}(t)=1,\; c_t=0$ and so we have $e_3=[t+k]\gen t + [t+x+b(u)]\gen x + c(u)\gen u$.

The residual equivalence group is $\mathscr{E}(e_1, e_2): t'=t+l,\; x'=x+Y(u),\; u'=U(u)$ with $U'(u)\neq 0$. Under such a transformation, $e_3$ is mapped to
\[
e'_3=[t'+k-l]\gen {t'} + [x'+t'+b(u) + c(u)Y'(u)-l-Y(u)]\gen {x'} + c(u)U'(u)\gen {u'}.
\]
We may always choose $l=k$ and for any given $c(u)\neq 0$ we may choose $Y(u)$ so that $b(u) + c(u)Y'(u)-l-Y(u)=0$ and $U(u)$ so that $c(u)U'(u)=U(u)$. Thus we obtain $e'_3=t'\gen {t'} + [t'+x']\gen {x'} + u'\gen {u'}$. So, $e_3=t\gen t + [t+x]\gen x + u\gen u$ in canonical form. We obtain the realization
\[
\mathsf{A}_{3.4}=\langle \gen x, \gen t, t\gen t + [t+x]\gen x + u\gen u\rangle.
\]

\medskip\noindent $\underline{\langle e_1, e_2\rangle=\langle \gen x, \gen u\rangle}.$ With $e_3=a(t)\gen t + b\gen x + c\gen u$ we must have $a(t)\neq 0$ since ${\rm rank}\,\mathsf{A}_{3.4}=3$ gives $e_1\wedge e_2\wedge e_3\neq 0$ and this means that $a\neq 0$. The commutation relations $[e_1, e_3]=e_1$ and $[e_2, e_3]=e_1+e_2$ give $b_x=1,\; c_x=0,\; b_u=1,\; c_u=1$ and so we have $e_3=a(t)\gen t + [x + u + b(t)]\gen x + [u+c(t)]\gen u$.

The residual equivalence group is $\mathscr{E}(e_1, e_2): t'=T(t),\; x'=x+Y(t),\; u'=u+U(t)$ with $\dot{T}(t)\neq 0,\; U'(u)\neq 0$. Under such a transformation, $e_3$ is mapped to
\[
e'_3=a(t)\dot{T}(t)\gen {t'} + [x'+u'+b(t)+a(t)\dot{Y}(t) - Y(t) - U(t)]\gen {x'} + [u'+c(t)-U(t)]\gen {u'}.
\]
We may always choose $Y(t),\; U(t)$ so that $c(t)-U(t)=0,\; b(t)+a(t)\dot{Y}(t) - Y(t) - U(t)=0$ and we choose $T(t)$ so that $a(t)\dot{T}(t)=T(t)$ and we then obtain $e_3=t\gen t + [x+u]\gen x + u\gen u$ in canonical form. We have the realization
\[
\mathsf{A}_{3.4}=\langle \gen x, \gen u, t\gen t + [x+u]\gen x + u\gen u \rangle.
\]

\subsection{The inequivalent non-linearizing realizations of $\mathsf{A}_{3.4}$:}
\begin{align*}
\mathsf{A}_{3.4}&=\langle\gen x, u\gen x, x\gen x - \gen u\rangle\\
\mathsf{A}_{3.4}&=\langle\gen x, u\gen x, t\gen t + x\gen x - \gen u\rangle\\
\mathsf{A}_{3.4}&=\langle \gen x, \gen t, t\gen t + [t+x]\gen x\rangle\\
\mathsf{A}_{3.4}&=\langle \gen x, \gen u, [x+u]\gen x, u\gen u\rangle\\
\mathsf{A}_{3.4}&=\langle \gen x, \gen t, t\gen t + [t+x]\gen x + u\gen u\rangle\\
\mathsf{A}_{3.4}&=\langle \gen x, \gen u, t\gen t + [x+u]\gen x + u\gen u \rangle.
\end{align*}

\section{Realizations of $\mathsf{A}_{3.5}$.} We have $\mathsf{A}_{3.5}=\langle e_1, e_2, e_3\rangle$ with the commutation relations $[e_1, e_2]=0,\; [e_1, e_3]=e_1,\; [e_2, e_3]=e_2$. Again, $\langle e_1, e_2\rangle$ is the two-dimensional abelian ideal of $\mathsf{A}_{3.5}$.

\medskip\noindent$\underline{{\rm rank}\,\mathsf{A}_{3.5}=1}:$ In this case ${\rm rank}\,\langle e_1, e_2\rangle=1$ so we have the canonical realization $\langle e_1, e_2\rangle=\langle \gen u, x\gen u\rangle$ and we must also have $e_3=c(t,x,u)\gen u$ because ${\rm rank}\,\mathsf{A}_{3.5}=1$. Since $[e_1, e_3]=e_1$ we have $c_u=1$, and $[e_2, e_3]=e_2$ gives $xc_u=x$ so that $c=u+c(t,x)$ and $e_3=[u+c(t,x)]\gen u$.

The residual equivalence group is $\mathscr{E}(e_1, e_2): t'=T(t),\; x'=x,\; u'=u+U(t,x)$. Under such a transformation, $e_3$ is mapped to
\[
e'_3=[u'+c(t,x)-U(t,x)]\gen {u'}.
\]
We choose $U(t,x)=c(t,x)$ and we obtain the canonical form $e_3=u\gen u$. We have the realization
\[
\mathsf{A}_{3.5}=\langle \gen u, x\gen u, u\gen u\rangle.
\]
An elementary calculation shows that this algebra linearizes the evolution equation: the function $F$ is independent of $(u, u_1, u_2)$ and $G$ is linear in $u_2$ and independent of $(u, u_1)$.

\medskip\noindent $\underline{{\rm rank}\,\mathsf{A}_{3.5}=2,\; {\rm rank}\,\langle e_1, e_2\rangle=1}:$ Again we have $\langle e_1, e_2\rangle=\langle \gen u, x\gen u\rangle$. Putting $e_3=a(t)\gen t + b\gen x + c\gen u$. From $[e_1, e_3]=e_1$ we find that $b_u=0,\; c_u=1$ so that $b=b(t,x),\; c=u+c(t,x)$. Then $[e_2, e_3]=e_2$ gives $xc_u-b=x$, from which we find that $b=0$. We also note that we must then have $a(t)\neq 0$ in order for us to have ${\rm rank}\,\mathsf{A}_{3.5}=2$. So $e_3=a(t)\gen t + [u+c(t,x)]\gen u$. The residual equivalence group is $\mathscr{E}(e_1, e_2): t'=T(t),\; x'=x,\; u'=u+U(t,x)$. Under such a transformation, $e_3$ is mapped to
\[
e'_3=a(t)\dot{T}(t)\gen {t'} + [u'+c(t,x)-U(t,x) + a(t)U_t(t,x)]\gen {u'}.
\]
Since $a(t)\neq 0$ we may choose $U(t,x)$ and $T(t)$ so that $c(t,x)-U(t,x) + a(t)U_t(t,x)=0$ and $a(t)\dot{T}(t)=T(t)$, giving $e'_3=t'\gen {t'} + u'\gen {u'}$, so that $e_3=t\gen t + u\gen u$ in canonical form, and we then have the realization
\[
\mathsf{A}_{3.5}=\langle \gen u, x\gen u, t\gen t + u\gen u\rangle.
\]

\medskip\noindent $\underline{{\rm rank}\,\mathsf{A}_{3.5}=2,\; {\rm rank}\,\langle e_1, e_2\rangle=2}:$

\medskip\noindent $\underline{\langle e_1, e_2\rangle=\langle \gen t, \gen x\rangle}:$ Because ${\rm rank}\,\mathsf{A}_{3.5}=2$ we have $e_3=a(t)\gen t + b\gen x$ (then $e_1\wedge e_2\wedge e_3=0$). from $[e_1, e_3]=e_1$ we have $a(t)=t+k,\, b_t=0$ and then $[e_2, e_3]=e_2$ gives $b_x=1$ so that $b=x+b(u)$ and $e_3=[t+k]\gen t + [x+b(u)]\gen x$. The residual equivalence group is $\mathscr{E}(e_1, e_2): t'=t+l,\; x'=x + Y(u),\; u'=U(u)$ with $U'(u)\neq 0$. Under such a transformation, $e_3$ is mapped to
\[
e'_3=[t'+l-k]\gen {t'} + [x'+b(u)-Y(u)]\gen {u'}.
\]
choose $l=k$ and $Y(u)$ so that $Y(u)=b(u)$ to give $e'_3=t'\gen {t'} + x'\gen {x'}$, so that $e_3=t\gen t + x\gen x$ in canonical form. We then have the realization
\[
\mathsf{A}_{3.5}=\langle \gen t, \gen x, t\gen t + x\gen x\rangle.
\]
We note that the other ordered choice $\langle e_1, e_2\rangle=\langle \gen x, \gen t\rangle$ gives the same algebra in canonical form.

\medskip\noindent $\underline{\langle e_1, e_2\rangle=\langle \gen x, \gen u\rangle}:$ Because ${\rm rank}\,\mathsf{A}_{3.5}=2$ we have $e_3=b\gen x + c\gen u$ (then $e_1\wedge e_2\wedge e_3=0$). From $[e_1, e_3]=e_1$ we obtain $e_3=[x+b(t,u)]\gen x + c(t,u)\gen u$ and from $[e_2, e_3]=e_2$ we have $e_3=[x+b(t)]\gen x + [u + c(t)]\gen u$.

The residual equivalence group is $\mathscr{E}: t'=T(t),\; x'=x+Y(t),\; u'=u+U(t)$ with $\dot{T}(t)\neq 0$. Under such a transformation, $e_3$ is mapped to
\[
e'_3=[x'+b(t)-Y(t)]\gen {x'} + [u'+c(t)-U(t)]\gen {u'}.
\]
We choose $Y(t)=b(t)$ and $U(t)=c(t)$ so that we have $e_3=x\gen x + u\gen u$ in canonical form, and we have the realization
\[
\mathsf{A}_{3.5}=\langle \gen x, \gen u, x\gen x + u\gen u\rangle.
\]

\medskip\noindent $\underline{{\rm rank}\,\mathsf{A}_{3.5}=3}:$ In this case we have ${\rm rank}\,\langle e_1, e_2\rangle=2$.

\medskip\noindent $\underline{\langle e_1, e_2\rangle=\langle \gen t, \gen x\rangle}:$ Put $e_3=a(t)\gen t + b\gen x + c\gen u$. ${\rm rank}\,\mathsf{A}_{3.5}=3$ so $e_1\wedge e_2\wedge e_3\neq 0$ so that $c(t,x,u)\neq 0$. From $[e_1, e_3]=e_1$ we have $\dot{a}(t)=1,\; b_t=c_t=0$ and from $[e_2, e_3]=e_2$ we have $b_x=1,\; c_x=0$ from which we have $e_3=[t+k]\gen t + [x+b(u)]\gen x + c(u)\gen u$.

The residual equivalence group is $\mathscr{E}(e_1, e_2): t'=t+l,\; x'=x + Y(u),\; u'=U(u)$ with $U'(u)\neq 0$. Under such a transformation, $e_3$ is mapped to
\[
e'_3=[t'+k-l]\gen {t'} + [x'+b(u)-Y(u)]\gen {x'} + c(u)U'(u)\gen {u'}.
\]
We choose $l=k$ and $Y(u)=b(u)$ and because $c(u)\neq 0$ we may choose $U(u)$ so that $c(u)U'(u)=-U(u)$ so that $e'_3=t'\gen {t'}+x'\gen {x'}+u'\gen {u'}$. This gives $e_3=t\gen t + x\gen x - u\gen u$ in canonical form and we have the realization
\[
\mathsf{A}_{3.5}=\langle \gen t, \gen x, t\gen t +  x\gen x + u\gen u\rangle.
\]
The same algebra is obtained with the ordered choice $\langle e_1, e_2\rangle=\langle \gen x, \gen t\rangle$.

\medskip\noindent $\underline{\langle e_1, e_2\rangle=\langle \gen x, \gen u\rangle}:$ Because ${\rm rank}\,\mathsf{A}_{3.5}=3$ we have $e_3=a(t)\gen t + b\gen x + c\gen u$ with $a(t)\neq 0$ (then $e_1\wedge e_2\wedge e_3\neq 0$). From $[e_1, e_3]=e_1$ we obtain $e_3=[x+b(t,u)]\gen x + c(t,u)\gen u$ and from $[e_2, e_3]=e_2$ we have $e_3=a(t)\gen t + [x+b(t)]\gen x + [u + c(t)]\gen u$.

The residual equivalence group is $\mathscr{E}: t'=T(t),\; x'=x+Y(t),\; u'=u+U(t)$ with $\dot{T}(t)\neq 0$. Under such a transformation, $e_3$ is mapped to
\[
e'_3=a(t)\dot{T}(t)\gen {t'} + [x' + b(t)-Y(t) + a(t)\dot{Y}(t)]\gen {t'} + [u' + c(t)-U(t) + a(t)\dot{U}(t)]\gen {u'}.
\]
For any choice of $a(t)\neq 0$ we may choose $Y(t),\; U(t)$ so that $b(t)-Y(t) + a(t)\dot{Y}(t)=0$ and $c(t)-U(t) + a(t)\dot{U}(t)=0$. We also choose $T(t)$ so that $a(t)\dot{T}(t)=T(t)$, so we have $e'_3=t'\gen {t'}+x'\gen {x'}+u'\gen {u'}$ and thus $e_3=t\gen t + x\gen x + u\gen u$ in canonical form. So we obtain the realization
\[
\mathsf{A}_{3.5}=\langle \gen x, \gen u, t\gen t +  x\gen x + u\gen u\rangle.
\]

\subsection{The inequivalent non-linearizing realizations of $\mathsf{A}_{3.5}$:}

\begin{align*}
\mathsf{A}_{3.5}&=\langle \gen u, x\gen u, t\gen t + u\gen u\rangle\\
\mathsf{A}_{3.5}&=\langle \gen t, \gen x, t\gen t + x\gen x\rangle\\
\mathsf{A}_{3.5}&=\langle \gen x, \gen u, x\gen x + u\gen u\rangle\\
\mathsf{A}_{3.5}&=\langle \gen t, \gen x, t\gen t +  x\gen x + u\gen u\rangle\\
\mathsf{A}_{3.5}&=\langle \gen x, \gen u, t\gen t +  x\gen x + u\gen u\rangle.
\end{align*}

\section{Realizations of $\mathsf{A}_{3.6}$.} We have $\mathsf{A}_{3.6}=\langle e_1, e_2, e_3\rangle$ with the commutation relations $[e_1, e_2]=0,\; [e_1, e_3]=e_1,\; [e_2, e_3]=-e_2$. Again, $\langle e_1, e_2\rangle$ is the two-dimensional abelian ideal of $\mathsf{A}_{3.6}$.

\medskip\noindent$\underline{{\rm rank}\,\mathsf{A}_{3.6}=1}:$ In this case ${\rm rank}\,\langle e_1, e_2\rangle=1$ so we have the canonical realization $\langle e_1, e_2\rangle=\langle \gen u, x\gen u\rangle$ and we must also have $e_3=c(t,x,u)\gen u$ because ${\rm rank}\,\mathsf{A}_{3.5}=1$. Since $[e_1, e_3]=e_1$ we have $c_u=1$, and $[e_2, e_3]=e_2$ gives $xc_u=-x$ so that we have no realization in this case.

\medskip\noindent $\underline{{\rm rank}\,\mathsf{A}_{3.6}=2,\; {\rm rank}\,\langle e_1, e_2\rangle=1}:$ Again we have $\langle e_1, e_2\rangle=\langle \gen u, x\gen u\rangle$. Putting $e_3=a(t)\gen t + b\gen x + c\gen u$. From $[e_1, e_3]=e_1$ we find that $b_u=0,\; c_u=1$ so that $b=b(t,x),\; c=u+c(t,x)$. Then $[e_2, e_3]=-e_2$ gives $xc_u-b=-x$, from which we find that $b=2x$.  So $e_3=a(t)\gen t + 2x\gen x + [u+c(t,x)]\gen u$. The residual equivalence group is $\mathscr{E}(e_1, e_2): t'=T(t),\; x'=x,\; u'=u+U(t,x)$. Under such a transformation, $e_3$ is mapped to
\[
e'_3=a(t)\dot{T}(t)\gen {t'} + 2x'\gen {x'} + [u'+c(t,x)-U(t,x) + a(t)U_t(t,x)]\gen {u'}.
\]
We may always choose $U(t,x)$ so that $c(t,x)-U(t,x) + a(t)U_t(t,x)=0$ giving $e'_3=a(t)\dot{T}(t)\gen {t'} + 2x'\gen {x'} + u'\gen {u'}$. If $a(t)=0$ we have $e_3=2x\gen x + u\gen u$ in canonical form. If $a(t)\neq 0$ we choose $T(t)$ so that $a(t)\dot{T}(t)=T(t)$ giving $e_3= t\gen t + 2x\gen x + u\gen u$ in canonical form. We then have the realizations
\begin{align*}
\mathsf{A}_{3.6}&=\langle \gen u, x\gen u, 2x\gen x + u\gen u\rangle\\
\mathsf{A}_{3.6}&=\langle \gen u, x\gen u, t\gen t + 2x\gen x + u\gen u\rangle
\end{align*}

\medskip\noindent $\underline{{\rm rank}\,\mathsf{A}_{3.6}=2,\; {\rm rank}\,\langle e_1, e_2\rangle=2}:$

\medskip\noindent $\underline{\langle e_1, e_2\rangle=\langle \gen t, \gen x\rangle}:$ Because ${\rm rank}\,\mathsf{A}_{3.5}=2$ we have $e_3=a(t)\gen t + b\gen x$ (then $e_1\wedge e_2\wedge e_3=0$). From $[e_1, e_3]=e_1$ we have $a(t)=t+k,\, b_t=0$ and then $[e_2, e_3]=-e_2$ gives $b_x=-1$ so that $b=-x+b(u)$ and $e_3=[t+k]\gen t + [-x+b(u)]\gen x$. The residual equivalence group is $\mathscr{E}(e_1, e_2): t'=t+l,\; x'=x + Y(u),\; u'=U(u)$ with $U'(u)\neq 0$. Under such a transformation, $e_3$ is mapped to
\[
e'_3=[t'+l-k]\gen {t'} + [-x'+b(u)+Y(u)]\gen {u'}.
\]
choose $l=k$ and $Y(u)$ so that $Y(u)+b(u)=0$ to give $e'_3=t'\gen {t'} - x'\gen {x'}$, so that $e_3=t\gen t + x\gen x$ in canonical form. We then have the realization
\[
\mathsf{A}_{3.6}=\langle \gen t, \gen x, t\gen t - x\gen x\rangle.
\]
We note that the other ordered choice $\langle e_1, e_2\rangle=\langle \gen x, \gen t\rangle$ gives the realization
\[
\mathsf{A}_{3.6}=\langle \gen x, \gen t, -t\gen t + x\gen x\rangle.
\]
which is the same algebra.

\medskip\noindent $\underline{\langle e_1, e_2\rangle=\langle \gen x, \gen u\rangle}:$ Because ${\rm rank}\,\mathsf{A}_{3.6}=2$ we have $e_3=b\gen x + c\gen u$ (then $e_1\wedge e_2\wedge e_3=0$). From $[e_1, e_3]=e_1$ we obtain $e_3=[x+b(t,u)]\gen x + c(t,u)\gen u$ and from $[e_2, e_3]=-e_2$ we have $e_3=[x+b(t)]\gen x + [-u + c(t)]\gen u$.

The residual equivalence group is $\mathscr{E}: t'=T(t),\; x'=x+Y(t),\; u'=u+U(t)$ with $\dot{T}(t)\neq 0$. Under such a transformation, $e_3$ is mapped to
\[
e'_3=[x'+b(t)-Y(t)]\gen {x'} + [-u'+c(t)+U(t)]\gen {u'}.
\]
We choose $Y(t)=b(t)$ and $U(t)+c(t)=0$ so that we have $e_3=x\gen x - u\gen u$ in canonical form, and we have the realization
\[
\mathsf{A}_{3.6}=\langle \gen x, \gen u, x\gen x - u\gen u\rangle.
\]

\medskip\noindent $\underline{{\rm rank}\,\mathsf{A}_{3.6}=3}:$ In this case we have ${\rm rank}\,\langle e_1, e_2\rangle=2$.

\medskip\noindent $\underline{\langle e_1, e_2\rangle=\langle \gen t, \gen x\rangle}:$ Put $e_3=a(t)\gen t + b\gen x + c\gen u$. Since ${\rm rank}\,\mathsf{A}_{3.3}=3$ we have $e_1\wedge e_2\wedge e_3\neq 0$ so that $c(t,x,u)\neq 0$. From $[e_1, e_3]=e_1$ we have $\dot{a}(t)=1,\; b_t=c_t=0$ and from $[e_2, e_3]=-e_2$ we have $b_x=-1,\; c_x=0$ from which we have $e_3=[t+k]\gen t + [-x+b(u)]\gen x + c(u)\gen u$.

The residual equivalence group is $\mathscr{E}(e_1, e_2): t'=t+l,\; x'=x + Y(u),\; u'=U(u)$ with $U'(u)\neq 0$. Under such a transformation, $e_3$ is mapped to
\[
e'_3=[t'+k-l]\gen {t'} + [-x'+b(u)-Y(u)]\gen {x'} + c(u)U'(u)\gen {u'}.
\]
We choose $l=k$ and $Y(u)=b(u)$ and because $c(u)\neq 0$ we may choose $U(u)$ so that $c(u)U'(u)=-U(u)$ so that $e'_3=t'\gen {t'}+x'\gen {x'}-u'\gen {u'}$. This gives $e_3=t\gen t - x\gen x + u\gen u$ in canonical form and we have the realization
\[
\mathsf{A}_{3.6}=\langle \gen t, \gen x, t\gen t -  x\gen x - u\gen u\rangle.
\]
The same algebra is obtained with the ordered choice $\langle e_1, e_2\rangle=\langle \gen x, \gen t\rangle$: we choose $U(u)$ so that $c(u)U(u)=U(u)$ and we obtain the canonical form $e_3=-t\gen t + x\gen x + u\gen u$ and this gives the realization $\mathsf{A}_{3.6}=\langle \gen t, \gen x, -t\gen t +  x\gen x + u\gen u\rangle$, which is the same algebra.

\medskip\noindent $\underline{\langle e_1, e_2\rangle=\langle \gen x, \gen u\rangle}:$ Because ${\rm rank}\,\mathsf{A}_{3.6}=3$ we have $e_3=a(t)\gen t + b\gen x + c\gen u$ with $a(t)\neq 0$ (then $e_1\wedge e_2\wedge e_3\neq 0$). From $[e_1, e_3]=e_1$ we obtain $e_3=[x+b(t,u)]\gen x + c(t,u)\gen u$ and from $[e_2, e_3]=-e_2$ we have $e_3=a(t)\gen t + [x+b(t)]\gen x + [-u + c(t)]\gen u$.

The residual equivalence group is $\mathscr{E}: t'=T(t),\; x'=x+Y(t),\; u'=u+U(t)$ with $\dot{T}(t)\neq 0$. Under such a transformation, $e_3$ is mapped to
\[
e'_3=a(t)\dot{T}(t)\gen {t'} + [x' + b(t)-Y(t) + a(t)\dot{Y}(t)]\gen {t'} + [-u' + c(t)+U(t) + a(t)\dot{U}(t)]\gen {u'}.
\]
For any choice of $a(t)\neq 0$ we may choose $Y(t),\; U(t)$ so that $b(t)-Y(t) + a(t)\dot{Y}(t)=0$ and $c(t)+U(t) + a(t)\dot{U}(t)=0$. We also choose $T(t)$ so that $a(t)\dot{T}(t)=T(t)$, so we have $e'_3=t'\gen {t'}+x'\gen {x'}-u'\gen {u'}$ and thus $e_3=t\gen t + x\gen x - u\gen u$ in canonical form. So we obtain the realization
\[
\mathsf{A}_{3.6}=\langle \gen x, \gen u, t\gen t +  x\gen x - u\gen u\rangle.
\]

\subsection{The inequivalent non-linearizing realizations of $\mathsf{A}_{3.6}$:}

\begin{align*}
\mathsf{A}_{3.6}&=\langle \gen u, x\gen u, 2x\gen x + u\gen u\rangle\\
\mathsf{A}_{3.6}&=\langle \gen u, x\gen u, t\gen t + 2x\gen x + u\gen u\rangle\\
\mathsf{A}_{3.6}&=\langle \gen t, \gen x, t\gen t - x\gen x\rangle\\
\mathsf{A}_{3.6}&=\langle \gen x, \gen u, x\gen x - u\gen u\rangle\\
\mathsf{A}_{3.6}&=\langle \gen t, \gen x, t\gen t -  x\gen x - u\gen u\rangle\\
\mathsf{A}_{3.6}&=\langle \gen x, \gen u, t\gen t +  x\gen x - u\gen u\rangle.
\end{align*}

\section{Realizations of $\mathsf{A}_{3.7}$.} We have $\mathsf{A}_{3.7}=\langle e_1, e_2, e_3\rangle$ with the commutation relations $[e_1, e_2]=0,\; [e_1, e_3]=e_1,\; [e_2, e_3]=qe_2$ with $0<|q|<1$. Again, $\langle e_1, e_2\rangle$ is the two-dimensional abelian ideal of $\mathsf{A}_{3.7}$.

\medskip\noindent$\underline{{\rm rank}\,\mathsf{A}_{3.7}=1}:$ In this case ${\rm rank}\,\langle e_1, e_2\rangle=1$ so we have the canonical realization $\langle e_1, e_2\rangle=\langle \gen u, x\gen u\rangle$ and we must also have $e_3=c(t,x,u)\gen u$ because ${\rm rank}\,\mathsf{A}_{3.7}=1$. Since $[e_1, e_3]=e_1$ we have $c_u=1$, and $[e_2, e_3]=qe_2$ gives $xc_u=qx$ which is a contradiction since $0<|q|<1$. So we have no realization in this case.

\medskip\noindent $\underline{{\rm rank}\,\mathsf{A}_{3.7}=2,\; {\rm rank}\,\langle e_1, e_2\rangle=1}:$ Again we have $\langle e_1, e_2\rangle=\langle \gen u, x\gen u\rangle$. Putting $e_3=a(t)\gen t + b\gen x + c\gen u$. From $[e_1, e_3]=e_1$ we find that $b_u=0,\; c_u=1$ so that $b=b(t,x),\; c=u+c(t,x)$. Then $[e_2, e_3]=qe_2$ gives $xc_u-b=qx$, from which we find that $b=(1-q)x$.  So $e_3=a(t)\gen t + (1-q)x\gen x + [u+c(t,x)]\gen u$. The residual equivalence group is $\mathscr{E}(e_1, e_2): t'=T(t),\; x'=x,\; u'=u+U(t,x)$. Under such a transformation, $e_3$ is mapped to
\[
e'_3=a(t)\dot{T}(t)\gen {t'} + (1-q)x'\gen {x'} + [u'+c(t,x)-U(t,x) + a(t)U_t + (1-q)xU_x]\gen {u'}.
\]
Since $1-q\neq 0$ we may choose $U(t,x)$ so that $c(t,x)-U(t,x) + a(t)U_t(t,x) + (1-q)xU_x(t,x)=0$. Then if $a(t)\neq 0$ we take $T(t)$ so that $a(t)\dot{T}(t)=T(t)$, and this gives $e'_3=t'\gen {t'} + (1-q)x'\gen {x'} + u'\gen {u'}$, so that $e_3=t\gen t + (1-q)x\gen x +  u\gen u$ in canonical form. If $a(t)=0$ then we have $e_3=(1-q)x\gen x +  u\gen u$ in canonical form and we then have the realizations
\begin{align*}
\mathsf{A}_{3.7}&=\langle \gen u, x\gen u, (1-q)x\gen x +  u\gen u\rangle\\
\mathsf{A}_{3.7}&=\langle \gen u, x\gen u, t\gen t + (1-q)x\gen x +  u\gen u\rangle.
\end{align*}

\medskip\noindent $\underline{{\rm rank}\,\mathsf{A}_{3.7}=2,\; {\rm rank}\,\langle e_1, e_2\rangle=2}:$

\medskip\noindent $\underline{\langle e_1, e_2\rangle=\langle \gen t, \gen x\rangle}:$ Because ${\rm rank}\,\mathsf{A}_{3.5}=2$ we have $e_3=a(t)\gen t + b\gen x$ (then $e_1\wedge e_2\wedge e_3=0$). From $[e_1, e_3]=e_1$ we have $a(t)=t+k,\, b_t=0$ and then $[e_2, e_3]=qe_2$ gives $b_x=q$ so that $b=qx+b(u)$ and $e_3=[t+k]\gen t + [qx+b(u)]\gen x$. The residual equivalence group is $\mathscr{E}(e_1, e_2): t'=t+l,\; x'=x + Y(u),\; u'=U(u)$ with $U'(u)\neq 0$. Under such a transformation, $e_3$ is mapped to
\[
e'_3=[t'+l-k]\gen {t'} + [qx'+b(u)-qY(u)]\gen {u'}.
\]
choose $l=k$ and $Y(u)$ so that $qY(u)=b(u)$ to give $e'_3=t'\gen {t'} + qx'\gen {x'}$, so that $e_3=t\gen t + qx\gen x$ in canonical form. We then have the realization
\[
\mathsf{A}_{3.7}=\langle \gen t, \gen x, t\gen t + qx\gen x\rangle.
\]
We note that the other ordered choice $\langle e_1, e_2\rangle=\langle \gen x, \gen t\rangle$ gives the algebra
\[
\mathsf{A}_{3.7}=\langle \gen x, \gen t, qt\gen t + x\gen x\rangle.
\]
as is easily computed using the above arguments. So we have the realizations
\begin{align*}
\mathsf{A}_{3.7}&=\langle \gen t, \gen x, t\gen t + qx\gen x\rangle\\
\mathsf{A}_{3.7}&=\langle \gen x, \gen t, qt\gen t + x\gen x\rangle.
\end{align*}

\medskip\noindent $\underline{\langle e_1, e_2\rangle=\langle \gen x, \gen u\rangle}:$ Because ${\rm rank}\,\mathsf{A}_{3.7}=2$ we have $e_3=b\gen x + c\gen u$ (then $e_1\wedge e_2\wedge e_3=0$). From $[e_1, e_3]=e_1$ we obtain
$e_3=[x+b(t,u)]\gen x + c(t,u)\gen u$ and from $[e_2, e_3]=qe_2$ we have $e_3=[x+b(t)]\gen x + [qu + c(t)]\gen u$.

The residual equivalence group is $\mathscr{E}: t'=T(t),\; x'=x+Y(t),\; u'=u+U(t)$ with $\dot{T}(t)\neq 0$. Under such a transformation, $e_3$ is mapped to
\[
e'_3=[x'+b(t)-Y(t)]\gen {x'} + [qu'+c(t)-qU(t)]\gen {u'}.
\]
We choose $Y(t)=b(t)$ and $U(t)$ so that  $qU(t)=c(t)$ giving $e_3=x\gen x + qu\gen u$ in canonical form, and we have the realization
\[
\mathsf{A}_{3.7}=\langle \gen x, \gen u, x\gen x + qu\gen u\rangle.
\]

\medskip\noindent $\underline{{\rm rank}\,\mathsf{A}_{3.7}=3}:$ In this case we have ${\rm rank}\,\langle e_1, e_2\rangle=2$.

\medskip\noindent $\underline{\langle e_1, e_2\rangle=\langle \gen t, \gen x\rangle}:$ Put $e_3=a(t)\gen t + b\gen x + c\gen u$. ${\rm rank}\,\mathsf{A}_{3.7}=3$ so $e_1\wedge e_2\wedge e_3\neq 0$ so that $c(t,x,u)\neq 0$. From $[e_1, e_3]=e_1$ we have $\dot{a}(t)=1,\; b_t=c_t=0$ and from $[e_2, e_3]=qe_2$ we have $b_x=q,\; c_x=0$ from which we have $e_3=[t+k]\gen t + [qx+b(u)]\gen x + c(u)\gen u$.

The residual equivalence group is $\mathscr{E}(e_1, e_2): t'=t+l,\; x'=x + Y(u),\; u'=U(u)$ with $U'(u)\neq 0$. Under such a transformation, $e_3$ is mapped to
\[
e'_3=[t'+k-l]\gen {t'} + [qx'+b(u)-qY(u) + c(u)Y'(u)]\gen {x'} + c(u)U'(u)\gen {u'}.
\]
We choose $l=k$ and $b(u)-qY(u) + c(u)Y'(u)=0$ and because $c(u)\neq 0$ we may choose $U(u)$ so that $c(u)U'(u)=U(u)$ so that $e'_3=t'\gen {t'}+qx'\gen {x'}+u'\gen {u'}$. This gives $e_3=t\gen t + qx\gen x + u\gen u$ in canonical form and we have the realization
\[
\mathsf{A}_{3.7}=\langle \gen t, \gen x, t\gen t +  qx\gen x + u\gen u\rangle.
\]
We note that with the ordered choice $\langle e_1, e_2\rangle=\langle \gen x, \gen t\rangle$ we obtain the realization
\[
\mathsf{A}_{3.7}=\langle \gen x, \gen t, qt\gen t +  x\gen x + u\gen u\rangle.
\]
So we have the realizations
\begin{align*}
\mathsf{A}_{3.7}&=\langle \gen t, \gen x, t\gen t +  qx\gen x + u\gen u\rangle\\
\mathsf{A}_{3.7}&=\langle \gen x, \gen t, qt\gen t +  x\gen x + u\gen u\rangle.
\end{align*}

\medskip\noindent $\underline{\langle e_1, e_2\rangle=\langle \gen x, \gen u\rangle}:$ Because ${\rm rank}\,\mathsf{A}_{3.7}=3$ we have $e_3=a(t)\gen t + b\gen x + c\gen u$ with $a(t)\neq 0$ (then $e_1\wedge e_2\wedge e_3\neq 0$). From $[e_1, e_3]=e_1$ we obtain $e_3=a(t)\gen t + [x+b(t,u)]\gen x + c(t,u)\gen u$ and from $[e_2, e_3]=qe_2$ we have $e_3=a(t)\gen t + [x+b(t)]\gen x + [qu + c(t)]\gen u$.

The residual equivalence group is $\mathscr{E}: t'=T(t),\; x'=x+Y(t),\; u'=u+U(t)$ with $\dot{T}(t)\neq 0$. Under such a transformation, $e_3$ is mapped to
\[
e'_3=a(t)\dot{T}(t)\gen {t'} + [x' + b(t)-Y(t) + a(t)\dot{Y}(t)]\gen {t'} + [qu' + c(t)-qU(t) + a(t)\dot{U}(t)]\gen {u'}.
\]
For any choice of $a(t)\neq 0$ we may choose $Y(t),\; U(t)$ so that $b(t)-Y(t) + a(t)\dot{Y}(t)=0$ and $c(t)-qU(t) + a(t)\dot{U}(t)=0$. We also choose $T(t)$ so that $a(t)\dot{T}(t)=T(t)$, so we have $e'_3=t'\gen {t'}+x'\gen {x'}+qu'\gen {u'}$ and thus $e_3=t\gen t + x\gen x + qu\gen u$ in canonical form. So we obtain the realization
\[
\mathsf{A}_{3.7}=\langle \gen x, \gen u, t\gen t +  x\gen x + qu\gen u\rangle.
\]

\subsection{The inequivalent non-linearizing realizations of $\mathsf{A}_{3.7}$:}

\begin{align*}
\mathsf{A}_{3.7}&=\langle \gen u, x\gen u, (1-q)x\gen x + u\gen u\rangle\\
\mathsf{A}_{3.7}&=\langle \gen u, x\gen u, t\gen t + (1-q)x\gen x + u\gen u\rangle\\
\mathsf{A}_{3.7}&=\langle \gen t, \gen x, t\gen t + qx\gen x\rangle\\
\mathsf{A}_{3.7}&=\langle \gen x, \gen t, qt\gen t + x\gen x\rangle\\
\mathsf{A}_{3.7}&=\langle \gen x, \gen u, x\gen x + qu\gen u\rangle\\
\mathsf{A}_{3.7}&=\langle \gen t, \gen x, t\gen t +  qx\gen x + u\gen u\rangle\\
\mathsf{A}_{3.7}&=\langle \gen x, \gen t, qt\gen t +  x\gen x + u\gen u\rangle\\
\mathsf{A}_{3.7}&=\langle \gen x, \gen u, t\gen t +  x\gen x + qu\gen u\rangle.
\end{align*}

\section{Realizations of $\mathsf{A}_{3.8}$.}  We have $\mathsf{A}_{3.8}=\langle e_1, e_2, e_3\rangle$ with $[e_1, e_2]=0,\; [e_1, e_2]=-e_2,\; [e_2, e_3]=e_1$. Note that $\langle e_1, e_2\rangle$ is the abelian ideal of $\mathsf{A}_{3.8}$.

\medskip\noindent $\underline{{\rm rank}\,\mathsf{A}_{3.8}=1}:$ In this case ${\rm rank}\,\langle e_1, e_2\rangle$ so that $\langle e_1, e_2\rangle=\langle \gen u, x\gen u\rangle$. Then $e_3=c(t,x,u)\gen u$ since ${\rm rank}\,\mathsf{A}_{3.8}=1$. First, $[e_1, e_3]=-e_2$ gives $c_u=-x$. Then $[e_2, e_3]=e_1$ gives $xc_u=1$ and this is a contradiction, so we have no realization in this case.

\medskip\noindent $\underline{{\rm rank}\,\mathsf{A}_{3.8}=2,\; {\rm rank}\,\langle e_1, e_2\rangle=1}:$ Again, we have $\langle e_1, e_2\rangle=\langle \gen u, x\gen u\rangle$ and we put $e_3=a(t)\gen t + b\gen x + c\gen u$ with $a(t)^2+b^2\neq 0$ since ${\rm rank}\,\mathsf{A}_{3.8}=2$. First, $[e_1, e_3]=-e_2$ gives $b_u=0,\; c_u=-x$ and then $[e_2, e_3]=e_1$ gives us $b=xc_u-1=-(1+x^2)$ so that
\[
e_3=a(t)\gen t -(1+x^2)\gen x + [-xu + c(t,x)]\gen u.
\]
The residual equivalence group is $\mathscr{E}: t'=T(t),\; x'=x,\; u'=u+U(t,x)$. Under such a transformation $e_3$ is mapped to
\[
e'_3=a(t)\dot{T}(t)\gen {t'} -(1+x'^2)\gen {x'}+[-x'u'+c(t,x)+xU+a(t)U_t - (1+x^2)U_x]\gen {u'}.
\]
We may always choose $U$ so that $c(t,x)+xU+a(t)U_t - (1+x^2)U_x=0$ whether $a(t)=0$ or $a(t)\neq 0$. So we then have
\[
e'_3=a(t)\dot{T}(t)\gen {t'} -(1+x'^2)\gen {x'} - x'u'\gen {u'}.
\]
If $a(t)=0$ we have $e'_3=-(1+x'^2)\gen {x'} - x'u'\gen {u'}$ giving $e_3=-(1+x^2)\gen x -xu\gen u$ in canonical form. If $a(t)\neq 0$ then we choose $T(t)$ so that $a(t)\dot{T}(t)=-T(t)$ so that $e'_3=-t'\gen {t'} -(1+x'^2)\gen {x'} - x'u'\gen {u'}$, giving $e_3=-t\gen t - (1+x^2)\gen x -xu \gen u$. Thus we have the two realizations
\begin{align*}
\mathsf{A}_{3.8}&=\langle \gen u, x\gen u, -(1+x^2)\gen x - xu\gen u\rangle\\
\mathsf{A}_{3.8}&=\langle \gen u, x\gen u, -t\gen t - (1+x^2)\gen x - xu\gen u\rangle.
\end{align*}

\medskip\noindent $\underline{{\rm rank}\,\mathsf{A}_{3.8}=2,\; {\rm rank}\,\langle e_1, e_2\rangle=2}:$

\medskip\noindent $\underline{\langle e_1, e_2\rangle=\langle \gen t, \gen x\rangle}:$ In this case, we have $e_3=a(t)\gen t + b\gen x$ since we must have $e_1\wedge e_2\wedge e_3=0$. Then $[e_2, e_3]=e_1$ gives $b_x\gen x=\gen t$ which is evidently impossible, so we have no realization in this case.

\medskip\noindent $\underline{\langle e_1, e_2\rangle=\langle \gen x, \gen t\rangle}:$ In this case, we again have $e_3=a(t)\gen t + b\gen x$ since we must have $e_1\wedge e_2\wedge e_3=0$. Then $[e_1, e_3]=-e_2$ gives $b_x\gen x=-\gen t$ which is evidently impossible, so we have no realization in this case.

\medskip\noindent $\underline{\langle e_1, e_2\rangle=\langle \gen x, \gen u\rangle}:$ In this case, we have $e_3= b\gen x + c\gen u$ since we must have $e_1\wedge e_2\wedge e_3=0$. Then $[e_1, e_3]=-e_2$ gives $b_x=0,\; c_x=-1$ and then $[e_2, e_3]=e_1$ gives $b_u=1,\; c_u=0$ so we have $e_3=[u+b(t)]\gen x + [-x+c(t)]\gen u$.

The residual equivalence group is $\mathscr{E}: t'=T(t),\; x'=x+Y(t),\; u'=u+U(t)$. Under such a transformation, $e_3$ is mapped to
\[
e'_3=[u'+b(t)-U(t)]\gen {x'} + [-x'+c(t)+Y(t)]\gen {u'}.
\]
Choose $Y(t),\, U(t)$ so that $b(t)-U(t)=0,\; c(t)+Y(t)=0$ and we then have $e'_3=u'\gen {x'}-x'\gen {u'}$, so that $e_3=u\gen x - x\gen u$ in canonical form. This gives the realization:
\[
\mathsf{A}_{3.8}=\langle \gen x, \gen u, u\gen x - x\gen u \rangle.
\]

\medskip\noindent $\underline{{\rm rank}\,\mathsf{A}_{3.8}=3}$. In this case we must have ${\rm rank}\,\langle e_1, e_2\rangle=2$.

\medskip\noindent $\underline{\langle e_1, e_2\rangle=\langle \gen t, \gen x\rangle}:$ In this case, we have $e_3=a(t)\gen t + b\gen x + c\gen u$ with $c(t,x,u)\neq 0$ since ${\rm rank}\,\mathsf{A}_{3.8}=3$ and so we must have $e_1\wedge e_2\wedge e_3\neq 0$. Then $[e_2, e_3]=e_1$ gives $b_x\gen x=\gen t$ which is evidently impossible, so we have no realization in this case.

\medskip\noindent $\underline{\langle e_1, e_2\rangle=\langle \gen x, \gen t\rangle}:$ In this case, we have $e_3=a(t)\gen t + b\gen x + c\gen u$ with $c(t,x,u)\neq 0$. Then $[e_1, e_3]=-e_2$ gives $b_x\gen x + c_x\gen u=-\gen t$ which is impossible, so we have no realization in this case.

\medskip\noindent $\underline{\langle e_1, e_2\rangle=\langle \gen x, \gen u\rangle}:$ In this case, we have $e_3= a(t)\gen t +b\gen x + c\gen u$ with $a(t)\neq 0$ since we must have $e_1\wedge e_2\wedge e_3\neq 0$. Then $[e_1, e_3]=-e_2$ gives $b_x=0,\; c_x=-1$ and then $[e_2, e_3]=e_1$ gives $b_u=1,\; c_u=0$ so we have $e_3=a(t)\gen t + [u+b(t)]\gen x + [-x+c(t)]\gen u$.

The residual equivalence group is $\mathscr{E}: t'=T(t),\; x'=x+Y(t),\; u'=u+U(t)$. Under such a transformation, $e_3$ is mapped to
\[
e'_3=a(t)\dot{T}(t)\gen {t'} + [u'+b(t)-U(t)+a(t)\dot{Y}(t)]\gen {x'} + [-x'+c(t)+Y(t) + a(t)\dot{U}(t)]\gen {u'}.
\]
For $a(t)\neq 0$ we may always solve the system of differential equations $a(t)\dot{Y}(t)=U(t)-b(t),\;  a(t)\dot{U}(t)=-Y(t)-c(t)$ for $Y(t)$ and $U(t)$, and we choose $T(t)$ so that $a(t)\dot{T}(t)=T(t)$ from which we have $e'_3=t'\gen {t'} + u'\gen {x'}-x'\gen {u'}$, so that $e_3=t\gen t + u\gen x - x\gen u$ in canonical form. This gives the realization:
\[
\mathsf{A}_{3.8}=\langle \gen x, \gen u, t\gen t + u\gen x - x\gen u \rangle.
\]

\subsection{The inequivalent non-linearizing realizations of $\mathsf{A}_{3.8}$:}

\begin{align*}
\mathsf{A}_{3.8}&=\langle \gen u, x\gen u, -(1+x^2)\gen x - xu\gen u\rangle\\
\mathsf{A}_{3.8}&=\langle \gen u, x\gen u, -t\gen t - (1+x^2)\gen x - xu\gen u\rangle\\
\mathsf{A}_{3.8}&=\langle \gen x, \gen u, u\gen x - x\gen u \rangle\\
\mathsf{A}_{3.8}&=\langle \gen x, \gen u, t\gen t + u\gen x - x\gen u \rangle.
\end{align*}

\section{Realizations of $\mathsf{A}_{3.9}$.}  We have $\mathsf{A}_{3.9}=\langle e_1, e_2, e_3\rangle$ with $[e_1, e_2]=0,\; [e_1, e_2]=qe_1-e_2,\; [e_2, e_3]=e_1+qe_2$ with $q>0$. Note that $\langle e_1, e_2\rangle$ is the abelian ideal of $\mathsf{A}_{3.9}$.

\medskip\noindent $\underline{{\rm rank}\,\mathsf{A}_{3.9}=1}:$ In this case ${\rm rank}\,\langle e_1, e_2\rangle$ so that $\langle e_1, e_2\rangle=\langle \gen u, x\gen u\rangle$. Then $e_3=c(t,x,u)\gen u$ since ${\rm rank}\,\mathsf{A}_{3.9}=1$. First, $[e_1, e_3]=qe_1-e_2$ gives $c_u=q-x$. Then $[e_2, e_3]=e_1+qe_2$ gives $xc_u=1+qx$ and so we have $x(q-x)=1+qx$ which gives $x^2+1=0$, a contradiction. Hence we have no realization in this case.

\medskip\noindent $\underline{{\rm rank}\,\mathsf{A}_{3.9}=2,\; {\rm rank}\,\langle e_1, e_2\rangle=1}:$ Again, we have $\langle e_1, e_2\rangle=\langle \gen u, x\gen u\rangle$ and we put $e_3=a(t)\gen t + b\gen x + c\gen u$ with $a(t)^2+b^2\neq 0$ since ${\rm rank}\,\mathsf{A}_{3.9}=2$. First, $[e_1, e_3]=qe_1-e_2$ gives $b_u=0,\; c_u=q-x$ and then $[e_2, e_3]=e_1+qe_2$ gives us $xc_u-b=1+qx$ from which we obtain $b=-(1+x^2),\; c=(q-x)u+c(t,x)$ and then
\[
e_3=a(t)\gen t -(1+x^2)\gen x + [(q-x)u + c(t,x)]\gen u.
\]
The residual equivalence group is $\mathscr{E}: t'=T(t),\; x'=x,\; u'=u+U(t,x)$. Under such a transformation $e_3$ is mapped to
\[
e'_3=a(t)\dot{T}(t)\gen {t'} -(1+x'^2)\gen {x'}+[(q-x')u'+c(t,x)+(q-x)U+a(t)U_t - (1+x^2)U_x]\gen {u'}.
\]
We may always choose $U$ so that $c(t,x)+(q-x)U+a(t)U_t - (1+x^2)U_x=0$ whether $a(t)=0$ or $a(t)\neq 0$. So we then have
\[
e'_3=a(t)\dot{T}(t)\gen {t'} -(1+x'^2)\gen {x'} + (q-x')u'\gen {u'}.
\]
If $a(t)=0$ we have $e'_3=-(1+x'^2)\gen {x'} + (q-x')u'\gen {u'}$ giving $e_3=-(1+x^2)\gen x + (q-x)u\gen u$ in canonical form. If $a(t)\neq 0$ then we choose $T(t)$ so that $a(t)\dot{T}(t)=-T(t)$ so that $e'_3=-t'\gen {t'} -(1+x'^2)\gen {x'} + (q-x')u'\gen {u'}$, giving $e_3=-t\gen t - (1+x^2)\gen x + (q-x)u\gen u$. Thus we have the two realizations
\begin{align*}
\mathsf{A}_{3.9}&=\langle \gen u, x\gen u, -(1+x^2)\gen x + (q-x)u\gen u\rangle\\
\mathsf{A}_{3.9}&=\langle \gen u, x\gen u, -t\gen t - (1+x^2)\gen x + (q-x)u\gen u\rangle.
\end{align*}

\medskip\noindent $\underline{{\rm rank}\,\mathsf{A}_{3.9}=2,\; {\rm rank}\,\langle e_1, e_2\rangle=2}:$

\medskip\noindent $\underline{\langle e_1, e_2\rangle=\langle \gen t, \gen x \rangle}:$ In this case, we have $e_3=a(t)\gen t + b\gen x$ since we must have $e_1\wedge e_2\wedge e_3=0$. Then $[e_2, e_3]=e_1+qe_2$ gives $b_x\gen x=\gen t+q\gen x$ which is evidently impossible, so we have no realization in this case.

\medskip\noindent $\underline{\langle e_1, e_2\rangle=\langle \gen x, \gen t \rangle}:$ In this case, we again have $e_3=a(t)\gen t + b\gen x$ since we must have $e_1\wedge e_2\wedge e_3=0$. Then $[e_1, e_3]=qe_1-e_2$ gives $b_x\gen x=q\gen x-\gen t$ which is evidently impossible, so we have no realization in this case.

\medskip\noindent $\underline{\langle e_1, e_2\rangle=\langle \gen x, \gen u \rangle}:$ In this case, we have $e_3= b\gen x + c\gen u$ since we must have $e_1\wedge e_2\wedge e_3=0$. Then $[e_1, e_3]=qe_1-e_2$ gives $b_x=q,\; c_x=-1$ and then $[e_2, e_3]=e_1+qe_2$ gives $b_u=1,\; c_u=q$ so we have $e_3=[qx+u+b(t)]\gen x + [qu-x+c(t)]\gen u$.

The residual equivalence group is $\mathscr{E}: t'=T(t),\; x'=x+Y(t),\; u'=u+U(t)$. Under such a transformation, $e_3$ is mapped to
\[
e'_3=[qx'+u'+b(t)-qY(t)-U(t)]\gen {x'} + [qu'-x+c(t)+Y(t)-qU(t)]\gen {u'}.
\]
Choose $Y(t),\, U(t)$ so that $b(t)-qY(t)-U(t)=0,\; c(t)+Y(t)-qU(t)=0$ and we then have $e'_3=[qx'+u']\gen {x'}+ [qu'-x']\gen {u'}$, so that $e_3=[qx+u]\gen x + [qu-x]\gen u$ in canonical form. This gives the realization:
\[
\mathsf{A}_{3.8}=\langle \gen x, \gen u, [qx+u]\gen x + [qu-x]\gen u \rangle.
\]

\medskip\noindent $\underline{{\rm rank}\,\mathsf{A}_{3.9}=3}$. In this case we must have ${\rm rank}\,\langle e_1, e_2\rangle=2$.

\medskip\noindent $\underline{\langle e_1, e_2\rangle=\langle \gen t, \gen x\rangle}:$ In this case, as above, we have no realization.

\medskip\noindent $\underline{\langle e_1, e_2\rangle=\langle \gen x, \gen t\rangle}:$ As above, we have no realization in this case.

\medskip\noindent $\underline{\langle e_1, e_2\rangle=\langle \gen x, \gen u\rangle}:$ In this case, we have $e_3=a(t)\gen t + b\gen x + c\gen u$ with $a(t)\neq 0$ since we must have $e_1\wedge e_2\wedge e_3\neq 0$. Then $[e_1, e_3]=qe_1-e_2$ gives $b_x=q,\; c_x=-1$ and then $[e_2, e_3]=e_1+qe_2$ gives $b_u=1,\; c_u=q$ so we have $e_3=a(t)\gen t + [qx+u+b(t)]\gen x + [qu-x+c(t)]\gen u$.

The residual equivalence group is $\mathscr{E}: t'=T(t),\; x'=x+Y(t),\; u'=u+U(t)$. Under such a transformation, $e_3$ is mapped to
\begin{multline*}
e'_3=a(t)\dot{T}(t)\gen {t'} + [qx'+u'+b(t)-qY(t)-U(t)+a(t)\dot{Y}(t)]\gen {x'}\\ + [qu'-x+c(t)+Y(t)-qU(t)+a(t)\dot{U}(t)]\gen {u'}.
\end{multline*}

Choose $Y(t),\, U(t)$ so that $b(t)-qY(t)-U(t)+a(t)\dot{Y}(t)=0,\; c(t)+Y(t)-qU(t)+a(t)\dot{U}(t)=0$ and take $T(t)$ so that $a(t)\dot{T}(t)=T(t)$ and we then have $e'_3=t'\gen {t'} +[qx'+u']\gen {x'}+ [qu'-x']\gen {u'}$, so that $e_3=t\gen t + [qx+u]\gen x + [qu-x]\gen u$ in canonical form. This gives the realization:
\[
\mathsf{A}_{3.9}=\langle \gen x, \gen u, t\gen t +[qx+u]\gen x + [qu-x]\gen u \rangle.
\]

\subsection{The inequivalent non-linearizing realizations of $\mathsf{A}_{3.9}$:}

\begin{align*}
\mathsf{A}_{3.9}&=\langle \gen u, x\gen u, -(1+x^2)\gen x + (q-x)u\gen u\rangle\\
\mathsf{A}_{3.9}&=\langle \gen u, x\gen u, -t\gen t - (1+x^2)\gen x + (q-x)u\gen u\rangle\\
\mathsf{A}_{3.9}&=\langle \gen x, \gen u, [qx+u]\gen x + [qu-x]\gen u \rangle\\
\mathsf{A}_{3.9}&=\langle \gen x, \gen u, t\gen t + [qx+u]\gen x + [qu-x]\gen u \rangle.
\end{align*}

\section{Realizations of four-dimensional solvable algebras.} First, we have the abelian algebra $\mathsf{A}=\langle e_1, e_2, e_3, e_4\rangle$. There are two realizations:
\begin{align*}
\mathsf{A} &=\langle \gen t, \gen u, x\gen u, c(x)\gen u \rangle,\; c''(x)\neq 0\\
\mathsf{A} &=\langle \gen u, x\gen u, c(t,x)\gen u, q(t,x)\gen u\rangle\\
& c_{xx}\neq 0,\, q_{xx}\neq 0.
\end{align*}
Both of these realizations lead to linearizable evolution equations.

\subsection{Decomposable four-dimensional solvable Lie algebras:} There are ten such types, including the abelian Lie algebra. The other nine are:
\[
\mathsf{A}_{2.1}\oplus\mathsf{A}_{2.2},\quad \mathsf{A}_{2.2}\oplus\mathsf{A}_{2.2},\quad \mathsf{A}_{3.k}\oplus\mathsf{A}_1\, (k=3,\dots, 9).
\]
Note that $\mathsf{A}_{2.1}\oplus\mathsf{A}_{2.2}$ may be considered as $\langle e_1, e_2, e_3\rangle \uplus\langle e_4\rangle$ with $\langle e_1, e_2, e_3\rangle$ as a three-dimensional abelian Lie algebra, and $e_4$ obeying the commutation relations $[e_1, e_4]=0,\; [e_2, e_4]=0,\; [e_3, e_4]=e_3$. This is how we shall classify the non-linearizing realizations. However, we may also consider $\mathsf{A}_{2.1}\oplus\mathsf{A}_{2.2}$ as $\mathsf{A}_{3.2}\oplus\mathsf{A}_1$ with $\langle e_1, e_2, e_3\rangle=\mathsf{A}_{3.2}$ and $\langle e_4\rangle=\mathsf{A}_1$. It turns out to be simpler to consider this algebra as $\langle e_1, e_2, e_3\rangle \uplus\langle e_4\rangle$ with $\langle e_1, e_2, e_3\rangle$ as a three-dimensional abelian Lie algebra, and $e_4$ obeying the commutation relations $[e_1, e_4]=0,\; [e_2, e_4]=0,\; [e_3, e_4]=e_3$.

\medskip\noindent $\underline{\mathsf{A}_{2.1}\oplus\mathsf{A}_{2.2}=\langle e_1, e_2, e_3\rangle \uplus\langle e_4\rangle}$. Here we have that $\langle e_1, e_2, e_3\rangle$ is an abelian Lie algebra, and thus we need only consider the rank-two and rank-three cases since rank-one realizations of three-dimensional abelian Lie algebras give us linearizable evolution equations.

\medskip\noindent $\underline{{\rm rank}\,\langle e_1, e_2, e_3\rangle=2}:$ In this case the canonical realization of $\langle e_1, e_2, e_3\rangle$ is $\langle e_1, e_2, e_3\rangle=\langle \gen t, \gen u, x\gen u\rangle$. Note that the order of this basis is important when we adjoin $\langle e_4\rangle$ with $[e_3, e_4]=e_3$. In fact, there is no equivalence transformation (of our type) which transforms $\langle \gen t, \gen u, x\gen u\rangle\uplus\langle e_4\rangle$ to $\langle \gen u, x\gen u, \gen t\rangle\uplus\langle e_4\rangle$, whereas $\langle \gen t, \gen u, x\gen u\rangle\uplus\langle e_4\rangle$ and $\langle \gen t, x\gen u, \gen u\rangle\uplus\langle e_4\rangle$ are equivalent under the equivalence transformation $t'=t,\; x'=1/x,\; u'=u/x$. Note also that $\langle \gen t, \gen u, x\gen u\rangle\uplus\langle e_4\rangle$ and $\langle \gen u, \gen t, x\gen u\rangle\uplus\langle e_4\rangle$ are indistinguishable because $[e_1, e_4]=0,\; [e_2, e_4]=0$. So we have the following two cases to consider: $\langle \gen t, \gen u, x\gen u\rangle\uplus\langle e_4\rangle$ and $\langle \gen u, x\gen u, \gen t\rangle\uplus\langle e_4\rangle$.

\medskip\noindent $\underline{\langle \gen t, x\gen u, \gen u\rangle\uplus\langle e_4\rangle}:$ Putting $e_4=a(t)\gen t + b\gen x + c\gen u$ we find from $[e_1, e_4]=0,\; [e_2, e_4]=0$ that $\dot{a}(t)=0,\; b_t=b_u=0,\; c_t=0,\, xc_u=b$. Then $[e_3, e_4]=e_3$ gives us $c_u=1$ and so $e_4=x\gen x + [u+c(x)]\gen u$.

\smallskip\noindent The residual equivalence group is $\mathscr{E}(e_1, e_2, e_3): t'=t+k,\; x'=x,\; u'=u+U(x)$. Under such a transformation $e_4$ is transformed to $e'_4=x'\gen {x'} + [u' + c(x) - U(x)]\gen {u'}$ and we choose $U(x)$ so that  $U(x)=c(x)$ and we find $e'_4=x'\gen {x'} + u'\gen {u'}$. So $e_4=x\gen x + u\gen u$ in canonical form and we have the realization
\[
\mathsf{A}_{2.1}\oplus\mathsf{A}_{2.2}=\langle\gen t, x\gen u, \gen u, x\gen x + u\gen u\rangle.
\]

\medskip\noindent $\underline{\langle \gen u, x\gen u, \gen t\rangle\uplus\langle e_4\rangle}:$ In this case, $[e_1, e_4]=0,\; [e_2, e_4]=0$ gives us $b=0,\; c_u=0$ so that $e_4=a(t)\gen t + c(t,x)\gen u$. Then $[e_3, e_4]=e_3$ gives $\dot{a}(t)=0,\; c_t=0$ and we than have $e_4=[t+l]\gen t + c(x)\gen u$. Using the residual equivalence group $\mathscr{E}(e_1, e_2, e_3): t'=t+k,\; x'=x,\; u'=u+U(x)$, then $e_4$ is transformed to $e'_4=[t'+l-k]\gen {t'} + c(x')\gen {u'}$ which gives $e_4=t\gen t + c(x)\gen u$ in canonical form on taking $k=l$. Note that we cannot have $c(x)=0$ since then $\gen t, \, t\gen t$ would be symmetries and then $F=0$, which is a contradiction. Also we must have $c''(x)\neq 0$ since otherwise the algebra would again have as basis $\gen t, t\gen t, \gen u, x\gen u$ and so we would have $F=0$. Thus we have $e_4=t\gen t + c(x)\gen u$ with $c''(x)\neq 0$ and we then obtain the realization:
\[
\mathsf{A}_{2.1}\oplus\mathsf{A}_{2.2}=\langle\gen u, x\gen u, \gen t, t\gen t + c(x)\gen u\rangle,\;\; c''(x)\neq 0.
\]

\medskip\noindent $\underline{{\rm rank}\,\langle e_1, e_2, e_3\rangle=3}:$ In this case the canonical realization of $\langle e_1, e_2, e_3\rangle$ is $\langle e_1, e_2, e_3\rangle=\langle \gen t, \gen x, \gen u\rangle$. There are two orderings to consider: $\langle e_1, e_2, e_3\rangle=\langle \gen t, \gen x, \gen u\rangle$ and $\langle e_1, e_2, e_3\rangle=\langle \gen x, \gen u, \gen t\rangle$, all other orderings being equivalent under the equivalence transformation $t'=t,\; x'=u,\; u'=x$.

\medskip\noindent $\underline{\langle \gen t, \gen x, \gen u\rangle\uplus\langle e_4\rangle}:$ First, $[e_1, e_4]=0,\; [e_2, e_4]=0$ gives us $e_4=a\gen t + b(u)\gen x + c(u)\gen u$ and we may assume $a=0$ without loss of generality, so that $e_4=b(u)\gen x + c(u)\gen u$. Then $[e_3, e_4]=e_3$ gives $b'(u)=0,\; c'(u)=u$ so that we may assume $b=0$ and so we find $e_4=[u+l]\gen u$. The residual equivalence group is $\mathscr{E}(e_1, e_2, e_3): t'=t+k, x'=x+p,\; u'=u+q$. Under such a transformation $e_4$ is mapped to $e'_4=[u'+l-q]\gen {u'}$, so that $e'_4=u'\gen {u'}$ on choosing $q=l$. Thus $e_4=u\gen u$ in canonical form and then we find the realization
\[
\mathsf{A}_{2.1}\oplus\mathsf{A}_{2.2}=\langle\gen t, \gen x, \gen u, u\gen u\rangle.
\]

\medskip\noindent $\underline{\langle \gen x, \gen u, \gen t\rangle\uplus\langle e_4\rangle}:$ First, $[e_1, e_4]=0,\; [e_2, e_4]=0$ gives us $e_4=a(t)\gen t + b(t)\gen x + c(t)\gen u $. Then $[e_3, e_4]=e_3$ gives $\dot{a}(t)=1,\;\dot{b}(t)=0,\; \dot{c}(t)=0$ so that we may assume $b=c=0$ and so we find $e_4=[t+l]\gen t$. However, this gives us $t\gen t$ as a symmetry and this is incompatible with $\gen t$ since these two together give $F=0$. Hence we have no realization in this case.

\subsection{The inequivalent non-linearizing realizations of $\mathsf{A}_{2.1}\oplus\mathsf{A}_{2.2}$:}

\begin{align*}
\mathsf{A}_{2.1}\oplus\mathsf{A}_{2.2}&=\langle\gen t, x\gen u, \gen u, x\gen x + u\gen u\rangle\\
\mathsf{A}_{2.1}\oplus\mathsf{A}_{2.2}&=\langle\gen u, x\gen u, \gen t, t\gen t + c(x)\gen u\rangle,\;\; c''(x)\neq 0\\
\mathsf{A}_{2.1}\oplus\mathsf{A}_{2.2}&=\langle\gen t, \gen x, \gen u, u\gen u\rangle.
\end{align*}

\section{Realizations of $\mathsf{A}_{2.2}\oplus\mathsf{A}_{2.2}$:} We may write
\[
\mathsf{A}_{2.2}\oplus\mathsf{A}_{2.2}=\langle e_1, e_2, e_3, e_4\rangle
\]
with $[e_1, e_2]=e_1,\; [e_3, e_4]=e_3,\; [e_1, e_3]=[e_1, e_4]=0,\; [e_2, e_3]=[e_2, e_4]=0$. We also note that $\langle e_1, e_2, e_3\rangle=\mathsf{A}_{3.2}$, so we look for one-dimensional extensions of $\mathsf{A}_{3.2}$ with $[e_3, e_4]=e_3$.

\medskip\noindent $\underline{\langle e_1, e_2, e_3\rangle=\langle \gen x, x\gen x, \gen t\rangle}:$ From the commutation relations $e_4=[t+l]\gen t + c(u)\gen u$ and we note that $c(u)\neq 0$ for otherwise we have no implementation of $\langle e_3, e_4\rangle$ as $\mathsf{A}_{2.2}$ with $[e_3, e_4]=e_3$. The residual symmetry group is $\mathscr{E}(e_1, e_2, e_3): t'=t+k,\, x'=x,\; u'=U(u)$ with $U'(u)\neq 0$. Under such a transformation $e_4$ is mapped to $e'_4=(t'+l-k)\gen {t'} + c(u)U'(u)\gen {u'}$ and we choose $k=l$ and $U(u)$ so that $c(u)U'(u)=U(u)$ giving $e_4=t\gen t + u\gen u$ in canonical form. We have the realization
\[
\mathsf{A}_{2.2}\oplus\mathsf{A}_{2.2}=\langle \gen x, x\gen x, \gen t, t\gen t + u\gen u\rangle.
\]

\medskip\noindent $\underline{\langle e_1, e_2, e_3\rangle=\langle \gen x, x\gen x, \gen u\rangle}:$ From the commutation relations $e_4=a(t)\gen t + [u + c(t)]\gen u$. The residual symmetry group is $\mathscr{E}(e_1, e_2, e_3): t'=T(t),\, x'=x,\; u'=u+ U(t)$. Under such a transformation $e_4$ is mapped to $e'_4=a(t)\dot{T}(t)\gen {t'} + [u'+a(t)\dot{c}(t) + c(t) - U(t)]\gen {u'}\gen {u'}$ and we may always choose $U(t)$ so that $a(t)\dot{c}(t)+ c(t) - U(t)=0$ giving $e'_4=a(t)\dot{T}(t)\gen {t'} + u'\gen {u'}$. If $a(t)=0$ we find $e_4=u\gen u$ in canonical form. If $a(t)\neq 0$ we choose $T(t)$ so that $a(t)\dot{T}(t)=T(t)$ giving $e_4=t\gen t + u\gen u$ in canonical form. We have the realizations
\begin{align*}
\mathsf{A}_{2.2}\oplus\mathsf{A}_{2.2}&=\langle \gen x, x\gen x, \gen u, u\gen u\rangle\\
\mathsf{A}_{2.2}\oplus\mathsf{A}_{2.2}&=\langle \gen x, x\gen x, \gen u, t\gen t + u\gen u\rangle.
\end{align*}

\medskip\noindent $\underline{\mathsf{A}_{3.2}=\langle \gen t, t\gen t + x\gen x, xu\gen x\rangle}:$ From the commutation relations $e_4=b(u)x\gen x-u\gen u$. The residual equivalence group is $\mathscr{E}(e_1, e_2, e_3): t'=t,\; x'=p(u)x,\, u'=u$ with $p(u)\neq 0$. Under such a transformation $e-4$ is mapped to $e'_4=[b(u)p(u)-up'(u)]x\gen {x'} -u'\gen {u'}$ and we may always choose $p(u)$ so that $b(u)p(u)-up'(u)=0$ so that $e_4=-u\gen u$ in canonical form. We have the realization
\[
\mathsf{A}_{2.2}\oplus\mathsf{A}_{2.2}=\langle \gen t, t\gen t + x\gen x, xu\gen x, -u\gen u\rangle.
\]

\medskip\noindent $\underline{\mathsf{A}_{3.2}=\langle \gen x, t\gen t + x\gen x, tu\gen x\rangle}:$ From the commutation relations we find that $e_4=at\gen t + tb(u)\gen x -(a+1)u\gen u$. The residual equivalence group is $\mathscr{E}(e_1, e_2, e_3): t'=k t,\; x'=x+p(u)t,\; u'=u/k$ with  $k\neq 0$, and under such a transformation $e_4$ is mapped to $e'_4=at'\gen {t'} + t[b(u)+ap(u)-(a+1)up'(u)]\gen {x'} - (a+1)u'\gen {u'}$. We may obviously always choose $p(u)$ so that $b(u)+ap(u)-(a+1)up'(u)=0$ and so we obtain $e'_4=at'\gen {t'}-(a+1)u'\gen {u'}$.  So we have the realization:
\[
\mathsf{A}_{2.2}\oplus\mathsf{A}_{2.2}=\langle \gen x, t\gen t + x\gen x, tu\gen x, at\gen t-(1+a)u\gen u\rangle,\quad a\in \mathbb{R}.
\]
We note that $a=-1$ gives us the algebra
\[
\mathsf{A}_{2.2}\oplus\mathsf{A}_{2.2}=\langle \gen x, t\gen t + x\gen x, tu\gen x, -t\gen t\rangle,
\]
which is equivalent to the algebra
\[
\mathsf{A}_{2.2}\oplus\mathsf{A}_{2.2}=\langle \gen u, t\gen t + x\gen x + u\gen u, x\gen u, -t\gen t-x\gen x\rangle
\]
under the equivalence transformation $t'=t,\, x'=tu,\; u'=x$. This algebra contains the subalgebra $\langle \gen u, x\gen u, u\gen u\rangle$ which linearizes the equation, giving $u_t=F(t,x)u_3+G(t,x)u_2$ as the evolution equation. Hence we have the non-linearizing realization
\[
\mathsf{A}_{2.2}\oplus\mathsf{A}_{2.2}=\langle \gen x, t\gen t + x\gen x, tu\gen x, at\gen t-(1+a)u\gen u\rangle,\quad a\neq -1.
\]

\medskip\noindent $\underline{\mathsf{A}_{3.2}=\langle \gen x, t\gen t + x\gen x, t\gen t\rangle}:$ From $[e_1, e_4]=[e_2, e_4]=0$ we have $e_4=at\gen t + tb(u)\gen x + c(u)\gen u$. Then $[e_3, e_4]=e_3$ is impossible to realize, so we have no realization in this case.

\medskip\noindent $\underline{\mathsf{A}_{3.2}=\langle \gen x, x\gen x + u\gen u, u\gen u\rangle}:$ Since we have a combination of two rank-two realizations of $\mathsf{A}_{2.2}$ then we have $e_4=a(t)\gen t + b\gen x + c\gen u$ with $a(t)^2+b^2\neq 0$. The commutation relations $[e_1, e_4]=[e_2, e_4]=0$ give us $e_4=a(t)\gen t + ub(t)\gen x + uc(t)\gen u$. Then $[e_3, e_4]=e_3$ gives $u\gen u=ub(t)\gen x$ which is impossible, so we have no realization in this case.

\medskip\noindent $\underline{\mathsf{A}_{3.2}=\langle \gen x, x\gen x + u\gen u, tu\gen u\rangle}:$ Following the calculations of the previous case with $e_3
=tu\gen u$, the relation $[e_3, e_4]=e_3$ gives us the contradiction $tu\gen u=ub(t)\gen x$, so we have no realization in this case.

\medskip\noindent $\underline{\mathsf{A}_{3.2}=\langle \gen x, x\gen x + u\gen u, u\gen x\rangle}:$ The commutation relations give $e_4=a(t)\gen t + ub(t)\gen x - u\gen u$. The residual equivalence group is $\mathscr{E}(e_1, e_2, e_3): t'=T(t),\; x'=x+p(t)u,\; u'=u$. Under such a transformation, $e_4$ is mapped to $e'_4=a(t)\dot{T}(t)\gen {t'} +[b(t)+a(t)\dot{p}(t)- p(t)]u\gen {x'} - u'\gen {u'}$ and we see that we may always choose $p(t)$ so that $b(t)+a(t)\dot{p}(t)- p(t)=0$, giving us $e'_4=a(t)\dot{T}(t)\gen {t'} - u'\gen {u'}$ and we have the canonical form $e_4=-u\gen u$ if $a(t)=0$ and, if $a(t)\neq 0$, we choose $T(t)$ so that $a(t)\dot{T}(t)=T(t)$ giving the canonical form $e_4=t\gen t-u\gen u$. We note that the realization $\mathsf{A}_{2.2}\oplus\mathsf{A}_{2.2}=\langle \gen x, x\gen x + u\gen u, u\gen x, -u\gen u\rangle$ contains the subalgebra $\langle \gen x, u\gen x, x\gen x\rangle$ which linearizes the evolution equation: in fact, $\langle \gen x, u\gen x, x\gen x\rangle$ is equivalent, under the equivalence transformation $t'=t,\; x'=u,\; u'=x$, to the algebra $\langle \gen u, x\gen u, u\gen u\rangle$ which gives the linear evolution equation $u_t=F(t,x)u_3 + G(t,x)u_2$. Thus we have the non-linearizing realization
\[
\mathsf{A}_{2.2}\oplus\mathsf{A}_{2.2}=\langle \gen x, x\gen x + u\gen u, u\gen x, t\gen t-u\gen u\rangle.
\]

\medskip\noindent $\underline{\mathsf{A}_{3.2}=\langle \gen t, t\gen t + x\gen x, \gen u\rangle}:$ The commutation relations $[e_1, e_4]=[e_2, e_4]=0,\; [e_3, e_4]=e_3$ give $e_4=qx\gen x + [u+l]\gen u$, with $q\in \mathbb{R}$. The residual equivalence group is $\mathscr{E}(e_1, e_2, e_3): t'=t,\; x'=px,\; u'=u+k$ with $p\neq 0$. Under such a transformation, $e_4$ is mapped to $e'_4=qx'\gen {x'} + [u'+l-k]\gen {u'}$ and we may take $k=l$ giving the canonical form $e_4=qx\gen x + u\gen u$. We have the realization
\[
\mathsf{A}_{2.2}\oplus\mathsf{A}_{2.2}=\langle \gen t, t\gen t + x\gen x, \gen u, qx\gen x + u\gen u\rangle,\quad q\in \mathbb{R}.
\]

\medskip\noindent $\underline{\mathsf{A}_{3.2}=\langle \gen x, t\gen t + x\gen x, \gen u\rangle}:$ The commutation relations give $e_4=at\gen t +qt\gen x + [u+l]\gen u$. The residual equivalence group is $\mathscr{E}(e_1, e_2, e_3): t'=\alpha t,\; x'=x + pt,\; u'=u+k$ with $\alpha\neq 0$. Under such a transformation, $e_4$ is mapped to $e'_4=at'\gen {t'} + t[q+ap]\gen {x'} + [u'+l-k]\gen {u'}$. We may take $k=l$; if $a\neq 0$ then we choose $p$ so that $q+ap=0$. Thus we have the canonical forms $e_4=qt\gen x + u\gen u$ if $a=0$ and $at\gen t + u\gen u$ if $a\neq 0$. We then have the realizations
\begin{align*}
\mathsf{A}_{2.2}\oplus\mathsf{A}_{2.2}&=\langle \gen x, t\gen t + x\gen x, \gen u, qt\gen x+u\gen u\rangle,\quad q\in \mathbb{R}\\
\mathsf{A}_{2.2}\oplus\mathsf{A}_{2.2}&=\langle \gen x, t\gen t + x\gen x, \gen u, at\gen t+u\gen u\rangle,\quad a\neq 0.
\end{align*}

\medskip\noindent $\underline{\mathsf{A}_{3.2}=\langle \gen x, t\gen t + x\gen x, t\gen t + u\gen u\rangle}:$ The commutation relations $[e_1, e_4]=[e_2, e_4]=0$ give us $e_4=at\gen t + tb(u)\gen x + c(u)\gen u$ and then we see that $[e_3, e_4]=e_3$ is impossible to realize, so we have no realization in this case.

\medskip\noindent $\underline{\mathsf{A}_{3.2}=\langle \gen x, x\gen x + u\gen u, \gen t\rangle}:$ The commutation relations $[e_1, e_4]=[e_2, e_4]=0,\; [e_3, e_4]=e_3$ give $e_4=[t+l]\gen t +pu\gen x + qu\gen u$. Note that $p^2+q^2\neq 0$ for other wise $\gen t, [t+l]\gen t$ would be symmetries, and these are incompatible since they give $F=0$ which is a contradiction. The residual equivalence group is $\mathscr{E}(e_1, e_2, e_3): t'=t+k,\; x'=x + \kappa u,\; u'=\alpha u$ with $\alpha\neq 0$. Under such a transformation, $e_4$ is mapped to $e'_4=[t'+l-k]\gen {t'} + u[p+q\kappa]\gen {x'} + qu'\gen {u'}$. We may take $k=l$; if $q\neq 0$ then choose $\kappa$ so that $p+q\kappa=0$; if $q=0$ then we have $p\neq 0$ and we choose $\alpha=p$. Thus we have the canonical forms $e_4=t\gen t + u\gen x$ if $q=0$ and $e_4=t\gen t + qu\gen u$ if $q\neq 0$. We then have the realizations
\begin{align*}
\mathsf{A}_{2.2}\oplus\mathsf{A}_{2.2}&=\langle \gen x, x\gen x + u\gen u, \gen t, t\gen t+qu\gen u\rangle,\quad q\neq 0\\
\mathsf{A}_{2.2}\oplus\mathsf{A}_{2.2}&=\langle \gen x, x\gen x + u\gen u, \gen t, t\gen t+u\gen x\rangle.
\end{align*}

\subsection{Inequivalent non-linearizing realizations of $\mathsf{A}_{2.2}\oplus\mathsf{A}_{2.2}$.}

\begin{align*}
\mathsf{A}_{2.2}\oplus\mathsf{A}_{2.2}&=\langle \gen x, x\gen x, \gen t, t\gen t + u\gen u\rangle\\
\mathsf{A}_{2.2}\oplus\mathsf{A}_{2.2}&=\langle \gen x, x\gen x, \gen u, u\gen u\rangle\\
\mathsf{A}_{2.2}\oplus\mathsf{A}_{2.2}&=\langle \gen x, x\gen x, \gen u, t\gen t + u\gen u\rangle\\
\mathsf{A}_{2.2}\oplus\mathsf{A}_{2.2}&=\langle \gen t, t\gen t + x\gen x, xu\gen x, -u\gen u\rangle\\
\mathsf{A}_{2.2}\oplus\mathsf{A}_{2.2}&=\langle \gen x, t\gen t + x\gen x, tu\gen x, at\gen t-(1+a)u\gen u\rangle,\quad a\neq -1\\
\mathsf{A}_{2.2}\oplus\mathsf{A}_{2.2}&=\langle \gen x, x\gen x + u\gen u, u\gen x, t\gen t-u\gen u\rangle\\
\mathsf{A}_{2.2}\oplus\mathsf{A}_{2.2}&=\langle \gen t, t\gen t + x\gen x, \gen u, qx\gen x + u\gen u\rangle,\quad q\in \mathbb{R}\\
\mathsf{A}_{2.2}\oplus\mathsf{A}_{2.2}&=\langle \gen x, t\gen t + x\gen x, \gen u, qt\gen x+u\gen u\rangle,\quad q\in \mathbb{R}\\
\mathsf{A}_{2.2}\oplus\mathsf{A}_{2.2}&=\langle \gen x, t\gen t + x\gen x, \gen u, at\gen t+u\gen u\rangle,\quad a\neq 0\\
\mathsf{A}_{2.2}\oplus\mathsf{A}_{2.2}&=\langle \gen x, x\gen x + u\gen u, \gen t, t\gen t+qu\gen u\rangle,\quad q\neq 0\\
\mathsf{A}_{2.2}\oplus\mathsf{A}_{2.2}&=\langle \gen x, x\gen x + u\gen u, \gen t, t\gen t+u\gen x\rangle.
\end{align*}

\section{Realizations of $\mathsf{A}_{3.3}\oplus\mathsf{A_1}$.}

\medskip\noindent $\underline{\mathsf{A}_{3.3}=\langle \gen u, x\gen u, -(\gen t + \gen x) \rangle}:$ With $e_4=a(t)\gen t + b\gen x + c\gen u$ we have $b=0,\; c_u=0$ from $[e_1, e_4]=[e_2, e_4]=0$. Then $[e_3, e_4]=0$ gives $\dot{a}(t)=0,\, c_t+c_x=0$. Thus we have $e_4=a\gen t + c(t-x)\gen u$. We note that we must have $a\neq 0$ since otherwise $e_1, e_2, e_4$ would be a rank-one three-dimensional abelian Lie algebra, which we know linearizes the evolution equation. So, dividing by $a\neq 0$ we have $e_4=\gen t + c(t-x)\gen u$. Note that the residual equivalence group is $\mathscr{E}(e_1, e_2, e_3): t'=t+k,\; x'=x+Y(t-x),\; u'=u + U(t-x)$. Under such a transformation, $e_4$ is mapped to $e'_4=\gen {t'} + [U_t(t-x)+c(t-x)]\gen {u'}$. Then we choose $U(t-x)$ so that $U_t(t-x)+c(t-x)=0$ and so we have $e'_4=\gen {t'}$, giving $e_4=\gen t$ in canonical form. We have the realization
\[
\mathsf{A}_{3.3}\oplus\mathsf{A}_1=\langle \gen u, x\gen u, -(\gen t + \gen x), \gen t \rangle.
\]

\medskip\noindent $\underline{\mathsf{A}_{3.3}=\langle \gen x, \gen t, (t+u)\gen x \rangle}:$ The commutation relations give $e_4=a(\gen t - \gen u)+b(u)\gen x$. Also, $a\neq 0$ for otherwise $e_1, e_2, e_4$ would be a rank-one three-dimensional abelian Lie algebra which linearizes the evolution equation. So we take $a=1$ without loss of generality (on dividing by $a\neq 0$) and then $e_4=\gen t -\gen u + b(u)\gen x$. The residual equivalence group is $\mathscr{E}(e_1, e_2, e_3): t'=t+k,\; x'=x+Y(u),\; u'=u-k$. Under such a transformation, $e_4$ is mapped to $e'_4=\gen {t'} - \gen {u'} + [b(u) - Y'(u)]\gen {x'}$. Then choose $Y(u)$ so that $Y(u)=b(u)$ so that $e'_4=\gen {t'}-\gen {u'}$. So $e_4=\gen t - \gen u$ in canonical form. We have the realization
\[
\mathsf{A}_{3.3}\oplus\mathsf{A}_1=\langle \gen x, \gen t, (t+u)\gen x, \gen t-\gen u \rangle.
\]

\medskip\noindent $\underline{\mathsf{A}_{3.3}=\langle \gen x, \gen u, u\gen x \rangle}:$ The commutation relations give us $e_4=a(t)\gen t + b(t)\gen x$. First note that $a(t)\neq 0$ for otherwise we would have $e_4=b(t)\gen x$. For $\dim \mathsf{A}_{3.3}\oplus\mathsf{A}_1=4$ we would need $\dot{b}(t)\neq 0$, but then $\gen x$ and $b(t)\gen x$ would be symmetries giving $\dot{b}(t)u_1=0$ in the equation for $G$, and this is a contradiction for $\dot{b}(t)\neq 0$.

The residual equivalence group is $\mathscr{E}(e_1, e_2, e_3): t'=T(t),\; x'=x+Y(t),\; u'=u$. Under such a transformation, $e_4$ is mapped to $e'_4=a(t)\dot{T}(t)\gen {t'} + [b(t)+a(t)\dot{Y}(t)]\gen {x'}$. Because $a(t)\neq 0$, we may choose $Y(t)$ so that $b(t)+a(t)\dot{Y}(t)=0$, and we choose $T(t)$ so that $a(t)\dot{T}(t)=1$, giving $e'_4=\gen {t'}$. So $e_4=\gen t$ in canonical form and we obtain the realization
\[
\mathsf{A}_{3.3}\oplus\mathsf{A}_1=\langle \gen x, \gen u, u\gen x, \gen t\rangle.
\]
Note that this algebra gives the same type of equation as
\[
\mathsf{A}_{3.3}\oplus\mathsf{A}_1=\langle \gen u, x\gen u, -(\gen t + \gen x), \gen t \rangle.
\]

\medskip\noindent $\underline{\mathsf{A}_{3.3}=\langle \gen u, x\gen u, -\gen x \rangle}:$ The commutation relations give us $e_4=a(t)\gen t + c(t)\gen u$. First note that $a(t)\neq 0$ for otherwise we would have $e_4=c(t)\gen u$. For $\dim \mathsf{A}_{3.3}\oplus\mathsf{A}_1=4$ we would need $\dot{c}(t)\neq 0$, but then $\gen u$ and $c(t)\gen u$ would be symmetries giving $\dot{c}(t)=0$ in the equation for $G$, and this is a contradiction for $\dot{c}(t)\neq 0$. The residual equivalence group is $\mathscr{E}(e_1, e_2, e_3): t'=T(t),\; x'=x,\; u'=u+U(t)$. Under such a transformation $e_4$ is mapped to $e'_4=a(t)\dot{T}(t)\gen {t'} + [c(t)+a(t)\dot{U}(t)]\gen {u'}$, and we may choose $U(t)$ so that $c(t)+a(t)\dot{U}(t)=0$ and we choose $T(t)$ so that $a(t)\dot{T}(t)=1$, giving $e_4=\gen t$ in canonical form. Thus we have the realization
\[
\mathsf{A}_{3.3}\oplus\mathsf{A}_1=\langle \gen u, x\gen u, -\gen x, \gen t \rangle.
\]

\medskip\noindent $\underline{\mathsf{A}_{3.3}=\langle \gen x, \gen u, u\gen x+t\gen u \rangle}:$ The commutation relations give us $e_4=b(t)\gen x$. Again, we cannot have $\dot{b}(t)\neq 0$ for then $\gen x$ together with $b(t)\gen x$ would give $\dot{b}(t)u_1=0$ in the equation for $G$. Thus $b(t)=\,{\rm constant}$, and then $\dim\mathsf{A}_{3.3}\oplus\mathsf{A}_1=3$, which is also a contradiction. Hence we have no realization in this case.

\medskip\noindent $\underline{\mathsf{A}_{3.3}=\langle \gen x, \gen t, t\gen x - \gen u \rangle}:$ The commutation relations give us $e_4=a\gen t + [p-au]\gen x + q\gen u$ with $p, q$ being constants. We must have $a^2+q^2\neq 0$ for if $a=q=0$ then $e_4=p\gen x$ and then we do not have $\dim\mathsf{A}_{3.3}\oplus\mathsf{A}_1=4$. The residual equivalence group is $\mathscr{E}(e_1, e_2, e_3): t'=t+k,\; x'=x - ku+l,\; u'=u+m$. Under such a transformation, $e_4$ is mapped to $e'_4=a\gen {t'} + [p+am-qk-au']\gen {x'} + q\gen {u'}$, and we may always arrange for $p+am-qk=0$ since $a^2+q^2\neq 0$. So we have $e_4=a(\gen t -u\gen x) + q\gen u$ in canonical form. If $a=0$ then $e_4=q\gen u$ and then $\dim\mathsf{A}_{3.3}\oplus\mathsf{A}_1$ is generated by the basis elements $\gen t, \gen x, \gen u, t\gen x$, and the two operators $\gen x, t\gen x$ gives us the $F=0$ which is a contradiction. Hence $a\neq 0$ and we may without loss of generality put $a=1$, giving $e_4=\gen t - u\gen x + q\gen u$. Putting these into the equation for $F$ we obtain $F=F(\tau)$ with $\tau=u_1^{-3}u_2$. Then the operators $\gen x, \gen x, t\gen x-\gen u$ in the equation for $G$ give $G_t=G_x=0$ and $G_u=u_1$ so that we have $G=uu_1+g(u_1, \tau)$ on changing variables $(t,x,u,u_1, u_2)\to (t,x,u,u_1, \tau)$. The operator $e_4=\gen t - u\gen x + q\gen u$ in the equation for $G$ gives
\[
qG_u+u_1^2G_{u_1}+3u_1u_2G_{u_2}=u_1G-3u_2^2F(\tau).
\]
Then putting $G=uu_1 + g(u_1, \tau)$ we find that
\[
u_1(q-g)+u_1^2g_{u_1}+3u_1^6\tau^2F(\tau)=0,
\]
which gives $\displaystyle G=uu_1 - \frac{3}{4}u_1^5\tau^2F(\tau) + q$.  Thus we have the realization
\[
\mathsf{A}_{3.3}\oplus\mathsf{A}_1=\langle \gen x, \gen t, t\gen x - \gen u, \gen t -u\gen x + q\gen u \rangle
\]
where $q\in \mathbb{R}$ is arbitrary.

\medskip\noindent $\underline{\mathsf{A}_{3.3}=\langle \gen x, \gen u, t\gen t+u\gen x \rangle}:$ The commutation relations give $e_4=mt\gen t + [q\ln t + p]\gen x + q\gen u$. If $m=0$, then $q\neq 0$ in order for $\dim\mathsf{A}_{3.3}\oplus\mathsf{A}_1=4$, and then the combination of symmetries $\gen x,\; \gen u$ and $m\gen t + [q\ln t + p]\gen x + q\gen u$ give us the contradiction $u_1=t$ in the equation for $G$. Thus, $m\neq 0$ and we may assume $m=1$ without loss of generality. So we have $e_4=t\gen t + q[\ln t\gen x + \gen u]$. The residual equivalence group is $\mathscr{E}(e_1, e_2, e_3): t'=kt,\; x'=x+Y(t),\; u'=u+U(t)$ with $U(t)=t\dot{Y}(t)$. Under such a transformation, $e_4$ is mapped to $e'_4=t'\gen {t'} + [t\dot{Y}(t)+q\ln t + p]\gen {x'}\ + [q+t\dot{U}(t)]\gen {u'}$. Choose $Y(t)$ so that $t\dot{Y}(t)+q\ln t + p=0$. Thus we have $U(t)+q\ln t + p=0$ from which we find that $t\dot{U}(t)+q=0$ and thus we find $e'_4=t'\gen {t'}$. So have $e_4=t\gen t$ in canonical form and we have the realization
\[
\mathsf{A}_{3.3}\oplus\mathsf{A}_1=\langle \gen x, \gen u, t\gen t + u\gen x, t\gen t\rangle.
\]

\subsection{Inequivalent non-linearizing realizations of $\mathsf{A}_{3.3}\oplus\mathsf{A_1}$.}

\begin{align*}
\mathsf{A}_{3.3}\oplus\mathsf{A}_1&=\langle \gen u, x\gen u, -(\gen t + \gen x), \gen t \rangle\\
\mathsf{A}_{3.3}\oplus\mathsf{A}_1&=\langle \gen u, x\gen u, -\gen x, \gen t \rangle\\
\mathsf{A}_{3.3}\oplus\mathsf{A}_1&=\langle \gen x, \gen t, (t+u)\gen x, \gen t-\gen u \rangle\\
\mathsf{A}_{3.3}\oplus\mathsf{A}_1&=\langle \gen x, \gen u, u\gen x, \gen t\rangle\\
\mathsf{A}_{3.3}\oplus\mathsf{A}_1&=\langle \gen x, \gen t, t\gen x - \gen u, \gen t -u\gen x + \kappa\gen u \rangle,\quad \kappa\in \mathbb{R}\\
\mathsf{A}_{3.3}\oplus\mathsf{A}_1&=\langle \gen x, \gen u, t\gen t + u\gen x, t\gen t\rangle.
\end{align*}

\section{Realizations of $\mathsf{A}_{3.4}\oplus\mathsf{A_1}$.}

\medskip\noindent $\underline{\mathsf{A}_{3.4}=\langle \gen x, u\gen x, x\gen x - \gen u \rangle}:$ From the commutation relations we find that $\displaystyle e_4=a(t)\gen t + q(t)e^{-u}\gen x$ and we need $a(t)\neq 0$ for otherwise we would have $e_1, e_2, e_4$ as a rank-one three-dimensional abelian Lie algebra which linearizes the evolution equation. The residual equivalence group is $\mathscr{E}(e_1, e_2, e_3): t'=T(t),\; x'=x+p(t)e^{-u},\; u'=u$. Under such a transformation, $e_4$ is mapped to $e'_4=a(t)\dot{T}(t)\gen {t'} + [a(t)\dot{p}(t)+q(t)]e^{-u}\gen {x'}$ and, since $a(t)\neq 0$ we choose $p(t)$ and $T(t)$ so that $a(t)\dot{p}(t)+q(t)=0$ and $a(t)\dot{T}(t)=1$ giving $e'_4=\gen {t'}$. So $e_4=\gen t$ in canonical form and we have the realization
\[
\mathsf{A}_{3.4}\oplus\mathsf{A_1}=\langle \gen x, u\gen x, x\gen x-\gen u, \gen t\rangle.
\]

\medskip\noindent $\underline{\mathsf{A}_{3.4}=\langle \gen x, u\gen x, t\gen t + x\gen x - \gen u \rangle}:$ In this case, we find from the commutation relations that $e_4=mt\gen t + t^{-1}b(\sigma)\gen x$ with $\sigma=te^{u}$. We also must have $m\neq 0$ since otherwise $e_1, e_2, e_4$ would form a rank-one three-dimensional abelian Lie algebra which linearizes the evolution equation. Hence, we may divide throughout by $m\neq 0$ and then we have $e_4=t\gen t + t^{-1}b(\sigma)\gen x$. The residual equivalence group is $\mathscr{E}(e_1, e_2, e_3): t'=kt,\; x'=x+t^{-1}Y(\sigma),\; u'=u$ with $k\neq 0$. Under such a transformation, $e_4$ is mapped to $e'_4=t'\gen {t'} + [b(\sigma)+Y(\sigma)]t^{-1}\gen {x'}$. We may choose $Y(\sigma)$ so that $b(\sigma)+Y(\sigma)=0$ and thus we obtain $e'_4=t'\gen {'}$ so that $e_4=t\gen t$ in canonical form. Thus we have the realization
\[
\mathsf{A}_{3.4}\oplus\mathsf{A_1}=\langle \gen x, u\gen x, t\gen t + x\gen x-\gen u, t\gen t\rangle.
\]

\medskip\noindent $\underline{\mathsf{A}_{3.4}=\langle \gen x, \gen t, t\gen t + [t+x]\gen x\rangle}:$ The commutation relations give $e_4=c(u)\gen u$ with $c(u)\neq 0$, by dimension. The residual equivalence group is $\mathscr{E}(e_1, e_2, e_3): t'=t,\; x'=x,\; u'=U(u)$ with $U'(u)\neq 0$. Under such a transformation, $e_4$ is mapped to $e'_4=c(u)U'(u)\gen {u'}$ and we then choose $U(u)$ so that $c(u)U'(u)=1$, and then $e'_4=\gen {u'}$. So, $e_4=\gen u$ in canonical form and we have the realization
\[
\mathsf{A}_{3.4}\oplus\mathsf{A_1}=\langle \gen x, \gen t, t\gen t + [t+x]\gen x, \gen u\rangle.
\]

\medskip\noindent $\underline{\mathsf{A}_{3.4}=\langle \gen x, \gen u, [x+u]\gen x+u\gen u\rangle}:$ The commutation relations give $e_4=a(t)\gen t$ and $a(t)\neq 0$. The residual equivalence group is $\mathscr{E}(e_1, e_2, e_3): t'=T(t),\; x'=x,\; u'=u$. Under such a transformation, $e_4$ is mapped to $e'_4=a(t)\dot{T}(t)\gen {t'}$ and we choose $T(t)$ so that $a(t)\dot{T}(t)=1$ giving $e_4=\gen t$ in canonical form. We have the realization
\[
\mathsf{A}_{3.4}\oplus\mathsf{A_1}=\langle \gen x, \gen u, [x+u]\gen t + u\gen u, \gen t\rangle.
\]

\medskip\noindent $\underline{\mathsf{A}_{3.4}=\langle \gen x, \gen t, t\gen t + [t+x]\gen x+u\gen u\rangle}:$ The commutation relations give $e_4=pu\gen x + qu\gen u$. If $q=0$ then we have $e_4=u\gen x$ in canonical form. If $q\neq 0$ then we have the residual equivalence group $\mathscr{E}(e_1, e_2, e_3): t'=t,\; x'=x+\kappa u,\; u'=\gamma u$ with $\gamma\neq 0$. Under such a transformation, $e_4$ is mapped to $e'_4=[p+\kappa q]u\gen {x'} + qu'\gen {u'}$ and for $q\neq 0$ we may choose $\kappa$ so that $p+\kappa q=0$ so that $e'_4=qu'\gen {u'}$ and hence $e_4=u\gen u$ in canonical form. We have the realizations
\begin{align*}
\mathsf{A}_{3.4}\oplus\mathsf{A_1}&=\langle \gen x, \gen t, t\gen t + [t+x]\gen x+u\gen u, u\gen x\rangle\\
\mathsf{A}_{3.4}\oplus\mathsf{A_1}&=\langle \gen x, \gen t, t\gen t + [t+x]\gen x+u\gen u, u\gen u\rangle.
\end{align*}

\medskip\noindent $\underline{\mathsf{A}_{3.4}=\langle \gen x, \gen u, t\gen t + [x+u]\gen x+u\gen u\rangle}:$ The commutation relations give $e_4=mt\gen t$ with $m\neq 0$ so we have the realization
\[
\mathsf{A}_{3.4}\oplus\mathsf{A_1}=\langle \gen x, \gen u, t\gen t + [x+u]\gen x+u\gen u, t\gen t\rangle.
\]

\subsection{Inequivalent non-linearizing realizations of $\mathsf{A}_{3.4}\oplus\mathsf{A}_1$.}

\begin{align*}
\mathsf{A}_{3.4}\oplus\mathsf{A_1}&=\langle \gen x, u\gen x, x\gen x-\gen u, \gen t\rangle\\
\mathsf{A}_{3.4}\oplus\mathsf{A_1}&=\langle \gen x, u\gen x, t\gen t + x\gen x-\gen u, t\gen t\rangle\\
\mathsf{A}_{3.4}\oplus\mathsf{A_1}&=\langle \gen x, \gen t, t\gen t + [t+x]\gen x, \gen u\rangle\\
\mathsf{A}_{3.4}\oplus\mathsf{A_1}&=\langle \gen x, \gen u, [x+u]\gen x + u\gen u, \gen t\rangle\\
\mathsf{A}_{3.4}\oplus\mathsf{A_1}&=\langle \gen x, \gen t, t\gen t + [t+x]\gen x+u\gen u, u\gen x\rangle\\
\mathsf{A}_{3.4}\oplus\mathsf{A_1}&=\langle \gen x, \gen t, t\gen t + [t+x]\gen x+u\gen u, u\gen u\rangle\\
\mathsf{A}_{3.4}\oplus\mathsf{A_1}&=\langle \gen x, \gen u, t\gen t + [x+u]\gen x+u\gen u, t\gen t\rangle.
\end{align*}

\section{Realizations of $\mathsf{A}_{3.5}\oplus\mathsf{A_1}$.}

\medskip\noindent $\underline{\mathsf{A}_{3.5}=\langle \gen u, x\gen u, t\gen t + u\gen u \rangle}:$ The commutation relations give $e_4=at\gen t + tc(x)\gen u$. If $a=0$ then we have $e_4=tc(x)\gen u$. If $c''(x)=0$ then $e_4=t(\alpha e_1+\beta e_2)$, and this together with $e_1, e_2$ as symmetries in the equation for $G$ give $-\alpha u_1 + \beta=0$ which is a contradiction for unless $\alpha=\beta=0$. Thus we must have $c''(x)\neq 0$. However, $e_1, e_2, e_4$ will then be a three-dimensional rank-one abelian Lie algebra which linearizes the evolution equation. Hence we have $a\neq 0$ and we may assume, without loss of generality, that $a=1$, so we take $e_4=t\gen t + tc(x)\gen u$. The residual equivalence group is $\mathscr{E}(e_1, e_2, e_3): t'=kt,\, x'=x,\; u'=u+tU(x)$. Under such a transformation $e_4$ is mapped to $e'_4=t'\gen {t'}+t[c(x)+U(x)]\gen {u'}$ and we may obviously choose $U(x)$ so that $c(x)+U(x)=0$ giving $e'_4=t'\gen {t'}$. Thus $e_4=t\gen t$ in canonical form and we have the realization
\[
\mathsf{A}_{3.5}\oplus\mathsf{A}_1=\langle \gen u, x\gen u, t\gen t + u\gen u, t\gen t\rangle.
\]
However, this algebra then contains the rank-one solvable Lie algebra $\langle \gen u, x\gen u, u\gen u\rangle$ which linearizes the evolution equation (we obtain $u_t=F(t,x)u_3+G(t,x)u_2$). Hence we have no non-linearizing realization in this case.

\medskip\noindent $\underline{\mathsf{A}_{3.5}=\langle \gen t, \gen x, t\gen t + x\gen x \rangle}:$ The commutation relations give $e_4=c(u)\gen u$. The residual equivalence group is $\mathscr{E}(e_1, e_2, e_3): t'=t,\, x'=x,u'=U(u)$ and under such a transformation $e_4$ is mapped to $e'_4=c(u)U'(u)\gen {u'}$. We may obviously choose $U(u)$ so that $c(u)U'(u)=1$ giving $e'_4=\gen {u'}$. Thus $e_4=\gen u$ in canonical form and we have the realization
\[
\mathsf{A}_{3.5}\oplus\mathsf{A}_1=\langle \gen t, \gen x, t\gen t + x\gen x, \gen u \rangle.
\]

\medskip\noindent $\underline{\mathsf{A}_{3.5}=\langle \gen x, \gen u, x\gen x + u\gen u \rangle}:$ The commutation relations give $e_4=a(t)\gen t$. The residual equivalence group is $\mathscr{E}(e_1, e_2, e_3): t'=T(t),\; x'=x,\; u'=u$. Under such a transformation $e_4$ is mapped to$e'_4=a(t)\dot{T}(t)\gen {t'}$ and we choose $T(t)$ so that $a(t)\dot{T}(t)=1$ so that $e'_4=\gen {t'}$ and $e_4=\gen t$ in canonical form and we have the realization
\[
\mathsf{A}_{3.5}\oplus\mathsf{A}_1=\langle \gen x, \gen u, x\gen x + u\gen u, \gen t \rangle.
\]

\medskip\noindent $\underline{\mathsf{A}_{3.5}=\langle \gen t, \gen x, t\gen t + x\gen x + u\gen u \rangle}:$ The commutation relations give $e_4=pu\gen x + qu\gen u$. The residual equivalence group is $\mathscr{E}(e_1, e_2, e_3): t'=t,\, x'=x + \alpha u,u'=\gamma u$ with $\gamma\neq 0$, and under such a transformation $e_4$ is mapped to $e'_4=(p+\alpha q)u\gen {x'} + qu'\gen {u'}$. If $q\neq 0$ then we choose $\alpha$ so that $p+\alpha q=0$ and so we find $e'_4=qu'\gen {u'}$, giving $e_4=qu\gen u$ in canonical form. If $q=0$ then $e_4=pu\gen x$. Thus we obtain the following realizations:
\begin{align*}
\mathsf{A}_{3.5}\oplus\mathsf{A}_1&=\langle \gen t, \gen x, t\gen t + x\gen x + u\gen u, u\gen x \rangle\\
\mathsf{A}_{3.5}\oplus\mathsf{A}_1&=\langle \gen t, \gen x, t\gen t + x\gen x + u\gen u, u\gen u \rangle.
\end{align*}
Note that the algebra $\langle \gen t, \gen x, t\gen t + x\gen x + u\gen u, u\gen u \rangle$ is the same as the algebra $\langle \gen t, \gen x, t\gen t + x\gen x, u\gen u \rangle$, and under the equivalence transformation $t'=t, x'=x, u'=\ln |u|$ this algebra is equivalent to $\langle \gen t, \gen x, t\gen t + x\gen x, \gen u \rangle$.

\medskip\noindent $\underline{\mathsf{A}_{3.5}=\langle \gen x, \gen u, t\gen t + x\gen x + u\gen u \rangle}:$ The commutation relations give $e_4=at\gen t + pt\gen x + qt\gen u$. We must have $a\neq 0$ for otherwise $e_4=pt\gen x+qt\gen u$ together with $\gen x, \gen u$ in the equation for $G$ give $-pu_1+q=0$, which is a contradiction unless $p=q=0$. So $a\neq 0$ and we may assume $a=1$, giving $e_4=t\gen t + pt\gen x+q\gen u$. The residual equivalence group is $\mathscr{E}(e_1, e_2, e_3): t'=kt,\; x'=x+\alpha t,\; u'=u+\beta t$. Under such a transformation $e_4$ is mapped to $e'_4=t'\gen {t'}+(p+\alpha)\gen {x'} + (q+\beta)\gen {u'}$ and we choose $\alpha, \beta$ so that $p+\alpha=q+\beta=0$ so that $e'_4=t'\gen {t'}$ and $e_4=t\gen t$ in canonical form and we have the realization
\[
\mathsf{A}_{3.5}\oplus\mathsf{A}_1=\langle \gen x, \gen u, t\gen t + x\gen x + u\gen u, t\gen t \rangle.
\]
Now this algebra is obviously the same as $\langle \gen x, \gen u, x\gen x + u\gen u, t\gen t\rangle$ and this in turn is equivalent to $\langle \gen x, \gen u, x\gen x + u\gen u, \gen t\rangle$ under the equivalence transformation $t'=\ln|t|,\; x'=x,\; u'=u$.

\subsection{Inequivalent non-linearizing realizations of $\mathsf{A}_{3.5}\oplus\mathsf{A}_1.$}

\begin{align*}
\mathsf{A}_{3.5}\oplus\mathsf{A_1}&=\langle \gen t, \gen x, t\gen t + x\gen x, \gen u \rangle\\
\mathsf{A}_{3.5}\oplus\mathsf{A_1}&=\langle \gen x, \gen u, x\gen x + u\gen u, \gen t \rangle\\
\mathsf{A}_{3.5}\oplus\mathsf{A_1}&=\langle \gen t, \gen x, t\gen t + x\gen x + u\gen u, u\gen x \rangle.
\end{align*}

\section{Realizations of $\mathsf{A}_{3.6}\oplus\mathsf{A_1}$.}

\medskip\noindent $\underline{\mathsf{A}_{3.6}=\langle \gen u, x\gen u, 2x\gen x + u\gen u \rangle}:$ From the commutation relations we find $e_4=a(t)\gen t + c(t)x^{1/2}\gen u$, and we note that $a(t)\neq 0$ for otherwise $e_4=c(t)x^{1/2}\gen u$ and together with $e_1, e_2$ this linearizes the evolution equation. The residual equivalence group is $\mathscr{E}(e_1, e_2, e_3): t'=T(t),\; x'=x,\, u'=u+x^{1/2}U(t)$. Under such a transformation $e_4$ is mapped to $e'_4=a(t)\dot{T}(t)\gen {t'} + x^{1/2}[a(t)\dot{U}(t)+c(t)]\gen {u'}$. We choose $U(t)$ so that $a(t)\dot{U}(t)+c(t)=0$ and $T(t)$ so that $a(t)\dot{T}(t)=1$. Thus $e'_4=\gen {t'}$, so that $e_4=\gen t$ in canonical form, and we have the realization
\[
\mathsf{A}_{3.6}\oplus\mathsf{A_1}= \langle \gen u, x\gen u, 2x\gen x + u\gen u, \gen t \rangle.
\]

\medskip\noindent $\underline{\mathsf{A}_{3.6}=\langle \gen u, x\gen u, t\gen t + 2x\gen x + u\gen u \rangle}:$ From the commutation relations we find $e_4=at\gen t + tc(\sigma)\gen u$ where $\sigma=x^2/t$ and we note that $a\neq 0$ for otherwise $e_1, e_2, e_4$ linearize the evolution equation. Thus we have $a\neq 0$ and we may assume $a=1$, so that $e_4=t\gen t + tc(\sigma)\gen u$. The residual equivalence group is $\mathscr{E}(e_1, e_2, e_3): t'=kt,\; x'=x,\, u'=u+tU(\sigma)$. Under such a transformation $e_4$ is mapped to $e'_4=t'\gen {t'} + t[c(\sigma) + U(\sigma)-\sigma U'(\sigma)]\gen {u'}$. We choose $U(\sigma)$ so that $c(\sigma) + U(\sigma)-\sigma U'(\sigma)$ and then $e'_4=t'\gen {t'}$, so that $e_4=t\gen t$ in canonical form, and we have the realization
\[
\mathsf{A}_{3.6}\oplus\mathsf{A_1}= \langle \gen u, x\gen u, t\gen t + 2x\gen x + u\gen u, t\gen t \rangle.
\]
However, this algebra is just the same as $\langle \gen u, x\gen u, 2x\gen x + u\gen u, t\gen t\rangle$ and this algebra in turn is equivalent to $\langle \gen u, x\gen u, 2x\gen x + u\gen u, \gen t\rangle$ under the equivalence transformation $t'=\ln |t|,\; x'=x,\; u'=u$. So we obtain the same algebra as in the first case.

\medskip\noindent $\underline{\mathsf{A}_{3.6}=\langle \gen t, \gen x, t\gen t -
x\gen x \rangle}:$ From the commutation relations we find $e_4=c(u)\gen u$. The residual equivalence group is $\mathscr{E}(e_1, e_2, e_3): t'=t,\; x'=x,\; u'=U(u)$. Under such a transformation $e_4$ is mapped to $e'_4=c(u)U'(u)\gen {u'}$. We choose $U(u)$ so that $c(u)U'(u)=1$ to give $e'_4=\gen {u'}$, and so $e_4=\gen u$ in canonical form. So we have the realization
\[
\mathsf{A}_{3.6}\oplus\mathsf{A}_1=\langle \gen t, \gen x, t\gen t-x\gen x, \gen u\rangle.
\]

\medskip\noindent $\underline{\mathsf{A}_{3.6}=\langle \gen x, \gen u, x\gen x -
u\gen u \rangle}:$ From the commutation relations we find $e_4=a(t)\gen t$. The residual equivalence group is $\mathscr{E}(e_1, e_2, e_3): t'=T(t),\; x'=x,\; u'=u$. Under such a transformation $e_4$ is mapped to $e'_4=a(t)\dot{T}(t)\gen {t'}$. We choose $T(t)$ so that $a(t)\dot{T}(t)=1$ to give $e'_4=\gen {t'}$, and so $e_4=\gen t$ in canonical form. So we have the realization
\[
\mathsf{A}_{3.6}\oplus\mathsf{A}_1=\langle \gen x, \gen u, x\gen x-u\gen u, \gen t\rangle.
\]

\medskip\noindent $\underline{\mathsf{A}_{3.6}=\langle \gen t, \gen x, t\gen t -
x\gen x - u\gen u \rangle}:$ From the commutation relations we find $e_4=pu\gen x + qu\gen u$. The residual equivalence group is $\mathscr{E}(e_1, e_2, e_3): t'=t,\; x'=x + \alpha u,\; u'=ku$. Under such a transformation $e_4$ is mapped to $e'_4=[p+\alpha q]u\gen {x'} + qu'\gen {u'}\gen {u'}$. If $q\neq 0$ then we choose $\alpha$ so that $p+\alpha q=0$ and then $e'_4=qu'\gen {u'}$, giving $e_4=u\gen u$ in canonical form.  If $q=0$ then we have $e_4=u\gen x$ in canonical form. So we have the realizations
\begin{align*}
\mathsf{A}_{3.6}\oplus\mathsf{A}_1&=\langle \gen t, \gen x, t\gen t-x\gen x-u\gen u, u\gen x\rangle\\
\mathsf{A}_{3.6}\oplus\mathsf{A}_1&=\langle \gen t, \gen x, t\gen t-x\gen x-u\gen u, u\gen u\rangle.
\end{align*}
Note however that $\langle \gen t, \gen x, t\gen t-x\gen x-u\gen u, u\gen u\rangle$ is the same algebra as $\langle \gen t, \gen x, t\gen t-x\gen x, u\gen u\rangle$, which is equivalent under $t'=t,\; x'=x,\; u'=\ln|u|$ to $\langle \gen t, \gen x, t\gen t-x\gen x, \gen u\rangle$.

\medskip\noindent $\underline{\mathsf{A}_{3.6}=\langle \gen x, \gen u, t\gen t + x\gen x -
u\gen u \rangle}:$ From the commutation relations we find $e_4=at\gen t+pt\gen x + qt\gen u$. We must have $a\neq 0$ for otherwise $e_4=pt\gen x+qt\gen u$ together with $e_1, e_2$ in the equation for $G$ give $q-pu_1=0$, which is a contradiction unless $p=q=0$. Thus $a\neq 0$ and we may assume $a=1$ and $e_4=t\gen t +pt\gen x+qt\gen u$. The residual equivalence group is $\mathscr{E}(e_1, e_2, e_3): t'=kt,\; x'=x+\alpha t,\; u'=u+\beta t$. Under such a transformation $e_4$ is mapped to $e'_4=t'\gen {t'}+t[p+\alpha]\gen {x'} + t[q+\beta]\gen {u'}$. We choose $\alpha,\; \beta$ so that $p+\alpha=q+\beta=0$ so that $e'_4=t'\gen {t'}$, and so $e_4=t\gen t$ in canonical form. So we have the realization
\[
\mathsf{A}_{3.6}\oplus\mathsf{A}_1=\langle \gen x, \gen u, t\gen t + x\gen x-u\gen u, t\gen t\rangle.
\]
Note that this algebra is just the same as $\langle \gen x, \gen u, x\gen x-u\gen u, t\gen t\rangle$ which is equivalent to $\langle \gen x, \gen u, x\gen x-u\gen u, \gen t\rangle$ under the equivalence transformation $t'=\ln|t|,\; x'=x,\; u'=u$.

\subsection{Inequivalent non-linearizing realizations of $\mathsf{A}_{3.6}\oplus\mathsf{A}_1.$}

\begin{align*}
\mathsf{A}_{3.6}\oplus\mathsf{A_1}&=\langle \gen u, x\gen u, 2x\gen x+u\gen u, \gen t \rangle\\
\mathsf{A}_{3.6}\oplus\mathsf{A_1}&=\langle \gen t, \gen x, t\gen t - x\gen x, \gen u \rangle\\
\mathsf{A}_{3.6}\oplus\mathsf{A_1}&=\langle \gen x, \gen u, x\gen x - u\gen u, \gen t \rangle\\
\mathsf{A}_{3.6}\oplus\mathsf{A_1}&=\langle \gen t, \gen x, t\gen t -x\gen x - u\gen u, u\gen x \rangle.
\end{align*}

\section{Realizations of $\mathsf{A}_{3.7}\oplus\mathsf{A_1}$.}

\medskip\noindent $\underline{\mathsf{A}_{3.7}=\langle \gen u, x\gen u, (1-q)x\gen x + u\gen u \rangle}:$ Note that $0<|q|<1$. From the commutation relations we find $e_4=a(t)\gen t + c(t)x^{1/(1-q)}\gen u$, and we note that $a(t)\neq 0$ for otherwise $e_1, e_2, e_4$ form a rank-one three-dimensional abelian Lie algebra which linearizes the evolution equation. The residual equivalence group is $\mathscr{E}(e_1, e_2, e_3): t'=T(t),\; x'=x,\, u'=u+x^{1/(1-q)}U(t)$. Under such a transformation $e_4$ is mapped to $e'_4=a(t)\dot{T}(t)\gen {t'} + x^{1/(1-q)}[a(t)\dot{U}(t)+c(t)]\gen {u'}$. We choose $U(t)$ so that $a(t)\dot{U}(t)+c(t)=0$ and $T(t)$ so that $a(t)\dot{T}(t)=1$. Thus $e'_4=\gen {t'}$, so that $e_4=\gen t$ in canonical form, and we have the realization
\[
\mathsf{A}_{3.7}\oplus\mathsf{A_1}= \langle \gen u, x\gen u, (1-q)x\gen x + u\gen u, \gen t \rangle.
\]

\medskip\noindent $\underline{\mathsf{A}_{3.7}=\langle \gen u, x\gen u, t\gen t + (1-q)x\gen x + u\gen u \rangle}:$ From the commutation relations we find $e_4=at\gen t + tc(\sigma)\gen u$ with $\sigma=t^{1-q}/x$. Again $a\neq 0$ for otherwise $e_1, e_2, e_3$ linearize the evolution equation. Thus we may assume that $a=1$ and we have $e_4=t\gen t + tc(\sigma)\gen u$. The residual equivalence group is $\mathscr{E}(e_1, e_2, e_3): t'=kt,\; x'=x,\; u'=u+tU(\sigma)$ and under such a transformation, $e_4$ is mapped to $e'_4=t'\gen {t'} + t[c(\sigma)+U(\sigma)+\sigma U'(\sigma)]\gen {u'}$ and we choose $U(\sigma)$ so that $c(\sigma)+U(\sigma)+\sigma U'(\sigma)=0$ giving $e'_4=t'\gen {t'}$ and so $e_4=t\gen t$ in canonical form. This gives the realization
\[
\mathsf{A}_{3.7}\oplus\mathsf{A_1}= \langle \gen u, x\gen u, t\gen t + (1-q)x\gen x + u\gen u, t\gen t \rangle.
\]
Now this is the same as the algebra $\langle \gen u, x\gen u, (1-q)x\gen x + u\gen u, t\gen t \rangle$, which in turn is equivalent to $\langle \gen u, x\gen u, (1-q)x\gen x + u\gen u, \gen t \rangle$ under the equivalence transformation $t'=\ln|t|,\; x'=x,\; u'=u$.

\medskip\noindent $\underline{\mathsf{A}_{3.7}=\langle \gen t, \gen x, t\gen t + qx\gen x\rangle}:$ From the commutation relations we have $e_4=c(u)\gen u$. The residual equivalence group is $\mathscr{E}(e_1, e_2, e_3): t'=t,\; x'=x+Y(u),\; u'=U(u)$, and under such a transformation $e_4$ is mapped to $e'_4=c(u)U'(u)\gen {u'}$. We choose $U(u)$ so that $c(u)U'(u)=1$ so $e'_4=\gen {u'}$ and thus $e_4=\gen u$ in canonical form. Then we have the realization
\[
\mathsf{A}_{3.7}\oplus\mathsf{A}_1=\langle \gen t, \gen x, t\gen t + qx\gen x, \gen u\rangle.
\]

\medskip\noindent $\underline{\mathsf{A}_{3.7}=\langle \gen x, \gen t, qt\gen t + x\gen x\rangle}:$ From the commutation relations we have $e_4=c(u)\gen u$ and then, arguing as in the previous case, we obtain $e_4=\gen u$ in canonical form. Then we have the realization
\[
\mathsf{A}_{3.7}\oplus\mathsf{A}_1=\langle \gen x, \gen t, qt\gen t + x\gen x, \gen u\rangle.
\]

\medskip\noindent $\underline{\mathsf{A}_{3.7}=\langle \gen x, \gen u, qx\gen x + u\gen u\rangle}:$ From the commutation relations we have $e_4=a(t)\gen t$ and then, arguing as in the previous case, with the residual equivalence group $\mathscr{E}(e_1, e_2, e_3): t'=T(t),\; x'=x,\; u'=u$ we obtain $e_4=\gen t$ in canonical form. Then we have the realization
\[
\mathsf{A}_{3.7}\oplus\mathsf{A}_1=\langle \gen x, \gen u, qx\gen x + u\gen u, \gen t\rangle.
\]

\medskip\noindent $\underline{\mathsf{A}_{3.7}=\langle \gen t, \gen x, t\gen t + qx\gen x+u\gen u\rangle}:$ From the commutation relations we have $e_4=\alpha u^q\gen x + \beta u\gen u$. If $\beta=0$ then we take $e_4=u^q\gen x$ in canonical form. The residual equivalence group is $\mathscr{E}(e_1, e_2, e_3): t'=t,\; x'=x+\lambda u^q,\; u'=ku$, and under such a transformation $e_4$ is mapped to $e'_4=[\alpha +\beta\lambda]u^q\gen {x'}+\beta u'\gen {u'}$. For $\beta\neq 0$ we choose $\lambda$ so that $\alpha +\beta\lambda=0$ and then we may take $e_4=u\gen u$ in canonical form. Then we have the realizations
\begin{align*}
\mathsf{A}_{3.7}\oplus\mathsf{A}_1&=\langle \gen t, \gen x, t\gen t + qx\gen x + u\gen u, u^q\gen x\rangle\\
\mathsf{A}_{3.7}\oplus\mathsf{A}_1&=\langle \gen t, \gen x, t\gen t + qx\gen x + u\gen u, u\gen u\rangle.
\end{align*}
Note that $\langle \gen t, \gen x, t\gen t + qx\gen x + u\gen u, u\gen u\rangle$ is the same algebra as $\langle \gen t, \gen x, t\gen t + qx\gen x, u\gen u\rangle$ which in turn is equivalent to $\langle \gen t, \gen x, t\gen t + qx\gen x, \gen u\rangle$ under the equivalence transformation $t'=t,\; x'=x,\; u'=\ln|u|$.

\medskip\noindent $\underline{\mathsf{A}_{3.7}=\langle \gen x, \gen t, qt\gen t + x\gen x+u\gen u\rangle}:$ From the commutation relations we have $e_4=\alpha u\gen x + \beta u\gen u$. Arguing as above, we obtain the two canonical forms for $e_4=u\gen x$ and $e_4=u\gen u$. Then we have the realizations
\begin{align*}
\mathsf{A}_{3.7}\oplus\mathsf{A}_1&=\langle \gen x, \gen t, qt\gen t + x\gen x + u\gen u, u\gen x\rangle\\
\mathsf{A}_{3.7}\oplus\mathsf{A}_1&=\langle \gen x, \gen t, qt\gen t + x\gen x + u\gen u, u\gen u\rangle.
\end{align*}
Note that $\langle \gen t, \gen t, qt\gen t + x\gen x + u\gen u, u\gen u\rangle$ is the same algebra as $\langle \gen x, \gen t, qt\gen t + x\gen x, u\gen u\rangle$ which in turn is equivalent to $\langle \gen x, \gen t, qt\gen t + x\gen x, \gen u\rangle$ under the equivalence transformation $t'=t,\; x'=x,\; u'=\ln|u|$.

\medskip\noindent $\underline{\mathsf{A}_{3.7}=\langle \gen x, \gen u, t\gen t + x\gen x + qu\gen u\rangle}:$ From the commutation relations we have $e_4=at\gen t$ so we take $e_4=t\gen t$ in canonical form, giving the realization
\[
\mathsf{A}_{3.7}\oplus\mathsf{A}_1=\langle \gen x, \gen u, t\gen t + x\gen x + qu\gen u, t\gen t\rangle.
\]
However, this is the same algebra as $\langle \gen x, \gen u, x\gen x + qu\gen u, t\gen t\rangle$ which is equivalent to $\langle \gen x, \gen u, x\gen x + qu\gen u, \gen t\rangle$ under the equivalence transformation $t'=\ln|t|,\; x'=x,\; u'=u$.

\subsection{Inequivalent non-linearizing realizations of $\mathsf{A}_{3.7}\oplus\mathsf{A}_1.$}

\begin{align*}
\mathsf{A}_{3.7}\oplus\mathsf{A_1}&=\langle \gen u, x\gen u, (1-q)x\gen x+u\gen u, \gen t \rangle\\
\mathsf{A}_{3.7}\oplus\mathsf{A_1}&=\langle \gen t, \gen x, t\gen t + qx\gen x, \gen u \rangle\\
\mathsf{A}_{3.7}\oplus\mathsf{A_1}&=\langle \gen x, \gen t, qt\gen t + x\gen x, \gen u \rangle\\
\mathsf{A}_{3.7}\oplus\mathsf{A_1}&=\langle \gen x, \gen u, x\gen x + qu\gen u, \gen t \rangle\\
\mathsf{A}_{3.7}\oplus\mathsf{A_1}&=\langle \gen t, \gen x, t\gen t + qx\gen x + u\gen u, u^q\gen x \rangle\\
\mathsf{A}_{3.7}\oplus\mathsf{A_1}&=\langle \gen x, \gen t, qt\gen t + x\gen x + u\gen u, u\gen x \rangle.
\end{align*}

\section{Realizations of $\mathsf{A}_{3.8}\oplus\mathsf{A_1}$.}

\medskip\noindent $\underline{\mathsf{A}_{3.8}=\langle \gen u, x\gen u, -(1+x^2)\gen x - xu\gen u\rangle}:$ We find that we have $e_4=a(t)\gen t + c(t)\sqrt{1+x^2}\gen u$ and we must have $a(t)\neq 0$ for otherwise $e_1, e_2, e_4$ will be a rank-one three-dimensional abelian Lie algebra which linearizes the evolution equation. The residual equivalence group is
$$\mathscr{E}(e_1, e_2, e_3):\quad  t'=T(t),\; x'=x,\; u'=u+\sqrt{1+x^2}U(t).$$ Under such a transformation, $e_4$ is mapped to
$$e'_4=a(t)\dot{T}(t)\gen {t'} + [c(t)+a(t)\dot{U}(t)]\sqrt{1+x^2}\gen {u'}$$ and we choose $U(t)$ so that $c(t)+a(t)\dot{U}(t)=0$ and $T(t)$ so that $a(t)\dot{T}(t)=1$ so that $e'_4=\gen {t'}$, giving $e_4=\gen t$ in canonical form. Thus we have the realization
\[
\mathsf{A}_{3.8}\oplus\mathsf{A_1}=\langle \gen u, x\gen u, -(1+x^2)\gen x - xu\gen u, \gen t\rangle
\]

\medskip\noindent $\underline{\mathsf{A}_{3.8}=\langle \gen u, x\gen u, -t\gen t -(1+x^2)\gen x - xu\gen u\rangle}:$ We find that we have $e_4=at\gen t + c(\sigma)(1+x^2)^{1/2}\gen u$ with $\sigma=t\exp(-\arctan x)$. Again we must have $a\neq 0$ for otherwise $e_1, e_2, e_4$ is a rank-one three-dimensional abelian Lie algebra which linearizes the evolution equation. We may take $a=1$ so that $e_4=t\gen t + c(\sigma)(1+x^2)^{1/2}\gen u$. The residual equivalence group is $$\mathscr{E}(e_1, e_2, e_3):\quad  t'=T(t),\quad x'=x,\quad u'=u+(1+x^2)^{1/2}U(\sigma).$$
Under such a transformation, $e_4$ is mapped to $$e'_4=t'\gen {t'} + [c(\sigma)+\sigma U'(\sigma)](1+x^2)^{1/2}\gen {u'}$$ and we choose $U(\sigma)$ so that $c(\sigma)+\sigma U'(\sigma)=0$ and $e'_4=t'\gen {t'}$, giving $e_4=t\gen t$ in canonical form. Thus we have the realization
\[
\mathsf{A}_{3.8}\oplus\mathsf{A_1}=\langle \gen u, x\gen u, -t\gen t-(1+x^2)\gen x - xu\gen u, t\gen t\rangle.
\]
This algebra is the same as $\langle \gen u, x\gen u, -(1+x^2)\gen x - xu\gen u, t\gen t\rangle$ and this in turn is equivalent to $\langle \gen u, x\gen u, -(1+x^2)\gen x - xu\gen u, \gen t\rangle$ under the equivalence transformation $t'=\ln|t|,\; x'=x,\; u'=u$.

\medskip\noindent $\underline{\mathsf{A}_{3.8}=\langle \gen x, \gen u, u\gen x-x\gen u\rangle}:$ The commutation relations give $e_4=a(t)\gen t$. The residual equivalence group is $\mathscr{E}(e_1, e_2, e_3): t'=T(t),\; x'=x,\; u'=u$. Under such a transformation, $e_4$ is mapped to $e'_4=a(t)\dot{T}(t)\gen {t'}$ and we choose $T(t)$ so that $a(t)\dot{T}(t)=1$ giving $e_4=\gen t$ in canonical form. We have the realization
\[
\mathsf{A}_{3.8}\oplus\mathsf{A_1}=\langle \gen x, \gen u, u\gen x - x\gen u, \gen t\rangle.
\]

\medskip\noindent $\underline{\mathsf{A}_{3.8}=\langle \gen x, \gen u, t\gen t + u\gen x-x\gen u\rangle}:$ The commutation relations give $e_4=at\gen t + b(t)\gen x + c(t)\gen u$ where $b(t)=\alpha \cos(\ln t) + \beta\sin(\ln t),\; c(t)=-\alpha \sin(\ln t) + \beta\cos(\ln t)$. We have $a\neq 0$ for otherwise $e_1, e_2, e_4$ give $-b'(t)u_1+c'(t)=0$ in the equation for $G$, and this is a contradiction unless $b(t)=c(t)=0$. So $a\neq 0$ and we may assume $a=1$. Thus we take $e_4=t\gen t + b(t)\gen x + c(t)\gen u$. The residual equivalence group is $\mathscr{E}(e_1, e_2, e_3): t'=kt,\; x'=x+Y(t),\; u'=u+U(t)$ where $Y(t)=p\cos(\ln t) + q\sin(\ln t)$ and $U(t)=-p\sin(\ln t) + q\cos(\ln t)$. Under such a transformation, $e_4$ is mapped to $e'_4=t'\gen {t'}+[b(_t)+U(t)]\gen {x'} + [c(t)-Y(t)]\gen {u'}$. We may choose $p, q$ so that $b(t)+U(t)=c(t)-Y(t)=0$, giving $e_4=t\gen t$ in canonical form.  We have the realization
\[
\mathsf{A}_{3.8}\oplus\mathsf{A_1}=\langle \gen x, \gen u, t\gen t +u\gen x - x\gen u, t\gen t\rangle.
\]
However, this is the same algebra as $\langle \gen x, \gen u, u\gen x - x\gen u, t\gen t\rangle$ which in turn is equivalent to $\langle \gen x, \gen u, u\gen x - x\gen u, \gen t\rangle$ under the equivalence transformation $t'=\ln |t|,\; x'=x,\; u'=u$.

\subsection{Inequivalent non-linearizing realizations of $\mathsf{A}_{3.8}\oplus\mathsf{A}_1.$}

\begin{align*}
\mathsf{A}_{3.8}\oplus\mathsf{A_1}&=\langle \gen u, x\gen u, -(1+x^2)\gen x - xu\gen u, \gen t\rangle\\
\mathsf{A}_{3.8}\oplus\mathsf{A_1}&=\langle \gen x, \gen u, u\gen x - x\gen u, \gen t\rangle.
\end{align*}

\section{Realizations of $\mathsf{A}_{3.9}\oplus\mathsf{A}_1.$}

\medskip\noindent $\underline{\mathsf{A}_{3.9}=\langle \gen u, x\gen u, -(1+x^2)\gen x+(q-x)u\gen u\rangle}:$ From the commutation relations we find that $e_4=a(t)\gen t + \kappa(t)(1+x^2)^{1/2}\exp(-q\arctan x)\gen u$, and $a(t)\neq 0$ for otherwise $e_1, e_2, e_4$ give a rank-one three-dimensional abelian Lie algebra which linearizes the evolution equation. The residual equivalence group is $\mathscr{E}(e_1, e_2, e_3): t'=T(t),\; x'=x,\; u'=u + \alpha(t)(1+x^2)^{1/2}\exp(-q\arctan x)$ and under such a transformation, $e_4$ is mapped to $e'_4=a(t)\dot{T}(t)\gen {t'} + [\kappa(t)+a(t)\dot{\alpha}(t)](1+x^2)^{1/2}\exp(-q\arctan x)\gen {u'}$. We choose $\alpha(t)$ so that $\kappa(t)+a(t)\dot{\alpha}(t)=0$ and $T(t)$ so that $a(t)\dot{T}(t)=1$ giving $e_4=\gen t$ in canonical form. We have the realization
\[
\mathsf{A}_{3.9}\oplus\mathsf{A}_1=\langle \gen u, x\gen u, -(1+x^2)\gen x +(q-x)u\gen u, \gen t\rangle.
\]

\medskip\noindent $\underline{\mathsf{A}_{3.9}=\langle \gen u, x\gen u, -t\gen t-(1+x^2)\gen x+(q-x)u\gen u\rangle}:$ From the commutation relations we find that $e_4=at\gen t + c(\sigma)(1+x^2)^{1/2}\exp(-q\arctan x)\gen u$ with $\sigma=t\exp(-\arctan x)$.Again $a\neq 0$ for otherwise $e_1, e_2, e_4$ give a rank-one three-dimensional abelian Lie algebra which linearizes the evolution equation. We may take $a=1$ so that $e_4=t\gen t + c(\sigma)(1+x^2)^{1/2}\exp(-q\arctan x)\gen u$. The residual equivalence group is $\mathscr{E}(e_1, e_2, e_3): t'=kt,\; x'=x,\; u'=u + U(\sigma)(1+x^2)^{1/2}\exp(-q\arctan x)$ and under such a transformation, $e_4$ is mapped to $e'_4=t'\gen {t'} + [c(\sigma)+\sigma U(\sigma)](1+x^2)^{1/2}\exp(-q\arctan x)\gen {u'}$. We choose $U(\sigma)$ so that $c(\sigma)+\sigma U(\sigma)=0$  giving $e_4=t\gen t$ in canonical form. We have the realization
\[
\mathsf{A}_{3.9}\oplus\mathsf{A}_1=\langle \gen u, x\gen u, t\gen t-(1+x^2)\gen x +(q-x)u\gen u, t\gen t\rangle.
\]
However, this algebra is the same as $\langle \gen u, x\gen u, -(1+x^2)\gen x +(q-x)u\gen u, t\gen t\rangle$ which in turn is equivalent to $\langle \gen u, x\gen u, -(1+x^2)\gen x +(q-x)u\gen u, \gen t\rangle$ under the equivalence transformation $t'=\ln |t|,\; x'=x,\; u'=u$.

\medskip\noindent $\underline{\mathsf{A}_{3.9}=\langle \gen x, \gen u, [qx+u]\gen x+(qu-x)\gen u\rangle}:$ From the commutation relations we find that $e_4=a(t)\gen t$. The residual equivalence group is $\mathscr{E}(e_1, e_2, e_3): t'=T(t),\; x'=x,\; u'=u$ and under such a transformation, $e_4$ is mapped to $e'_4=a(t)\dot{T}(t)\gen {t'}$. We choose $T(t)$ so that $a(t)\dot{T}(t)=1$  giving $e_4=\gen t$ in canonical form. We have the realization
\[
\mathsf{A}_{3.9}\oplus\mathsf{A}_1=\langle \gen x, \gen u, [qx+u]\gen x+(qu-x)\gen u, \gen t\rangle.
\]

\medskip\noindent $\underline{\mathsf{A}_{3.9}=\langle \gen x, \gen u, t\gen t + [qx+u]\gen x+(qu-x)\gen u\rangle}:$ From the commutation relations we find that $e_4=at\gen t + b(t)\gen x + c(t)\gen u$ where $b(t)=\alpha\cos(\ln |t|)+\beta\sin(\ln |t|),\; c(t)=-\alpha\sin(\ln |t|)+\beta\cos(\ln |t|)$. We have $a\neq 0$ for otherwise $e_1, e_2, e_4$ in the equation for $G$ give $-\dot{b}(t)u_1+\dot{c}(t)=0$ which is a contradiction unless $b(t)=c(t)=0$. So $a\neq 0$ and we may take $a=1$ so that $e_4=t\gen t + b(t)\gen x + c(t)\gen u$. The residual equivalence group is $\mathscr{E}(e_1, e_2, e_3): t'=kt,\; x'=x + Y(t),\; u'=u+ U(t)$ where $k\neq 0$ and $Y(t)=t^q[\lambda\cos(\ln|t|) + \mu\sin(\ln|t|)]$ and $U(t)=t^q[-\lambda\sin(\ln|t|) + \mu\cos(\ln|t|)]$ and under such a transformation, $e_4$ is mapped to $e'_4=t'\gen {t'} + t^q[b(t)+qY(t)+t\dot{Y}(t)]\gen {x'} + t^q[c(t)+qU(t)+t\dot{U}(t)]\gen {u'}$. We may choose $Y(t),\; U(t)$ so that $b(t)+qY(t)+t\dot{Y}(t)=c(t)+qU(t)+t\dot{U}(t)=0$ giving $e_4=t\gen t$ in canonical form. We have the realization
\[
\mathsf{A}_{3.9}\oplus\mathsf{A}_1=\langle \gen x, \gen u, t\gen t + [qx+u]\gen x+(qu-x)\gen u, t\gen t\rangle.
\]
However, this algebra is the same as $\langle \gen x, \gen u, [qx+u]\gen x+(qu-x)\gen u, t\gen t\rangle$ which is in turn equivalent to the algebra $\langle \gen x, \gen u,[qx+u]\gen x+(qu-x)\gen u, \gen t\rangle$ under the equivalence transformation $t'=\ln|t|,\; x'=x,\; u'=u$.

\subsection{Inequivalent non-linearizing realizations of $\mathsf{A}_{3.9}\oplus\mathsf{A}_1.$}

\begin{align*}
\mathsf{A}_{3.9}\oplus\mathsf{A}_1&=\langle \gen u, x\gen u, -(1+x^2)\gen x +(q-x)u\gen u, \gen t\rangle\\
\mathsf{A}_{3.9}\oplus\mathsf{A}_1&=\langle \gen x, \gen u, [qx+u]\gen x+(qu-x)\gen u, \gen t\rangle.
\end{align*}

\section{Realizations of non-decomposable four-dimensional solvable Lie algebras.}

\medskip\noindent $\underline{\mathsf{A}_{4.1}}:$ Here we have $\mathsf{A}_{4.1}=\langle e_1, e_2, e_3, e_4\rangle$ where $\langle e_1, e_2, e_3\rangle$ is a three-dimensional abelian Lie algebra and $[e_1, e_4]=0,\; [e_2, e_4]=e_1,\; [e_3, e_4]=e_2$. We have ${\rm rank}\,\langle e_1, e_2, e_3\rangle=2$ or ${\rm rank}\,\langle e_1, e_2, e_3\rangle=3$ since if $\langle e_1, e_2, e_3\rangle$ is a rank-one realization,  the evolution equation is linearizable.

\smallskip\noindent $\underline{{\rm rank}\,\langle e_1, e_2, e_3\rangle=3}:$ In this case we have three inequivalent possibilities:

$\langle \gen t, \gen x, \gen u \rangle$, $\langle \gen x, \gen t, \gen u \rangle$ or $\langle \gen x, \gen u, \gen t \rangle$.
The first two are impossible to implement since if $e_1=\gen t,\; e_2=\gen x$ and $e_4=a(t)\gen t + b\gen x + c\gen u$ then $[e_2, e_4]=e_1$ is not possible; if $e_2=\gen t, e_3=\gen u$ then $[e_3, e_4]=e_2$ is impossible. Thus we have the possibility $\langle e_1, e_2, e_3\rangle=\langle \gen x, \gen u, \gen t \rangle$ and the commutation relations give $e_4=a\gen t + [u+\beta]\gen x + [t+\gamma]\gen u$. The residual equivalence group is $\mathscr{E}(e_1, e_2, e_3): t'=t+k,\; x'=x+l,\, u'=u+m$ and under such a transformation $e_4$ is mapped to $e'_4=a\gen {t'} + [u'+\beta - m]\gen {x'} + [t'+\gamma-k]\gen {u'}$ and we choose $k=\gamma,\, m=\beta$ so we obtain $e_4=a\gen t + u\gen x+t\gen u$ in canonical form. We have the realization
\[
\mathsf{A}_{4.1}=\langle \gen x, \gen u, \gen t, u\gen x + t\gen u\rangle.
\]

\smallskip\noindent $\underline{{\rm rank}\,\langle e_1, e_2, e_3\rangle=2}:$ In this case we have the three possibilities $\langle \gen t , \gen u, x\gen u\rangle$, $\langle \gen u , \gen t, x\gen u\rangle$, or $\langle\gen u, x\gen u, \gen t\rangle$. The first two are not implementable since we either have $[\gen u, e_4]=\gen t$ or $[x\gen u, e_4]=\gen t$ which are clearly impossible with our type of symmetry vector field. For $\langle\gen u, x\gen u, \gen t\rangle$ the commutation relations give $e_4=a\gen t - \gen x + [tx + c(x)]\gen u$. The residual equivalence group is $\mathscr{E}(e_1, e_2, e_3): t'=t+k,\; x'=x,\; u'=u+U(x)$ and under such a transformation $e_4$ is mapped to $e'_4=a\gen {t'} -\gen {x'} + [t'x'+c(x)-kx-U'(x)]\gen {u'}$ and we choose $U(x)$ so that $c(x)-kx-U'(x)=0$ so that $e_4=a\gen t -\gen x + tx\gen u$ in canonical form. We obtain the realization
\[
\mathsf{A}_{4.1}=\langle \gen u, x\gen u, \gen t, -\gen x + tx\gen u\rangle.
\]

\bigskip\noindent{\bf Realizations of $\mathsf{A}_{4.1}$:}

\begin{align*}
\mathsf{A}_{4.1}&=\langle \gen x, \gen u, \gen t, u\gen x + t\gen u\rangle\\
\mathsf{A}_{4.1}&=\langle \gen u, x\gen u, \gen t, -\gen x + tx\gen u\rangle.
\end{align*}

\bigskip\noindent $\underline{\mathsf{A}_{4.2}}:$ In this case $\langle e_1, e_2, e_3\rangle$ is an abelian ideal and $[e_1, e_4]=qe_1,\; [e_2, e_4]=e_2,\; [e_3, e_4]=e_2+e_3$ and $q\neq 0$. Again we consider only the rank-two and rank-three realizations of $\langle e_1, e_2, e_3\rangle$ since the rank-one realization linearizes the evolution equation.

\smallskip\noindent $\underline{\langle e_1, e_2, e_3\rangle=\langle \gen t, \gen x, \gen u \rangle}:$ The commutation relations give $e_4=[qt+\alpha]\gen t + [x+u+\beta]\gen x + [u+\gamma]\gen u$. The residual equivalence group is $\mathscr{E}(e_1, e_2, e_3): t'=t+k,\; x'=x+l,\; u'=u+m$. Under such a transformation $e_4$ is mapped to $e'_4=[qt'+\alpha - qk]\gen {t'} + [x'+u' + \beta-l-m]\gen {x'} + [u'+\gamma - m]\gen {u'}$ and we may choose $k, l, m$ so that $\alpha - qk=\beta-l-m=\gamma - m=0$ since $q\neq 0$, so that $e_4=qt\gen t + [x+u]\gen x + u\gen u$ in canonical form. We have the realization
\[
\mathsf{A}_{4.2}=\langle \gen t, \gen x, \gen u, qt\gen t + [u+x]\gen x + u\gen u\rangle.
\]

\smallskip\noindent $\underline{\langle e_1, e_2, e_3\rangle=\langle \gen x, \gen t, \gen u \rangle}:$ This is impossible since we require $[e_3, e_4]=[\gen u, e_4]=\gen t+\gen u$ which is impossible with $e_4=a(t)\gen t + b\gen x + c\gen u$.

\smallskip\noindent $\underline{\langle e_1, e_2, e_3\rangle=\langle \gen x, \gen u, \gen t \rangle}:$ The commutation relations give $e_4=[t+l]\gen t + [qx+\beta]\gen x + [t+u+\gamma]\gen u$ and arguing as above with the residual equivalence group  $\mathscr{E}(e_1, e_2, e_3)$ we find $e_4=t\gen t + qx\gen x + [t+u]\gen u$ in canonical form. We have the realization
\[
\mathsf{A}_{4.2}=\langle \gen x, \gen u, \gen t, t\gen t + qx\gen x + [t+u]\gen u\rangle.
\]

\smallskip\noindent $\underline{\langle e_1, e_2, e_3\rangle=\langle \gen t, \gen u, x\gen u \rangle}:$ The commutation relations give $e_4=[qt+l]\gen t-\gen x + [u+c(x)]\gen u$. The residual equivalence group is   $\mathscr{E}(e_1, e_2, e_3): t'=t+k,\; x'=x,\; u'=u+U(x)$ and under such a transformation $e_4$ is mapped to $e'_4=[qt+l-qk]\gen {t'} -\gen {x'} + [u'+c(x)-U'(x)]\gen {u'}$ and we choose $U(x)$ so that $c(x)-U'(x)=0$ and $k$ such that $l-qk=0$ since $q\neq 0$ so that $e_4=qt\gen t -\gen x + u\gen u$ in canonical form. We obtain the realization
\[
\mathsf{A}_{4.2}=\langle \gen t, \gen u, x\gen u, qt\gen t - \gen x + u\gen u\rangle.
\]

\smallskip\noindent $\underline{\langle e_1, e_2, e_3\rangle=\langle \gen u, \gen t, x\gen u \rangle}:$ This is impossible to implement since with $e_4=a(t)\gen t + b\gen x + c\gen u$ we require $[x\gen u, e_4]=\gen t + x\gen u$ which is impossible.

\smallskip\noindent $\underline{\langle e_1, e_2, e_3\rangle=\langle \gen u, x\gen u, \gen t \rangle}:$ The commutation relations give $e_4=[t+l]\gen t + (q-1)x\gen x + [qu+tx+c(x)]\gen u$. Using the residual group $\mathscr{E}(e_1, e_2, e_3): t'=t+k,\; x'=x,\; u'=u+U(x)$ as above, we find that we may choose $e_4=t\gen t + (q-1)x\gen x + [qu + tx]\gen u$ in canonical form. We have the realization
\[
\mathsf{A}_{4.2}=\langle \gen u, x\gen u, \gen t, t\gen t + (q-1)x\gen x + [qu+tx]\gen u\rangle.
\]

\bigskip\noindent{\bf Realizations of $\mathsf{A}_{4.2}$:}

\begin{align*}
\mathsf{A}_{4.2}&=\langle \gen t, \gen x, \gen u, qt\gen t + [u+x]\gen x + u\gen u\rangle\\
\mathsf{A}_{4.2}&=\langle \gen x, \gen u, \gen t, t\gen t + qx\gen x + [t+u]\gen u\rangle\\
\mathsf{A}_{4.2}&=\langle \gen t, \gen u, x\gen u, qt\gen t - \gen x + u\gen u\rangle\\
\mathsf{A}_{4.2}&=\langle \gen u, x\gen u, \gen t, t\gen t + (q-1)x\gen x + [qu+tx]\gen u\rangle.
\end{align*}

\bigskip\noindent $\underline{\mathsf{A}_{4.3}}:$ In this case $\langle e_1, e_2, e_3\rangle$ is an abelian ideal and $[e_1, e_4]=e_1,\; [e_2, e_4]=0,\; [e_3, e_4]=e_2$ and $q\neq 0$. Again we consider only the rank-two and rank-three realizations of $\langle e_1, e_2, e_3\rangle$ since the rank-one realization linearizes the evolution equation.

\smallskip\noindent $\underline{\langle e_1, e_2, e_3\rangle=\langle \gen t, \gen x, \gen u \rangle}:$ The commutation relations give $e_4=[t+l]\gen t + [u+m]\gen x$. Using the residual equivalence group $\mathscr{E}(e_1, e_2, e_3): t'=t+k,\; x'=x+\beta,\, u'=u+\gamma$ as before, we find that we may choose $e_4=t\gen t + u\gen x$ in canonical form. We have the realization
\[
\mathsf{A}_{4.3}=\langle \gen t, \gen x, \gen u, t\gen t + u\gen x\rangle.
\]

\smallskip\noindent $\underline{\langle e_1, e_2, e_3\rangle=\langle \gen x, \gen t, \gen u \rangle}:$ In this case the relation $[e_3, e_4]=e_2$ is impossible to implement.

\smallskip\noindent $\underline{\langle e_1, e_2, e_3\rangle=\langle \gen x, \gen u, \gen t \rangle}:$ The commutation relations together with the residual equivalence group $\mathscr{E}(e_1, e_2, e_3)$ give $e_4=[t+x]\gen x$ in canonical form. We have the realization
\[
\mathsf{A}_{4.3}=\langle \gen x, \gen u, \gen t, [t+x]\gen x\rangle.
\]

\smallskip\noindent $\underline{\langle e_1, e_2, e_3\rangle=\langle \gen t, \gen u, x\gen u \rangle}:$ The commutation relations give $e_4=[t+l]\gen t - \gen x + c(x)\gen u$. Using the residual equivalence group $\mathscr{E}(e_1, e_2, e_3): t'=t+k,\; x'=x+\beta,\; u'=u+U(x)$ we find that $e_4$ is mapped to $e'_4=[t'+l-k]\gen {t'} - \gen {x'} + [c(x)-U'(x)]\gen {u'}$ and we may always choose $k,\; U(x)$ so that $k=l,\; U'(x)=c(x)$ giving $e_4=t\gen t - \gen x$ in canonical form. We have the realization
\[
\mathsf{A}_{4.3}=\langle \gen t, \gen u, x\gen u, t\gen t - \gen x\rangle.
\]

\smallskip\noindent $\underline{\langle e_1, e_2, e_3\rangle=\langle \gen u, \gen t, x\gen u \rangle}:$ This is impossible to implement since $[e_3, e_4]=e_2$ gives $[x\gen u, e_4]=\gen t$, which is not possible.

\smallskip\noindent $\underline{\langle e_1, e_2, e_3\rangle=\langle \gen u, x\gen u, \gen t \rangle}:$ The commutation relations give $e_4=a\gen t + x\gen x + [u+tx + c(x)]\gen u$. Using the residual equivalence group $\mathscr{E}(e_1, e_2, e_3)$ we find that we may choose $e_4=a\gen t +x\gen x + [u+tx]\gen u$ in canonical form. We obtain the realization
\[
\mathsf{A}_{4.3}=\langle \gen u, x\gen u, \gen t, x\gen x + [u+tx]\gen u\rangle.
\]

\bigskip\noindent{\bf Realizations of $\mathsf{A}_{4.3}$:}

\begin{align*}
\mathsf{A}_{4.3}&=\langle \gen t, \gen x, \gen u, t\gen t + u\gen x\rangle\\
\mathsf{A}_{4.3}&=\langle \gen x, \gen u, \gen t, [t+x]\gen x\rangle\\
\mathsf{A}_{4.3}&=\langle \gen t, \gen u, x\gen u, t\gen t - \gen x\rangle\\
\mathsf{A}_{4.3}&=\langle \gen u, x\gen u, \gen t, x\gen x + [u+tx]\gen u\rangle.
\end{align*}

\bigskip\noindent $\underline{\mathsf{A}_{4.4}}:$ In this case $\langle e_1, e_2, e_3\rangle$ is an abelian ideal and $[e_1, e_4]=e_1,\; [e_2, e_4]=e_1+e_2,\; [e_3, e_4]=e_2+e_3$ and $q\neq 0$. Again we consider only the rank-two and rank-three realizations of $\langle e_1, e_2, e_3\rangle$ since the rank-one realization linearizes the evolution equation.

We note that the realizations $\langle e_1, e_2, e_3\rangle=\langle \gen t, \gen x, \gen u\rangle$ and $\langle e_1, e_2, e_3\rangle=\langle \gen x, \gen t, \gen u\rangle$ are impossible to implement since $[e_2, e_4]=e_1+e_2$ is not possible in the first case, and $[e_3, e_4]=e_2+e_3$ is not possible in the second case, with $e_4=a(t)\gen t + b\gen x + c\gen u$. Also $\langle e_1, e_2, e_3\rangle=\langle \gen t, \gen u, x\gen u\rangle$ and $\langle e_1, e_2, e_3\rangle=\langle \gen u, \gen t, x\gen u\rangle$ are impossible to realize for the same reasons. So we are left with two possibilities for $\langle e_1, e_2, e_3\rangle$.

\smallskip\noindent $\underline{\langle e_1, e_2, e_3\rangle=\langle \gen x, \gen u, \gen t\rangle}:$ In this case the commutation relations give $e_4=[t+l]\gen t + [x+u+\beta]\gen x + [u+t+\gamma]\gen u$. Using the residual equivalence group $\mathscr{E}(e_1, e_2, e_3)$ we may choose $e_4=t\gen t + [x+u]\gen x + [u+t]\gen u$ in canonical form. So we have the realization
\[
\mathsf{A}_{4.4}=\langle \gen x, \gen u, \gen t, t\gen t + [x+u]\gen x + [u+t]\gen u\rangle.
\]

\smallskip\noindent $\underline{\langle e_1, e_2, e_3\rangle=\langle \gen u, x\gen u, \gen t\rangle}:$ In this case the commutation relations give $e_4=[t+l]\gen t - x\gen x + [u+t+c(x)]\gen u$. Using the residual equivalence group $\mathscr{E}(e_1, e_2, e_3)$ we may choose $e_4=t\gen t - x+\gen x + [u+t]\gen u$ in canonical form. So we have the realization
\[
\mathsf{A}_{4.4}=\langle \gen u, x\gen u, \gen t, t\gen t - x\gen x + [u+t]\gen u\rangle.
\]

\bigskip\noindent{\bf Realizations of $\mathsf{A}_{4.4}$:}

\begin{align*}
\mathsf{A}_{4.4}&=\langle \gen x, \gen u, \gen t, t\gen t + [x+u]\gen x + [u+t]\gen u\rangle\\
\mathsf{A}_{4.4}&=\langle \gen u, x\gen u, \gen t, t\gen t - x\gen x + [u+t]\gen u\rangle.
\end{align*}

\bigskip\noindent $\underline{\mathsf{A}_{4.5}}:$ In this case $\langle e_1, e_2, e_3\rangle$ is an abelian ideal and $[e_1, e_4]=e_1,\; [e_2, e_4]=qe_2,\; [e_3, e_4]=pe_3$ and $-1\leq p\leq q\leq 1,\; pq\neq 0$. Again we consider only the rank-two and rank-three realizations of $\langle e_1, e_2, e_3\rangle$.

\medskip\noindent $\underline{\langle e_1, e_2, e_3\rangle=\langle \gen t, \gen x, \gen u\rangle}:$ The commutation relations give $e_4=(t+l)\gen t + (qx+\beta)\gen x + (pu+\gamma)\gen u$. Using the residual equivalence group $\mathscr{E}(e_1, e_2, e_3)$ and standard arguments, we may choose $e_4=t\gen t + qx\gen x + pu\gen u$. We have the realization
\[
\mathsf{A}_{4.5}=\langle \gen t, \gen x, \gen u, t\gen t + qx\gen x + pu\gen u\rangle.
\]

\medskip\noindent $\underline{\langle e_1, e_2, e_3\rangle=\langle \gen x, \gen t, \gen u\rangle}:$ The commutation relations together with the residual equivalence group $\mathscr{E}(e_1, e_2, e_3)$ we find that we may take $e_4=qt\gen t + x\gen x + pu\gen u$. We have the realization
\[
\mathsf{A}_{4.5}=\langle \gen x, \gen t, \gen u, qt\gen t + x\gen x + pu\gen u\rangle.
\]

\medskip\noindent $\underline{\langle e_1, e_2, e_3\rangle=\langle \gen x, \gen u, \gen t\rangle}:$ The commutation relations together with the residual equivalence group $\mathscr{E}(e_1, e_2, e_3)$ we find that we may take $e_4=pt\gen t + x\gen x + qu\gen u$. We have the realization
\[
\mathsf{A}_{4.5}=\langle \gen x, \gen u, \gen t, pt\gen t + x\gen x + qu\gen u\rangle.
\]

\medskip\noindent $\underline{\langle e_1, e_2, e_3\rangle=\langle \gen t, \gen u, x\gen u\rangle}:$ The commutation relations together with the residual equivalence group $\mathscr{E}(e_1, e_2, e_3)$ and standard arguments, we find that we may take $e_4=t\gen t + (q-p)x\gen x + qu\gen u$. We have the realization
\[
\mathsf{A}_{4.5}=\langle \gen t, \gen u, x\gen u, t\gen t + (q-p)x\gen x + qu\gen u\rangle.
\]

\medskip\noindent $\underline{\langle e_1, e_2, e_3\rangle=\langle \gen u, \gen t, x\gen u\rangle}:$ The commutation relations together with the residual equivalence group $\mathscr{E}(e_1, e_2, e_3)$ and standard arguments, we find that we may take $e_4=qt\gen t + (1-p)x\gen x + u\gen u$. We have the realization
\[
\mathsf{A}_{4.5}=\langle \gen u, \gen t, x\gen u, qt\gen t + (1-p)x\gen x + u\gen u\rangle.
\]

\medskip\noindent $\underline{\langle e_1, e_2, e_3\rangle=\langle \gen u, \gen t, x\gen u\rangle}:$ The commutation relations together with the residual equivalence group $\mathscr{E}(e_1, e_2, e_3)$ and standard arguments, we find that we may take $e_4=pt\gen t + (1-q)x\gen x + u\gen u$. We have the realization
\[
\mathsf{A}_{4.5}=\langle \gen u, x\gen u, \gen t, pt\gen t + (1-q)x\gen x + u\gen u\rangle.
\]

\bigskip\noindent{\bf Realizations of $\mathsf{A}_{4.5}$:}

\begin{align*}
\mathsf{A}_{4.5}&=\langle \gen t, \gen x, \gen u, t\gen t + qx\gen x + pu\gen u\rangle\\
\mathsf{A}_{4.5}&=\langle \gen x, \gen t, \gen u, qt\gen t + x\gen x + pu\gen u\rangle\\
\mathsf{A}_{4.5}&=\langle \gen x, \gen u, \gen t, pt\gen t + x\gen x + qu\gen u\rangle\\
\mathsf{A}_{4.5}&=\langle \gen t, \gen u, x\gen u, t\gen t + (q-p)x\gen x + qu\gen u\rangle\\
\mathsf{A}_{4.5}&=\langle \gen u, \gen t, x\gen u, qt\gen t + (1-p)x\gen x + u\gen u\rangle\\
\mathsf{A}_{4.5}&=\langle \gen u, x\gen u, \gen t, pt\gen t + (1-q)x\gen x + u\gen u\rangle.
\end{align*}

\bigskip\noindent $\underline{\mathsf{A}_{4.6}}:$ In this case $\langle e_1, e_2, e_3\rangle$ is an abelian ideal and $[e_1, e_4]=qe_1,\; [e_2, e_4]=pe_2-e_3,\; [e_3, e_4]=e_2+pe_3$ and $q\neq 0,\; p\geq 0$. Again we consider only the rank-two and rank-three realizations of $\langle e_1, e_2, e_3\rangle$.

\medskip\noindent $\langle e_1, e_2, e_3\rangle=\langle \gen x, \gen t, \gen u\rangle$ and $\langle e_1, e_2, e_3\rangle=\langle \gen x, \gen u, \gen t\rangle$ are impossible to implement since in the first case $[e_3, e_4]=e_2+pe_3$ is impossible to realize, and in the second case $[e_2, e_4]=pe_2-e_3$ is impossible to realize. The same applies to $\langle e_1, e_2, e_3\rangle=\langle \gen u, \gen t, x\gen u\rangle$ and $\langle e_1, e_2, e_3\rangle=\langle \gen u, x\gen u, \gen t\rangle$.

\medskip\noindent $\underline{\langle e_1, e_2, e_3\rangle=\langle \gen t, \gen x, \gen u\rangle}:$ The commutation relations together with the residual equivalence group $\mathscr{E}(e_1, e_2, e_3)$ and standard arguments, we find that we may take $e_4=qt\gen t + (px+u)\gen x + (pu-x)\gen u$. We have the realization
\[
\mathsf{A}_{4.6}=\langle \gen t, \gen x, \gen u, qt\gen t + (px+u)\gen x + (pu-x)\gen u\rangle.
\]

\medskip\noindent $\underline{\langle e_1, e_2, e_3\rangle=\langle \gen t, \gen u, x\gen u\rangle}:$ The commutation relations give $e_4=(qt+l)\gen t -(1+x^2)\gen x + [(p-x)u+c(x)]\gen u$. The residual equivalence group is $\mathscr{E}(e_1, e_2, e_3): t'=t+k,\; x'=x,\; u'=U(x)$ and under such a transformation $e_4$ is mapped to $e'_4=[qt'+l-qk]\gen {t'} -(1+x'^2)\gen {x'} + [(p-x')u'+c(x)-(p-x)U(x)-(1+x^2)U'(x)]\gen {u'}$ and we choose $k$ so that $l-qk=0$ and $U(x)$ so that $c(x)-(p-x)U(x)-(1+x^2)U'(x)=0$. Thus we have $e_4=qt\gen t - (1+x^2)\gen x + (p-x)u\gen u$ in canonical form. We have the realization
\[
\mathsf{A}_{4.6}=\langle \gen t, \gen u, x\gen u,  qt\gen t - (1+x^2)\gen x + (p-x)u\gen u\rangle.
\]

\bigskip\noindent{\bf Realizations of $\mathsf{A}_{4.6}$:}

\begin{align*}
\mathsf{A}_{4.6}&=\langle \gen t, \gen x, \gen u, qt\gen t + (px+u)\gen x + (pu-x)\gen u\rangle\\
\mathsf{A}_{4.6}&=\langle \gen t, \gen u, x\gen u qt\gen t - (1+x^2)\gen x + (p-x)u\gen u\rangle.
\end{align*}

\bigskip\noindent $\underline{\mathsf{A}_{4.7}}:$ In this case $\langle e_1, e_2, e_3\rangle=\mathsf{A}_{3.3}$ That is, $[e_1, e_2]=[e_1, e_3]=0,\; [e_2, e_3]=e_1$ and we also have $[e_1, e_4]=2e_1,\; [e_2, e_4]=e_2,\; [e_3, e_4]=e_2+e_3$.

We find that we have no realizations in the following cases:
\begin{align*}
\mathsf{A}_{3.3}&=\langle \gen x, \gen t, (t+u)\gen x\rangle\\
\mathsf{A}_{3.3}&=\langle \gen x, \gen u, u\gen x\rangle\\
\mathsf{A}_{3.3}&=\langle \gen x, \gen t, t\gen x-\gen u\rangle\\
\mathsf{A}_{3.3}&=\langle \gen x, \gen u, t\gen t +u\gen x\rangle,
\end{align*}
as follows from the commutation relations for $\mathsf{A}_{4.7}$. We are then left with two cases of $\mathsf{A}_{3.3}$.

\medskip\noindent$\underline{\langle e_1, e_2, e_3\rangle=\langle \gen u, x\gen u, -\gen x\rangle}:$ The commutation rules give
\[
e_4=a(t)\gen t + x\gen x + [2u-\frac{x^2}{2}+c(t)]\gen u.
\]
The residual equivalence group is $\mathscr{E}: t'=T(t),\; x'=x,\; u'=u + U(t)$ and under such a transformation $e_4$ is mapped to
\[
e'_4=a(t)\dot{T}(t)\gen {t'} + x'\gen {x'} + [2u'-\frac{x'^2}{2}+c(t)-2U(t)+a(t)\dot{U}(t)]\gen {u'}
\]
and we may always choose $U(t)$ so that $c(t)-2U(t)+a(t)\dot{U}(t)=0$ so that
\[
e'_4=a(t)\dot{T}(t)\gen {t'} + x'\gen {x'} + [2u'-\frac{x'^2}{2}]\gen {u'}
\]
If $a(t)=0$ we have $\displaystyle e_4= x\gen x + [2u-\frac{x^2}{2}]\gen u$ in canonical form. If $a(t)\neq 0$ we choose $T(t)$ so that $a(t)\dot{T}(t)=T(t)$ so that  we obtain $\displaystyle e_4= t\gen t + x\gen x + [2u-\frac{x^2}{2}]\gen u$ in canonical form. We have the realizations
\begin{align*}
\mathsf{A}_{4.7}&=\langle \gen u, x\gen u, -\gen x, x\gen x + [2u-\frac{x^2}{2}]\gen u\rangle\\
\mathsf{A}_{4.7}&=\langle \gen u, x\gen u, -\gen x, t\gen t + x\gen x + [2u-\frac{x^2}{2}]\gen u\rangle.
\end{align*}

\medskip\noindent$\underline{\langle e_1, e_2, e_3\rangle=\langle \gen u, x\gen u, -(\gen t+\gen x)\rangle}:$ The commutation rules give
\[
e_4=(t+l)\gen t + x\gen x + [2u-\frac{x^2}{2}+c(t-x)]\gen u.
\]
The equivalence group is $\mathscr{E}(e_1, e_2, e_3): t'=t+k,\; x'=x,\; u'=u + U(t-x)$. Under such a transformation $e_4$ is mapped to $e'_4$ given by
\[
e'_4=(t'+l-k)\gen {t'} + x'\gen {x'} + [2u'-\frac{x'^2}{2}+c(t-x)-2U(t-x)+(t-x+l)U'(t-x)]\gen {u'},
\]
and we choose $k=l$ and $U(t-x)$ so that $c(t-x)-2U(t-x)+(t-x+l)U'(t-x)=0$ and then we obtain
\[
e_4=t\gen t + x\gen x + [2u-\frac{x^2}{2}]\gen u
\]
in canonical form. We have the realization
\[
\mathsf{A}_{4.7}=\langle \gen u, x\gen u, -(\gen t+\gen x), t\gen t + x\gen x + [2u-\frac{x^2}{2}]\gen u\rangle.
\]

\medskip\noindent$\underline{\langle e_1, e_2, e_3\rangle=\langle \gen x, \gen u, u\gen x+t\gen u\rangle}:$ The commutation rules give $e_4=-\gen t + [2x + b(t)]\gen x + u\gen u$. The residual equivalence group is $\mathscr{E}(e_1, e_2, e_3): t'=t,\; x'=x+Y(t),\; u'=u$. Under such a transformation $e_4$ is mapped to $e'_4=-\gen {t'} + [2x'+b(t)-2Y(t)-\dot{Y}(t)]\gen {x'} + u'\gen {u'}$ and we choose $Y(t)$ so that $b(t)-2Y(t)-\dot{Y}(t)=0$ so that $e_4=-\gen t + 2x\gen x + u\gen u$ in canonical form. We have the realization
\[
\mathsf{A}_{4.7}=\langle \gen x, \gen u, u\gen x+t\gen u, -\gen t + 2x\gen x + u\gen u\rangle.
\]

\bigskip\noindent{\bf Realizations of $\mathsf{A}_{4.7}$:}

\begin{align*}
\mathsf{A}_{4.7}&=\langle \gen u, x\gen u, -\gen x, x\gen x + [2u-\frac{x^2}{2}]\gen u\rangle\\
\mathsf{A}_{4.7}&=\langle \gen u, x\gen u, -\gen x, t\gen t + x\gen x + [2u-\frac{x^2}{2}]\gen u\rangle\\
\mathsf{A}_{4.7}&=\langle \gen u, x\gen u, -(\gen t+\gen x), t\gen t + x\gen x + [2u-\frac{x^2}{2}]\gen u\rangle\\
\mathsf{A}_{4.7}&=\langle \gen x, \gen u, u\gen x+t\gen u, -\gen t + 2x\gen x + u\gen u\rangle.
\end{align*}

\bigskip\noindent $\underline{\mathsf{A}_{4.8}}:$ In this case $\langle e_1, e_2, e_3\rangle=\mathsf{A}_{3.3}$ That is, $[e_1, e_2]=[e_1, e_3]=0,\; [e_2, e_3]=e_1$ and we also have $[e_1, e_4]=(1+q)e_1,\; [e_2, e_4]=e_2,\; [e_3, e_4]=qe_3$ with $|q|\leq 1$.

\medskip\noindent$\underline{\langle e_1, e_2, e_3\rangle=\langle \gen u, x\gen u, -\gen x\rangle}:$ The commutation relations give
$$e_4=a(t)\gen t + qx\gen x + [(1+q)u+c(t)]\gen u.$$
The residual equivalence group is $\mathscr{E}(e_1, e_2, e_3): t'=T(t),\; x'=x,\; u'=u+U(t)$ and under such a transformation $e_4$ is mapped to $$e'_4=a(t)\dot{T}(t)\gen {t'} + qx'\gen {x'} + [(1+q)u' + c(t)-(1+q)U(t)+a(t)\dot{U}(t)]\gen {u'}.$$
For $a(t)\neq 0$ we may always choose $U(t)$ so that $c(t)-(1+q)U(t)+a(t)\dot{U}(t)=0$, and  we may choose $T(t)$ so that $a(t)\dot{T}(t)=T(t)$ so that
$$e_4=t\gen t + qx\gen x + (1+q)u\gen u$$
in canonical form. If $a(t)=0$ then for $-1<q\leq 1$ we may always choose $U(t)$ so that $c(t)-(1+q)U(t)=0$ so that we obtain $e_4=qx\gen x + (1+q)u\gen u$ in canonical form.  However, if $q=-1,\; a(t)=0$ then we have  $e_4=-x\gen x + c(t)\gen u$. In this case, if $\dot{c}(t)=0$ then we have $c(t)=\,{\rm const}$ and we have $e_4=-x\gen x + c\gen u$ with $c\in \mathbb{R}$. If $\dot{c}(t)\neq 0$ then $e_4$ is mapped by an equivalence transformation to $e'_4=-x'\gen {x'} + c(t)\gen {u'}$ and we may take $c(t)=T(t)$ in this case since $\dot{c}(t)\neq 0$ so that $e_4=-x\gen x + t\gen u$ in canonical form. We have the following realizations:
\begin{align*}
\mathsf{A}_{4.8}&=\langle \gen u, x\gen u, -\gen x, t\gen t + qx\gen x + (1+q)u\gen u\rangle,\\
\mathsf{A}_{4.8}&=\langle \gen u, x\gen u, -\gen x, qx\gen x + (1+q)u\gen u\rangle,\;\; -1<q\leq 1\\
\mathsf{A}_{4.8}&=\langle \gen u, x\gen u, -\gen x, -x\gen x + c\gen u\rangle,\;\,c\in \mathbb{R},\;\; q=-1\\
\mathsf{A}_{4.8}&=\langle \gen u, x\gen u, -\gen x, -x\gen x + t\gen u\rangle,\;\; q=-1.
\end{align*}

\medskip\noindent$\underline{\langle e_1, e_2, e_3\rangle=\langle \gen u, x\gen u, -(\gen t + \gen x)\rangle}:$ The commutation relations give
$$e_4=[qt+a]\gen t + qx\gen x + [(1+q)u + c(t-x)]\gen u.$$
The residual equivalence group is
$$\mathscr{E}(e_1, e_2, e_3): t'=t+k,\; x'=x,\; u'=u + U(t-x).$$ Under such a transformation $e_4$ is mapped to
\begin{equation*}
  \begin{split}
      e'_4=&[qt'+a-qk]\gen {t'} + qx'\gen {x'} + [(1+q)u' + c(t-x) - (1+q)U(t-x)+ \\
      & [q(t-x) + a]U'(t-x)]\gen {u'}.
  \end{split}
\end{equation*}
Now for any $-1\leq q\leq 1$ we may choose $U(t-x)$ so that $c(t-x) - (1+q)U(t-x) + [q(t-x) + a]U'(t-x)=0$, so we then have
$$e'_4=[qt'+a-qk]\gen {t'} + qx'\gen {x'} + (1+q)u'\gen {u'}.$$
If $q\neq 0$ then we may choose $k$ so that $a-qk=0$ giving $e_4=qt\gen t + qx\gen x + (1+q)u\gen u$ in canonical form. If $q=0$ we have $e_4=a\gen t +u\gen u$. However, the three operators $\gen u, x\gen u, u\gen u$ linearize the evolution equation, so we must have $a\neq 0$. We have the following realizations
\begin{align*}
\mathsf{A}_{4.8}&=\langle \gen u, x\gen u, -(\gen t + \gen x), qt\gen t + qx\gen x + (1+q)u\gen u\rangle,\;\; q\neq 0\\
\mathsf{A}_{4.8}&=\langle \gen u, x\gen u, -(\gen t + \gen x), a\gen t + u\gen u\rangle,\;\; a\neq 0,\; q=0.
\end{align*}

\medskip\noindent$\underline{\langle e_1, e_2, e_3\rangle=\langle \gen x, \gen t, (t + u)\gen x\rangle}:$ The commutation relations give $e_4=(t+l)\gen t + [(1+q)x+b(u)]\gen x + (u-l)\gen u$. The residual equivalence group is $\mathscr{E}(e_1, e_2, e_3): t'=t+k,\; x'=x+Y(u),\; u'=u-k$. Under such a transformation $e_4$ is mapped to $e'_4=(t'+l-k)\gen {t'} + [(1+q)x'+b(u)-(1+q)Y(u)+(u-l)Y'(u)]\gen {u'} + (u'+k-l)\gen {u'}$. We may always choose $Y(u)$ so that $b(u)-(1+q)Y(u)+(u-l)Y'(u)=0$ and we choose $k=l$ so that we obtain $e_4=t\gen t + (1+q)x\gen x + u\gen u$ in canonical form. Thus we have the realization
\[
\mathsf{A}_{4.8}=\langle \gen x, \gen t, (t+u)\gen x, t\gen t + (1+q)x\gen x + u\gen u\rangle.
\]

\medskip\noindent$\underline{\langle e_1, e_2, e_3\rangle=\langle \gen x, \gen u, u\gen x\rangle}:$ The commutation relations give $e_4=a(t)\gen t + [(1+q)x + b(t)]\gen x+ u\gen u$. The residual equivalence group is $\mathscr{E}(e_1, e_2, e_3): t'=T(t),\; x'=x+Y(t),\; u'=u$. Under such a transformation, $e_4$ is mapped to $e'_4=a(t)\dot{T}(t)\gen {t'} + [(1+q)x' + b(t)-(1+q)Y(t)+a(t)\dot{Y}(t)]\gen {x'} + u'\gen u$. If $a(t)\neq 0$ we may always choose $Y(t)$ so that $b(t)-(1+q)Y(t)+a(t)\dot{Y}(t)=0$ and $T(t)$ so that $a(t)\dot{T}(t)=T(t)$ giving $e_4=t\gen t + (1+q)x\gen x + u\gen u$ in canonical form. However, if $a(t)=0$ we have $e'_4=[(1+q)x' + b(t)-(1+q)Y(t)]\gen {x'} + u'\gen u$, and for $q\neq -1$ we may choose $Y(t)$ so that $b(t)-(1+q)Y(t)=0$ so that $e_4=(1+q)x\gen x + u\gen u$ in canonical form. If $a(t)=0,\; q=-1$ we have $e_4=b(t)\gen x + u\gen u$. If $b(t)={\rm const}$, then we have $e_4=p\gen t + u\gen u$ with $p\in \mathbb{R}$. If $\dot{b}(t)\neq 0$ we have that $e'_4=b(t)\gen {x'}+u'\gen {u'}$ and we may choose $b(t)=T(t)$ so that $t\gen x + u\gen u$ in canonical form. We have the realizations:
\begin{align*}
\mathsf{A}_{4.8}&=\langle \gen x, \gen u, u\gen x, t\gen t + (1+q)x\gen x + u\gen u\rangle\\
\mathsf{A}_{4.8}&=\langle \gen x, \gen u, u\gen x, (1+q)x\gen x + u\gen u\rangle,\;\; q\neq -1.\\
\mathsf{A}_{4.8}&=\langle \gen x, \gen u, u\gen x, p\gen x + u\gen u\rangle,\;\; p\in \mathbb{R},\; q=-1\\
\mathsf{A}_{4.8}&=\langle \gen x, \gen u, u\gen x, t\gen x + u\gen u\rangle, \; q=-1.
\end{align*}

\medskip\noindent$\underline{\langle e_1, e_2, e_3\rangle=\langle \gen x, \gen u, u\gen x+t\gen u\rangle}:$ The commutation relations give $e_4=(1-q)t\gen t + [(1+q)x + b(t)]\gen x + u\gen u$. The residual equivalence group is $\mathscr{E}(e_1, e_2, e_3): t'=t,\; x'=x+Y(t),\; u'=u$. Under such a transformation $e_4$ is mapped to $e'_4=(1-q)t'\gen {t'} + [(1+q)x' + b(t)-(1+q)Y(t) + (1-q)t\dot{Y}(t)]\gen {x'} + u'\gen {u'}$. We may always choose $Y(t)$ so that $b(t)-(1+q)Y(t) + (1-q)t\dot{Y}(t)=0$ and so we find that $e_4=(1-q)t\gen t + (1+q)x\gen x + u\gen u$ in canonical form. We have the realization
\[
\mathsf{A}_{4.8}=\langle \gen x, \gen u, u\gen x+t\gen u, (1-q)t\gen t + (1+q)x\gen x + u\gen u\rangle.
\]

\medskip\noindent$\underline{\langle e_1, e_2, e_3\rangle=\langle \gen x, \gen t, t\gen x - \gen u\rangle}:$ The commutation relations give $e_4=(t+l)\gen t + [(1+q)x -lu + \beta]\gen x + [qu + \gamma]\gen u$. The residual equivalence group is $\mathscr{E}(e_1, e_2, e_3): t'=t+k,\; x'=x-ku + \sigma,\; u'=u+\lambda$. Under such a transformation $e_4$ is mapped to $e'_4=(t'+l-k)\gen {t'} + [(1+q)x'+ (k-l)u+\beta-k\gamma-(1+q)\sigma]\gen {x'} + [qu'+\gamma-q\lambda]\gen {u'}$. We may always choose $k=l$. For $q\neq 0, -1$ we may choose $\sigma,\; \lambda$ so that $\beta-k\gamma-(1+q)\sigma=\gamma-q\lambda=0$ giving $e_4=t\gen t + (1+q)x\gen x + qu\gen u $ in canonical form. If $q=0$ then we may still choose $\sigma$ so that $\beta-k\gamma-\sigma=0$ which then gives $e_4=t\gen t + x\gen x +c\gen u$ in canonical form, with $c\in \mathbb{R}$. If $q=-1$ then we may choose $\lambda$ so that $\lambda+\gamma=0$ so that we have $e_4=t\gen t + \beta\gen x -u\gen u$ in canonical form, with $\beta\in \mathbb{R}$. Note $\gen x$ already belongs to the algebra. Thus we have the realizations
\begin{align*}
\mathsf{A}_{4.8}&=\langle \gen x, \gen t, t\gen x-\gen u, t\gen t + (1+q)x\gen x + qu\gen u\rangle, \;\; q\neq 0\\
\mathsf{A}_{4.8}&=\langle \gen x, \gen t, t\gen x-\gen u, t\gen t + x\gen x + c\gen u\rangle,\;c\in \mathbb{R},\;\; q=0.
\end{align*}

\medskip\noindent$\underline{\langle e_1, e_2, e_3\rangle=\langle \gen x, \gen u, t\gen t + u\gen x\rangle}:$ The commutation relations give $e_4=t(q\ln |t|+\kappa)\gen t + [(1+q)x + p\ln |t| + \lambda]\gen x + [u + p]\gen u$, with $p,\, \kappa,\, \lambda \in \mathbb{R}$. The residual equivalence group is $\mathscr{E}(e_1, e_2, e_3): t'=kt,\; x'=x+Y(t),\; u'=u+ U(t)$ with $U(t)=t\dot {Y}(t)$. Under such a transformation $e_4$ is mapped to $e'_4=t'(q\ln |t'| + \kappa +q\ln |k|)\gen {t'} + [(1+q)x' + p\ln t - (1+q)Y(t) + (q\ln t + \kappa)t\dot{Y}(t) + \lambda]\gen {x'} + [u'+ (q\ln t + \kappa)t\dot{U}(t)-U(t) + p]\gen {u'}$. We can always arrange for $Y(t)$ to be such that $p\ln t - (1+q)Y(t) + (q\ln t + \kappa)t\dot{Y}(t) + \lambda=0$ for $|q|\leq 1$. Then, differentiating this with respect to $t$ and then multiplying by $t$ and noting that $U(t)=t\dot{Y}(t)$, we find that we automatically have $(q\ln t + \kappa)t\dot{U}(t)-U(t) + p=0$. Thus we may choose $Y(t)$ so that $e'_4=t'(q\ln |t'| + \kappa +q\ln |k|)\gen {t'} + (1+q)x'\gen {x'} + u'\gen {u'}$. For $q\neq 0$ we choose $k$ so that $\kappa +q\ln |k|=0$ and so we find that $e_4=q\ln |t|\gen t + (1+q)x\gen x + u\gen u$ in canonical form. If $q=0$ then we have $e_4=\alpha t\gen t + x\gen x + u\gen u$ in canonical form. Thus we have the realizations
\begin{align*}
\mathsf{A}_{4.8}&=\langle \gen x, \gen u, t\gen t+u\gen x, qt\ln |t|\gen t + (1+q)x\gen x + u\gen u\rangle, \;\; q\neq 0\\
\mathsf{A}_{4.8}&=\langle \gen x, \gen u, t\gen t+u\gen x, \alpha t\gen t + x\gen x + u\gen u\rangle, \;\;\alpha\in \mathbb{R},\;\; q=0.
\end{align*}

\bigskip\noindent{\bf Realizations of $\mathsf{A}_{4.8}$:}

\begin{align*}
\mathsf{A}_{4.8}&=\langle \gen u, x\gen u, -\gen x, t\gen t + qx\gen x + (1+q)u\gen u\rangle,\\
\mathsf{A}_{4.8}&=\langle \gen u, x\gen u, -\gen x, qx\gen x + (1+q)u\gen u\rangle,\;\; -1<q\leq 1\\
\mathsf{A}_{4.8}&=\langle \gen u, x\gen u, -\gen x, -x\gen x + c\gen u\rangle,\;\,c\in \mathbb{R},\;\; q=-1\\
\mathsf{A}_{4.8}&=\langle \gen u, x\gen u, -\gen x, -x\gen x + t\gen u\rangle,\;\; q=-1\\
\mathsf{A}_{4.8}&=\langle \gen u, x\gen u, -(\gen t + \gen x), qt\gen t + qx\gen x + (1+q)u\gen u\rangle,\;\; q\neq 0\\
\mathsf{A}_{4.8}&=\langle \gen u, x\gen u, -(\gen t + \gen x), a\gen t + u\gen u\rangle,\;\; a\neq 0,\; q=0\\
\mathsf{A}_{4.8}&=\langle \gen x, \gen t, (t+u)\gen x, t\gen t + (1+q)x\gen x + u\gen u\rangle\\
\mathsf{A}_{4.8}&=\langle \gen x, \gen u, u\gen x, t\gen t + (1+q)x\gen x + u\gen u\rangle\\
\mathsf{A}_{4.8}&=\langle \gen x, \gen u, u\gen x, (1+q)x\gen x + u\gen u\rangle,\;\; q\neq -1.\\
\mathsf{A}_{4.8}&=\langle \gen x, \gen u, u\gen x, p\gen x + u\gen u\rangle,\;\; p\in \mathbb{R},\; q=-1\\
\mathsf{A}_{4.8}&=\langle \gen x, \gen u, u\gen x, t\gen x + u\gen u\rangle, q=-1.\\
\mathsf{A}_{4.8}&=\langle \gen x, \gen u, u\gen x+t\gen u, (1-q)t\gen t + (1+q)x\gen x + u\gen u\rangle\\
\mathsf{A}_{4.8}&=\langle \gen x, \gen t, t\gen x-\gen u, t\gen t + (1+q)x\gen x + qu\gen u\rangle, \;\; q\neq 0\\
\mathsf{A}_{4.8}&=\langle \gen x, \gen t, t\gen x-\gen u, t\gen t + x\gen x + c\gen u\rangle,\;\, b,\, c\in \mathbb{R},\; q=0\\
\mathsf{A}_{4.8}&=\langle \gen x, \gen u, t\gen t+u\gen x, qt\ln |t|\gen t + (1+q)x\gen x + u\gen u\rangle, \;\; q\neq 0\\
\mathsf{A}_{4.8}&=\langle \gen x, \gen u, t\gen t+u\gen x, \alpha t\gen t + x\gen x + u\gen u\rangle, \;\;\alpha\in \mathbb{R},\;\; q=0.
\end{align*}

\bigskip\noindent $\underline{\mathsf{A}_{4.9}}:$ In this case $\langle e_1, e_2, e_3\rangle=\mathsf{A}_{3.3}$ That is, $[e_1, e_2]=[e_1, e_3]=0,\; [e_2, e_3]=e_1$ and we also have $[e_1, e_4]=2qe_1,\; [e_2, e_4]=qe_2-e_3,\; [e_3, e_4]=e_2+qe_3$ with $q\geq 0$.

We have the following realizations of $\mathsf{A}_{4.9}$:

It is impossible to realize $\mathsf{A}_{4.9}$ with $\langle e_1, e_2, e_3\rangle$ equal to the following realizations:

\begin{align*}
\mathsf{A}_{3.3}&=\langle \gen u, x\gen u, -\gen x\rangle\\
\mathsf{A}_{3.3}&=\langle \gen u, x\gen u, -\gen t -\gen x\rangle\\
\mathsf{A}_{3.3}&=\langle \gen x, \gen t, (t+u)\gen x\rangle\\
\mathsf{A}_{3.3}&=\langle \gen x, \gen u, u\gen x\rangle\\
\mathsf{A}_{3.3}&=\langle \gen x, \gen t, t\gen x - \gen u\rangle\\
\mathsf{A}_{3.3}&=\langle \gen x, \gen u, t\gen t + u\gen x\rangle.
\end{align*}

\medskip\noindent $\underline{\langle e_1, e_2, e_3\rangle=\langle \gen u, x\gen u, -\gen x\rangle}:$ The commutation relation $[e_1, e_4]=2qe_1$ gives $e_4=a(t)\gen t + b(t,x)\gen x + [2qu+c(t,x)]\gen u$ and then  $[e_2, e_4]=qe_2-e_3$ gives $[2qx + c_x(t,x)-b(t, x)]\gen u= qx\gen u + \gen x$ which is impossible.

\medskip\noindent $\underline{\langle e_1, e_2, e_3\rangle=\langle \gen u, x\gen u, -(\gen t + \gen x)\rangle}:$ We have $[e_2, e_4]=[2qx + c_x(t,x)-b(t, x)]\gen u=qx\gen u + \gen t + \gen x$ which is impossible.

\medskip\noindent $\underline{\langle e_1, e_2, e_3\rangle=\langle \gen x, \gen t, (t+u)\gen x\rangle}:$ Here we cannot possibly implement $[e_3, e_4]=e_2+qe_3$.

\medskip\noindent $\underline{\langle e_1, e_2, e_3\rangle=\langle \gen x, \gen u, u\gen x\rangle}:$ From the commutation relations
$$[e_1, e_4]=2qe_1,\quad [e_2, e_4]=qe_2-e_3$$ we find $\displaystyle e_4=a(t)\gen t + [2qx - \frac{u^2}{2}+b(t)]\gen x + [qu+c(t)]\gen u$ then we find $[e_3, e_4]=e_2+qe_3$ gives $[qu-c(t)]\gen x=\gen u + qu\gen x$ which is impossible.

\medskip\noindent $\underline{\langle e_1, e_2, e_3\rangle=\langle \gen x, \gen t, t\gen x-\gen u\rangle}:$ From the commutation relations
$$[e_1, e_4]=2qe_1,\quad [e_2, e_4]=qe_2-e_3$$ we find $\displaystyle e_4=(qt+l)\gen t + [2qx - \frac{t^2}{2}+b(u)]\gen x + [t+c(u)]\gen u$ then we find $[e_3, e_4]=e_2+qe_3$ gives $[qt-b'(u)-l]\gen x-c'(u)\gen u=\gen t + qt\gen x-q\gen u$ which is impossible.

\medskip\noindent $\underline{\langle e_1, e_2, e_3\rangle=\langle \gen x, \gen u, t\gen t+ u\gen x \rangle}:$ The relation $[e_2, e_4]=qe_2-e_3$ is impossible to implement.

\medskip\noindent $\underline{\langle e_1, e_2, e_3\rangle=\langle \gen x, \gen u, u\gen x+t\gen u\rangle}:$ From the commutation relations we find $\displaystyle e_4=-(1+t^2)\gen t + [2qx - \frac{u^2}{2}+b(t)]\gen x + (q-t)u\gen u$. The residual equivalence group is $\mathscr{E}(e_1, e_2, e_3): t'=t,\; x'=x + Y(t),\; u'=u$. Under such a transformation $e_4$ is mapped to $\displaystyle e'_4=-(1+t'^2)\gen {t'} + [2qx'-\frac{u'^2}{2} + b(t) - 2qY(t)-(1+t^2)\dot{Y}(t)]\gen {x'} + (q-t')u'\gen {u'}$. we may always choose $Y(t)$ so that $b(t) - 2qY(t)-(1+t^2)\dot{Y}(t)=0$ and thus we have $\displaystyle e_4=-(1+t^2)\gen t + [2qx - \frac{u^2}{2}]\gen x + (q-t)u\gen u$ in canonical form. We have the realization
\[
\mathsf{A}_{4.9}=\langle \gen x, \gen u, u\gen x + t\gen u, -(1+t^2)\gen t + [2qx - \frac{u^2}{2}]\gen x + (q-t)u\gen u\rangle.
\]

\bigskip\noindent $\underline{\mathsf{A}_{4.10}}:$ In this case $\langle e_1, e_2, e_3\rangle=\mathsf{A}_{3.5}$ That is, $[e_1, e_3]=e_1,\; [e_1, e_2]=0,\; [e_2, e_3]=e_2$ and we also have $[e_1, e_4]=-e_2,\; [e_2, e_4]=e_1,\; [e_3, e_4]=0$.

\medskip\noindent The following realizations of $\langle e_1, e_2, e_3\rangle=\mathsf{A}_{3.5}$ do not allow an extension to $\mathsf{A}_{4.10}$:

\begin{align*}
\mathsf{A}_{3.5}&=\langle \gen x, \gen t, t\gen t + x\gen x\rangle\\
\mathsf{A}_{3.5}&=\langle \gen t, \gen x, t\gen t + x\gen x\rangle\\
\mathsf{A}_{3.5}&=\langle \gen t, \gen x, t\gen t +  x\gen x + u\gen u\rangle\\
\mathsf{A}_{3.5}&=\langle \gen x, \gen t, t\gen t +  x\gen x + u\gen u\rangle.
\end{align*}
The reason for this is that if $e_1=\gen t,\; e_2=\gen x$ then $[e_2, e_4]=e_1$ is impossible to implement since $e_4=a(t)\gen t + b\gen x + c\gen u$. The same applies to $e_1=\gen x,\, e_2=\gen t$: then we cannot implement $[e_1, e_4]=-e_2$.

\medskip\noindent $\underline{\langle e_1, e_2, e_3\rangle=\langle \gen u, x\gen u, t\gen t+ u\gen u \rangle}:$ The commutation relations give $e_4=at\gen t -(1+x^2)\gen x + [tc(x)-xu]\gen u$. The residual equivalence group is $\mathscr{E}(e_1, e_2, e_3): t'=kt,\; x'=x,\; u'=u+tU(x)$. Under such a transformation, $e_4$ is mapped to $e'_4=at'\gen {t'} - (1+x'^2)\gen {x'} +[-x'u'+t[aU(x)-(1+x^2)U'(x)+c(x)-xU(x)]]\gen {u'}$ and we may always choose $U(x)$ so that $aU(x)-(1+x^2)U'(x)+c(x)-xU(x)=0$ giving $e_4=at\gen t - (1+x^2)\gen x - xu\gen u$ in canonical form. We have the realization
\[
\mathsf{A}_{4.10}=\langle \gen u, x\gen u, t\gen t+ u\gen u, at\gen t - (1+x^2)\gen x - xu\gen u \rangle,\;\; a\in \mathbb{R}.
\]

\medskip\noindent $\underline{\langle e_1, e_2, e_3\rangle=\langle \gen x, \gen u, x\gen x+ u\gen u \rangle}:$ The commutation relations give $e_4=a(t)\gen t + u\gen x - x\gen u$. The residual equivalence group is $\mathscr{E}(e_1, e_2, e_3): t'=T(t),\; x'=x,\; u'=u$. Under such a transformation, $e_4$ is mapped to $e'_4=a(t)\dot{T}(t)\gen {t'} + u'\gen {x'} - x'\gen {u'}$ and if $a(t)\neq 0$ we may choose $T(t)$ so that $a(t)\dot{T}(t)=T(t)$ giving $e_4=t\gen t + u\gen x-x\gen u$ in canonical form. We have $e_4=u\gen x-x\gen u$ if $a(t)=0$. We have the realizations
\begin{align*}
\mathsf{A}_{4.10}&=\langle \gen x, \gen u, x\gen x+ u\gen u, u\gen x - x\gen u \rangle\\
\mathsf{A}_{4.10}&=\langle \gen x, \gen u, t\gen t + x\gen x+ u\gen u, u\gen x - x\gen u \rangle.
\end{align*}

\medskip\noindent $\underline{\langle e_1, e_2, e_3\rangle=\langle \gen x, \gen u, t\gen t + x\gen x+ u\gen u \rangle}:$ The commutation relations give $e_4=at\gen t + [u + t\beta]\gen x - [x + t\gamma]\gen u$. The residual equivalence group is $\mathscr{E}(e_1, e_2, e_3): t'=kt,\; x'=x + Y(t),\; u'=u + U(t)$ with $Y(t)=pt$ and $U(t)=qt$. Under such a transformation, $e_4$ is mapped to $e'_4=at'\gen {t'} + [u'+t(\beta-q+ap)]\gen {x'} + [-x'+ t(\gamma + p + aq)]\gen {u'}$ and we may always choose $p,q$ so that $\beta-q+ap=\gamma + p + aq=0$, giving $e_4=at\gen t + u\gen x-x\gen u$ in canonical form. We have the realization
\[
\mathsf{A}_{4.10}=\langle \gen x, \gen u, t\gen t + x\gen x+ u\gen u, at\gen t + u\gen x - x\gen u \rangle,\;\; a\in \mathbb{R}.
\]

\subsection{Non-linearizing realizations of admissible non-decomposable four-dimensional solvable algebras.}

\bigskip\noindent $\mathsf{A}_{4.1}$:

\begin{align*}
\mathsf{A}_{4.1}&=\langle \gen x, \gen u, \gen t, u\gen x + t\gen u\rangle\\
\mathsf{A}_{4.1}&=\langle \gen u, x\gen u, \gen t, -\gen x + tx\gen u\rangle.
\end{align*}

\medskip\noindent $\mathsf{A}_{4.2}$:

\begin{align*}
\mathsf{A}_{4.2}&=\langle \gen t, \gen x, \gen u, qt\gen t + [u+x]\gen x + u\gen u\rangle\\
\mathsf{A}_{4.2}&=\langle \gen x, \gen u, \gen t, t\gen t + qx\gen x + [t+u]\gen u\rangle\\
\mathsf{A}_{4.2}&=\langle \gen t, \gen u, x\gen u, qt\gen t - \gen x + u\gen u\rangle\\
\mathsf{A}_{4.2}&=\langle \gen u, x\gen u, \gen t, t\gen t + (q-1)x\gen x + [qu+tx]\gen u\rangle.
\end{align*}

\bigskip\noindent $\mathsf{A}_{4.3}$:

\begin{align*}
\mathsf{A}_{4.3}&=\langle \gen t, \gen x, \gen u, t\gen t + u\gen x\rangle\\
\mathsf{A}_{4.3}&=\langle \gen x, \gen u, \gen t, [t+x]\gen x\rangle\\
\mathsf{A}_{4.3}&=\langle \gen t, \gen u, x\gen u, t\gen t - \gen x\rangle\\
\mathsf{A}_{4.3}&=\langle \gen u, x\gen u, \gen t, x\gen x + [u+tx]\gen u\rangle.
\end{align*}

\bigskip\noindent $\mathsf{A}_{4.4}$:

\begin{align*}
\mathsf{A}_{4.4}&=\langle \gen x, \gen u, \gen t, t\gen t + [x+u]\gen x + [u+t]\gen u\rangle\\
\mathsf{A}_{4.4}&=\langle \gen u, x\gen u, \gen t, t\gen t - x\gen x + [u+t]\gen u\rangle.
\end{align*}

\bigskip\noindent $\mathsf{A}_{4.5}$:

\begin{align*}
\mathsf{A}_{4.5}&=\langle \gen t, \gen x, \gen u, t\gen t + qx\gen x + pu\gen u\rangle\\
\mathsf{A}_{4.5}&=\langle \gen x, \gen t, \gen u, qt\gen t + x\gen x + pu\gen u\rangle\\
\mathsf{A}_{4.5}&=\langle \gen x, \gen u, \gen t, pt\gen t + x\gen x + qu\gen u\rangle\\
\mathsf{A}_{4.5}&=\langle \gen t, \gen u, x\gen u, t\gen t + (q-p)x\gen x + qu\gen u\rangle\\
\mathsf{A}_{4.5}&=\langle \gen u, \gen t, x\gen u, qt\gen t + (1-p)x\gen x + u\gen u\rangle\\
\mathsf{A}_{4.5}&=\langle \gen u, x\gen u, \gen t, pt\gen t + (1-q)x\gen x + u\gen u\rangle.
\end{align*}

\bigskip\noindent $\mathsf{A}_{4.6}$:

\begin{align*}
\mathsf{A}_{4.6}&=\langle \gen t, \gen x, \gen u, qt\gen t + (px+u)\gen x + (pu-x)\gen u\rangle\\
\mathsf{A}_{4.6}&=\langle \gen t, \gen u, x\gen u qt\gen t - (1+x^2)\gen x + (p-x)u\gen u\rangle.
\end{align*}

\bigskip\noindent $\mathsf{A}_{4.7}$:

\begin{align*}
\mathsf{A}_{4.7}&=\langle \gen u, x\gen u, -\gen x, x\gen x + [2u-\frac{x^2}{2}]\gen u\rangle\\
\mathsf{A}_{4.7}&=\langle \gen u, x\gen u, -\gen x, t\gen t + x\gen x + [2u-\frac{x^2}{2}]\gen u\rangle\\
\mathsf{A}_{4.7}&=\langle \gen u, x\gen u, -(\gen t+\gen x), t\gen t + x\gen x + [2u-\frac{x^2}{2}]\gen u\rangle\\
\mathsf{A}_{4.7}&=\langle \gen x, \gen u, u\gen x+t\gen u, -\gen t + 2x\gen x + u\gen u\rangle.
\end{align*}

\bigskip\noindent $\mathsf{A}_{4.8}$:

\begin{align*}
\mathsf{A}_{4.8}&=\langle \gen u, x\gen u, -\gen x, t\gen t + qx\gen x + (1+q)u\gen u\rangle,\\
\mathsf{A}_{4.8}&=\langle \gen u, x\gen u, -\gen x, qx\gen x + (1+q)u\gen u\rangle,\;\; -1<q\leq 1\\
\mathsf{A}_{4.8}&=\langle \gen u, x\gen u, -\gen x, -x\gen x + c\gen u\rangle,\;\,c\in \mathbb{R},\;\; q=-1\\
\mathsf{A}_{4.8}&=\langle \gen u, x\gen u, -\gen x, -x\gen x + t\gen u\rangle,\;\; q=-1\\
\mathsf{A}_{4.8}&=\langle \gen u, x\gen u, -(\gen t + \gen x), qt\gen t + qx\gen x + (1+q)u\gen u\rangle,\;\; q\neq 0\\
\mathsf{A}_{4.8}&=\langle \gen u, x\gen u, -(\gen t + \gen x), a\gen t + u\gen u\rangle,\;\; a\neq 0,\; q=0\\
\mathsf{A}_{4.8}&=\langle \gen x, \gen t, (t+u)\gen x, t\gen t + (1+q)x\gen x + u\gen u\rangle\\
\mathsf{A}_{4.8}&=\langle \gen x, \gen u, u\gen x, t\gen t + (1+q)x\gen x + u\gen u\rangle\\
\mathsf{A}_{4.8}&=\langle \gen x, \gen u, u\gen x, (1+q)x\gen x + u\gen u\rangle,\;\; q\neq -1.\\
\mathsf{A}_{4.8}&=\langle \gen x, \gen u, u\gen x, p\gen x + u\gen u\rangle,\;\; p\in \mathbb{R},\; q=-1\\
\mathsf{A}_{4.8}&=\langle \gen x, \gen u, u\gen x, t\gen x + u\gen u\rangle, q=-1.\\
\mathsf{A}_{4.8}&=\langle \gen x, \gen u, u\gen x+t\gen u, (1-q)t\gen t + (1+q)x\gen x + u\gen u\rangle\\
\mathsf{A}_{4.8}&=\langle \gen x, \gen t, t\gen x-\gen u, t\gen t + (1+q)x\gen x + qu\gen u\rangle, \;\; q\neq 0\\
\mathsf{A}_{4.8}&=\langle \gen x, \gen t, t\gen x-\gen u, t\gen t + x\gen x + c\gen u\rangle,\; c\in \mathbb{R},\; q=0\\
\mathsf{A}_{4.8}&=\langle \gen x, \gen u, t\gen t+u\gen x, qt\ln |t|\gen t + (1+q)x\gen x + u\gen u\rangle, \;\; q\neq 0\\
\mathsf{A}_{4.8}&=\langle \gen x, \gen u, t\gen t+u\gen x, \alpha t\gen t + x\gen x + u\gen u\rangle, \;\;\alpha\in \mathbb{R},\;\; q=0.
\end{align*}

\bigskip\noindent $\mathsf{A}_{4.9}$:

\[
\mathsf{A}_{4.9}=\langle \gen x, \gen u, u\gen x + t\gen u, -(1+t^2)\gen t + [2qx - \frac{u^2}{2}]\gen x + (q-t)u\gen u\rangle.
\]

\bigskip\noindent $\mathsf{A}_{4.10}$:

\begin{align*}
\mathsf{A}_{4.10}&=\langle \gen u, x\gen u, t\gen t+ u\gen u, at\gen t - (1+x^2)\gen x - xu\gen u \rangle,\;\; a\in \mathbb{R}\\
\mathsf{A}_{4.10}&=\langle \gen x, \gen u, x\gen x+ u\gen u, u\gen x - x\gen u \rangle\\
\mathsf{A}_{4.10}&=\langle \gen x, \gen u, t\gen t + x\gen x+ u\gen u, u\gen x - x\gen u \rangle\\
\mathsf{A}_{4.10}&=\langle \gen x, \gen u, t\gen t + x\gen x+ u\gen u, at\gen t + u\gen x - x\gen u \rangle,\;\; a\in \mathbb{R}.
\end{align*}

\section{Admissible five-dimensional solvable Lie algebras.} In this section we give a list of admissible five-dimensional solvable Lie algebras. We list only those realizations which do not automatically give linear (or linearizable) evolution equations. At the basis of this is the result that any evolution equation (quasi-linear, of third order) which admits an abelian symmetry algebra of dimension four or greater, or a rank-one three-dimensional abelian symmetry algebra, is linear or linearizable by a point transformation (which preserves evolution equations). For a classification of linear evolution equations we refer the reader to \cite{GungorLahnoZhdanov2004, BasarabHorwathGuengoer2017}.

\medskip\noindent The five-dimensional solvable Lie algebras are of two types: decomposable and non-decomposable. There are twenty-seven decomposable types: $5\mathsf{A}_1,\; \mathsf{A}_{2.2}\oplus 3\mathsf{A}_{1},\; 2\mathsf{A}_{2.2}\oplus\mathsf{A}_{1},\; \mathsf{A}_{3.k}\oplus 2\mathsf{A}_{1}\;\; (k=3,\dots, 9),\; \mathsf{A}_{3.k}\oplus\mathsf{A}_{2.2}\;\; (k=3,\dots, 9),\; \mathsf{A}_{4.k}\oplus\mathsf{A}_{1}\;\; (k=1,\dots, 10)$. Then there are thirty-nine non-decomposable types, labeled $\mathsf{A}_{5.k},\;\; k=1,\dots, 39$.

\medskip\noindent Note that the three-dimensional solvable Lie algebras $\mathsf{A}_{3.k}\;\; (k=3,\dots, 9)$ all contain a two-dimensional abelian ideal, so that $\mathsf{A}_{3.k}\oplus 2\mathsf{A}_{1}\;\; (k=3,\dots, 9)$ will always contain a four-dimensional abelian subalgebra, and so these algebras will either linearize the evolution equation or give no realization. The same applies to  $\mathsf{A}_{2.2}\oplus 3\mathsf{A}_{1}$: this contains a four-dimensional abelian Lie algebra and thus linearizes the equation or gives no realization. The five-dimensional abelian Lie algebra $5\mathsf{A}_1$ linearizes the evolution equation. Thus we have only to consider the decomposable Lie algebras $2\mathsf{A}_{2.2}\oplus\mathsf{A}_{1},\; \mathsf{A}_{3.k}\oplus\mathsf{A}_{2.2}\;\; (k=3,\dots, 9),\; \mathsf{A}_{4.k}\oplus\mathsf{A}_{1}\;\; (k=1,\dots, 10)$.

\medskip\noindent For the non-decomposable solvable Lie algebras, $\mathsf{A}_{5.1},\ \mathsf{A}_{5.2}$ and $\mathsf{A}_{5.7}$ through to $\mathsf{A}_{5.18}$ all contain a four-dimensional abelian Lie algebra, and so linearize the evolution equation or give no realization. Thus we need only look at $\mathsf{A}_{5.3}, \mathsf{A}_{5.4},\; \mathsf{A}_{5.5},\;\mathsf{A}_{5.6}$ and $\mathsf{A}_{5.19}$ through to $\mathsf{A}_{5.39}$.

\subsection{$2\mathsf{A}_{2.2}\oplus\mathsf{A}_{1}$:}

\smallskip\noindent{$\underline{2\mathsf{A}_{2.2}=\langle \gen x, x\gen x, \gen t, t\gen t+u\gen u\rangle}$:} From the commutation relations we find $e_5=u\gen u$ in canonical form. Thus $2\mathsf{A}_{2.2}\oplus\mathsf{A}_{1}$ will contain $\gen t$ and $t\gen t$ which give $F=0$ in the equation for $F$. This is a contradiction, so we have no realization in this case.

\smallskip\noindent{$\underline{2\mathsf{A}_{2.2}=\langle \gen x, x\gen x, \gen u, u\gen u\rangle}$:} From the commutation relations we find $e_5=a(t)\gen t$. The residual equivalence group is $\mathscr{E}(e_1, e_2, e_3, e_4): t'=T(t),\; x'=x,\; u'=u$. Under such a transformation $e_5$ is mapped to $e'_5=a(t)\dot{T}(t)\gen {t'}$ and we may choose $T(t)$ so that $a(t)\dot{T}=1$ giving $e_5=\gen t$ in canonical form. Thus we have the realization
\[
2\mathsf{A}_{2.2}\oplus\mathsf{A}_{1}=\langle \gen x, x\gen x, \gen u, u\gen u, \gen t\rangle.
\]

\smallskip\noindent{$\underline{2\mathsf{A}_{2.2}=\langle \gen x, x\gen x, \gen u, t\gen t + u\gen u\rangle}$:} From the commutation relations we find $e_5=at\gen t+qt\gen u$. Note that we cannot have $a=0$ for then we would have $\gen u,\; t\gen u$ as symmetries, and these give the contradiction $1=0$ in the equation for $G$. So we may assume $a=1$ and we then take $e_5=t\gen t + qt\gen u$. The residual equivalence group is $\mathscr{E}(e_1, e_2, e_3, e_4): t'=kt,\; (k\neq 0),\; x'=x,\; u'=u+\kappa t$. Under such a transformation $e_5$ is mapped to $e'_5=t'\gen {t'}+[q+\kappa]t\gen {u'}$ and we choose $\kappa$ so that $q+\kappa=0$ giving $e_5=t\gen t$ in canonical form. Thus we have the realization
\[
2\mathsf{A}_{2.2}\oplus\mathsf{A}_{1}=\langle \gen x, x\gen x, \gen u, t\gen t + u\gen u, t\gen t\rangle.
\]
Clearly this algebra is equivalent to
\[
2\mathsf{A}_{2.2}\oplus\mathsf{A}_{1}=\langle \gen x, x\gen x, \gen u, u\gen u, \gen t\rangle
\]
on transforming by the equivalence transformation $t'=\ln |t|,\; x'=x,\; u'=u$.

\smallskip\noindent $\underline{2\mathsf{A}_{2.2}=\langle \gen t, t\gen t + x\gen x, xu\gen x, -u\gen u\rangle}$: From the commutation relations we find $e_5=pxu\gen x$ so that $\dim 2\mathsf{A}_{2.2}\oplus\mathsf{A}_{1}=4$ which is a contradiction. Hence we have no realization in this case.

\smallskip\noindent $\underline{2\mathsf{A}_{2.2}=\langle \gen x, t\gen t + x\gen x, tu\gen x, at\gen t - (1+a)u\gen u\rangle}$ with $a\neq -1$: The commutation relations give $e_5=\alpha(t\gen t -u\gen u) + tb(u)\gen x$ with $(1+a)ub'(u)=ab(u)$ so that we have $b(u)=\kappa u^{a/(1+a)}$ with $\kappa\in\mathbb{R}$.  When $\alpha=0$ we have $e_5=tb(u)\gen x$, which, together with $\gen x,\; tu\gen x$ give a rank-one three-dimensional abelian Lie algebra which linearizes the evolution equation for $a\neq 0$, and gives $e_5=t\gen x$ when $a=0$ and this gives together with $\gen x$, the contradiction $u_1=0$ in the equation for $G$. Thus we must have $\alpha\neq 0$ and we may then take $\alpha=1$ and $e_5=t\gen t - u\gen u + \kappa tu^{a/(1+a)}\gen x$. Furthermore, we have the residual equivalence group $\displaystyle\mathscr{E}(e_1, e_2, e_3, e_4): t'=kt,\; x'=x+\lambda tu^{a/(1+a)},\; u'=u/k$. Under such a transformation, $e_5$ is mapped to $\displaystyle e'_5=t'\gen {t'} - u'\gen {u'} +tu^{a/(1+a)}[\kappa +\lambda/(1+a)]\gen {x'}$ and we can clearly choose $\lambda$ so that $\kappa +\lambda/(1+a)=0$, thus giving $e_5=t\gen t - u\gen u$ in canonical form. We then have the following  realization:
\[
2\mathsf{A}_{2.2}\oplus\mathsf{A}_1=\langle \gen x, t\gen t + x\gen x, tu\gen x, at\gen t-(1+a)u\gen u, t\gen t - u\gen u\rangle,
\]
with $a\neq -1$, which is the same algebra as
\[
2\mathsf{A}_{2.2}\oplus\mathsf{A}_1=\langle \gen x, t\gen t + x\gen x, tu\gen x, -u\gen u, t\gen t - u\gen u \rangle.
\]
Now this algebra contains the subalgebra $\langle \gen x, tu\gen x, x\gen x\rangle$ and, as already remarked, this linearizes the evolution equation: under the equivalence transformation $t'=t,\; x'=tu,\; u'=x$ we find that $\langle \gen x, tu\gen x, x\gen x\rangle$ is equivalent to $\langle \gen u, x\gen u, u\gen u\rangle$ which gives the evolution equation $u_t=F(t,x)u_3 + G(t,x)u_2$. Thus we have no non-linearizing realization in this case.

\medskip\noindent $\underline{2\mathsf{A}_{2.2}=\langle \gen x, x\gen x + u\gen u, u\gen x, t\gen t - u\gen u\rangle}:$ The commutation relations give $e_5=at\gen t + ptu\gen x$. Clearly, $a\not=0$ since otherwise $u\gen x,\; tu\gen x$ would be symmetries giving the contradiction $u_1=0$ in the equation for $G$. Thus we may take $a=1$ and $e_5=t\gen t + ptu\gen x$. The residual equivalence group is $\mathscr{E}(e_1, e_2, e_3, e_4): t'=kt,\; x'=x+qtu,\; u'=u$. Under such a transformation $e_5$ is mapped to $e'_5=t'\gen t' + [p+q]tu\gen {x'}$. We can always choose $q$ so that $p+q=0$ so that $e_5=t\gen t$ in canonical form. This gives us the realization
\[
2\mathsf{A}_{2.2}\oplus\mathsf{A}_1=\langle \gen x, x\gen x + u\gen u, u\gen x, t\gen t - u\gen u, t\gen t\rangle.
\]
However, this algebra contains the subalgebra $\langle \gen x, u\gen x, x\gen x\rangle$ and as we have seen above this linearizes the evolution equation. Thus we have no non-linearizing realization in this case.

\medskip\noindent $\underline{2\mathsf{A}_{2.2}=\langle \gen t, t\gen t + x\gen x, \gen u, qx\gen x + u\gen u\rangle}:$ The commutation relations give $e_5=x\gen x$, so that the resulting algebra contains the operators $\gen t\; t\gen t$ and these are incompatible since they give $F=0$. Hence we have no realization in this case.

\medskip\noindent $\underline{2\mathsf{A}_{2.2}=\langle \gen x, t\gen t + x\gen x, \gen u, qt\gen x + u\gen u\rangle}:$ The commutation relations give $e_5=at\gen t + pt\gen x$ and $aq=0$. Obviously, $a\neq 0$ since otherwise we have $e_5=t\gen x$ and this together with $\gen x$ in the equation for $G$ gives $u_1=0$ which a contradiction. Thus we must have $q=0$ and we then we may take $a=1$ so that we have the algebra $\langle \gen x, t\gen t+x\gen x, \gen u, u\gen u, t\gen t + pt\gen x\rangle$. We have the residual equivalence group $\mathscr{E}(e_1, e_2, e_3, e_4): t'=kt,\; x'=x+\lambda t,\; u'=u$. Under such a transformation $e_5$ is mapped to $e'_5=t'\gen {t'} + [p+\lambda]t\gen {x'}$ and we may choose $\lambda=-p$ so that $e_5=t\gen t$ in canonical form. Thus we have the realization
\[
2\mathsf{A}_{2.2}\oplus\mathsf{A}_1=\langle \gen x, t\gen t+x\gen x, \gen u, u\gen u, t\gen t\rangle.
\]
This algebra is equivalent to
\[
2\mathsf{A}_{2.2}\oplus\mathsf{A}_1=\langle \gen x, x\gen x, \gen u, u\gen u, \gen t\rangle
\]
under the equivalence transformation $t'=\ln|t|,\; x'=x,\, u'=u$.

\medskip\noindent $\underline{2\mathsf{A}_{2.2}=\langle \gen x, t\gen t + x\gen x, \gen u, at\gen t + u\gen u\rangle}:$ The commutation relations give $e_5=t\gen t$ and then we have the realization
\[
2\mathsf{A}_{2.2}\oplus\mathsf{A}_1=\langle \gen x, t\gen t+x\gen x, \gen u, at\gen t + u\gen u, t\gen t\rangle,
\]
which is obviously the same as the algebra
\[
2\mathsf{A}_{2.2}\oplus\mathsf{A}_1=\langle \gen x, t\gen t+x\gen x, \gen u, u\gen u, t\gen t\rangle,
\]
and this is, in turn, equivalent to
\[
2\mathsf{A}_{2.2}\oplus\mathsf{A}_1=\langle \gen x, x\gen x, \gen u, u\gen u, \gen t\rangle
\]
as we have already seen.

\medskip\noindent $\underline{2\mathsf{A}_{2.2}=\langle \gen x, x\gen x + u\gen u, \gen t, t\gen t + qu\gen u\rangle}$ with $q\neq 0$: The commutation relations give $e_5=u\gen u$ and the resulting algebra will then contain the two operators $\gen t,\, t\gen t$ which give $F=0$ in the equation for $F$, which is a contradiction. Thus we have no realization in this case.

\medskip\noindent $\underline{2\mathsf{A}_{2.2}=\langle \gen x, x\gen x + u\gen u, \gen t, t\gen t + u\gen x\rangle}:$ the commutation relations give $e_5=u\gen x$. The resulting algebra then contains the operators $\gen t,\, t\gen t$ giving $F=0$, which is a contradiction. Thus we have no realization in this case.

\subsection{Non-linearizing realizations of $2\mathsf{A}_{2.2}\oplus\mathsf{A}_1$:}

\[
2\mathsf{A}_{2.2}\oplus\mathsf{A}_1=\langle \gen x, x\gen x, \gen u, u\gen u, \gen t\rangle.
\]

\section{Realizations of $\mathsf{A}_{3.k}\oplus\mathsf{A}_{2.2}$ for $k=3,\dots, 9$:} We have $\mathsf{A}_{3.k}\oplus\mathsf{A}_{2.2}=\langle e_1, e_2, e_3\rangle\oplus\langle e_4, e_5\rangle$ with $\mathsf{A}_{3.k}=\langle e_1, e_2, e_3\rangle$ and $\mathsf{A}_{2.2}=\langle e_1, e_2\rangle$. We note that $\langle e_1, e_2, e_3, e_4\rangle=\mathsf{A}_{3.k}\oplus\mathsf{A}_1$ as well as $\langle e_1, e_2, e_3, e_5\rangle=\mathsf{A}_{3.k}\oplus\mathsf{A}_1$. Thus we consider the admissible extensions $\mathsf{A}_{3.k}\oplus\mathsf{A}_1$ with $\mathsf{A}_1=\langle e_4\rangle$ or $\mathsf{A}_1=\langle e_5\rangle$. We can always assume that we $[e_4, e_5]=e_4$ in our realization of $\mathsf{A}_{2.2}$. Finally, we may use the residual equivalence group $\mathscr{E}(e_1, e_2, e_3)$ to give, first, a canonical form for $e_4$ and then to exploit the commutation relation $[e_4, e_5]=e_4$ together with the residual equivalence group $\mathscr{E}(e_1, e_2, e_3, e_4)$ in order to fix a canonical form for $e_5$. From this, it follows that we may always take $\mathsf{A}_{3.k}\oplus\mathsf{A}_1=\langle e_1, e_2, e_3, e_4\rangle$ and then proceed to find a canonical form for $e_5$.

\subsection{$\mathsf{A}_{3.3}\oplus\mathsf{A}_{2.2}$:} In the following we show that there are no non-linearizing realizations, or that there are no possible admissible realizations.

\medskip\noindent $\underline{\mathsf{A}_{3.3}\oplus\mathsf{A}_1=\langle \gen u, x\gen u, -(\gen t + \gen x), \gen t\rangle}:$ If $e_4=\gen t$ then $e_5=a\gen t + c(t-x)\gen u$ and so $c(t-x)\gen u$ belongs to the algebra. Hence the algebra will contain  the subalgebra $\langle\gen u, x\gen u, c(t-x)\gen u\rangle$ which is a rank-one three-dimensional abelian Lie algebra which linearizes the evolution equation. So we have no non-linearizing realization. The same is true if we take $e_5=\gen t$ in canonical form.

\medskip\noindent $\underline{\mathsf{A}_{3.3}\oplus\mathsf{A}_1=\langle \gen u, x\gen u, -\gen x, \gen t\rangle}:$ If $e_4=\gen t$ then $e_5=a(t)\gen t + c(t)\gen u$. Note that if $\dot{a}(t)=0$ then the algebra contains the operator $c(t)\gen u$ which together with $\gen u$ gives $\dot{c}(t)=0$ in the equation for $G$. This contravenes the requirements of dimension. Thus $\dot{a}(t)\neq 0$. Then $[e_4, e_5]=e_4$ gives $e_5=(t+l)\gen t + \gamma \gen u$. Thus the algebra contains $\gen t$ and $t\gen t$. These give $F=0$ in the equation for $F$. Thus we have no realization in this case. The same discussion holds if $e_5=\gen t$ in canonical form.

\medskip\noindent $\underline{\mathsf{A}_{3.3}\oplus\mathsf{A}_1=\langle \gen x, \gen t, [t+u]\gen x, \gen t - \gen u\rangle}:$ If $e_4=\gen t - \gen u$ then $e_5=ae_4 + b(u)\gen x$. Then the algebra will contain the rank-one three-dimensional abelian subalgebra $\langle \gen x, (t+u)\gen x, b(u)\gen x\rangle$ which linearizes the evolution equation. The same holds for $e_5=\gen t - \gen u$ in canonical form. We have no non-linearizing realization in this case.

\medskip\noindent $\underline{\mathsf{A}_{3.3}\oplus\mathsf{A}_1=\langle \gen x, \gen u, u\gen x, \gen t\rangle}:$ With $e_4=\gen t$ in canonical form, $e_5=a(t)\gen t + b(t)\gen x$. Then $[e_4, e_5]=e_4$ gives $e_5=t\gen t + le_4+\beta e_1$ so that the algebra contains $\gen t,\; t\gen t$ which give $F=0$ in the equation for $F$. Thus we have no realization in this case. Taking $e_5=\gen t$ in canonical form we have $e_4=a(t)\gen t + b(t)\gen x$ and $[e_4, e_5]=e_4$ gives $e_4=\alpha e^{-t}\gen t + \beta e^{-t}\gen x$ and this is equivalent to $t\gen t + t\gen x$ under the equivalence transformation $t'=e^{-t},\; x'=x,\, u'=u$, then $e_5$ becomes $t\gen t$ and the algebra then contains the operators $\gen x, t\gen x$ so that the equation for $G$ gives $u_1=0$, which is a contradiction.

\medskip\noindent $\underline{\mathsf{A}_{3.3}\oplus\mathsf{A}_1=\langle \gen x, \gen t, t\gen x-\gen u, \gen t-u\gen x + \kappa \gen u\rangle}:$ With $e_4=\gen t-u\gen x + \kappa \gen u$, we have $e_5=a(\gen t - u\gen x) + \beta e_1 + \gamma \gen u=ae_4+\beta e_1+(\gamma - a\kappa)\gen u$. Clearly then $\gen u$ is also a symmetry if $\gamma - a\kappa\neq 0$, from which it follows that $\gen x, t\gen x$ are symmetries. However, they give the contradiction $u_1=0$ in the equation for $G$. Thus $\gamma - a\kappa=0$ and so $e_5=ae_4+\beta e_1$ which contradicts the requirements of dimension. the same is true if we take $e_5=\gen t-u\gen x + \kappa \gen u$ in canonical form.

\medskip\noindent $\underline{\mathsf{A}_{3.3}\oplus\mathsf{A}_1=\langle \gen x, \gen u, t\gen t+ u\gen x, t\gen t\rangle}:$ If $e_4=t\gen t$ then $e_5=ae_4+\gamma e_2 + \beta e_1 + \gamma\ln|t|\gen x$. If $\gamma\neq 0$ then $\ln|t|\gen x$ is a symmetry together with $\gen x$ and these two give the contradiction $u_1=0$ in the equation for $G$. The same occurs if we take $e_5=t\gen t$ in canonical form. Thus we have no realization in this case.

\subsection{$\mathsf{A}_{3.4}\oplus\mathsf{A}_{2.2}$:}

\medskip\noindent $\underline{\mathsf{A}_{3.4}\oplus\mathsf{A}_1=\langle \gen x, u\gen x, x\gen x- \gen u, \gen t\rangle}:$ With $e_4=\gen t$ we have $e_5=a(t)\gen t + b(t)e^{-u}\gen x$. We cannot have $a(t)=0$ for then we would have $b(t)e^{-u}\gen x$ as a symmetry and then the algebra would contain the subalgebra $\langle \gen x, u\gen x, b(t)e^{-u}\gen x\rangle$ which linearizes the equation. From $[e_4, e_5]=e_4$ we have $e_5=(t+l)\gen t + \kappa e^{-u}\gen x$. We note that $\kappa\neq 0$ for otherwise $t\gen t$ is a symmetry, together with $\gen t$ and they give $F=0$, which is a contradiction. So we find that $e_5=t\gen t + \kappa e^{-u}\gen x$ with $\kappa\neq 0$. Putting these operators as symmetries into the equation for $F$ gives $F=0$. Taking $e_5=\gen t$ in canonical form, we have $e_4=-(t\gen t + \kappa e^{-u}\gen x)$ and the same result occurs. Thus we have no realization in this case.

\medskip\noindent $\underline{\mathsf{A}_{3.4}\oplus\mathsf{A}_1=\langle \gen x, u\gen x, t\gen t + x\gen x- \gen u, t\gen t\rangle}:$ With $e_4=t\gen t$ in canonical form, $e_5=at\gen t + tb(\tau)\gen x$ where $\tau=te^{u}$, and $[e_4, e_5]=t[b(\tau)+\tau b'(\tau)]\gen x$ which is also to be a symmetry. Thus the Lie algebra will contain a rank-one three-dimensional abelian subalgebra of the form $\langle \gen x, u\gen x, tb(\tau)\gen x\rangle$ which linearizes the evolution equation or gives a contradiction. Thus we have no non-linearizing realization in this case.

\medskip\noindent $\underline{\mathsf{A}_{3.4}\oplus\mathsf{A}_1=\langle \gen x, \gen t, t\gen t + [t+x]\gen x, \gen u\rangle}:$ If $e_4=\gen u$ then we find $e_5=c(u)\gen u$ and $[e_4, e_5]=e_4$ yields $e_5=u\gen u$ in canonical form. If $e_5=\gen u$ then $e_4=c(u)\gen u$ and $[e_4, e_5]=e_4$ gives $e_4=-u\gen u$ in canonical form. Thus we have the realization
\[
\mathsf{A}_{3.4}\oplus\mathsf{A}_{2.2}=\langle \gen x, \gen t, t\gen t + [t+x]\gen x, \gen u, u\gen u\rangle.
\]

\medskip\noindent $\underline{\mathsf{A}_{3.4}\oplus\mathsf{A}_1=\langle \gen x, \gen t, t\gen t + [t+x]\gen x + u\gen u, u\gen u\rangle}:$ If $e_4=u\gen u$ then $e_5=ae_2+ce_4 + bu\gen x$ and $[e_4, e_5]=e_4$ is impossible. If $e_5=u\gen u$ then $e_4=ae_2+ce_4+bu\gen x$ and $[e_4, e_5]=-bu\gen x$. Thus we may put $e_4=u\gen x,\; e_5=-u\gen u$ and we obtain the realization
\[
\mathsf{A}_{3.4}\oplus\mathsf{A}_{2.2}=\langle \gen x, \gen t, t\gen t + [t+x]\gen x + u\gen u, u\gen x, -u\gen u\rangle.
\]

\medskip\noindent $\underline{\mathsf{A}_{3.4}\oplus\mathsf{A}_1=\langle \gen x, \gen t, t\gen t + [t+x]\gen x + u\gen u, u\gen x\rangle}:$ With $e_4=u\gen x$ we have $e_5=ae_2 + be_4+cu\gen u$ and then $[e_4, e_5]=-ce_4$ so that $e_5=-u\gen u$ in canonical form. If $e_5=u\gen x$ then we are unable to implement $[e_4, e_5]=e_4$. Thus we have the same realization as in the previous case.

\medskip\noindent $\underline{\mathsf{A}_{3.4}\oplus\mathsf{A}_1=\langle \gen x, \gen u, t\gen t + [t+x]\gen x + u\gen u, t\gen t\rangle}:$ With $e_4=t\gen t$ we find $e_5=at\gen t + [qt\ln|t|+\kappa t]\gen x + qt\gen u$. Consequently, $[qt\ln|t|+\kappa t]\gen x + qt\gen u$ would also be a symmetry, and this, together with $\gen x,\; \gen u$ as symmetries in the equation for $G$, gives $[q\ln|t| + q+\kappa]u_1-q=0$ which is a contradiction unless $q=\kappa=0$. Then we have $e_5=ae_4$ and this is not allowed because of the requirements of dimension. Hence we have no realization in this case.

We therefore have two non-linearizing realizations:
\begin{align*}
\mathsf{A}_{3.4}\oplus\mathsf{A}_{2.2}&=\langle \gen x, \gen t, t\gen t + [t+x]\gen x, \gen u, u\gen u\rangle\\
\mathsf{A}_{3.4}\oplus\mathsf{A}_{2.2}&=\langle \gen x, \gen t, t\gen t + [t+x]\gen x + u\gen u, u\gen x, -u\gen u\rangle.
\end{align*}

\subsection{Realizations of $\mathsf{A}_{3.5}\oplus\mathsf{A}_{2.2}$:}

\medskip\noindent $\underline{\mathsf{A}_{3.5}\oplus\mathsf{A}_1=\langle \gen t, \gen x, t\gen t + x\gen x, \gen u\rangle}:$ With $e_4=\gen u$ then $e_5=c(u)\gen u$ and we then find $e_5=u\gen u$ in canonical form if $[e_4, e_5]=e_4$. If $e_5=\gen u$ then we have $e_4=c(u)\gen u$ and then $[e_4, e_5]=e_4$ gives $e_4=-u\gen u$. Thus we have the realization
\[
\mathsf{A}_{3.5}\oplus\mathsf{A}_{2.2}=\langle \gen t, \gen x, t\gen t + x\gen x, \gen u, u\gen u\rangle.
\]

\medskip\noindent $\underline{\mathsf{A}_{3.5}\oplus\mathsf{A}_1=\langle \gen x, \gen u, x\gen x + u\gen u, \gen t\rangle}:$ If $e_4=\gen t$ then $e_5=a(t)\gen t$ so that  $[e_4, e_5]=e_4$ gives $e_5=t\gen t$ in canonical form. Thus we have $\gen t$ and $t\gen t$ as symmetries, and they give $F=0$ in the equation for $F$. If $e_5=\gen t$ then $e_4=a(t)\gen t$ and we again obtain $\gen t$ and $t\gen t$ as symmetries, giving $F=0$. Thus we have no realization in this case.

\medskip\noindent $\underline{\mathsf{A}_{3.5}\oplus\mathsf{A}_1=\langle \gen t, \gen x, t\gen t + x\gen x + u\gen u, u\gen x\rangle}:$ If $e_4=u\gen x$ then $e_5=pu\gen x + qu\gen u$. Then $[e_4, e_5]=-qe_4$ so we have $e_5=-u\gen u$ in canonical form if $[e_4, e_5]=e_4$. If $e_5=u\gen x$ then we have $e_4=pu\gen x + qu\gen u$ and we can only have $[e_4, e_5]=qe_5$. We obtain the realization
\[
\mathsf{A}_{3.5}\oplus\mathsf{A}_{2.2}=\langle \gen t, \gen x, t\gen t + x\gen x + u\gen u, u\gen x, -u\gen u\rangle.
\]

\medskip\noindent $\underline{\mathsf{A}_{3.5}\oplus\mathsf{A}_1=\langle \gen t, \gen x, t\gen t + x\gen x + u\gen u, u\gen u\rangle}:$ A similar discussion to the previous case leads to the realization

\[
\mathsf{A}_{3.5}\oplus\mathsf{A}_{2.2}=\langle \gen t, \gen x, t\gen t + x\gen x + u\gen u, u\gen x, -u\gen u\rangle,
\]
so we have no new realization.

We then have the following two non-linearizing realizations:
\begin{align*}
\mathsf{A}_{3.5}\oplus\mathsf{A}_{2.2}&=\langle \gen t, \gen x, t\gen t + x\gen x, \gen u, u\gen u\rangle\\
\mathsf{A}_{3.5}\oplus\mathsf{A}_{2.2}&=\langle \gen t, \gen x, t\gen t + x\gen x + u\gen u, u\gen x, -u\gen u\rangle.
\end{align*}

\subsection{Realizations of $\mathsf{A}_{3.6}\oplus\mathsf{A}_{2.2}$:}

\medskip\noindent $\underline{\mathsf{A}_{3.6}\oplus\mathsf{A}_1=\langle \gen u, x\gen u, 2x\gen x + u\gen u, \gen t\rangle}:$ If $e_4=\gen t$ then $e_5=a(t)\gen t + c(t)x^{1/2}\gen u$. We want $a(t)\neq 0$ for otherwise the algebra will contain $\langle \gen u, x\gen u, c(t)x^{1/2}\gen u\rangle$ which is a rank-one three-dimensional abelian Lie algebra which linearizes the evolution equation. In this case, $[e_4, e_5]=e_4$ gives $e_5=t\gen t + \kappa x^{1/2}\gen u$ in canonical form. Note that now we require $\kappa\neq 0$ for otherwise we would have $\gen t,\; t\gen t$ as symmetries, giving $F=0$ in the equation for $F$, which is a contradiction. Then we have the realization $\langle \gen u, x\gen u, 2x\gen x + u\gen u, \gen t, t\gen t + \kappa x^{1/2}\gen u\rangle$ which gives $F=0$ in the equation for $F$. If on the other hand we have $e_5=\gen t$ then $e_4=a(t)\gen t + c(t)x^{1/2}\gen u$, then $[e_4, e_5]=e_4$ gives $e_4=ae^{-t}\gen t + \kappa e^{-t}x^ {1/2}\gen u$ and we may take $a=1$ since we require $a\neq 0$ and we want $\kappa\neq 0$ for otherwise the operators $\gen t, e^{-t}\gen t$ give $F=0$ in the equation for $F$. Again, the resulting realization gives $F=0$. Thus we have no realization in this case.

\medskip\noindent $\underline{\mathsf{A}_{3.6}\oplus\mathsf{A}_1=\langle \gen u, x\gen u, t\gen t + 2x\gen x + u\gen u, t\gen t\rangle}:$ if $e_4=t\gen t$ then $e_5=at\gen t + tc(\sigma)\gen u$ with $\sigma=x/t^2$. Consequently, the algebra contains the operator $tc(\sigma)\gen u$ and hence it contains the rank-one three-dimensional abelian Lie algebra
$$\langle \gen u, x\gen u, tc(\sigma)\gen u\rangle$$ which linearizes the evolution equation. The same holds for $e_5=t\gen t$ with $e_4=at\gen t + tc(\sigma)\gen u$. So we have no realization in this case.

\medskip\noindent $\underline{\mathsf{A}_{3.6}\oplus\mathsf{A}_1=\langle \gen t, \gen x, t\gen t - x\gen x, \gen u\rangle}:$ With $e_4=\gen u$ we have $e_5=c(u)\gen u$ and then $[e_4, e_5]=e_4$ gives $e_5=u\gen u$ in canonical form. If $e_5=\gen u$ then $[e_4, e_5]=e_4$ gives $e_4=e^{-u}\gen u$. Then with the equivalence transformation $t'=t,\; x'=x,\, u'=e^{u}$ we find that $e'_4=\gen {u'},\; e'_5=u'\gen {u'}$ so we have only one inequivalent realization:
\[
\mathsf{A}_{3.6}\oplus\mathsf{A}_{2.2}=\langle \gen t, \gen x, t\gen t - x\gen x, \gen u, u\gen u\rangle.
\]

\medskip\noindent $\underline{\mathsf{A}_{3.6}\oplus\mathsf{A}_1=\langle \gen t, \gen x, t\gen t - x\gen x-u\gen u, u\gen x\rangle}:$ With $e_4=u\gen x$ we have $e_5=pu\gen x+qu\gen u$. Then we have $[e_4, pu\gen x+qu\gen u]=-qu\gen x$ so we take $e_4=u\gen x$ and $e_5=-u\gen u$ and then we have $[e_4, e_5]=e_4$ so we have the realization:
\[
\mathsf{A}_{3.6}\oplus\mathsf{A}_{2.2}=\langle \gen t, \gen x, t\gen t - x\gen x-u\gen u, u\gen x, -u\gen u\rangle.
\]

\medskip\noindent $\underline{\mathsf{A}_{3.6}\oplus\mathsf{A}_1=\langle \gen t, \gen x, t\gen t - x\gen x-u\gen u, u\gen u\rangle}:$ With $e_4=u\gen u$ we have $e_5=pu\gen x+qu\gen u$. Then $[u\gen u, pu\gen x+qu\gen u]=pu\gen x$ so we relabel and put $e_4=u\gen x$ and $e_5=-u\gen u$ and then we have $[e_4, e_5]=e_4$ so we have only one inequivalent realization:
\[
\mathsf{A}_{3.6}\oplus\mathsf{A}_{2.2}=\langle \gen t, \gen x, t\gen t - x\gen x-u\gen u, u\gen x, -u\gen u\rangle.
\]

\medskip\noindent $\underline{\mathsf{A}_{3.6}\oplus\mathsf{A}_1=\langle \gen x, \gen u, t\gen t + x\gen x-u\gen u, t\gen t\rangle}:$ With $e_4=t\gen t$ we have $e_5=at\gen t + pt\gen x+qt\gen u$. It follows that the algebra will contain the operator $pt\gen x+qt\gen u$ which, together with $\gen x,\; \gen u$ in the equation for $G$ give $pu_1-q=0$ and this can only be true if $p=q=0$. Hence $e_5=ae_4$ and this contravenes the requirements of dimension. Hence we have no realization in this case.

We then have the following non-linearizing realizations:
\begin{align*}
\mathsf{A}_{3.6}\oplus\mathsf{A}_{2.2}&=\langle \gen t, \gen x, t\gen t - x\gen x, \gen u, u\gen u\rangle\\
\mathsf{A}_{3.6}\oplus\mathsf{A}_{2.2}&=\langle \gen t, \gen x, t\gen t - x\gen x-u\gen u, u\gen x, -u\gen u\rangle.
\end{align*}

\subsection{Realizations of $\mathsf{A}_{3.7}\oplus\mathsf{A}_{2.2}$:}

\medskip\noindent $\underline{\mathsf{A}_{3.7}\oplus\mathsf{A}_1=\langle \gen u, x\gen u, (1-q)x\gen x + u\gen u, \gen t\rangle}:$ We have $e_5=a(t)\gen t + c(t)x^{1/(1-q)}\gen u$. Then the relation $[e_4, e_5]=e_4$ gives $e_5=t\gen t + \kappa x^{1/(1-q)}\gen u$ in canonical form. We require $\kappa\neq 0$ for otherwise the algebra will contain the operators $\gen t,\; t\gen t$ which give $F=0$ in the equation for $F$, a contradiction. We have the realization:
\[
\mathsf{A}_{3.7}\oplus\mathsf{A}_{2.2}=\langle \gen u, x\gen u, (1-q)x\gen x + u\gen u, \gen t, t\gen t + \kappa x^{1/(1-q)}\gen u\rangle.
\]

\medskip\noindent $\underline{\mathsf{A}_{3.7}\oplus\mathsf{A}_1=\langle \gen x, x\gen u, t\gen t + (1-q)x\gen x + u\gen u, t\gen t\rangle}:$ In this case $e_5=at\gen t +ct\gen u$. From this we see that the algebra will contain the subalgebra $\langle \gen u, x\gen u, ct\gen u\rangle$ which linearizes the evolution equation since it is a rank-one three-dimensional abelian Lie algebra. So we have no non-linearizing realization in this case.

\medskip\noindent $\underline{\mathsf{A}_{3.7}\oplus\mathsf{A}_1=\langle \gen t, \gen x, t\gen t + qx\gen x, \gen u\rangle}:$ We have $e_5=c(u)\gen u$ and then $[e_4, e_5]=e_4$ gives $e_5=u\gen u$ in canonical form. Thus we have the realization
\[
\mathsf{A}_{3.7}\oplus\mathsf{A}_{2.2}=\langle \gen t, \gen x, t\gen t + qx\gen x, \gen u, u\gen u\rangle.
\]

\medskip\noindent $\underline{\mathsf{A}_{3.7}\oplus\mathsf{A}_1=\langle \gen x, \gen t, qt\gen t + x\gen x, \gen u\rangle}:$ We have $e_5=c(u)\gen u$ and then $[e_4, e_5]=e_4$ gives $e_5=u\gen u$ in canonical form. Thus we have the realization
\[
\mathsf{A}_{3.7}\oplus\mathsf{A}_{2.2}=\langle \gen x, \gen t, qt\gen t + x\gen x, \gen u, u\gen u\rangle.
\]

\medskip\noindent $\underline{\mathsf{A}_{3.7}\oplus\mathsf{A}_1=\langle \gen x, \gen u, qx\gen x + u\gen u, \gen t\rangle}:$ We have $e_5=a(t)\gen t$ and then $[e_4, e_5]=e_4$ gives $e_5=t\gen t$ in canonical form. The operators $\gen t,\; t\gen t$ give $F=0$ in the equation for $F$, so we have no realization in this case.

\medskip\noindent $\underline{\mathsf{A}_{3.7}\oplus\mathsf{A}_1=\langle \gen t, \gen x, t\gen t + qx\gen x + u\gen u, u\gen u\rangle}:$ We have $e_5=\alpha u^q\gen x + \beta u\gen u$ so that we need $\alpha\neq 0$ on grounds of dimension, and then $[e_4, e_5]=\alpha qu^q\gen x\neq e_4$ . Thus we have no realization in this case.

\medskip\noindent $\underline{\mathsf{A}_{3.7}\oplus\mathsf{A}_1=\langle \gen t, \gen x, t\gen t + qx\gen x + u\gen u, u^q\gen x\rangle}:$ We have $e_5=\alpha u^q\gen x + \beta u\gen u$ so that we need $\beta\neq 0$ on grounds of dimension, and then $[e_4, e_5]=-\beta qu^q\gen x=e_4$ if $\beta=-1/q$. Thus we have the realization:
\[
\mathsf{A}_{3.7}\oplus\mathsf{A}_{2.2}=\langle \gen t, \gen x, t\gen t + qx\gen x + u\gen u, u^q\gen x, -\frac{1}{q}u\gen u\rangle.
\]

\medskip\noindent $\underline{\mathsf{A}_{3.7}\oplus\mathsf{A}_1=\langle \gen x, \gen t, qt\gen t + x\gen x + u\gen u, u\gen x\rangle}:$ We have $e_5=\alpha u\gen x + \beta u\gen u$ so that we need $\beta\neq 0$ on grounds of dimension, and then $[e_4, e_5]=-\beta u\gen x=e_4$ if $\beta=-1$. Thus we have the realization:
\[
\mathsf{A}_{3.7}\oplus\mathsf{A}_{2.2}=\langle \gen x, \gen t, qt\gen t + x\gen x + u\gen u, u\gen x, -u\gen u\rangle.
\]

\medskip\noindent $\underline{\mathsf{A}_{3.7}\oplus\mathsf{A}_1=\langle \gen x, \gen t, qt\gen t + x\gen x + u\gen u, u\gen u\rangle}:$ We have $e_5=\alpha u\gen x + \beta u\gen u$ and then $[e_4, e_5]=-\beta u\gen x=e_4$ is not possible to implement. Thus we have no realization in this case.

\medskip\noindent $\underline{\mathsf{A}_{3.7}\oplus\mathsf{A}_1=\langle \gen x, \gen u, t\gen t + x\gen x + qu\gen u, t\gen t\rangle}:$ We have $e_5=\alpha t\gen t $ and is not possible on grounds of dimension. Thus we have no realization in this case.

\subsection{Realizations of $\mathsf{A}_{3.8}\oplus\mathsf{A}_{2.2}$:}

\medskip\noindent $\underline{\mathsf{A}_{3.8}\oplus\mathsf{A}_1=\langle \gen u, x\gen u, -(1+x^2)\gen x - xu\gen u, \gen t\rangle}:$ we find $e_5=a(t)\gen t + c(t)\sqrt{1+x^2}\gen u$. We require $a(t)\neq 0$ for otherwise the algebra would contain $\langle \gen u, x\gen u, c(t)\sqrt{1+x^2}\gen u\rangle$ which linearizes the evolution equation because it is a rank-one three-dimensional abelian Lie algebra. Then $[e_4, e_5]=e_4$ gives $e_5=t\gen t + \kappa\sqrt{1+x^2}\gen u$ and we also require $\kappa\neq 0$ for otherwise the algebra would contain the operators $\gen t,\; t\gen t$ which give $F=0$ in the equation for $F$. We have the realization
\[
\mathsf{A}_{3.8}\oplus\mathsf{A}_{2.2}=\langle \gen u, x\gen u, -(1+x^2)\gen x - xu\gen u, \gen t, t\gen t + \kappa\sqrt{1+x^2}\gen u\rangle.
\]

\medskip\noindent $\underline{\mathsf{A}_{3.8}\oplus\mathsf{A}_1=\langle \gen u, x\gen u, -t\gen t -(1+x^2)\gen x - xu\gen u, t\gen t\rangle}:$ We have $$e_5=at\gen t + c(\sigma)\sqrt{1+x^2}\gen u$$ with $\sigma=t\exp(-\arctan x)$. From this it follows that the algebra contains the subalgebra $\langle \gen u, x\gen u, c(\sigma)\sqrt{1+x^2}\gen u\rangle$ which linearizes the evolution equation because it is a rank-one three-dimensional abelian Lie algebra. Thus we have no non-linearizing realization in this case.

\medskip\noindent $\underline{\mathsf{A}_{3.8}\oplus\mathsf{A}_1=\langle \gen x, \gen u, u\gen x - x\gen u, \gen t\rangle}:$ We find $e_5=a(t)\gen t$ and then $[e_4, e_5]=e_4$ gives $e_5=t\gen t$ in canonical form. However, $e_4, e_5$ then give $F=0$ in the equation for $F$. We have no realization in this case.

\medskip\noindent $\underline{\mathsf{A}_{3.8}\oplus\mathsf{A}_1=\langle \gen x, \gen u, t\gen t + u\gen x - x\gen u, \gen t\rangle}:$ We find $e_5=at\gen t + b(t)\gen x + c(t)\gen u$ with $b(t)=\alpha\cos(\ln|t|) + \beta\sin(\ln|t|)$ and $c(t)=-\alpha\sin(\ln|t|) + \beta\cos(\ln|t|)$. Thus the algebra contains the operator $b(t)\gen x + c(t)\gen u$ and this together with $\gen x,\; \gen u$ give $\dot{b}(t)u_1-\dot{c}(t)=0$ in the equation for $G$, which is impossible unless $b(t)=c(t)=0$. Thus we have no realization in this case.

\subsection{Realizations of $\mathsf{A}_{3.9}\oplus\mathsf{A}_{2.2}$:}

\medskip\noindent $\underline{\mathsf{A}_{3.9}\oplus\mathsf{A}_1=\langle \gen u, x\gen u, -(1+x^2)\gen x + (q-x)u\gen u, \gen t\rangle}:$ We have $e_5=a(t)\gen t + c(t)\sqrt{1+x^2}\exp(-q\arctan x)\gen u$. Then $[e_4, e_5]=e_4$ gives
$$e_5=t\gen t + \kappa \sqrt{1+x^2}\exp(-q\arctan x)\gen u$$
in canonical form. We require $\kappa\neq 0$ since otherwise we have $\gen t$ and $t\gen t$ as symmetries and they give $F=0$ in the equation for $F$. Putting the operators of this algebra into the equation for $F$ gives $F=0$ so we have no realization in this case.

\medskip\noindent $\underline{\mathsf{A}_{3.9}\oplus\mathsf{A}_1=\langle \gen u, x\gen u, -t\gen t -(1+x^2)\gen x + (q-x)u\gen u, t\gen t\rangle}:$ In this case we have $e_5=at\gen t + c(\sigma)\sqrt{1+x^2}\exp(-q\arctan x)\gen u$ with $\sigma=t\exp(-q\arctan x)$. We see that the algebra will contain the subalgebra
$$\langle \gen u, x\gen u, c(\sigma)\sqrt{1+x^2}\exp(-q\arctan x)\rangle$$
which linearizes the evolution equation since it is a rank-one three-dimensional abelian Lie algebra. So we have no non-linearizing realization in this case.

\medskip\noindent $\underline{\mathsf{A}_{3.9}\oplus\mathsf{A}_1=\langle \gen x, \gen u, [qx+u]\gen x + (qu-x)u\gen u, \gen t\rangle}:$ We have $e_5=a(t)\gen t$ and this together with $[e_4, e_5]=0$ gives $F=0$ in the equation for $F$. Thus we have no realization in this case.

\medskip\noindent $\underline{\mathsf{A}_{3.9}\oplus\mathsf{A}_1=\langle \gen x, \gen u, t\gen t + [qx+u]\gen x + (qu-x)u\gen u, t\gen t\rangle}:$ In this case we find $e_5=at\gen t + b(t)\gen x +c(t)\gen u$ with $b(t)=\alpha\cos(\ln|t|) + \beta\sin(\ln|t|)$ and $c(t)=-\alpha\sin(\ln|t|) + \beta\cos(\ln|t|)$. Consequently the algebra contains the operator $b(t)\gen x +c(t)\gen u$ which, together with $\gen x,\, \gen u$ in the equation for $G$ gives $\dot{b}(t)u_1-\dot{c}(t)=0$. This is a contradiction unless $b(t)=c(t)=0$ and then $e_5=ae_4$ which contradicts the requirement of dimension. Thus we have no realization in this case.

\subsection{Non-linearizing realizations of $\mathsf{A}_{3.k}\oplus\mathsf{A}_{2.2}$:}

\medskip\noindent $\underline{\mathsf{A}_{3.4}\oplus\mathsf{A}_{2.2}}$:

\begin{align*}
\mathsf{A}_{3.4}\oplus\mathsf{A}_{2.2}&=\langle \gen x, \gen t, t\gen t + [t+x]\gen x, \gen u, u\gen u\rangle\\
\mathsf{A}_{3.4}\oplus\mathsf{A}_{2.2}&=\langle \gen x, \gen t, t\gen t + [t+x]\gen x + u\gen u, u\gen x, -u\gen u\rangle.
\end{align*}

\medskip\noindent $\underline{\mathsf{A}_{3.5}\oplus\mathsf{A}_{2.2}}$:
\begin{align*}
\mathsf{A}_{3.5}\oplus\mathsf{A}_{2.2}&=\langle \gen t, \gen x, t\gen t + x\gen x, \gen u, u\gen u\rangle\\
\mathsf{A}_{3.5}\oplus\mathsf{A}_{2.2}&=\langle \gen t, \gen x, t\gen t + x\gen x + u\gen u, u\gen x, -u\gen u\rangle.
\end{align*}

\medskip\noindent $\underline{\mathsf{A}_{3.6}\oplus\mathsf{A}_{2.2}}$:

\begin{align*}
\mathsf{A}_{3.6}\oplus\mathsf{A}_{2.2}&=\langle \gen t, \gen x, t\gen t - x\gen x, \gen u, u\gen u\rangle\\
\mathsf{A}_{3.6}\oplus\mathsf{A}_{2.2}&=\langle \gen t, \gen x, t\gen t - x\gen x-u\gen u, u\gen x, -u\gen u\rangle.
\end{align*}

\medskip\noindent $\underline{\mathsf{A}_{3.7}\oplus\mathsf{A}_{2.2}}$:

\begin{align*}
\mathsf{A}_{3.7}\oplus\mathsf{A}_{2.2}&=\langle \gen u, x\gen u, (1-q)x\gen x + u\gen u, \gen t, t\gen t + \kappa x^{1/(1-q)}\gen u\rangle\\
\mathsf{A}_{3.7}\oplus\mathsf{A}_{2.2}&=\langle \gen t, \gen x, t\gen t + qx\gen x, \gen u, u\gen u\rangle\\
\mathsf{A}_{3.7}\oplus\mathsf{A}_{2.2}&=\langle \gen x, \gen t, qt\gen t + x\gen x, \gen u, u\gen u\rangle\\
\mathsf{A}_{3.7}\oplus\mathsf{A}_{2.2}&=\langle \gen t, \gen x, t\gen t + qx\gen x + u\gen u, u^q\gen x, -\frac{1}{q}u\gen u\rangle\\
\mathsf{A}_{3.7}\oplus\mathsf{A}_{2.2}&=\langle \gen x, \gen t, qt\gen t + x\gen x + u\gen u, u\gen x, -u\gen u\rangle.
\end{align*}

\medskip\noindent $\underline{\mathsf{A}_{3.8}\oplus\mathsf{A}_{2.2}}$:

\begin{align*}
\mathsf{A}_{3.8}\oplus\mathsf{A}_{2.2}&=\langle \gen u, x\gen u, -(1+x^2)\gen x - xu\gen u, \gen t, t\gen t + \kappa\sqrt{1+x^2}\gen u\rangle.
\end{align*}

\section{Realizations of $\mathsf{A}_{4.k}\oplus\mathsf{A}_1$:}

\medskip\noindent $\underline{\mathsf{A}_{4.k}\oplus\mathsf{A}_1,\quad k=1,\dots 6:}$ The algebras $\mathsf{A}_{4.k}$, with $k=1,\dots, 6$ all contain a three-dimensional abelian ideal, so that in these cases $\mathsf{A}_{4.k}\oplus\mathsf{A}_1$ contains a four-dimensional abelian Lie subalgebra and so will either linearize the evolution equation or not be admissible. Hence we need consider only $\mathsf{A}_{4.7}\oplus\mathsf{A}_1,\; \mathsf{A}_{4.8}\oplus\mathsf{A}_1,\; \mathsf{A}_{4.9}\oplus\mathsf{A}_1,\; \mathsf{A}_{4.10}\oplus\mathsf{A}_1$

\medskip\noindent $\underline{\mathsf{A}_{4.7}\oplus\mathsf{A}_1:}$

\medskip\noindent $\underline{\mathsf{A}_{4.7}=\langle \gen u, x\gen u, -\gen x, x\gen x + \left[2u-\frac{x^2}{2}\right]\gen u \rangle}:$ in this case $e_5=a(t)\gen t$. The residual equivalence group is $\mathscr{E}(e_1, e_2, e_3, e_4): t'=T(t),\; x'=x,\; u'=u$ and under such a transformation, $e_5$ is mapped to $e'_5=a(t)\dot{T}(t)\gen {t'}$. Choose $T(t)$ so that $a(t)\dot{T}=1$ and then $e_5=\gen t$ in canonical form. We have the realization
\[
\mathsf{A}_{4.7}\oplus\mathsf{A}_1= \langle \gen u, x\gen u, -\gen x, x\gen x + \left[2u-\frac{x^2}{2}\right]\gen u, \gen t \rangle.
\]

\medskip\noindent $\underline{\mathsf{A}_{4.7}=\langle \gen u, x\gen u, -\gen x, t\gen t + x\gen x + \left[2u-\frac{x^2}{2}\right]\gen u \rangle}:$ in this case $e_5=at\gen t + qt^2\gen t$. We must have $a\neq 0$ for otherwise $t^2\gen u$ will belong to the algebra as well as $\gen u$ and these two give the contradiction $t=0$ in the equation for $G$. So we may assume $a=1$. The residual equivalence group is $\mathscr{E}(e_1, e_2, e_3, e_4): t'=kt,\; x'=x,\; u'=u+\kappa t^2$ and under such a transformation, $e_5=\gen t + qt^2\gen t$ is mapped to $e'_5=t'\gen {t'} + [q+2\kappa]t^2\gen {u'}$. Choose $\kappa$ so that $q+2\kappa=0$ and then $e_5=t\gen t$ in canonical form. We have the realization
\[
\mathsf{A}_{4.7}\oplus\mathsf{A}_1= \langle \gen u, x\gen u, -\gen x, t\gen t + x\gen x + \left[2u-\frac{x^2}{2}\right]\gen u, t\gen t \rangle.
\]

\medskip\noindent $\underline{\mathsf{A}_{4.7}=\langle \gen u, x\gen u, -(\gen t + \gen x), t\gen t + x\gen x + \left[2u-\frac{x^2}{2}\right]\gen u \rangle}:$ We find that $e_5=(t-x)^2\gen u$ so that $\mathsf{A}_{4.7}\oplus\mathsf{A}_1$ contains the subalgebra $\langle \gen u, x\gen u, (t-x)^2\gen u\rangle$ which linearizes the evolution equation since it is a rank-one three-dimensional abelian Lie algebra. So we have no non-linearizing realization in this case.

\medskip\noindent $\underline{\mathsf{A}_{4.7}=\langle \gen x, \gen u, u\gen x + t\gen u, -\gen t + 2x\gen x + u\gen u \rangle}:$ We find that $e_5=e^{-2t}\gen x$ which together with the operator $\gen x$ gives the contradiction $u_1=0$ in the equation for $G$. So we have no realization in this case.

\subsection{Non-linearizing realizations of $\mathsf{A}_{4.7}\oplus\mathsf{A}_1$:}

\begin{align*}
\mathsf{A}_{4.7}\oplus\mathsf{A}_1&= \langle \gen u, x\gen u, -\gen x, x\gen x + \left[2u-\frac{x^2}{2}\right]\gen u, \gen t \rangle\\
\mathsf{A}_{4.7}\oplus\mathsf{A}_1&= \langle \gen u, x\gen u, -\gen x, t\gen t + x\gen x + \left[2u-\frac{x^2}{2}\right]\gen u, t\gen t \rangle.
\end{align*}

\medskip\noindent $\underline{\mathsf{A}_{4.8}\oplus\mathsf{A}_1:}$

\medskip\noindent $\underline{\mathsf{A}_{4.8}=\langle \gen u, x\gen u, -\gen x, t\gen t + qx\gen x + (1+q)u\gen u\rangle}:$ We find $e_5=at\gen t + \kappa t^{1+q}\gen u$. For $q=-1$ we may take $e_5=t\gen t$ in canonical form since then $e_5=at\gen t + \kappa \gen u$ and because $\gen u\in \mathsf{A}_{4.8}$. For $-1<q\leq 1$ we must have $a\neq 0$ for otherwise we would have $t^{1+q}\gen u$ and $\gen u$ belonging to the algebra and these give the contradiction $t^q=0$ in the equation for $G$. Hence we may assume $a=1$ without loss of generality. We have the residual equivalence group $\mathscr{E}(e_1, e_2, e_3, e_4): t'=kt,\; x'=x,\; u'=u+\lambda t^{1+q}$. Under such a transformation $e_5$ is mapped to $e'_5=t'\gen {t'} + [\kappa + (1+q)\lambda]t^{q+1}\gen {u'}$, and we see that we may choose $\lambda$ so that $\kappa + (1+q)\lambda=0$ giving $e_5=t\gen t$ in canonical form. Thus we have the realization
\[
\mathsf{A}_{4.8}\oplus\mathsf{A}_1=\langle \gen u, x\gen u, -\gen x, t\gen t + qx\gen x + (1+q)u\gen u, t\gen t\rangle.
\]

\medskip\noindent $\underline{\mathsf{A}_{4.8}=\langle \gen u, x\gen u, -\gen x, qx\gen x + (1+q)u\gen u\rangle,\; -1< q\leq 1}:$  We find $e_5=a(t)\gen t $. We have the residual equivalence group $\mathscr{E}(e_1, e_2, e_3, e_4): t'=T(t),\; x'=x,\; u'=u$. Under such a transformation $e_5$ is mapped to $e'_5=a(t)\dot{T}(t)\gen {t'}$, and we see that we may choose $T(t)$ so that $a(t)\dot{T}(t)=1$ giving $e_5=\gen t$ in canonical form. Thus we have the realization
\[
\mathsf{A}_{4.8}\oplus\mathsf{A}_1=\langle \gen u, x\gen u, -\gen x, qx\gen x + (1+q)u\gen u, \gen t\rangle,\;\; -1<q\leq 1.
\]

\medskip\noindent $\underline{\mathsf{A}_{4.8}=\langle \gen u, x\gen u, -\gen x, -x\gen x + c\gen u\rangle,\; q=-1,\;\; c\in \mathbb{R}}:$ In this case we may take $\mathsf{A}_{4.8}=\langle \gen u, x\gen u, -\gen x, -x\gen x\rangle$ and then we find that $e_5=a(t)\gen t$. The same analysis as in the previous case for  $\mathsf{A}_{4.8}=\langle \gen u, x\gen u, -\gen x, qx\gen x + (1+q)u\gen u\rangle,\; -1< q\leq 1$ gives $e_5=\gen t$ in canonical form. Thus we have the realization
\[
\mathsf{A}_{4.8}\oplus\mathsf{A}_1=\langle \gen u, x\gen u, -\gen x, qx\gen x + (1+q)u\gen u, \gen t\rangle,\;\; -1\leq q\leq 1.
\]

\medskip\noindent $\underline{\mathsf{A}_{4.8}=\langle \gen u, x\gen u, -\gen x, -x\gen x + t\gen u\rangle,\;\; q=-1}:$ In this case $e_5=c(t)\gen u$ and this together with the operator $\gen u$ gives $\dot{c}(t)=0$ in the equation for $G$. Thus we have $e_5=\gen u$ and this contradicts the requirement of dimension. Hence we have no realization in this case.

\medskip\noindent $\underline{\mathsf{A}_{4.8}=\langle \gen u, x\gen u, -(\gen t + \gen x), qt\gen t + qx\gen x + (1+q)u\gen u\rangle,\;\; q\neq 0}:$ In this case we find $e_5=at\gen t + \kappa(t-x)^{(q+1)/q}\gen u$. We need $a\neq 0$ for otherwise the algebra will contain the subalgebra $\langle \gen u, x\gen u, (t-x)^{(q+1)/q}\gen u\rangle$ which linearizes the evolution equation since it is a rank-one three-dimensional abelian Lie algebra. So we take $a=1$ (we divide throughout by $a$) and then $e_5=t\gen t + \kappa(t-x)^{(q+1)/q}\gen u$. We have the residual equivalence group $\mathscr{E}(e_1, e_2, e_3, e_4): t'=t,\; x'=x,\; u'=u+\lambda(t-x)^{(q+1)/q}$ and under such a transformation $e_5$ is mapped to $\displaystyle e'_5=t'\gen {t'} + [\kappa + \lambda(1+\frac{1}{q})](t-x)^{(q+1)/q}\gen {u'}$. We choose $\lambda$ so that $\displaystyle \kappa + \lambda(1+\frac{1}{q})=0$ for $q\neq -1$ giving $e_5=t\gen t$ in canonical form. If $q=-1$ then $e_5=t\gen t + \kappa \gen u$ and even in this case we may take $e_5=t\gen t$ in canonical form since $\gen u$ belongs to $\mathsf{A}_{4.8}$. Thus we have the realization
\[
\mathsf{A}_{4.8}\oplus\mathsf{A}_1=\langle \gen u, x\gen u, -(\gen t + \gen x), qt\gen t + qx\gen x + (1+q)u\gen u, t\gen t\rangle, \;\; q\neq 0.
\]

\medskip\noindent $\underline{\mathsf{A}_{4.8}=\langle \gen u, x\gen u, -(\gen t + \gen x), a\gen t + u\gen u\rangle},\; a\neq 0:$ In this case $\displaystyle e_5=\kappa \gen t + \lambda e^{(t-x)/a}\gen u$. The linear combination $ae_5-\kappa e_4$ gives the operator $c(t,u)\gen u$ where $\displaystyle c(t,u)=a\lambda e^{(t-x)/a}-\kappa u$ and then the algebra contains the rank-one three-dimensional solvable Lie subalgebra $\langle \gen u, x\gen u, c(t,u)\gen u\rangle$ which linearizes the evolution equation (any rank-one three-dimensional solvable Lie algebra with an abelian ideal of dimension two linearizes the evolution equation). So we have no non-linearizing realization in this case.

\medskip\noindent $\underline{\mathsf{A}_{4.8}=\langle \gen x, \gen t, (t+u)\gen x, t\gen t + (1+q)x\gen x + u\gen u\rangle}:$ We find that $e_5=u^{q+1}\gen x$. For $q=-1$ we see that $e_5=e_1$, so we have no realization in this case. For $-1<q\leq 1$ we see that the algebra contains the subalgebra $\langle \gen x, (t+u)\gen x, u^{q+1}\gen x\rangle$ which linearizes the evolution equation since it is a rank-one three-dimensional abelian Lie algebra. So we have no non-linearizing realization in this case.

\medskip\noindent $\underline{\mathsf{A}_{4.8}=\langle \gen x, \gen u, u\gen x, t\gen t + (1+q)x\gen x + u\gen u\rangle}:$ We have $e_5=a\gen t +\beta t^{q+1}\gen x$. We must have $a\neq 0$ for otherwise the algebra contains the operator $e_5=t^{q+1}\gen x$; for $q=-1$ this contravenes the requirements of dimension, for $-1<q\leq 1$ the operator $t^{q+1}\gen x$ together with $\gen x$ gives the contradiction $u_1=0$ in the equation for $G$. So we take $a=1$ and then $e_5=\gen t +\beta t^{q+1}\gen x$.  We have the residual equivalence group $\mathscr{E}(e_1, e_2, e_3, e_4): t'=kt,\; x'=x+\kappa t^{(q+1)},\; u'=u$ and under such a transformation $e_5$ is mapped to $e'_5=t'\gen {t'} + [\beta + (1+q)\kappa]t^{q+1}\gen {x'}$. We choose $\kappa$ so that $\beta + (1+q)\kappa=0$ for $q\neq -1$ and thus we have $e_5=t\gen t$ in canonical form; for $q=-1$ we have $e_5=t\gen t + \beta \gen x=t\gen t + \beta e_1$ and in this case we may take $e_5=t\gen t$ in canonical form. Thus we have the realization
\[
\mathsf{A}_{4.8}\oplus\mathsf{A}_1=\langle \gen x, \gen u, u\gen x, t\gen t + (1+q)x\gen x + u\gen u, t\gen t\rangle.
\]

\medskip\noindent $\underline{\mathsf{A}_{4.8}=\langle \gen x, \gen u, u\gen x, (1+q)x\gen x + u\gen u\rangle,\;\; -1<q\leq 1}:$ In this case $e_5=a(t)\gen t$. We have the residual equivalence group $\mathscr{E}(e_1, e_2, e_3, e_4): t'=T(t),\, x'=x,\; u'=u$ and under such a transformation $e_5$ is mapped to $e'_5=a(t)\dot{T}(t)\gen {t'}$. We choose $T(t)$ so that $a(t)\dot{T}(t)=1$ giving $e_5=\gen t$ in canonical form. Thus we have the realization
\[
\mathsf{A}_{4.8}\oplus\mathsf{A}_1=\langle \gen x, \gen u, u\gen x, (1+q)x\gen x + u\gen u, \gen t\rangle, \; -1< q\leq 1.
\]

\medskip\noindent $\underline{\mathsf{A}_{4.8}=\langle \gen x, \gen u, u\gen x, p\gen x + u\gen u\rangle,\; p\in \mathbb{R},\;  q=-1}:$ In this case we have $e_5=a(t)\gen t + c(t)\gen x$. We have $a(t)\neq 0$ for otherwise we have $b(t)\gen x$ as a symmetry which, together with $\gen x$ gives $\dot{b}(t)u_1=0$ in the equation for $G$. This can be true only if $b(t)=$ constant, and this then gives $e_5=be_1$, contravening the requirements of dimension. So $a(t)\neq 0$. We have the residual equivalence group $\mathscr{E}(e_1, e_2, e_3, e_4): t'=T(t),\, x'=x+Y(t),\; u'=u$ and under such a transformation $e_5$ is mapped to $e'_5=a(t)\dot{T}(t)\gen {t'} + [a(t)\dot{Y}(t) + b(t)]\gen {x'}$. We choose $T(t)$ so that $a(t)\dot{T}(t)=1$ and $Y(t)$ so that $a(t)\dot{Y}(t) + b(t)=0$ giving $e_5=\gen t$ in canonical form. Thus we have the realization
\[
\mathsf{A}_{4.8}\oplus\mathsf{A}_1=\langle \gen x, \gen u, u\gen x, u\gen u, \gen t\rangle.
\]
Putting this realization together with the previous one we have the realization
\[
\mathsf{A}_{4.8}\oplus\mathsf{A}_1=\langle \gen x, \gen u, u\gen x, (1+q)x\gen x + u\gen u, \gen t\rangle, \; -1\leq q\leq 1.
\]

\medskip\noindent $\underline{\mathsf{A}_{4.8}=\langle \gen x, \gen u, u\gen x, t\gen x + u\gen u\rangle,\;  q=-1}:$ Here we find $e_5=b(t)\gen x$ which together with $\gen x$ gives $\dot{b}(t)u_1=0$ in the equation for $G$. This is only possible for $b(t)=$constant, and this contravenes the requirements of dimension. So we have no realization in this case.

\medskip\noindent $\underline{\mathsf{A}_{4.8}=\langle \gen x, \gen u, u\gen x + t\gen u, (1-q)t\gen t + (1+q)x\gen x + u\gen u\rangle}:$ Here we find $e_5=b(t)\gen x$ which together with $\gen x$ gives $\dot{b}(t)u_1=0$ in the equation for $G$. This is only possible for $b(t)=$constant, and this contravenes the requirements of dimension. So we have no realization in this case.

\medskip\noindent $\underline{\mathsf{A}_{4.8}=\langle \gen x, \gen t, t\gen x - \gen u, t\gen t + (1+q)x\gen x + qu\gen u\rangle,\; q\neq 0}:$ Here we find $e_5=0$ if $q\neq -1$ and  $e_5=\beta \gen x$ if $q=-1$, and this contravenes the requirements of dimension. So we have no realization in this case.

\medskip\noindent $\underline{\mathsf{A}_{4.8}=\langle \gen x, \gen t, t\gen x - \gen u, t\gen t + [x+\beta]\gen x + \gamma \gen u\rangle,\; q= 0}:$ Here we find $e_5=\gen u$ and then $\gen x$ and $t\gen x$ belong to the algebra, giving $u_1=0$ in the equation for $G$. Thus we have no realization in this case.

\medskip\noindent $\underline{\mathsf{A}_{4.8}=\langle \gen x, \gen u, t\gen t + u\gen x, qt\ln|t|\gen t + (1+q)x\gen x + u\gen u\rangle,\; q\neq 0}:$ Here we find $e_5=0$ if $q\neq -1$ and  $e_5=\beta \gen x$ if $q=-1$, and this contravenes the requirements of dimension. So we have no realization in this case.

\medskip\noindent $\underline{\mathsf{A}_{4.8}=\langle \gen x, \gen u, t\gen t + u\gen x, \alpha t\gen t + x\gen x + u\gen u\rangle,\; q=0}:$ Here we find $e_5=at\gen t$. Hence we have the realization
\[
\mathsf{A}_{4.8}\oplus\mathsf{A}_1=\langle \gen x, \gen u, t\gen t + u\gen x, x\gen x + u\gen u, t\gen t\rangle.
\]

\subsection{Non-linearizing realizations of $\mathsf{A}_{4.8}\oplus\mathsf{A}_1$:}

\begin{align*}
\mathsf{A}_{4.8}\oplus\mathsf{A}_1&=\langle \gen u, x\gen u, -\gen x, t\gen t + qx\gen x + (1+q)u\gen u, t\gen t\rangle\\
\mathsf{A}_{4.8}\oplus\mathsf{A}_1&=\langle \gen u, x\gen u, -\gen x, qx\gen x + (1+q)u\gen u, \gen t\rangle\\
\mathsf{A}_{4.8}\oplus\mathsf{A}_1&=\langle \gen u, x\gen u, -(\gen t + \gen x), qt\gen t + qx\gen x + (1+q)u\gen u, t\gen t\rangle, \;\; q\neq 0\\
\mathsf{A}_{4.8}\oplus\mathsf{A}_1&=\langle \gen x, \gen u, u\gen x, t\gen t + (1+q)x\gen x + u\gen u, t\gen t\rangle\\
\mathsf{A}_{4.8}\oplus\mathsf{A}_1&=\langle \gen x, \gen u, u\gen x, (1+q)x\gen x + u\gen u, \gen t\rangle\\
\mathsf{A}_{4.8}\oplus\mathsf{A}_1&=\langle \gen x, \gen u, t\gen t + u\gen x, x\gen x + u\gen u, t\gen t\rangle,\;\; q=0.
\end{align*}

\subsection{Realizations of $\mathsf{A}_{4.9}\oplus\mathsf{A}_1$:} We have only one realization of $\mathsf{A}_{4.9}$: that is
\[
\mathsf{A}_{4.9}=\langle \gen x, \gen u, u\gen x + t\gen u, -(1+t^2)\gen t+ [2qx-\frac{u^2}{2}]\gen x + (q-t)u\gen u\rangle,\;\; q\geq 0.
\]
Then we find that $e_5=b(t)\gen x$ and we see that this, together with $\gen x$ gives $\dot{b}(t)u_1=0$ which requires $b(t)=$constant for consistency. This however means that $e_5\in \mathsf{A}_{4.9}$ and so we have no realization of $\mathsf{A}_{4.9}\oplus\mathsf{A}_1$.

\subsection{Realizations of $\mathsf{A}_{4.10}\oplus\mathsf{A}_1$:}

\medskip\noindent $\underline{\mathsf{A}_{4.10}=\langle \gen u, x\gen u, t\gen t + u\gen u, at\gen t-(1+x^2)\gen x - xu\gen u\rangle}:$ Here we find $e_5=\alpha t\gen t + \kappa t\sqrt{1+x^2}\exp(a\arctan x)\gen u$. Clearly $\alpha\neq 0$ for otherwise the algebra
\[
\langle \gen u, x\gen u, t\sqrt{1+x^2}\exp(a\arctan x)\gen u\rangle
\]
linearizes the evolution equation since it is a rank-one three-dimensional abelian Lie algebra. So we may assume $\alpha=1$ and then
$$e_5=t\gen t + \kappa t\sqrt{1+x^2}\exp(a\arctan x)\gen u.$$
We have the residual equivalence group $\mathscr{E}(e_1, e_2, e_3, e_4): t'=kt,\; x'=x,\; u'=u+\lambda t\sqrt{1+x^2}\exp(a\arctan x)$ and under such a transformation $e_5$ is mapped to $e'_5=t'\gen {t'} + t[\kappa + \lambda]\sqrt{1+x^2}\exp(a\arctan x)\gen {u'}$. We may choose $\lambda$ so that $\kappa + \lambda=0$ and then $e_5=t\gen t$ in canonical form. The algebra then contains the subalgebra $\langle \gen u, x\gen u, u\gen u\rangle$ which linearizes the evolution equation since it is a rank-one three-dimensional solvable Lie algebra with two-dimensional abelian ideal. Thus we have no non-linearizing realization in this case.

\medskip\noindent $\underline{\mathsf{A}_{4.10}=\langle \gen x, \gen u, x\gen x + u\gen u, u\gen x - x\gen u\rangle}:$ Here we find $e_5=a(t)\gen t$. We have the residual equivalence group $\mathscr{E}(e_1, e_2, e_3, e_4): t'=T(t),\; x'=x,\; u'=u$ and under such a transformation $e_5$ is mapped to $e'_5=a(t)\dot{T}(t)\gen {t'}$. We choose $T(t)$ so that $a(t)\dot{T}(t)=1$ and so $e_5=\gen t$ in canonical form. Hence we have the realization
\[
\mathsf{A}_{4.10}\oplus\mathsf{A}_1=\langle \gen x, \gen u, x\gen x + u\gen u, u\gen x-x\gen u, \gen t\rangle.
\]

\medskip\noindent $\underline{\mathsf{A}_{4.10}=\langle \gen x, \gen u, x\gen x + u\gen u, t\gen t + u\gen x - x\gen u\rangle}:$ Here we find $e_5=at\gen t$. Hence we have the realization

\[
\mathsf{A}_{4.10}\oplus\mathsf{A}_1=\langle \gen x, \gen u, x\gen x + u\gen u, t\gen t +u\gen x-x\gen u, t\gen t\rangle.
\]
However, this is equivalent to
\[
\mathsf{A}_{4.10}\oplus\mathsf{A}_1=\langle \gen x, \gen u, x\gen x + u\gen u, u\gen x-x\gen u, \gen t\rangle
\]
under the equivalence transformation $t'=\ln|t|,\; x'=x,\; u'=u$.

\medskip\noindent $\underline{\mathsf{A}_{4.10}=\langle \gen x, \gen u, t\gen t + x\gen x + u\gen u, at\gen t + u\gen x - x\gen u\rangle}:$ Here we find $e_5=\alpha t\gen t$. Hence we have the realization

\[
\mathsf{A}_{4.10}\oplus\mathsf{A}_1=\langle \gen x, \gen u, x\gen x + u\gen u, t\gen t +u\gen x-x\gen u, t\gen t\rangle.
\]
However, as we have already remarked, this is equivalent to
\[
\mathsf{A}_{4.10}\oplus\mathsf{A}_1=\langle \gen x, \gen u, x\gen x + u\gen u, u\gen x-x\gen u, \gen t\rangle
\]
under the equivalence transformation $t'=\ln|t|,\; x'=x,\; u'=u$. Consequently we have only one inequivalent non-linearizing realization of $\mathsf{A}_{4.10}\oplus\mathsf{A}_1$:
\[
\mathsf{A}_{4.10}\oplus\mathsf{A}_1=\langle \gen x, \gen u, x\gen x + u\gen u, u\gen x-x\gen u, \gen t\rangle.
\]

\section{Realizations of the non-decomposable five-dimensional solvable Lie algebras.} As we have already noted, there are thirty nine such types, labelled $\mathsf{A}_{5.k}$ with $k=1,\dots 39$. The following give either a linear (or linearizable) evolution equation or are inadmissible: $\mathsf{A}_{5.1}$ and $\mathsf{A}_{5.2}$ and $\mathsf{A}_{5.7}$ through to $\mathsf{A}_{5.18}$, by virtue of them containing a four-dimensional abelian Lie algebra, and we know that any evolution equation admitting such a symmetry algebra is linear (or linearizable).

We also note that $\mathsf{A}_{5.3},\; \mathsf{A}_{5.4}\; \mathsf{A}_{5.5},\; \mathsf{A}_{5.6}$ all contain an abelian Lie subalgebra (in fact, it is always $\langle e_1, e_2, e_3\rangle$) and we know that rank-one three-dimensional abelian Lie algebras linearize our evolution equation (when they are admitted as symmetries). Thus we look only at the cases of rank-two and rank-three three-dimensional abelian Lie algebras. These have the inequivalent canonical forms
\[
\langle \gen t, \gen u, x\gen u\rangle,\; \langle \gen u, \gen t, x\gen u\rangle,\; \langle \gen u, x\gen u, \gen t\rangle
\]
for rank-two realizations, and
\[
\langle \gen t, \gen x, \gen u\rangle,\; \langle \gen x, \gen t, \gen u\rangle,\; \langle \gen x, \gen u, \gen t\rangle
\]
for rank-three realizations. Here we note that $\gen u,\; x\gen u$ are interchangeable under the equivalence transformation $t'=t,\; x'=1/x,\; u'=u/x$, and $\gen x$ is interchangeable with $\gen u$ under the equivalence transformation $t'=t,\; x'=u,\; u'=x$.

\medskip\noindent $\underline{\mathsf{A}_{5.3}}:$ Here we have that $\langle e_1, e_2, e_3\rangle$ is an abelian Lie algebra and $[e_2, e_4]=e_3,\; [e_2, e_5]=e_1,\; [e_4, e_5]=e_2$, all other commutators being zero.

\medskip\noindent $\underline{\langle e_1, e_2, e_3\rangle=\langle \gen t, \gen x, \gen u\rangle}:$ This is impossible to realize since $e_5=a(t)\gen t + b(t,x,u)\gen x + c(t,x,u)\gen u$ and $[e_2, e_5]=e_1$ requires $[\gen x, e_5]=\gen t$ which is clearly impossible.

\medskip\noindent $\underline{\langle e_1, e_2, e_3\rangle=\langle \gen x, \gen u, \gen t\rangle}:$ This is impossible to realize since $e_4=a(t)\gen t + b(t,x,u)\gen x + c(t,x,u)\gen u$ and $[e_2, e_4]=e_1$ requires $[\gen x, e_5]=\gen t$ which is clearly impossible.

\medskip\noindent $\underline{\langle e_1, e_2, e_3\rangle=\langle \gen t, \gen u, x\gen u\rangle}:$ This is impossible to realize since $e_5=a(t)\gen t + b(t,x,u)\gen x + c(t,x,u)\gen u$ and $[e_2, e_5]=e_1$ requires $[\gen u, e_5]=\gen t$ which is clearly impossible.

\medskip\noindent $\underline{\langle e_1, e_2, e_3\rangle=\langle \gen u, x\gen u, \gen t\rangle}:$ This is impossible to realize since $e_4=a(t)\gen t + b(t,x,u)\gen x + c(t,x,u)\gen u$ and $[e_2, e_4]=e_1$ requires $[x\gen u, e_5]=\gen t$ which is clearly impossible.

\medskip\noindent $\underline{\langle e_1, e_2, e_3\rangle=\langle \gen x, \gen t, \gen u\rangle}:$ From the commutation relations $e_4=\alpha \gen t+\beta \gen x + [t+l]\gen u$ so that we would have $\gen u$ and $t\gen u$ as symmetries, giving the contradiction $1=0$ in the equation for $G$. So we have no realization in this case.

\medskip\noindent $\underline{\langle e_1, e_2, e_3\rangle=\langle \gen u, \gen t, x\gen u\rangle}:$ From the commutation relations $e_4=\alpha \gen t+ [tx+c(x)]\gen u$ and $e_5=\alpha \gen t + [t+\gamma(x)]\gen u$. Then $[e_4, e_5]=e_2$ is impossible. So we have no realization in this case.

So we have no non-linearizing realization for $\mathsf{A}_{5.3}$.

\medskip\noindent $\underline{\mathsf{A}_{5.4}}:$ Here we have that $\langle e_1, e_2, e_3\rangle$ is an abelian Lie algebra and $[e_2, e_4]=e_1,\; [e_3, e_5]=e_1$, all other commutators being zero.

We note immediately that $\langle \gen t, \gen u, x\gen u\rangle$ and $\langle \gen t, \gen x, \gen u\rangle$ are impossible to realize since we require $[e_2, e_4]=e_1$ and with $e_4=a(t)\gen t + b(t,x,u)\gen x + c(t,x,u)\gen u$ this is impossible for $e_2=\gen x$ or $e_2=\gen u$.

\medskip\noindent $\underline{\langle e_1, e_2, e_3\rangle=\langle \gen x, \gen t, \gen u\rangle}:$ From the commutation relations $e_4=a\gen t + [t+l]\gen x + c\gen u$ so that the algebra will contain the operators $\gen x,\, t\gen x$ which give $F=0$ in the equation for $F$, which is a contradiction. Thus we have no realization.

\medskip\noindent $\underline{\langle e_1, e_2, e_3\rangle=\langle \gen x, \gen u, \gen t\rangle}:$ From the commutation relations $e_5=a\gen t + [t+l]\gen x + c\gen u$ so that the algebra will contain the operators $\gen x,\, t\gen x$ which give $F=0$ in the equation for $F$, which is a contradiction. Thus we have no realization.

\medskip\noindent $\underline{\langle e_1, e_2, e_3\rangle=\langle \gen u, \gen t, x\gen u\rangle}:$ From the commutation relations $[e_1, e_4]=0,\, [e_2, e_4]=e_1,\; [e_3, e_4]=0$ we find that $e_4=\alpha \gen t + [t+\gamma(x)]\gen u$. We also find  $e_5=a\gen t -\gen x + c(x)\gen u$ from $[e_1, e_5]=[e_2, e_5]=0,\, [e_3, e_5]=e_1$. We have the residual equivalence group  $\mathscr{E}(e_1, e_2, e_3): t'=t+k,\, x'=x,\, u'=u + U(x)$ and under such a transformation $e_5$ is mapped to $e'_5=a\gen {t'} -\gen {x'} + [c(x)-U'(x)]\gen {u'}$ and we choose $U(x)$ so that $c(x)=U'(x)$ giving $e_5=a\gen t - \gen x$ in canonical form. Then $[e_4, e_5]=0$ gives $e_4=\alpha \gen t + [t+ax+l]\gen u$ so that $e_4=ae_2 + le_1 + ae_3+t\gen u$. Hence we have the operators $\gen u$ and $ t\gen u$ belonging to the algebra, and these give the contradiction $1=0$ in the equation for $G$, so we have no realization here.

\medskip\noindent $\underline{\langle e_1, e_2, e_3\rangle=\langle \gen u, x\gen u, \gen t\rangle}:$ From the commutation relations $[e_1, e_4]=0,\, [e_2, e_4]=e_1,\; [e_3, e_4]=0$ we find that $e_4=\alpha \gen t -\gen x +  \gamma(x)\gen u$. We also find  $e_5=a\gen t + [t+c(x)]\gen u$ from $[e_1, e_5]=[e_2, e_5]=0,\, [e_3, e_5]=e_1$. We have the residual equivalence group  $\mathscr{E}(e_1, e_2, e_3): t'=t+k,\, x'=x,\, u'=u + U(x)$ and under such a transformation $e_4$ is mapped to $e'_4=a\gen {t'} -\gen {x'} + [\gamma(x)-U'(x)]\gen {u'}$ and we choose $U(x)$ so that $\gamma(x)=U'(x)$ giving $e_4=\alpha\gen t - \gen x$ in canonical form. Then $[e_4, e_5]=0$ gives $e_5= a\gen t + [t+\alpha x+l]\gen u$ so that $e_5=ae_2 + le_1 + \alpha e_2+t\gen u$. Hence we have the operators $\gen u$ and $ t\gen u$ belonging to the algebra, and these give the contradiction $1=0$ in the equation for $G$, so we have no realization here.

\medskip\noindent $\underline{\mathsf{A}_{5.5}}:$ Here we have that $\langle e_1, e_2, e_3\rangle$ is an abelian Lie algebra and $[e_2, e_5]=e_1,\; [e_3, e_4]=e_1,\; [e_3, e_5]=e_2$, all other commutators being zero.

We note immediately that $\langle \gen t, \gen u, x\gen u\rangle$ and $\langle \gen t, \gen x, \gen u\rangle$ are impossible to realize since we require $[e_3, e_4]=e_1$ and with $e_4=a(t)\gen t + b(t,x,u)\gen x + c(t,x,u)\gen u$ this is impossible for $e_2=\gen x$ or $e_2=\gen u$. Also $\langle \gen u, \gen t, x\gen u\rangle$ and $\langle \gen x, \gen t, \gen u\rangle$ are impossible to realize since we require $[e_3, e_5]=e_2$ and with $e_5=a(t)\gen t + b(t,x,u)\gen x + c(t,x,u)\gen u$ this is impossible for $e_3=\gen u$ or $e_3=x\gen u$. Thus we must have $\langle e_1, e_2, e_3\rangle=\langle \gen x, \gen u, \gen t\rangle$ or $\langle e_1, e_2, e_3\rangle=\langle \gen u, x\gen u, \gen t\rangle$.

\medskip\noindent $\underline{\langle e_1, e_2, e_3\rangle=\langle \gen u, x\gen u, \gen t\rangle}:$ In this case $e_4=a(t)\gen t + c(t,x)\gen u$ from $[e_1, e_4]=[e_2, e_4]=0$. Then $[e_3, e_4]=e_1$ gives $e_4=a\gen t + [t+c(x)]\gen u=ae_3 + [t+c(x)]\gen u$ so that the algebra contains the subalgebra $\langle \gen u, x\gen u, [t+c(x)]\gen u\rangle$ which either linearizes the evolution equation since it is a rank-one three-dimensional abelian Lie algebra, or it is inadmissible. Thus we have no non-linearizing realization in this case.

\medskip\noindent $\underline{\langle e_1, e_2, e_3\rangle=\langle \gen x, \gen u, \gen t\rangle}:$ From the commutation relations we obtain $e_4=a\gen t + [t+l]\gen x + c\gen u$ and from this it follows that the algebra contains the operators $\gen x,\; t\gen x$ which are incompatible because they give $F=0$ in the equation for $F$. Hence we have no realization in this case.

\medskip\noindent $\underline{\mathsf{A}_{5.6}}:$ Here we have that $\langle e_1, e_2, e_3\rangle$ is an abelian Lie algebra and $[e_2, e_5]=e_1,\; [e_3, e_4]=e_1,\; [e_3, e_5]=e_2,\; [e_4, e_5]=e_3$, all other commutators being zero.

We note that $\langle e_1, e_2, e_3\rangle=\langle \gen t, \gen x, \gen u\rangle$ is impossible, since we need $[e_2, e_5]=e_1$ and this is not possible with $e_5=a(t)\gen t + b(t,x,u)\gen x + c(t,x,u)\gen u$. Similarly, $\langle e_1, e_2, e_3\rangle=\langle \gen x, \gen t, \gen u\rangle$ is impossible, since we need $[e_3, e_5]=e_2$. For the same reasons, $\langle e_1, e_2, e_3\rangle=\langle \gen t, \gen u, x\gen u\rangle$ and $\langle e_1, e_2, e_3\rangle=\langle \gen u, \gen t, x\gen u\rangle$ are impossible.

\medskip\noindent $\underline{\langle e_1, e_2, e_3\rangle=\langle \gen x, \gen u, \gen t\rangle}:$ From the commutation relations, $e_4=a\gen t + [t+l]\gen x + c\gen u$ so that the algebra will contain the operators $\gen x,\; t\gen x$ which are incompatible since they give $F=0$.

\medskip\noindent $\underline{\langle e_1, e_2, e_3\rangle=\langle \gen u, x\gen u, \gen t\rangle}:$ From the commutation relations, $e_4=a\gen t + [t+c(x)]\gen u$ and $e_5=\alpha \gen t - \gen x + [tx + \gamma(x)]\gen u$. Then $[e_4, e_5]=e_3$ is impossible to realize. Hence we have no realization in this case.

\medskip\noindent\underline{The algebras $\mathsf{A}_{5.7}$ to $\mathsf{A}_{5.18}$:} These all contain the four-dimensional abelian Lie algebra $\langle e_1, e_2, e_3, e_4\rangle$, and any evolution equation (of our type) admitting such an algebra is linearizable, so we have no non-linearizing realizations of these algebras.

\medskip\noindent $\underline{\mathsf{A}_{5.19}}:$ In this case we have $\langle e_1, e_2, e_3, e_4\rangle=\mathsf{A}_{3.3}\oplus\mathsf{A}_1$ and we also have the commutation relations
\[
[e_1, e_5]=(p+1)e_1,\; [e_2, e_5]=e_2,\; [e_3, e_5]=pe_3,\; [e_4, e_5]=qe_4
\]
with $p\in \mathbb{R},\; q\neq 0$.

\medskip\noindent $\underline{\langle e_1, e_2, e_3, e_4\rangle=\langle \gen u, x\gen u, -(\gen t+\gen x), \gen t \rangle}:$ The commutation relations give $e_5=[qt+l]\gen t + qx\gen x + [(q+1)u+\gamma]\gen u$ and we must have $p=q$. It is clear that we may take $e_5=qt\gen t + qx\gen x + (q+1)u\gen u$ since $\gen t,\; \gen u\in \langle e_1, e_2, e_3, e_4\rangle$. Thus we have the realization
\[
\mathsf{A}_{5.19}=\langle \gen u, x\gen u, -(\gen t+\gen x), \gen t, qt\gen t + qx\gen x + (q+1)u\gen u\rangle.
\]

\medskip\noindent $\underline{\langle e_1, e_2, e_3, e_4\rangle=\langle \gen u, x\gen u, -\gen x, \gen t \rangle}:$ The commutation relations give $e_5=[qt+l]\gen t + px\gen x + [(p+1)u+\gamma]\gen u$. It is clear that we may take $e_5=qt\gen t + qx\gen x + (q+1)u\gen u$ since $\gen t,\; \gen u\in \langle e_1, e_2, e_3, e_4\rangle$. Thus we have the realization
\[
\mathsf{A}_{5.19}=\langle \gen u, x\gen u, -\gen x, \gen t, qt\gen t + px\gen x + (p+1)u\gen u\rangle.
\]
Note that the algebra $\langle \gen u, x\gen u, -\gen x, \gen t\rangle$ gives the same equation as $\langle \gen u, x\gen u, -(\gen t + \gen x), \gen t\rangle$ and thus we see that
\[
\mathsf{A}_{5.19}=\langle \gen u, x\gen u, -\gen x, \gen t, qt\gen t + px\gen x + (p+1)u\gen u\rangle,\;\; q\neq 0,\;\; p\in \mathbb{R}.
\]
is a generalization of
\[
\mathsf{A}_{5.19}=\langle \gen u, x\gen u, -(\gen t+\gen x), \gen t, qt\gen t + qx\gen x + (q+1)u\gen u\rangle
\]
since it allows for the cases where $p\neq q$. Hence we have only one inequivalent realization
\[
\mathsf{A}_{5.19}=\langle \gen u, x\gen u, -\gen x, \gen t, qt\gen t + px\gen x + (p+1)u\gen u\rangle, \;\; q\neq 0,\;\; p\in \mathbb{R}.
\]
which subsumes
\[
\mathsf{A}_{5.19}=\langle \gen u, x\gen u, -(\gen t+\gen x), \gen t, qt\gen t + qx\gen x + (q+1)u\gen u\rangle.
\]

\medskip\noindent $\underline{\langle e_1, e_2, e_3, e_4\rangle=\langle \gen x, \gen t, (t+u)\gen x, \gen t-\gen u \rangle}:$ The commutation relations give $q=1$ and $e_5=[t+l]\gen t + [(1+p)x + r]\gen x + [u-l]\gen u$, so we have $e_5=le_4+re_1 + t\gen t + (1+p)\gen x + u\gen u$ and thus we may take $e_5=t\gen t + (1+p)x\gen x + u\gen u$. Thus we have the realization
\[
\mathsf{A}_{5.19}=\langle \gen x, \gen t, (t+u)\gen x, \gen t-\gen u, t\gen t + (1+p)x\gen x + u\gen u\rangle,\;\; q=1,\;\, p\in \mathbb{R}.
\]

\medskip\noindent $\underline{\langle e_1, e_2, e_3, e_4\rangle=\langle \gen x, \gen u, u\gen x, \gen t\rangle}:$ The commutation relations give $e_5=[qt+l]\gen t + [(1+p)x + r]\gen x + u\gen u$, so we have $e_5=le_4+re_1 + qt\gen t + (1+p)x\gen x + u\gen u$ and thus we may take $e_5=qt\gen t + (1+p)x\gen x + u\gen u$. Thus we have the realization
\[
\mathsf{A}_{5.19}=\langle \gen x, \gen u, u\gen x, \gen t, qt\gen t + (1+p)x\gen x + u\gen u\rangle.
\]

\medskip\noindent $\underline{\langle e_1, e_2, e_3, e_4\rangle=\langle \gen x, \gen t, t\gen x-\gen u, \gen t-u\gen x + \kappa \gen u\rangle,\;\; \kappa\in \mathbb{R}}:$ The commutation relations give $q=1$ and  $e_5=[t+l]\gen t + [(1+p)x -lu + r]\gen x + [pu + l\kappa]\gen u$, so we have $e_5=le_4+re_1 + t\gen t + (1+p)x\gen x + pu\gen u$ and thus we may take $e_5=t\gen t + (1+p)x\gen x + pu\gen u$. Thus we have the realization
\[
\mathsf{A}_{5.19}=\langle \gen x, \gen t, t\gen x-\gen u, \gen t-u\gen x + \kappa \gen u, t\gen t + (1+p)x\gen x + pu\gen u\rangle,
\]
with $q=1$ and  $p,\kappa\in \mathbb{R}$.

\medskip\noindent $\underline{\langle e_1, e_2, e_3, e_4\rangle=\langle \gen x, \gen u, t\gen t+u\gen x, t\gen t\rangle}:$ The commutation relations give $p=q$ and   $e_5=[qt\ln|t|+lt]\gen t + [(1+q)x + r]\gen x + u]\gen u$, so we have $e_5=le_4+re_1 + qt\ln|t|\gen t + (1+q)x\gen x + u\gen u$ and thus we may take $e_5=qt\ln|t|\gen t + (1+q)x\gen x + u\gen u$. Thus we have the realization
\[
\mathsf{A}_{5.19}=\langle \gen x, \gen u, t\gen t+u\gen x, t\gen t, qt\ln|t|\gen t + (1+q)x\gen x + u\gen u\rangle,\;\; q=p.
\]

\subsection{Inequivalent non-linearizing realizations of $\mathsf{A}_{5.19}$:}

\begin{align*}
\mathsf{A}_{5.19}&=\langle \gen u, x\gen u, -\gen x, \gen t, qt\gen t + px\gen x + (p+1)u\gen u\rangle,\\
\mathsf{A}_{5.19}&=\langle \gen x, \gen t, (t+u)\gen x, \gen t-\gen u, t\gen t + (1+p)x\gen x + u\gen u\rangle,\;\; q=1,\\
\mathsf{A}_{5.19}&=\langle \gen x, \gen u, u\gen x, \gen t, qt\gen t + (1+p)x\gen x + u\gen u\rangle\\
\mathsf{A}_{5.19}&=\langle \gen x, \gen t, t\gen x-\gen u, \gen t-u\gen x + \kappa \gen u, t\gen t + (1+p)x\gen x + pu\gen u\rangle,\;\;q=1\\
\mathsf{A}_{5.19}&=\langle \gen x, \gen u, t\gen t+u\gen x, t\gen t, qt\ln|t|\gen t + (1+q)x\gen x + u\gen u\rangle,\;\; p=q\neq 0.
\end{align*}
In the above $p, \kappa\in \mathbb{R}$ and $q\neq 0$ unless otherwise indicated.

\bigskip\noindent \underline{$\mathsf{A}_{5.20}$:} In this case, $\langle e_1, e_2, e_3, e_4\rangle=\mathsf{A}_{3.3}\oplus\mathsf{A}_1$ as well as $[e_1, e_5]=(p+1)e_1,\; [e_2, e_5]=e_2,\; [e_3, e_5]=pe_3,\; [e_4, e_5]=e_1+(p+1)e_4$ with $p\in \mathbb{R}$.

\medskip\noindent \underline{$\langle e_1, e_2, e_3, e_4\rangle=\langle \gen u, x\gen u, -(\gen t + \gen x), \gen t\rangle$:} With $e_5=a(t)\gen t + b(t,x,u)\gen x + c(t,x,u)\gen u$ we have from $[e_3, e_5]=pe_3$ that $\dot{a}(t)=p$ and from $[e_4, e_5]=e_1+(p+1)e_4$ we have $\dot{a}(t)=(p+1)$ and this is a contradiction, so we have no realization in this case.

\medskip\noindent \underline{$\langle e_1, e_2, e_3, e_4\rangle=\langle \gen u, x\gen u, -\gen x, \gen t\rangle$:} From $[e_1, e_5]=(p+1)e_1$,\; $[e_2, e_5]=e_2$,\; $[e_3, e_5]=pe_3$ and $[e_4, e_5]=e_1+(p+1)e_4$, we obtain $e_5=[(p+1)t + l]\gen t + px\gen x + [(p+1)u+t+m]\gen u$. The residual equivalence group $\mathscr{E}(e_1, e_2, e_3, e_4)$ consists of the transformations $t'=t+k$, $x'=x$, $u'=u+\gamma$, and under such a transformation $e_5$ is transformed to
\[
e'_5=[(1+p)t'+l-(1+p)k]\gen {t'} + px'\gen {x'} + [(1+p)u'+t'+m-k-(1+p)\gamma]\gen {u'}.
\]
It is now clear that, for $p+1\neq 0$ we may choose $k=l$ and $\gamma$ so that $m-k-(1+p)\gamma=0$, giving
\[
e_5=(1+p)t\gen {t} + px\gen {x} + [(1+p)u+t]\gen {u}
\]
in canonical form. If $p=-1$ we have $e'_5=l\gen {t'} - x'\gen {x'} + [t'+m-k]\gen {u'}$ then we choose $k=m$ so that $e_5=l\gen {t} - x\gen {x} + t\gen {u}$ in canonical form. Thus we have the realization
\[
\mathsf{A}_{5.20}=\langle \gen u, x\gen u, -\gen x, \gen t, (1+p)t\gen {t} + px\gen {x} + [(1+p)u+t]\gen {u}\rangle.
\]

\medskip\noindent \underline{$\langle e_1, e_2, e_3, e_4\rangle=\langle \gen x, \gen t, (t+u)\gen x, \gen t-\gen u\rangle$:} From $[e_2, e_5]=e_2$ we obtain $\dot{a}(t)=1$ and  $[e_4, e_5]=e_1+(p+1)e_4$ gives $\dot{a}(t)=1+p$ which gives us $p=0$. hence we have the commutation relations $[e_1, e_5]=e_1$, $[e_2, e_5]=e_2$, $[e_3, e_5]=0$ and $[e_4, e_5]=e_1+e_4$. These then give us $e_5=(t+l)\gen t + [x-u+m]\gen x + (u-l)\gen u$. The residual equivalence group $\mathscr{E}(e_1, e_2, e_3)$ is given by transformations of the form $t'=t+k$, $x'=x+\beta$, $u'=u-k$, and under such a transformation $e_5$ is transformed to

\[
e'_5=[t'+l-k]\gen {t'} + [x'-u'+m + k - \beta]\gen {x'} + [u'+k-l]\gen {u'}.
\]
We see that we may choose $k=l$ and $\beta=m+k$, giving
\[
e_5=t\gen t + (x-u)\gen x + u\gen u
\]
in canonical form. Thus we have the realization
\[
\mathsf{A}_{5.20}=\langle \gen x, \gen t, (t+u)\gen x, \gen t-\gen u, t\gen t + (x-u)\gen x + u\gen u\rangle.
\]

\medskip\noindent \underline{$\langle e_1, e_2, e_3, e_4\rangle=\langle \gen x, \gen u, u\gen x, \gen t\rangle$:} The commutation relations $[e_1, e_5]=(p+1)e_1$, $[e_2, e_5]=e_2$, $[e_3, e_5]=pe_3$ and $[e_4, e_5]=e_1+(1+p)e_4$ gives $e_5=[(1+p)t+l]\gen t + [(1+p)x+t+m]\gen x + u\gen u$. The residual equivalence group $\mathscr{E}(e_1, e_2, e_3)$ is given by transformations of the form $t'=t+k$, $x'=x+\beta$, $u'=u$, and under such a transformation $e_5$ is transformed to
\[
e'_5=[(1+p)t'+l-(1+p)k]\gen {t'} + [(1+p)x'+t'+m-k-(1+p)\beta]\gen {x'} + u'\gen {u'}.
\]
We see that, for $1+p\neq 0$, we may choose $k$ so that $(1+p)k=l$ and $\beta$ so that $m-k-(1+p)\beta=0$, giving
\[
e_5=(1+p)t\gen t + [(1+p)x+t]\gen x + u\gen u
\]
in canonical form. For $p=-1$ we have
\[
e'_5=l\gen {t'} + [t'+m-k]\gen {x'} + u'\gen {u'}
\]
and here we may choose $k=m$ so that we obtain $e_5=l\gen t + t\gen x + u\gen u$ in canonical form. From this we obtain the realization
\[
\mathsf{A}_{5.20}=\langle\gen x, \gen u, u\gen x, \gen t, (1+p)t\gen t + [(1+p)x+t]\gen x + u\gen u\rangle.
\]

\medskip\noindent \underline{$\langle e_1, e_2, e_3, e_4\rangle=\langle  \gen x, \gen t, t\gen x-\gen u, \gen t-u\gen x + \kappa\gen u\rangle,\;\; \kappa\in \mathbb{R}$:} From $[e_2, e_5]=e_2$ we find that $\dot{a}(t)=1$ and from $[e_4, e_5]=e_1+(1+p)e_4$ we obtain $\dot{a}(t)=1+p$, so we must have  $p=0$. Thus we have the commutation relations $[e_1, e_5]=e_1$, $[e_2, e_5]=e_2$, $[e_3, e_5]=0$ and $[e_4, e_5]=e_1+e_4$. From these we obtain $\kappa=0$ and $e_5=[t+l]\gen t + [x-lu+m]\gen x + \gen u$. Note that  $e_4=\gen t -u\gen x$. The residual equivalence group $\mathscr{E}(e_1, e_2, e_3)$ is given by transformations of the form $t'=t+k$, $x'=x-ku + \beta$, $u'=u$, and under such a transformation $e_5$ is transformed to
\[
e'_5=[t'+l-k]\gen {t'} + [x'+(k-l)u+m-\beta]\gen {x'} + \gen {u'}.
\]
We choose $k=l$ and $\beta=m$ so that $e_5=t\gen t + x\gen x + \gen u$ in canonical form. Thus we obtain the realization
\[
\mathsf{A}_{5.20}=\langle \gen x, \gen t, t\gen x-\gen u, \gen t-u\gen x, t\gen t + x\gen x + \gen u\rangle.
\]

\medskip\noindent \underline{$\langle e_1, e_2, e_3, e_4\rangle=\langle \gen x, \gen u, t\gen t+u\gen x, t\gen t\rangle$:} From $[e_4, e_5]=e_1+(p+1)e_4$ we have $[t\gen t, e_5]=(p+1)t\gen t + \gen x$ and $[e_3, e_4]=pe_3$ gives $(p+1)t\gen t + \gen x + [u\gen x, e_5]=pt\gen t+pu\gen x$ so that $[u\gen x, e_5]=-\gen t + (pu-1)\gen x$ and this is clearly impossible with $e_5=a(t)\gen t + b(t,x,u)\gen x + c(t,x,u)\gen u$. Thus we have no realization in this case.

\subsection{Inequivalent non-linearizing realizations of $\mathsf{A}_{5.20}$:}

\begin{align*}
\mathsf{A}_{5.20}&=\langle \gen u, x\gen u, -\gen x, \gen t, (1+p)t\gen {t} + px\gen {x} + [(1+p)u+t]\gen {u}\rangle,\quad p\in \mathbb{R}\\
\mathsf{A}_{5.20}&=\langle \gen x, \gen t, (t+u)\gen x, \gen t-\gen u, t\gen t + (x-u)\gen x + u\gen u\rangle,\quad p=0\\
\mathsf{A}_{5.20}&=\langle\gen x, \gen u, u\gen x, \gen t, (1+p)t\gen t + [(1+p)x+t]\gen x + u\gen u\rangle, \quad p\in \mathbb{R}\\
\mathsf{A}_{5.20}&=\langle \gen x, \gen t, t\gen x-\gen u, \gen t-u\gen x, t\gen t + x\gen x + \gen u\rangle,\quad p=0.
\end{align*}

\bigskip\noindent \underline{$\mathsf{A}_{5.21}$:} In this case, $\langle e_1, e_2, e_3, e_4\rangle=\mathsf{A}_{3.3}\oplus\mathsf{A}_1$ as well as $[e_1, e_5]=2e_1,\; [e_2, e_5]=e_2+e_3,\; [e_3, e_5]=e_3+e_4,\; [e_4, e_5]=e_4$.

\medskip\noindent We have no realizations for this algebra:

\medskip\noindent \underline{$\langle e_1, e_2, e_3, e_4\rangle=\langle \gen u, x\gen u, -(\gen t + \gen x), \gen t\rangle$:} Here $[e_2, e_5]=e_2+e_3$ requires $[x\gen u, e_5]=x\gen u -(\gen t+\gen x)$ which is clearly not possible.

\medskip\noindent \underline{$\langle e_1, e_2, e_3, e_4\rangle=\langle \gen u, x\gen u, -\gen x, \gen t\rangle$:} Here $[e_3, e_5]=e_3+e_4$ requires $[\gen x, e_5]=\gen x -\gen t$ which is clearly not possible.

\medskip\noindent \underline{$\langle e_1, e_2, e_3, e_4\rangle=\langle \gen x, \gen t, (t+u)\gen x, \gen t-\gen u\rangle$:} Here $[e_3, e_5]=e_3+e_4$ requires $[(t+u)\gen x, e_5]=(t+u)\gen x +\gen t-\gen u$ which is clearly not possible.

\medskip\noindent \underline{$\langle e_1, e_2, e_3, e_4\rangle=\langle \gen x, \gen u, u\gen x, \gen t\rangle$:} Here $[e_3, e_5]=e_3+e_4$ requires $[u\gen x, e_5]=u\gen x +\gen t$ which is clearly not possible.

\bigskip\noindent \underline{$\mathsf{A}_{5.22}$:} In this case, $\langle e_1, e_2, e_3, e_4\rangle=\mathsf{A}_{3.3}\oplus\mathsf{A}_1$ as well as $[e_1, e_5]=0,\; [e_2, e_5]=e_3,\; [e_3, e_5]=0,\; [e_4, e_5]=e_4$.

\medskip\noindent We have only one admissible realization for this algebra:

\medskip\noindent \underline{$\langle e_1, e_2, e_3, e_4\rangle=\langle \gen u, x\gen u, -(\gen t + \gen x), \gen t\rangle$:} Here $[e_2, e_5]=e_3$ requires $[x\gen u, e_5]=-(\gen t+\gen x)$ which is clearly not possible.

\medskip\noindent \underline{$\langle e_1, e_2, e_3, e_4\rangle=\langle \gen u, x\gen u, -\gen x, \gen t\rangle$:} Here $[e_1, e_5]=0$ gives $e_5=a(t)\gen t + b(t,x)\gen x + c(t,x)\gen u$ and then $[e_2, e_5]=e_3$ gives $b(t,x)\gen u=\gen x$ which is clearly not possible.

\medskip\noindent \underline{$\langle e_1, e_2, e_3, e_4\rangle=\langle \gen x, \gen t, (t+u)\gen x, \gen t-\gen u\rangle$:} The commutation relations $[e_1, e_5]=0,\; [e_2, e_5]=e_3,\; [e_3, e_5]=0$ give $e_5=ae_4+[\frac{u^2}{2}+tu+b(u)]\gen x$. Then $[e_4, e_5]=e_4$ is impossible to implement. Hence we have no realization in this case.

\medskip\noindent \underline{$\langle e_1, e_2, e_3, e_4\rangle=\langle \gen x, \gen u, u\gen x, \gen t\rangle$:} The commutation relations give $e_5=t\gen t+\frac{u^2}{2}\gen x + le_4+\kappa e_1$. Thus we have the realization
\[
\mathsf{A}_{5.22}=\langle \gen x, \gen u, u\gen x, \gen t, t\gen t + \frac{u^2}{2}\gen x\rangle.
\]

\medskip\noindent \underline{$\langle e_1, e_2, e_3, e_4\rangle=\langle \gen x, \gen t, t\gen x-\gen u, \gen t-u\gen x + \kappa \gen u\rangle$:} The commutation relations  $[e_1, e_5]=0,\; [e_2, e_5]=e_3,\; [e_3, e_5]=0$ give $e_5=ae_4+\beta\gen x-t\gen u + [\gamma-a\kappa]\gen u$. Then $[e_4, e_5]=e_4$ is impossible to implement. Hence we have no realization in this case.

\medskip\noindent \underline{$\langle e_1, e_2, e_3, e_4\rangle=\langle \gen x, \gen u, t\gen t+u\gen x, t\gen t\rangle$:} Here $[e_2, e_5]=e_3$ requires $[\gen u, e_5]=t\gen t+u\gen x$ which is clearly not possible.

\subsection{Realizations of $\mathsf{A}_{5.22}$:}

\[
\mathsf{A}_{5.22}=\langle \gen x, \gen u, u\gen x, \gen t, t\gen t + \frac{u^2}{2}\gen x\rangle.
\]

\bigskip\noindent \underline{$\mathsf{A}_{5.23}$:} In this case, $\langle e_1, e_2, e_3, e_4\rangle=\mathsf{A}_{3.3}\oplus\mathsf{A}_1$ as well as $[e_1, e_5]=2e_1,\; [e_2, e_5]=e_2+e_3,\; [e_3, e_5]=e_3,\; [e_4, e_5]=pe_4,\; p\neq 0$.

\medskip\noindent \underline{$\langle e_1, e_2, e_3, e_4\rangle=\langle \gen u, x\gen u, -(\gen t + \gen x), \gen t\rangle$:} Here $[e_2, e_5]=e_2+e_3$ requires $[x\gen u, e_5]=x\gen u-(\gen t+\gen x)$ which is clearly not possible.

\medskip\noindent \underline{$\langle e_1, e_2, e_3, e_4\rangle=\langle \gen u, x\gen u, -\gen x, \gen t\rangle$:} Here $[e_1, e_5]=2e_1$ gives $e_5=a(t)\gen t + b(t,x)\gen x + [2u+c(t,x)]\gen u$ and then $[e_2, e_5]=e_2+e_3$ gives $[2x-b(t,x)]\gen u=-\gen x+x\gen u$ which is clearly not possible.

\medskip\noindent \underline{$\langle e_1, e_2, e_3, e_4\rangle=\langle \gen x, \gen t, (t+u)\gen x, \gen t-\gen u\rangle$:} The commutation relations give  $e_5=t\gen t + [2x+\frac{(t+u)^2}{2}]\gen x + u\gen u + le_4+\kappa e_1$ and we have only $p=1$.  Hence we have the realization
\[
\mathsf{A}_{5.23}=\langle \gen x, \gen t, (t+u)\gen x, \gen t-\gen u, t\gen t + [2x+\frac{(t+u)^2}{2}]\gen x + u\gen u\rangle.
\]

\medskip\noindent \underline{$\langle e_1, e_2, e_3, e_4\rangle=\langle \gen x, \gen u, u\gen x, \gen t\rangle$:} The commutation relations give $e_5=pt\gen t+[2x+\frac{u^2}{2}]\gen x+ u\gen u + le_4+\kappa e_1$. Thus we have the realization
\[
\mathsf{A}_{5.23}=\langle \gen x, \gen u, u\gen x, \gen t, pt\gen t+[2x+\frac{u^2}{2}]\gen x+ u\gen u\rangle.
\]

\medskip\noindent \underline{$\langle e_1, e_2, e_3, e_4\rangle=\langle \gen x, \gen t, t\gen x-\gen u, \gen t-u\gen x + \kappa \gen u\rangle$:} The commutation relations  $[e_1, e_5]=2e_1,\; [e_2, e_5]=e_2+e_3,\; [e_3, e_5]=e_3$ give $e_5=t\gen t + [2x+\frac{t^2}{2}\gen x + [u-t]\gen u + \beta e_1 + le_4 + (\gamma-l\kappa)\gen u$. Then $[e_4, e_5]=pe_4$ leads to the equation $(\gamma - l\kappa)\gen x + \gen t-u\gen x + [\kappa - 1] \gen u=p[\gen t-u\gen x + \kappa \gen u]$. This is not possible to implement. Hence we have no realization in this case.

\medskip\noindent \underline{$\langle e_1, e_2, e_3, e_4\rangle=\langle \gen x, \gen u, t\gen t+u\gen x, t\gen t\rangle$:} Here $[e_2, e_5]=e_2+e_3$ requires $[\gen u, e_5]=x\gen u + t\gen t+u\gen x$ which is clearly not possible.

\subsection{Realizations of $\mathsf{A}_{5.23}$:}

\begin{align*}
\mathsf{A}_{5.23}&=\langle \gen x, \gen t, (t+u)\gen x, \gen t-\gen u, t\gen t + [2x+\frac{(t+u)^2}{2}]\gen x + u\gen u\rangle\\
\mathsf{A}_{5.23}&=\langle \gen x, \gen u, u\gen x, \gen t, pt\gen t+[2x+\frac{u^2}{2}]\gen x+ u\gen u\rangle.
\end{align*}

\bigskip\noindent \underline{$\mathsf{A}_{5.24}$:} In this case, $\langle e_1, e_2, e_3, e_4\rangle=\mathsf{A}_{3.3}\oplus\mathsf{A}_1$ as well as $[e_1, e_5]=2e_1,\; [e_2, e_5]=e_2+e_3,\; [e_3, e_5]=e_3,\; [e_4, e_5]=\epsilon e_1+2e_4,\; \epsilon=\pm 1$.

\medskip\noindent \underline{$\langle e_1, e_2, e_3, e_4\rangle=\langle \gen u, x\gen u, -(\gen t + \gen x), \gen t\rangle$:} Here $[e_2, e_5]=e_2+e_3$ requires $[x\gen u, e_5]=x\gen u-(\gen t+\gen x)$ which is clearly not possible.

\medskip\noindent \underline{$\langle e_1, e_2, e_3, e_4\rangle=\langle \gen u, x\gen u, -\gen x, \gen t\rangle$:} Here $[e_1, e_5]=2e_1$ gives $e_5=a(t)\gen t + b(t,x)\gen x + [2u+c(t,x)]\gen u$ and then $[e_2, e_5]=e_2+e_3$ gives $[2x-b(t,x)]\gen u=-\gen x+x\gen u$ which is clearly not possible.

\medskip\noindent \underline{$\langle e_1, e_2, e_3, e_4\rangle=\langle \gen x, \gen t, (t+u)\gen x, \gen t-\gen u\rangle$:} The commutation relations $[e_1, e_5]=2e_1,\, [e_2, e_5]=e_2+e_3,\; [e_3, e_5]=e_3$ give  $e_5=t\gen t + [2x+\frac{t^2}{2}+b(u)]\gen x + u\gen u + le_4$. Then $[e_4, e_5]=\epsilon e_1+2e_4$ gives $[u-b'(u)]\gen x=\epsilon\gen x + \gen t-\gen u$ and this is clearly not possible. Hence we have no realization in this case.

\medskip\noindent \underline{$\langle e_1, e_2, e_3, e_4\rangle=\langle \gen x, \gen u, u\gen x, \gen t\rangle$:} The commutation relations give $e_5=2t\gen t+[2x+\frac{u^2}{2}+\epsilon t]\gen x+ u\gen u + le_4+\kappa e_1$. Thus we have the realization
\[
\mathsf{A}_{5.24}=\langle \gen x, \gen u, u\gen x, \gen t, 2t\gen t+[2x+\frac{u^2}{2}+\epsilon t]\gen x+ u\gen u\rangle.
\]

\medskip\noindent \underline{$\langle e_1, e_2, e_3, e_4\rangle=\langle \gen x, \gen t, t\gen x-\gen u, \gen t-u\gen x + \kappa \gen u\rangle$:} The commutation relations  $[e_1, e_5]=2e_1,\; [e_2, e_5]=e_2+e_3,\; [e_3, e_5]=e_3$ give $e_5=t\gen t + [2x+\frac{t^2}{2}\gen x + [u-t]\gen u + \beta e_1 + le_4 + (\gamma-l\kappa)\gen u$. Then $[e_4, e_5]=pe_4$ leads to the equation $(\gamma - l\kappa)\gen x + \gen t-u\gen x + [\kappa - 1] \gen u=\epsilon \gen x +2[\gen t-u\gen x + \kappa \gen u]$. This is not possible to implement. Hence we have no realization in this case.

\medskip\noindent \underline{$\langle e_1, e_2, e_3, e_4\rangle=\langle \gen x, \gen u, t\gen t+u\gen x, t\gen t\rangle$:} Here $[e_2, e_5]=e_2+e_3$ requires $[\gen u, e_5]=x\gen u + t\gen t+u\gen x$ which is clearly not possible.

\subsection{Realizations of $\mathsf{A}_{5.24}$:}

\[
\mathsf{A}_{5.24}=\langle \gen x, \gen u, u\gen x, \gen t, 2t\gen t+[2x+\frac{u^2}{2}+\epsilon t]\gen x+ u\gen u\rangle,\;\; \epsilon=\pm 1.
\]

\bigskip\noindent \underline{$\mathsf{A}_{5.25}$:} In this case, $\langle e_1, e_2, e_3, e_4\rangle=\mathsf{A}_{3.3}\oplus\mathsf{A}_1$ as well as $[e_1, e_5]=2pe_1,\; [e_2, e_5]=pe_2+e_3,\; [e_3, e_5]=-e_2+pe_3,\; [e_4, e_5]=qe_4,\; p\in \mathbb{R},\; q\neq 0$.

\medskip\noindent \underline{$\langle e_1, e_2, e_3, e_4\rangle=\langle \gen u, x\gen u, -(\gen t + \gen x), \gen t\rangle$:} Here $[e_2, e_5]=pe_2+e_3$ requires $[x\gen u, e_5]=px\gen u-(\gen t+\gen x)$ which is clearly not possible.

\medskip\noindent \underline{$\langle e_1, e_2, e_3, e_4\rangle=\langle \gen u, x\gen u, -\gen x, \gen t\rangle$:} Here $[e_1, e_5]=2pe_1$ gives $e_5=a(t)\gen t + b(t,x)\gen x + [2pu+c(t,x)]\gen u$ and then $[e_2, e_5]=pe_2+e_3$ gives $[2px-b(t,x)]\gen u=-\gen x+px\gen u$ which is clearly not possible.

\medskip\noindent \underline{$\langle e_1, e_2, e_3, e_4\rangle=\langle \gen x, \gen t, (t+u)\gen x, \gen t-\gen u\rangle$:} The relation $[e_3, e_5]=-e_2+pe_3$ requires $[(t+u)\gen x, e_5]=-\gen t + p(t+u)\gen x$ and this is clearly not possible. Hence we have no realization in this case.

\medskip\noindent \underline{$\langle e_1, e_2, e_3, e_4\rangle=\langle \gen x, \gen u, u\gen x, \gen t\rangle$:} The commutation relations $[e_1, e_5]=2pe_1,\; [e_2, e_5]=pe_2+e_3$ give $e_5=a(t)\gen t+[2px+\frac{u^2}{2}+b(t)]\gen x+ [pu+c(t)]\gen u$. Then $[e_3, e_5]=-e_2+pe_3$ gives $[pu-c(t)]\gen x=-\gen u + pu\gen x$ which is not possible. Thus we have nos realization.

\medskip\noindent \underline{$\langle e_1, e_2, e_3, e_4\rangle=\langle \gen x, \gen t, t\gen x-\gen u, \gen t-u\gen x + \kappa \gen u\rangle$:} The relation $[e_3, e_5]=-e_2+pe_3$ requires $[t\gen x-\gen u, e_5]=-\gen t + pt\gen x -p\gen u$ which is clearly impossible. Hence we have no realization in this case.

\medskip\noindent \underline{$\langle e_1, e_2, e_3, e_4\rangle=\langle \gen x, \gen u, t\gen t+u\gen x, t\gen t\rangle$:} Here $[e_2, e_5]=pe_2+e_3$ requires $[\gen u, e_5]=p\gen u + t\gen t+u\gen x$ which is clearly not possible. Hence we have no realization in this case.

\bigskip\noindent \underline{$\mathsf{A}_{5.26}$:} In this case, $\langle e_1, e_2, e_3, e_4\rangle=\mathsf{A}_{3.3}\oplus\mathsf{A}_1$ as well as $[e_1, e_5]=2pe_1,\; [e_2, e_5]=pe_2+e_3,\; [e_3, e_5]=-e_2+pe_3,\; [e_4, e_5]=\epsilon e_1 + 2pe_4,\; p\in \mathbb{R},\; \epsilon =\pm 1$.

\medskip\noindent \underline{$\langle e_1, e_2, e_3, e_4\rangle=\langle \gen u, x\gen u, -(\gen t + \gen x), \gen t\rangle$:} Here $[e_2, e_5]=pe_2+e_3$ requires $[x\gen u, e_5]=px\gen u-(\gen t+\gen x)$ which is clearly not possible.

\medskip\noindent \underline{$\langle e_1, e_2, e_3, e_4\rangle=\langle \gen u, x\gen u, -\gen x, \gen t\rangle$:} Here $[e_1, e_5]=2pe_1$ gives $e_5=a(t)\gen t + b(t,x)\gen x + [2pu+c(t,x)]\gen u$ and then $[e_2, e_5]=pe_2+e_3$ gives $[2px-b(t,x)]\gen u=-\gen x+px\gen u$ which is clearly not possible.

\medskip\noindent \underline{$\langle e_1, e_2, e_3, e_4\rangle=\langle \gen x, \gen t, (t+u)\gen x, \gen t-\gen u\rangle$:} The relation $[e_3, e_5]=-e_2+pe_3$ requires $[(t+u)\gen x, e_5]=-\gen t + p(t+u)\gen x$ and this is clearly not possible. Hence we have no realization in this case.

\medskip\noindent \underline{$\langle e_1, e_2, e_3, e_4\rangle=\langle \gen x, \gen u, u\gen x, \gen t\rangle$:} The commutation relations $[e_1, e_5]=2pe_1,\; [e_2, e_5]=pe_2+e_3$ give $e_5=a(t)\gen t+[2px+\frac{u^2}{2}+b(t)]\gen x+ [pu+c(t)]\gen u$. Then $[e_3, e_5]=-e_2+pe_3$ gives $[pu-c(t)]\gen x=-\gen u + pu\gen x$ which is not possible. Thus we have nos realization.

\medskip\noindent \underline{$\langle e_1, e_2, e_3, e_4\rangle=\langle \gen x, \gen t, t\gen x-\gen u, \gen t-u\gen x + \kappa \gen u\rangle$:} The relation $[e_3, e_5]=-e_2+pe_3$ requires $[t\gen x-\gen u, e_5]=-\gen t + pt\gen x -p\gen u$ which is clearly impossible. Hence we have no realization in this case.

\medskip\noindent \underline{$\langle e_1, e_2, e_3, e_4\rangle=\langle \gen x, \gen u, t\gen t+u\gen x, t\gen t\rangle$:} Here $[e_2, e_5]=pe_2+e_3$ requires $[\gen u, e_5]=p\gen u + t\gen t+u\gen x$ which is clearly not possible. Hence we have no realization in this case.

\bigskip\noindent \underline{$\mathsf{A}_{5.27}$:} In this case, $\langle e_1, e_2, e_3, e_4\rangle=\mathsf{A}_{3.3}\oplus\mathsf{A}_1$ as well as $[e_1, e_5]=e_1,\; [e_2, e_5]=0,\; [e_3, e_5]=e_3+e_4,\; [e_4, e_5]=e_1 + e_4$

\medskip\noindent \underline{$\langle e_1, e_2, e_3, e_4\rangle=\langle \gen u, x\gen u, -(\gen t + \gen x), \gen t\rangle$:} From $[e_1, e_5]=e_1,\; [e_2, e_5]=0,\; [e_3, e_5]=e_3+e_4$ we obtain $e_5=a\gen t + x\gen x + [u+c(t-x)]\gen u$ and then $[e_4, e_5]=e_1+e_4$ gives $c'(t-x)\gen u=\gen t+\gen u$ which is clearly impossible. Thus we have no realization in this case.

\medskip\noindent \underline{$\langle e_1, e_2, e_3, e_4\rangle=\langle \gen u, x\gen u, -\gen x, \gen t\rangle$:} Here $[e_3, e_5]=e_3+e_4$ requires $[\gen x,  e_5]=-\gen t + \gen x$ which is clearly not possible. Hence we have no realization in this case.

\medskip\noindent \underline{$\langle e_1, e_2, e_3, e_4\rangle=\langle \gen x, \gen t, (t+u)\gen x, \gen t-\gen u\rangle$:} Here $[e_3, e_5]=e_3+e_4$ requires $[(t+u)\gen x,  e_5]=\gen t + (t+u)\gen x-\gen u$ which is clearly not possible. Hence we have no realization in this case

\medskip\noindent \underline{$\langle e_1, e_2, e_3, e_4\rangle=\langle \gen x, \gen u, u\gen x, \gen t\rangle$:} Here $[e_3, e_5]=e_3+e_4$ requires $[u\gen x,  e_5]=\gen t + u\gen x$ which is clearly not possible. Hence we have no realization in this case

\medskip\noindent \underline{$\langle e_1, e_2, e_3, e_4\rangle=\langle \gen x, \gen t, t\gen x-\gen u, \gen t-u\gen x + \kappa \gen u\rangle$:} Here $[e_3, e_5]=e_3+e_4$ requires $[t\gen x-\gen u,  e_5]=\gen t + (t-u)\gen x+(\kappa-1)\gen u$ which is clearly not possible. Hence we have no realization in this case

\medskip\noindent \underline{$\langle e_1, e_2, e_3, e_4\rangle=\langle \gen x, \gen u, t\gen t+u\gen x, t\gen t\rangle$:} We have $[e_4, e_5]=e_1+e_4$ which gives $[t\gen t, e_5]=t\gen t +\gen x$ and we also have $[e_3, e_5]=e_3+e_4$ which gives $[t\gen t + u\gen x, e_5]=2t\gen t+ u\gen x$. Combining these two results we find that $[u\gen x, e_5]=t\gen t+(u-1)\gen x$ which is clearly impossible. Hence we have no realizations in this case.

\bigskip\noindent \underline{$\mathsf{A}_{5.28}$:} In this case, $\langle e_1, e_2, e_3, e_4\rangle=\mathsf{A}_{3.3}\oplus\mathsf{A}_1$ as well as $[e_1, e_5]=(p+1)e_1,\; [e_2, e_5]=pe_2,\; [e_3, e_5]=e_3+e_4,\; [e_4, e_5]=e_4,\; p\in \mathbb{R}$.

\medskip\noindent \underline{$\langle e_1, e_2, e_3, e_4\rangle=\langle \gen u, x\gen u, -(\gen t + \gen x), \gen t\rangle$:} From $[e_4, e_5]=e_4$ yields $e_5=(t+l)\gen t + b(x,u)\gen x + c(x,u)\gen u$ and then $[e_3, e_5]=e_3+e_4$ gives $\gen x=\gen t+b_x\gen x + c_x\gen u$ which is clearly impossible. Thus we have no realization in this case.

\medskip\noindent \underline{$\langle e_1, e_2, e_3, e_4\rangle=\langle \gen u, x\gen u, -\gen x, \gen t\rangle$:} Here $[e_3, e_5]=e_3+e_4$ requires $[\gen x,  e_5]=-\gen t + \gen x$ which is clearly not possible. Hence we have no realization in this case.

\medskip\noindent \underline{$\langle e_1, e_2, e_3, e_4\rangle=\langle \gen x, \gen t, (t+u)\gen x, \gen t-\gen u\rangle$:} Here $[e_3, e_5]=e_3+e_4$ requires $[(t+u)\gen x,  e_5]=\gen t + (t+u)\gen x-\gen u$ which is clearly not possible. Hence we have no realization in this case

\medskip\noindent \underline{$\langle e_1, e_2, e_3, e_4\rangle=\langle \gen x, \gen u, u\gen x, \gen t\rangle$:} Here $[e_3, e_5]=e_3+e_4$ requires $[u\gen x,  e_5]=\gen t + u\gen x$ which is clearly not possible. Hence we have no realization in this case

\medskip\noindent \underline{$\langle e_1, e_2, e_3, e_4\rangle=\langle \gen x, \gen t, t\gen x-\gen u, \gen t-u\gen x + \kappa \gen u\rangle$:} Here $[e_3, e_5]=e_3+e_4$ requires $[t\gen x-\gen u,  e_5]=\gen t + (t-u)\gen x+(\kappa-1)\gen u$ which is clearly not possible. Hence we have no realization in this case

\medskip\noindent \underline{$\langle e_1, e_2, e_3, e_4\rangle=\langle \gen x, \gen u, t\gen t+u\gen x, t\gen t\rangle$:} We have $[e_4, e_5]=e_4$ gives $[t\gen t, e_5]=t\gen t$ and we also have $[e_3, e_5]=e_3+e_4$ which gives $[t\gen t + u\gen x, e_5]=2t\gen t+ u\gen x$. Combining these two results we find that $[u\gen x, e_5]=t\gen t+u\gen x$ which is clearly impossible. Hence we have no realizations in this case.

\bigskip\noindent \underline{$\mathsf{A}_{5.29}$:} In this case, $\langle e_1, e_2, e_3, e_4\rangle=\mathsf{A}_{3.3}\oplus\mathsf{A}_1$ as well as $[e_1, e_5]=e_1,\; [e_2, e_5]=e_2,\; [e_3, e_5]=e_4,\; [e_4, e_5]=0$.

\medskip\noindent \underline{$\langle e_1, e_2, e_3, e_4\rangle=\langle \gen u, x\gen u, -(\gen t + \gen x), \gen t\rangle$:} From $[e_4, e_5]=e_4$ we have $[\gen t, e_5]=0$ and then we have $[e_3, e_5]=e_4$ which gives $[\gen t + \gen x, e_5]=-\gen t$. combining these two relations, we have $[\gen x, e_5]=-\gen t$ which is impossible. Hence we have no realization in this case.

\medskip\noindent \underline{$\langle e_1, e_2, e_3, e_4\rangle=\langle \gen x, \gen t, (t+u)\gen x, \gen t-\gen u\rangle$:} Here $[e_3, e_5]=e_4$ requires $[(t+u)\gen x,  e_5]=\gen t$ which is clearly not possible. Hence we have no realization in this case

\medskip\noindent \underline{$\langle e_1, e_2, e_3, e_4\rangle=\langle \gen x, \gen u, u\gen x, \gen t\rangle$:} Here $[e_3, e_5]=e_4$ requires $[u\gen x,  e_5]=\gen t$ which is clearly not possible. Hence we have no realization in this case

\medskip\noindent \underline{$\langle e_1, e_2, e_3, e_4\rangle=\langle \gen x, \gen t, t\gen x-\gen u, \gen t-u\gen x + \kappa \gen u\rangle$:} Here $[e_3, e_5]=e_4$ requires $[t\gen x-\gen u,  e_5]=\gen t-u\gen x + \kappa\gen u$ which is clearly not possible. Hence we have no realization in this case.

\medskip\noindent \underline{$\langle e_1, e_2, e_3, e_4\rangle=\langle \gen x, \gen u, t\gen t+u\gen x, t\gen t\rangle$:} We have $[e_4, e_5]=0$ gives $[t\gen t, e_5]=0$ and we also have $[e_3, e_5]=e_4$ which gives $[t\gen t + u\gen x, e_5]=t\gen t$. Combining these two results we find that $[u\gen x, e_5]=t\gen t$ which is clearly impossible. Hence we have no realizations in this case.

\bigskip\noindent\underline{$\mathsf{A}_{5.30}$:} For this algebra we have $\langle e_1, e_2, e_3, e_4\rangle=\mathsf{A}_{4.1}$ and $[e_1, e_5]=(p+2)e_1,\, [e_2, e_5]=(p+1)e_2,\; [e_3, e_5]=pe_3,\; [e_4, e_5]=e_4,\; p\in \mathbb{R}$.

\medskip\noindent $\langle e_1, e_2, e_3, e_4\rangle=\langle \gen x, \gen u , \gen t, u\gen x + t\gen u\rangle:$ The commutation relations give $e_5=pt\gen t+(p+2)x\gen x + (p+1)u\gen u + \kappa e_1$ so we have the realization
\[
\mathsf{A}_{5.30}=\langle \gen x, \gen u , \gen t, u\gen x + t\gen u, pt\gen t+(p+2)x\gen x + (p+1)u\gen u\rangle.
\]

\medskip\noindent $\langle e_1, e_2, e_3, e_4\rangle=\langle \gen u, x\gen u , \gen t, -\gen x + tx\gen u\rangle:$ The commutation relations give $e_5=(pt+l)\gen t+x\gen x + [(p+2)u -\frac{lx^2}{2} + \gamma]\gen u$. The residual equivalence group is $\mathscr{E}(e_1, e_2, e_3, e_4): t'=t+k,\; x'=x,\, u'=u-\frac{kx^2}{2} + \lambda$ and under such a transformation $e_5$ is mapped to $e'_5=[pt'+l-pk]\gen {t'} + x'\gen {x'} + [(p+2)u'+(pk-l)\frac{x^2}{2} + \gamma - (p+2)\lambda]\gen {u'}$. For $p\neq 0$ we choose $k$ so that $pk-l=0$ and then $e'_5=pt'\gen {t'} + x'\gen {x'} + [(p+2)u'+\gamma - (p+2)\lambda]\gen {u'}$, giving $e_5=pt\gen t + x\gen x + (p+2)u\gen u + \alpha e_1$ with $\alpha=\gamma-(p+2)\lambda$, so that we have the realization
\[
\mathsf{A}_{5.30}=\langle \gen u, x\gen u , \gen t, -\gen x + tx\gen u, pt\gen t + x\gen x + (p+2)u\gen u\rangle,\; p\neq 0.
\]
If $p=0$ then we have $e_5=l\gen t + x\gen x + [2u - \frac{lx^2}{2} + \gamma]\gen u$ and we find that we must have $l\neq 0$ since $l=0$ gives $F=0$ in the equation for $F$. Thus, writing $l=2\alpha$, we have $e_5=x\gen x + 2u\gen u + \alpha(2\gen t -x^2\gen u)$ with $\alpha\neq 0$ when $p=0$, and we have the realization
\[
\mathsf{A}_{5.30}=\langle \gen u, x\gen u , \gen t, -\gen x + tx\gen u, x\gen x + 2u\gen u + \alpha(2\gen t -x^2\gen u)\rangle,\; \alpha\neq 0.
\]

\subsection{Realizations of $\mathsf{A}_{5.30}$:}

\begin{align*}
\mathsf{A}_{5.30}&=\langle \gen x, \gen u , \gen t, u\gen x + t\gen u, pt\gen t+(p+2)x\gen x + (p+1)u\gen u\rangle\\
\mathsf{A}_{5.30}&=\langle \gen u, x\gen u , \gen t, -\gen x + tx\gen u, pt\gen t + x\gen x + (p+2)u\gen u\rangle,\; p\neq 0\\
\mathsf{A}_{5.30}&=\langle \gen u, x\gen u , \gen t, -\gen x + tx\gen u, x\gen x + 2u\gen u + \alpha(2\gen t -x^2\gen u)\rangle,\; \alpha\neq 0.
\end{align*}

\bigskip\noindent\underline{$\mathsf{A}_{5.31}$:} For this algebra we have $\langle e_1, e_2, e_3, e_4\rangle=\mathsf{A}_{4.1}$ and $[e_1, e_5]=3e_1,\, [e_2, e_5]=2e_2,\; [e_3, e_5]=e_3,\; [e_4, e_5]=e_3+e_4$.

\medskip\noindent $\langle e_1, e_2, e_3, e_4\rangle=\langle \gen x, \gen u , \gen t, u\gen x + t\gen u\rangle:$ The commutation relations give $[e_4, e_5]=e_3+e_4$ requires $[u\gen x + t\gen u, e_5]=\gen t + u\gen x + t\gen u$ which is clearly impossible. hence we have no realization in this case.

\medskip\noindent $\langle e_1, e_2, e_3, e_4\rangle=\langle \gen u, x\gen u , \gen t, -\gen x + tx\gen u\rangle:$ The commutation relations give $[e_4, e_5]=e_3+e_4$ requires $[-\gen x + tx\gen u, e_5]=\gen t - \gen x + tx\gen u$ which is clearly impossible. hence we have no realization in this case.

\bigskip\noindent\underline{$\mathsf{A}_{5.32}$:} For this algebra we have $\langle e_1, e_2, e_3, e_4\rangle=\mathsf{A}_{4.1}$ and $[e_1, e_5]=e_1,\, [e_2, e_5]=e_2,\; [e_3, e_5]=pe_1+e_3,\; [e_4, e_5]=0$.

\medskip\noindent $\langle e_1, e_2, e_3\rangle=\langle \gen t, \gen x ,\gen u\rangle:$ The commutation relations give $[e_4, e_5]=e_3+e_4$ requires $[-\gen x + tx\gen u, e_5]=\gen t - \gen x + tx\gen u$ which is clearly impossible. hence we have no realization in this case.

\medskip\noindent\underline{$\langle e_1, e_2, e_3, e_4\rangle=\langle \gen u, x\gen u, \gen t, -\gen x + tx\gen u\rangle$:} From $[e_1, e_5]=e_1,\, [e_2, e_5]=e_2,\; [e_3, e_5]=pe_1+e_3$ we find $e_5=(t+l)\gen t + 2x\gen x + [u+pt+c(x)]\gen u$ and then $[e_4, e_5]=0$ gives $2\gen x + [2tx+c'(x) + lx]\gen u=0$ which is impossible, so we have no realization in this case.

\medskip\noindent\underline{$\langle e_1, e_2, e_3, e_4\rangle=\langle \gen x, \gen u, \gen t, u\gen x + t\gen u\rangle$:} From the commutation relations we find $e_5=t\gen t + [x+pt]\gen x + u\gen u +\beta e_1$ and we then have the realization
\[
\mathsf{A}_{5.32}=\langle \gen x, \gen u, \gen t, u\gen x + t\gen u, t\gen t + [x+pt]\gen x + u\gen u\rangle
\]

\bigskip\noindent\underline{$\mathsf{A}_{5.33}$:} Here we have $\langle e_1, e_2, e_3\rangle =3\mathsf{A}_1$ and $[e_1, e_4]=e_1,\, [e_2, e_4]=0,\, [e_3, e_4]=pe_3$ as well as $[e_1, e_5]=0,\; [e_2, e_5]=e_2, \; [e_3, e_5]=qe_3,\, [e_4, e_5]=0$ and $p^2+q^2\neq 0$.

\medskip\noindent\underline{$\langle e_1, e_2, e_3\rangle=\langle \gen t, \gen x, \gen u\rangle$:} From the commutation relations $[e_1, e_4]=e_1,\, [e_2, e_4]=0,\, [e_3, e_4]=pe_3$ as well as $[e_1, e_5]=0,\; [e_2, e_5]=e_2, \; [e_3, e_5]=qe_3$ we find that $e_4=(t+l)\gen t + \beta \gen x + [pu+\gamma]\gen u$ and $e_5=\alpha \gen t + [x+\kappa]\gen x + [qu+\lambda]\gen u$. Then $[e_4, e_5]=0$ gives $\alpha=\beta=0$ and $p\lambda-q\gamma=0$ from which we see that we must have $(\gamma, \lambda)=\sigma(p,q)$ and thus $e_4=(t+l)\gen t + p[u+\sigma]\gen u$ and $e_5=[x+\kappa]\gen x + q[u+\sigma]\gen u$. We have the residual equivalence group $\mathscr{E}(e_1, e_2, e_3): t'=t+k,\; x'=x+\mu,\; u'=u+\nu$ and under such a transformation $e_4$ is mapped to $e'_4=[t'+l-k]\gen {t'} + p[u'+\sigma - \nu]\gen {u'}$ and $e_5$ is mapped to $e'_5=[x'+\kappa-\mu]\gen {x'} + q[u'+\sigma-\nu]\gen {u'}$. Then we choose $k=l,\; \mu=\kappa,\; \mu=\sigma$ giving $e_4=t\gen t + p\gen u$ and $e_5=x\gen x + qu\gen u$ in canonical form. Hence we have the realization

\[
\mathsf{A}_{5.33}=\langle \gen t, \gen x, \gen u, t\gen t + p\gen u, x\gen x + qu\gen u\rangle,\; p\neq 0.
\]
Note that $p\neq 0$ for otherwise $\gen t$ and $t\gen t$ would be symmetries and together they give $F=0$, a contradiction. In a similar way we also have the realizations
\[
\mathsf{A}_{5.33}=\langle \gen x, \gen t, \gen u, x\gen x + pu\gen u, t\gen t + qu\gen u\rangle,\; q\neq 0.
\]
and
\[
\mathsf{A}_{5.33}=\langle \gen x, \gen u, \gen t, pt\gen t + x\gen x, qt\gen t + u\gen u\rangle.
\]
So we have the following realizations:
\begin{align*}
\mathsf{A}_{5.33}&=\langle \gen t, \gen x, \gen u, t\gen t + pu\gen u, x\gen x + qu\gen u\rangle,\; p\neq 0\\
\mathsf{A}_{5.33}&=\langle \gen x, \gen t, \gen u, x\gen x + pu\gen u, t\gen t + qu\gen u\rangle, \; q\neq 0\\
\mathsf{A}_{5.33}&=\langle \gen x, \gen u, \gen t, pt\gen t + x\gen x, qt\gen t + u\gen u\rangle.
\end{align*}

\bigskip\noindent\underline{$\langle e_1, e_2, e_3\rangle=\langle \gen t, \gen u, x\gen u\rangle$:} From the commutation relations we find $e_4=(t+l)\gen t -px\gen x + c(x)\gen u$ and $e_5=(1-q)x\gen x + [u+\gamma(x)]\gen u$ with $(1-q)xc'(x)-c(x)=-px\gamma'(x)$. We have the residual equivalence group $\mathscr{E}(e_1, e_2, e_3): t'=t+k,\; x'=x,\; u'=u+U(x)$. Under such a transformation $e_5$ is mapped to $e'_5=(1-q)x'\gen {x'}+ [u'+\gamma(x)-U(x)+(1-q)xU'(x)]\gen {u'}$ and $e_4$ is mapped to $e'_4=[t'+l-k]\gen {t'}-px'\gen {x'} + [c(x)-pxU'(x)]\gen {u'}$. We may choose $k=l$ and we may always choose $U(x)$ so that $\gamma(x)-U(x)+(1-q)xU'(x)=0$, and then the canonical forms are $e_4=t\gen t -px\gen x + c(x)\gen u$ and $e_5=(1-q)x\gen x + u\gen u$. Then the condition $[e_4, e_5]=0$ gives $(1-q)xc'(x)-c(x)=0$ and this gives $c(x)=\kappa x^{1/(1-q)}$ when $q\neq 1$. In fact, we consider only $q\neq 1$ for if $q=1$ then the algebra contains the solvable subalgebra $\langle \gen u, x\gen u, u\gen u\rangle$ which linearizes the evolution equation. Thus $q\neq 1$ and we have $e_4=t\gen t - px\gen x + \kappa x^{1/(1-q)}\gen u$ and $e_5=(1-q)x\gen x + u\gen u$. We have the residual equivalence group $\mathscr{E}(e_1, e_2, e_3, e_5): t'=t+k,\; x'=x,\; u'=u+\lambda x^{1/(1-q)}$ and under such a transformation $e_4$ is mapped to $e'_4=t'\gen {t'}-px'\gen {x'} + [\kappa -\frac{p\lambda}{(1-q)}]x^{1/(1-q)}\gen {u'}$ on taking $k=0$. For $p\neq 0$ we may always choose $\lambda$ so that $\kappa -\frac{p\lambda}{(1-q)}=0$ giving $e_4=t\gen t -px\gen x$. For $p=0$ we must have $\kappa\neq 0$ for otherwise we will have $\gen t$ and $t\gen t$ as symmetries, which give $F=0$, a contradiction. Thus we have the realizations
\begin{align*}
\mathsf{A}_{5.33}&=\langle \gen t, \gen u, x\gen u, t\gen t - px\gen x, (1-q)x\gen x + u\gen u\rangle,\; q\neq 1,\; p\neq 0\\
\mathsf{A}_{5.33}&=\langle \gen t, \gen u, x\gen u, t\gen t + \kappa x^{1/(1-q)}\gen u, (1-q)x\gen x + u\gen u\rangle,\; \kappa\neq 0,\; q\neq 1.
\end{align*}
In a similar way we obtain the following realizations for $\langle e_1, e_2, e_3\rangle=\langle \gen u, \gen t, x\gen u\rangle$ and $\langle e_1, e_2, e_3\rangle=\langle \gen u, x\gen u, \gen t\rangle$:

\begin{align*}
\mathsf{A}_{5.33}&=\langle \gen u, \gen t, x\gen u, (1-p)x\gen x + u\gen u, t\gen t - qx\gen x\rangle,\; p\neq 1,\; q\neq 0\\
\mathsf{A}_{5.33}&=\langle \gen u, \gen t, x\gen u, (1-p)x\gen x + u\gen u, t\gen t + \kappa x^{1/(1-p)}\gen u\rangle,\; \kappa\neq 0,\; p\neq 1\\
\mathsf{A}_{5.33}&=\langle \gen u, x\gen u, \gen t, pt\gen t + x\gen x + u\gen u, qt\gen t - x\gen x\rangle.
\end{align*}

\subsection{Realizations of $\mathsf{A}_{5.33}$:}

\begin{align*}
\mathsf{A}_{5.33}&=\langle \gen t, \gen x, \gen u, t\gen t + pu\gen u, x\gen x + qu\gen u\rangle,\; p\neq 0\\
\mathsf{A}_{5.33}&=\langle \gen x, \gen t, \gen u, x\gen x + pu\gen u, t\gen t + qu\gen u\rangle,\; q\neq 0\\
\mathsf{A}_{5.33}&=\langle \gen x, \gen u, \gen t, pt\gen t + x\gen x, qt\gen t + u\gen u\rangle\\
\mathsf{A}_{5.33}&=\langle \gen t, \gen u, x\gen u, t\gen t - px\gen x, (1-q)x\gen x + u\gen u\rangle,\; q\neq 1,\; p\neq 0\\
\mathsf{A}_{5.33}&=\langle \gen t, \gen u, x\gen u, t\gen t + \kappa x^{1/(1-q)}\gen u, (1-q)x\gen x + u\gen u\rangle,\; \kappa\neq 0,\; q\neq 1\\
\mathsf{A}_{5.33}&=\langle \gen u, \gen t, x\gen u, (1-p)x\gen x + u\gen u, t\gen t - qx\gen x\rangle,\; p\neq 1,\; q\neq 0\\
\mathsf{A}_{5.33}&=\langle \gen u, \gen t, x\gen u, (1-p)x\gen x + u\gen u, t\gen t + \kappa x^{1/(1-p)}\gen u\rangle,\; \kappa\neq 0,\; p\neq 1\\
\mathsf{A}_{5.33}&=\langle \gen u, x\gen u, \gen t, pt\gen t + x\gen x + u\gen u, qt\gen t - x\gen x\rangle.
\end{align*}

\bigskip\noindent \underline{$\mathsf{A}_{5.34}$:} In this case, $\langle e_1, e_2, e_3\rangle=3\mathsf{A}_1$ as well as $[e_1, e_4]=pe_1,\;  [e_2, e_4]=e_2,\; [e_3, e_4]=e_3,\; [e_1, e_5]=e_1,\; [e_2, e_5]=0,\;[e_3, e_5]=e_2,\; [e_4, e_5]=0$.

\medskip\noindent Note that since $[e_3, e_5]=e_2$ the cases $\langle e_1, e_2, e_3\rangle=\langle \gen x, \gen t, \gen u\rangle$ and $\langle e_1, e_2, e_3\rangle=\langle \gen u, \gen t, x\gen u\rangle$ are inadmissible.

\medskip\noindent\underline{$\langle e_1, e_2, e_3\rangle=\langle \gen t, \gen x, \gen u\rangle$:} From the commutation relations $[e_1, e_4]=pe_1,\;  [e_2, e_4]=e_2,\; [e_3, e_4]=e_3,\; [e_1, e_5]=e_1,\; [e_2, e_5]=0,\;[e_3, e_5]=e_2$ we obtain $e_4=(pt+l)\gen t + (x+\beta)\gen x + (u+\gamma)\gen u$ and $e_5=(t+\alpha)\gen t + (u+\lambda)\gen x + \mu\gen u$. Then $[e_4, e_5]=0$ gives $\mu=0,\; \lambda=\gamma$ and $l=p\alpha$. Then $e_4=p(t+\alpha)\gen t + (x+\beta)\gen x + (u+\gamma)\gen u$ and $e_5=(t+\alpha)\gen t + (u+\gamma)\gen x$. We have the residual equivalence group $\mathscr{E}(e_1, e_2, e_3): t'=t+k,\; x'=x + \kappa,\; u'=u + \nu$ and under such a transformation $e_4$ and $e_5$ are mapped to $e'_4=p(t'+\alpha-k)\gen {t'} + (x'+\beta-\kappa)\gen {x'} + (u'+\gamma-\nu)\gen {u'}$ and $e'_5=(t'+\alpha-k)\gen {t'} + (u'+\gamma-\nu)\gen {x'}$. Choose $k=\alpha$ and $\nu=\gamma$ so that $e_4=pt\gen t + x\gen x + u\gen u$ and $e_5=t\gen t + u\gen x$ in canonical form. Thus we have the realization
\[
\mathsf{A}_{5.34}=\langle \gen t, \gen x, \gen u, pt\gen t + x\gen x + u\gen u, t\gen t + u\gen x\rangle.
\]

\medskip\noindent\underline{$\langle e_1, e_2, e_3\rangle=\langle \gen x, \gen u, \gen t\rangle$:} From the commutation relations $[e_1, e_4]=pe_1,\;  [e_2, e_4]=e_2,\; [e_3, e_4]=e_3,\; [e_1, e_5]=e_1,\; [e_2, e_5]=0,\;[e_3, e_5]=e_2$ we obtain $e_4=(t+l)\gen t + (px+\beta)\gen x + (u+\gamma)\gen u$ and $e_5=\alpha\gen t + (x+\lambda)\gen x + (t+\mu)\gen u$. Then $[e_4, e_5]=0$ gives $\alpha=0,\; \mu=l,\; p\lambda=\beta$ . Then $e_4=(t+l)\gen t + p(x+\lambda)\gen x + (u+\gamma)\gen u$ and $e_5=(x+\lambda)\gen x + (t+l)\gen u$. We have the residual equivalence group $\mathscr{E}(e_1, e_2, e_3): t'=t+k,\; x'=x + \kappa,\; u'=u + \nu$ and under such a transformation $e_4$ and $e_5$ are mapped to $e'_4=(t'+l-k)\gen {t'} + p(x'+\lambda-\kappa)\gen {x'} + (u'+\gamma-\nu))\gen {u'}$ and $e'_5=(t'+l-k)\gen {t'} + (x'+\lambda-\kappa)\gen {x'} +(t'+l-k)\gen {u'}$. Choose $k=l,\; \nu=\gamma$ and $\kappa=\lambda$ so that $e_4=t\gen t + px\gen x + u\gen u$ and $e_5=x\gen x + t\gen u$ in canonical form. Thus we have the realization
\[
\mathsf{A}_{5.34}=\langle \gen x, \gen u, \gen t, t\gen t + px\gen x + u\gen u, x\gen x + t\gen u\rangle.
\]

\medskip\noindent\underline{$\langle e_1, e_2, e_3\rangle=\langle \gen t, \gen u, x\gen u\rangle$:} From the commutation relations $[e_1, e_4]=pe_1,\;  [e_2, e_4]=e_2,\; [e_3, e_4]=e_3,\; [e_1, e_5]=e_1,\; [e_2, e_5]=0,\;[e_3, e_5]=e_2$ we obtain $e_4=(pt+l)\gen t + u\gen u$ and $e_5=(t+\alpha)\gen t - \gen x + c(x)\gen u$. Then $[e_4, e_5]=0$ gives $c(x)=0$ and $l=p\alpha$, so that $e_4=p(t+\alpha)\gen t + u\gen u$ and $e_5=(t+\alpha)\gen t - \gen x$. Since $\gen t$ is already an element of the algebra we obtain the realization
\[
\mathsf{A}_{5.34}=\langle \gen t, \gen u, x\gen u, pt\gen t + u\gen u, t\gen t - \gen x\rangle.
\]

\medskip\noindent\underline{$\langle e_1, e_2, e_3\rangle=\langle \gen u, x\gen u, \gen t\rangle$:} From the commutation relations $[e_1, e_4]=pe_1,\;  [e_2, e_4]=e_2,\; [e_3, e_4]=e_3,\; [e_1, e_5]=e_1,\; [e_2, e_5]=0,\;[e_3, e_5]=e_2$ we obtain $e_4=(t+l)\gen t + (p-1)x\gen x + [pu+c(x)]\gen u$ and $e_5=\alpha \gen t + x\gen x + [u+tx+\gamma(x)]\gen u$. Then $[e_4, e_5]=0$ gives $\alpha=0$ and $xc'(x)-c(x)+p\gamma(x)-(p-1)x\gamma'(x)=lx$. Thus we have $e_5=x\gen x + [u+tx+\gamma(x)]\gen u$. We have the residual equivalence group $\mathscr{E}(e_1, e_2, e_3): t'=t+k,\; x'=x,\, u'=u+U(x)$. Then under such a transformation $e_5$ is mapped to $e'_5=x'\gen {x'} + [u'+t'x'+\gamma(x)-kx-U(x)+ xU'(x)]\gen {u'}$ and we may always choose $U(x)$ so that $\gamma(x)-kx-U(x)+ xU'(x)=0$ giving $e_5=x\gen x + [u+tx]\gen u$ in canonical form. We still have $e_4=(t+l)\gen t + (p-1)x\gen x + [pu+c(x)]\gen u$ and now we have $xc'(x)-c(x)=lx$ from $[e_4, e_5]=0$, so that $c(x)=lx\ln|x|+\sigma x$. We now have the residual equivalence group $\mathscr{E}(e_1, e_2, e_3, e_5): t'=t+k,\; x'=x,\, u'=u+kx\ln|x|+\lambda x$ and under such a transformation $e_4$ is mapped to $e'_4=(t'+l-k)\gen {t'}+ (p-1)x'\gen {x'} + [pu'+(l-k)x\ln|x|+\sigma x+(p-1)kx-\lambda x]\gen {u'}$. We choose $k=l$ and $\lambda=(p-1)k+\sigma$ so that $e_4=t\gen t + (p-1)x\gen x + pu\gen u$ in canonical form. We then have the realization
\[
\mathsf{A}_{5.34}=\langle \gen u, x\gen u, \gen t, t\gen t + (p-1)x\gen x + pu\gen u, x\gen x + [u+tx]\gen u\rangle.
\]

\subsection{Realizations of $\mathsf{A}_{5.34}$:}

\begin{align*}
\mathsf{A}_{5.34}&=\langle \gen t, \gen x, \gen u, pt\gen t + x\gen x + u\gen u, t\gen t + u\gen x\rangle\\
\mathsf{A}_{5.34}&=\langle \gen x, \gen u, \gen t, t\gen t + px\gen x + u\gen u, x\gen x + t\gen u\rangle\\
\mathsf{A}_{5.34}&=\langle \gen t, \gen u, x\gen u, pt\gen t + u\gen u, t\gen t - \gen x\rangle\\
\mathsf{A}_{5.34}&=\langle \gen u, x\gen u, \gen t, t\gen t + (p-1)x\gen x + pu\gen u, x\gen x + [u+tx]\gen u\rangle.
\end{align*}

\bigskip\noindent \underline{$\mathsf{A}_{5.35}$:} In this case, $\langle e_1, e_2, e_3\rangle=3\mathsf{A}_1$ as well as $[e_1, e_4]=pe_1,\; [e_2, e_4]=e_2,\; [e_3, e_4]=e_3,\; [e_1, e_5]=qe_1,\; [e_2, e_5]=-e_3,\; [e_3, e_5]=e_2,\; [e_4, e_5]=0 p, q\in \mathbb{R},\; p^2+q^2\neq 0$.

\medskip\noindent Note that $\langle e_1, e_2, e_3\rangle=\langle \gen x, \gen t, \gen u\rangle$ and $\langle e_1, e_2, e_3\rangle=\langle \gen x, \gen u, \gen t\rangle$ are not possible since we require $[e_3, e_5]=e_2$ and $[e_2, e_5]=-e_3$ and these lead to $[\gen u, e_5]=-\gen t$ in the first case, and $[\gen u, e_5]=\gen t$ in the second case, these are impossible to implement with $e_5=a(t)\gen t + b(t,x,u)\gen x + c(t,x,u)\gen u$. the same remarks apply to $\langle e_1, e_2, e_3\rangle=\langle \gen u, \gen t, x\gen u\rangle$ and $\langle e_1, e_2, e_3\rangle=\langle \gen u, x\gen u, \gen t\rangle$.

\medskip\noindent\underline{$\langle e_1, e_2, e_3\rangle=\langle \gen t, \gen x, \gen u\rangle$:} The commutation relations $[e_1, e_4]=pe_1,\; [e_2, e_4]=e_2,\; [e_3, e_4]=e_3,\; [e_1, e_5]=qe_1,\; [e_2, e_5]=-e_3,\; [e_3, e_5]=e_2$ yield $e_4=(pt+l)\gen t + (x+\beta)\gen x + (u+\mu)\gen u$ and $e_5=(qt+m)\gen t + (u+\nu)\gen x + (-x+\kappa)\gen u$. Then $[e_4, e_5]=0$ gives $pm-ql=0$ as well as $\mu=\nu,\, \kappa=-\beta$. Thus $(l,m)=\sigma(p,q)$ since $p^2+q^2\neq 0$ and we then have $e_4=p(t+\sigma)\gen t + (x+\beta)\gen x + (u+\mu)\gen u$ as well as $e_5=q(t+\sigma)\gen t + (u+\mu)\gen x - (x+\beta)\gen u$. We have the residual equivalence group $\mathscr{E}(e_1, e_2, e_3): t'=t+k,\; x'=x + \lambda,\; u'=u + \tau$ and under such a transformation $e_4$ and $e_5$ are mapped to $e'_4=p(t'+\sigma-k)\gen {t'} + (x'+\beta-\lambda)\gen {x'} + (u'+\mu-\tau)\gen {u'}$ and $e_5=q(t'+\sigma-k)\gen {t'} + (u'+\mu-\tau)\gen {x'} - (x'+\beta-\lambda)\gen {u'}$. Then we choose $k=\sigma,\, \lambda=\beta,\; \tau=\mu$ to give $e_4=pt\gen t + x\gen x + u\gen u$ and $e_5=qt\gen t + u\gen x-x\gen u$ in canonical form. Hence we have the realization
\[
\mathsf{A}_{5.35}=\langle \gen t, \gen x, \gen u, pt\gen t + x\gen x + u\gen u, qt\gen t + u\gen x-x\gen u\rangle.
\]

\medskip\noindent\underline{$\langle e_1, e_2, e_3\rangle=\langle \gen t, \gen u, x\gen u\rangle$:} The commutation relations $[e_1, e_4]=pe_1,\; [e_2, e_4]=e_2,\; [e_3, e_4]=e_3,\; [e_1, e_5]=qe_1,\; [e_2, e_5]=-e_3,\; [e_3, e_5]=e_2$ yield $e_4=(pt+l)\gen t + [u+c(x)]\gen u$ and $e_5=(qt+m)\gen t -(1+x^2)\gen x + [\gamma(x)-xu]\gen u$. Then $[e_4, e_5]=0$ gives $ql-pm=0$ and $(1+x^2)c'(x)-xc(x)=\gamma(x)$. So we have $(l,m)=\sigma(p,q)$ since $p^2+q^2\neq 0$. Thus $e_4=p(t+\sigma)\gen t + [u+c(x)]\gen u$ and $e_5=q(t+\sigma)\gen t -(1+x^2)\gen x + [\gamma(x)-xu]\gen u$.  We have the residual equivalence group $\mathscr{E}(e_1, e_2, e_3): t'=t+k,\; x'=x,\; u'=u + U(x)$ and under such a transformation $e_4$ and $e_5$ are mapped to $e'_4=p(t'+\sigma-k)\gen {t'} + [u'+c(x)-U(x)]\gen {u'}$ and $e'_5=q(t'+\sigma-k)\gen {t'} -(1+x'^2)\gen {x'} + [\gamma(x)+xU(x)-(1+x^2)U'(x)-x'u']\gen {u'}$. We choose $k=\sigma$ and $U(x)=c(x)$ so that we have $\gamma(x)+xU(x)-(1+x^2)U'(x)=0$ since $(1+x^2)c'(x)-xc(x)=\gamma(x)$. hence we have $e_4=p\gen t + u\gen u$ and $e_5=qt\gen t - (1+x^2)\gen x -xu\gen u$ in canonical form. We then have the realization
\[
\mathsf{A}_{5.35}=\langle \gen t, \gen u, x\gen u, pt\gen t + u\gen u, qt\gen t - (1+x^2)\gen x-xu\gen u\rangle.
\]

\subsection{Realizations of $\mathsf{A}_{5.35}$:}

\begin{align*}
\mathsf{A}_{5.35}&=\langle \gen t, \gen x, \gen u, pt\gen t + x\gen x + u\gen u, qt\gen t + u\gen x-x\gen u\rangle\\
\mathsf{A}_{5.35}&=\langle \gen t, \gen u, x\gen u, pt\gen t + u\gen u, qt\gen t - (1+x^2)\gen x-xu\gen u\rangle.
\end{align*}

\bigskip\noindent \underline{$\mathsf{A}_{5.36}$:} In this case, $\langle e_1, e_2, e_3, e_4\rangle=3\mathsf{A}_{4.8}$ with $q=0$ as well as $[e_1, e_5]=0,\; [e_2, e_5]=-e_2,\; [e_3, e_5]=e_3,\; [e_4, e_5]=0$.

\medskip\noindent\underline{$\langle e_1, e_2, e_3, e_4\rangle=\langle \gen u, x\gen u, -\gen x, t\gen t + u\gen u\rangle$:} From the commutation relations we have $e_5=pt\gen t + x\gen x + \gamma t\gen u$. We have the residual equivalence group $\mathscr{E}(e_1, e_2, e_3, e_4): t'=kt,\; x'=x,\; u'=u+\lambda t$. Under such a transformation $e_5$ is mapped to $e'_5=pt'\gen {t'} + x'\gen {x'} + [\gamma+p\lambda]t\gen {u'}$. if $p\neq 0$ then we may choose $\lambda$ so that $\gamma + p\lambda=0$ giving $e_5=pt\gen t + x\gen x$ in canonical form. If $p=0$ then either $\gamma=0$ giving $e_5=x\gen x$ or $\gamma\neq 0$ and then $e'_5=x'\gen {x'}+\frac{\gamma}{k}t'\gen {u'}$ and here we choose $k=\gamma$ so that we have $e_5=x\gen x + t\gen u$ in canonical form. Thus we have the realizations
\begin{align*}
\mathsf{A}_{5.36}&=\langle\gen u, x\gen u, -\gen x, t\gen t + u\gen u, pt\gen t + x\gen x\rangle,\; p\in \mathbb{R}\\
\mathsf{A}_{5.36}&=\langle\gen u, x\gen u, -\gen x, t\gen t + u\gen u, x\gen x + t\gen u\rangle.
\end{align*}

\medskip\noindent\underline{$\langle e_1, e_2, e_3, e_4\rangle=\langle \gen u, x\gen u, -\gen x, u\gen u\rangle$:} This linearizes the equation.

\medskip\noindent\underline{$\langle e_1, e_2, e_3, e_4\rangle=\langle \gen u, x\gen u, -(\gen t + \gen x), a\gen t + u\gen u\rangle,\; a\neq 0$:} The commutation relations $[e_1, e_5]=0,\; [e_2, e_5]=-e_2,\; [e_3, e_5]=e_3$ give $e_5=(t+l)\gen t + x\gen x + c(\xi)\gen u$ where $\xi=t-x$. Then $[e_4, e_5]=0$ gives $a=0$ which is a contradiction. Hence we have no realization in this case.

\medskip\noindent\underline{$\langle e_1, e_2, e_3, e_4\rangle=\langle \gen x, \gen t, (t+u)\gen x, t\gen t + x\gen x + u\gen u\rangle$:} The commutation relations give $e_5=-t\gen t + \kappa u\gen x -u\gen u$. We have the residual equivalence group $\mathscr{E}(e_1, e_2, e_3, e_4): t'=t,\; x'=x+\lambda u,\, u'=u$ and under such a transformation $e_5$ is mapped to $e'_5=-t'\gen {t'}+[\kappa+\lambda]u\gen x - u'\gen {u'}$ and on choosing $\lambda=-\kappa$ we may take $e_5=-t\gen t-u\gen u$ in canonical form. Thus we have the realization
\[
\mathsf{A}_{5.36}=\langle\gen x, \gen t, (t+u)\gen x, t\gen t + x\gen x + u\gen u, -t\gen t - u\gen u\rangle.
\]

\medskip\noindent\underline{$\langle e_1, e_2, e_3, e_4\rangle=\langle \gen x, \gen u, u\gen x, t\gen t + x\gen x + u\gen u\rangle$:} From the commutation relations we have $e_5=\alpha t\gen t + \beta t\gen x -u\gen u$. We have the residual equivalence group $\mathscr{E}(e_1, e_2, e_3, e_4): t'=kt,\; x'=x+\lambda t,\, u'=u$ and under such a transformation $e_5$ is mapped to $e'_5=\alpha t'\gen {t'} + [\beta + \alpha\lambda]t\gen {x'} - u'\gen {u'}$. If $\alpha\neq 0$ we choose $\lambda$ so that $\beta+\alpha\lambda=0$ and then $e_5=\alpha t\gen t -u\gen u$ in canonical form. If $\alpha=0$ then $e_5=\beta t\gen x -u\gen u$. For $\beta\neq 0$ we then have $e'_5=\frac{\beta}{k}t'\gen {x'} - u'\gen {u'}$ and we take $k=\beta$ so that $e_5=t\gen x-u\gen u$ in canonical form. If $\alpha=0$ and $\beta=0$ we have $e_5=-u\gen u$. Thus we have the realizations
\begin{align*}
\mathsf{A}_{5.36}&=\langle\gen x, \gen u, u\gen x, t\gen t + x\gen x +u\gen u, \alpha t\gen t - u\gen u\rangle,\; \alpha\in \mathbb{R}\\
\mathsf{A}_{5.36}&=\langle\gen x, \gen u, u\gen x, t\gen t + x\gen x +u\gen u, t\gen x - u\gen u\rangle.
\end{align*}

\medskip\noindent\underline{$\langle e_1, e_2, e_3, e_4\rangle=\langle \gen x, \gen u, u\gen x+t\gen u, t\gen t + x\gen x + u\gen u\rangle$:} The commutation relations give $e_5=-2t\gen t + \beta t\gen x - u\gen u$. We have the residual equivalence group $\mathscr{E}(e_1, e_2, e_3, e_4): t'=t,\; x'=x+\lambda t,\, u'=u$ and under such a transformation $e_5$ is mapped to $e'_5=-2t'\gen {t'} + [\beta -2\lambda]t\gen {x'} - u'\gen {u'}$. We choose $\lambda$ so that $\beta-2\lambda=0$ and this gives $e_5=-2t\gen t-u\gen u$ in canonical form. hence we have the realization
\[
\mathsf{A}_{5.36}=\langle \gen x, \gen u, u\gen x+t\gen u, t\gen t + x\gen x +u\gen u, -2t\gen t - u\gen u\rangle.
\]

\medskip\noindent\underline{$\langle e_1, e_2, e_3, e_4\rangle=\langle \gen x, \gen t, t\gen x-\gen u, t\gen t + x\gen x + c\gen u\rangle$:} The commutation relations give $c=0$ and $e_5=-t\gen t + [u+\gamma]\gen u$. We have the residual equivalence group $\mathscr{E}(e_1, e_2, e_3, e_4): t'=t,\; x'=x,\, u'=u+\sigma$ and under such a transformation $e_5$ is mapped to $e'_5=-t'\gen {t'} + [u'+\gamma-\sigma]\gen {u'}$, so choosing $\sigma=-\gamma$ gives $e_5=-t\gen t+u\gen u$ in canonical form. Hence we have the realization
\[
\mathsf{A}_{5.36}=\langle \gen x, \gen t, t\gen x-\gen u, t\gen t + x\gen x, -t\gen t + u\gen u\rangle.
\]

\medskip\noindent\underline{$\langle e_1, e_2, e_3, e_4\rangle=\langle \gen x, \gen u, t\gen t+u\gen x, \alpha t\gen t + x\gen x + u\gen u\rangle$:} The commutation relations give $\alpha=0$ and $e_5=[t\ln|t|+\kappa t]\gen t -u\gen u$. We have the residual equivalence group $\mathscr{E}(e_1, e_2, e_3, e_4). t'=kt,\, x'=x,\; u'=u$ and under such a transformation $e_5$ is mapped to $e'_5=t'[\ln|t'|+(\kappa-k)]\gen {t'}-u'\gen {u'}$. we choose $k=\kappa$ and so $e_5=t\ln|t|\gen t - u\gen u$ in canonical form. Thus we have the realization
\[
\mathsf{A}_{5.36}=\langle \gen x, \gen u, t\gen t+u\gen x, x\gen x + u\gen u, t\ln|t|\gen t -u\gen u\rangle.
\]
This realization is equivalent to
\[
\mathsf{A}_{5.36}=\langle \gen x, \gen u, \gen t+u\gen x, x\gen x + u\gen u, t\gen t -u\gen u\rangle
\]
under the equivalence transformation $t'=\ln|t|,\, x'=x,\, u'=u$.

\subsection{Realizations of $\mathsf{A}_{5.36}$:}

\begin{align*}
\mathsf{A}_{5.36}&=\langle\gen u, x\gen u, -\gen x, t\gen t + u\gen u, pt\gen t + x\gen x\rangle,\; p\in \mathbb{R}\\
\mathsf{A}_{5.36}&=\langle\gen u, x\gen u, -\gen x, t\gen t + u\gen u, x\gen x + t\gen u\rangle\\
\mathsf{A}_{5.36}&=\langle\gen x, \gen t, (t+u)\gen x, t\gen t + x\gen x + u\gen u, -t\gen t - u\gen u\rangle\\
\mathsf{A}_{5.36}&=\langle\gen x, \gen u, u\gen x, t\gen t + x\gen x +u\gen u, \alpha t\gen t - u\gen u\rangle,\; \alpha\in \mathbb{R}\\
\mathsf{A}_{5.36}&=\langle\gen x, \gen u, u\gen x, t\gen t + x\gen x +u\gen u, t\gen x - u\gen u\rangle.
\\\mathsf{A}_{5.36}&=\langle \gen x, \gen u, u\gen x+t\gen u, t\gen t + x\gen x +u\gen u, -2t\gen t - u\gen u\rangle\\
\mathsf{A}_{5.36}&=\langle \gen x, \gen t, t\gen x-\gen u, t\gen t + x\gen x, -t\gen t + u\gen u\rangle\\
\mathsf{A}_{5.36}&=\langle \gen x, \gen u, \gen t+u\gen x, x\gen x + u\gen u, t\gen t -u\gen u\rangle
\end{align*}

\bigskip\noindent \underline{$\mathsf{A}_{5.37}$:} In this case, $\langle e_1, e_2, e_3\rangle=\mathsf{A}_{4.8}$ with $q=1$ as well as $[e_1, e_5]=0,\; [e_2, e_5]=-e_3,\;[e_3, e_5]=e_2,\; [e_4, e_5]=0$.

\medskip\noindent Because we require $[e_2, e_5]=-e_3$ and $[e_3, e_5]=e_2$ it is impossible to implement the following cases of $\mathsf{A}_{5.37}$:

\begin{align*}
&\langle \gen u, x\gen u, -(\gen t+\gen x), qt\gen t + qx\gen x+(1+q)u\gen u\rangle\\
&\langle \gen x, \gen t, (t+u)\gen x, t\gen t + (1+q)x\gen x+u\gen u\rangle\\
&\langle \gen x, \gen t, t\gen x-\gen u, t\gen t + (1+q)x\gen x+qu\gen u\rangle\\
&\langle \gen x, \gen u, t\gen t+u\gen x, qt\ln|t|\gen t + (1+q)x\gen x+u\gen u\rangle.
\end{align*}

\medskip\noindent We now examine each remaining realization of $\mathsf{A}_{4.8}$ where $q=1$ is allowed.

\medskip\noindent\underline{$\langle e_1, e_2, e_3, e_4\rangle=\langle \gen u, x\gen u, -\gen x, t\gen t + x\gen x+2u\gen u\rangle$:} From $[e_1, e_5]=0$ we have $e_5=a(t)\gen t + b(t,x)\gen x + c(t,x)\gen u$ and then $[e_2, e_5]=-e_3$ gives $b(t,x)\gen u=\gen x$ which is impossible. Hence we have no realization in this case.

\medskip\noindent\underline{$\langle e_1, e_2, e_3, e_4\rangle=\langle \gen u, x\gen u, -\gen x, x\gen x+2u\gen u\rangle$:} From $[e_1, e_5]=0$ we have $e_5=a(t)\gen t + b(t,x)\gen x + c(t,x)\gen u$ and then $[e_2, e_5]=-e_3$ gives $b(t,x)\gen u=\gen x$ which is impossible. Hence we have no realization in this case.

\medskip\noindent\underline{$\langle e_1, e_2, e_3, e_4\rangle=\langle \gen x, \gen u, u\gen x, t\gen t + 2x\gen x+u\gen u\rangle$:} From $[e_1, e_5]=0$ and $[e_2, e_5]=-e_3$ we have $\displaystyle e_5=a(t)\gen t + [b(t)-\frac{u^2}{2}]\gen x + c(t)\gen u$ and then $[e_3, e_5]=e_2$ gives $-c(t)\gen x=\gen u$ which is impossible. Hence we have no realization in this case.

\medskip\noindent\underline{$\langle e_1, e_2, e_3, e_4\rangle=\langle \gen x, \gen u, u\gen x, 2x\gen x+u\gen u\rangle$:} From $[e_1, e_5]=0$ and $[e_2, e_5]=-e_3$ we have $\displaystyle e_5=a(t)\gen t + [b(t)-\frac{u^2}{2}]\gen x + c(t)\gen u$ and then $[e_3, e_5]=e_2$ gives $-c(t)\gen x=\gen u$ which is impossible. Hence we have no realization in this case.

\medskip\noindent\underline{$\langle e_1, e_2, e_3, e_4\rangle=\langle \gen x, \gen u, u\gen x+t\gen u, 2x\gen x+u\gen u\rangle$:} This realization of $\mathsf{A}_{4.8}$ gives $F=0$ which is a contradiction. Hence we have no realization in this case.

\bigskip\noindent \underline{$\mathsf{A}_{5.38}$:} In this case we have $\langle e_1, e_2, e_3\rangle=3\mathsf{A}_1$ together with the commutation relations $[e_1, e_4]=e_1,\; [e_2, e_4]=[e_3, e_4]=0$ and $[e_1, e_5]=0,\; [e_2, e_5]=e_2,\;[e_3, e_5]=0,\; [e_4, e_5]=e_3$. We see that $\langle e_1, e_2, e_3, e_4\rangle=\mathsf{A}_{2.1}\oplus\mathsf{A}_{2.2}$ with $\langle e_1, e_4\rangle=\mathsf{A}_{2.2}$. Thus we have the following cases for $\langle e_1, e_2, e_3, e_4\rangle$:

\begin{align*}
&\langle \gen u, x\gen u, \gen t, x\gen x+u\gen u\rangle\\
&\langle \gen u, \gen t, x\gen u, x\gen x+u\gen u\rangle\\
&\langle \gen t, x\gen u, \gen u, t\gen t + c(x)\gen u\rangle,\; c''(x)\neq 0\\
&\langle \gen u, \gen x, \gen t, u\gen u\rangle\\
&\langle \gen u, \gen t, \gen x, u\gen u\rangle.
\end{align*}

\medskip\noindent Because we require $[e_4, e_5]=e_3$ we see that we have no realizations in the cases $\langle \gen u, x\gen u, \gen t, x\gen x+u\gen u\rangle$ and $\langle \gen u, \gen x, \gen t, u\gen u\rangle$.

\medskip\noindent\underline{$\langle e_1, e_2, e_3, e_4\rangle=\langle \gen u, \gen t, x\gen u, x\gen x+u\gen u\rangle$:} From the commutation relations we obtain $e_5=(t+l)\gen t + [x\ln|x|+\beta x]\gen u$. Since $\gen t$ and $x\gen u$ belong to the algebra, we may take $e_5=t\gen t + x\ln|x|\gen u$ in canonical form. Thus we have the realization
\[
\mathsf{A}_{5.38}=\langle \gen u, \gen t, x\gen u, x\gen x+u\gen u, t\gen t + x\ln|x|\gen u\rangle
\]

\medskip\noindent\underline{$\langle e_1, e_2, e_3, e_4\rangle=\langle \gen t, x\gen u, \gen u, t\gen t + c(x)\gen u\rangle$ with $c''(x)\neq 0$:} From the commutation relations we find that we must have $xc'(x)=1$ so that $e_4=t\gen t + [\ln|x|+\kappa]\gen u$ and then $e_5=-x\gen x + \gamma(x)\gen u$. We have the residual group $\mathscr{E}(e_1, e_2, e_3, e_4): t'=t,\; x'=x,\; u'=u+U(x)$. Under such a transformation $e_5$ is mapped to $e'_5=-x'\gen {x'} + [\gamma(x)-xU'(x)]\gen {u'}$. We choose $U(x)$ such that $xU'(x)=\gamma(x)$ so that $e_5=-x\gen x$ in canonical form. We note that $\gen u$ is in the centre of the algebra so we may take $e_4=t\gen t + \ln|x|\gen u$. Thus we have the realization
\[
\mathsf{A}_{5.38}=\langle \gen t, x\gen u, \gen u, t\gen t + \ln|x|\gen u, -x\gen x\rangle.
\]

\medskip\noindent\underline{$\langle e_1, e_2, e_3, e_4\rangle=\langle \gen u, \gen t, \gen x, u\gen u\rangle$:} The commutation relations give $e_5=(t+l)\gen t +\beta\gen x + \gamma\gen u$ from which it now follows that the algebra contains the two operators $\gen t$ and $t\gen t$. These give $F=0$ in the equation for $F$ and so we have no realization in this case.

\subsection{Realizations of $\mathsf{A}_{5.38}$:}

\begin{align*}
\mathsf{A}_{5.38}&=\langle \gen u, \gen t, x\gen u, x\gen x+u\gen u, t\gen t + x\ln|x|\gen u\rangle\\
\mathsf{A}_{5.38}&=\langle \gen u, \gen t, x\gen u, t\gen t + \ln|x|\gen u, -x\gen x\gen u\rangle.
\end{align*}

\bigskip\noindent \underline{$\mathsf{A}_{5.39}$:} In this case we have $\langle e_1, e_2, e_3\rangle=3\mathsf{A}_1$ together with the commutation relations $[e_1, e_4]=e_1,\; [e_2, e_4]=e_2,\; [e_3, e_4]=0$ and $[e_1, e_5]=-e_2,\; [e_2, e_5]=e_1,\;[e_3, e_5]=0,\; [e_4, e_5]=e_3$.

\medskip\noindent We first note that the following cases cannot be implemented:

\begin{align*}
&\langle e_1, e_2, e_3\rangle=\langle \gen t, \gen x, \gen u\rangle\\
&\langle e_1, e_2, e_3\rangle=\langle \gen x, \gen t, \gen u\rangle\\
&\langle e_1, e_2, e_3\rangle=\langle \gen t, \gen u, x\gen u\rangle\\
&\langle e_1, e_2, e_3\rangle=\langle \gen u, \gen t, x\gen u\rangle.
\end{align*}
since we require $[e_1, e_5]=-e_2,\; [e_2, e_5]=e_1$, and it is easy to calculate that these are impossible when $e_5=a(t)\gen t + b(t,x,u)\gen x + c(t,x,u)\gen u$.

\medskip\noindent\underline{$\langle e_1, e_2, e_3\rangle=\langle \gen x, \gen u, \gen t\rangle$:} From $[e_3, e_4]=0$ we find that $e_4=a\gen t + b(x,u)\gen x + c(x,u)\gen u$ and that $e_5=\alpha \gen t + \beta(x,u)\gen x + \gamma(x,u)\gen u$ from $[e_3, e_5]=0$. Then it is impossible to implement $[e_4, e_5]=e_3=\gen t$. Thus we have no realization in this case.

\medskip\noindent\underline{$\langle e_1, e_2, e_3\rangle=\langle \gen u, x\gen u, \gen t\rangle$:} From $[e_3, e_4]=0$ we find that $e_4=a\gen t + b(x,u)\gen x + c(x,u)\gen u$ and that $e_5=\alpha \gen t + \beta(x,u)\gen x + \gamma(x,u)\gen u$ from $[e_3, e_5]=0$. Then it is impossible to implement $[e_4, e_5]=e_3=\gen t$. Thus we have no realization in this case.

\section{Realizations of five-dimensional solvable Lie algebras.}

\subsection{Decomposable algebras:}

\bigskip\noindent{\bf Non-linearizing realizations of $2\mathsf{A}_{2.2}\oplus\mathsf{A}_1$:}

\[
2\mathsf{A}_{2.2}\oplus\mathsf{A}_1=\langle \gen x, x\gen x, \gen u, u\gen u, \gen t\rangle.
\]

\bigskip\noindent{\bf Non-linearizing realizations of $\mathsf{A}_{3.k}\oplus\mathsf{A}_{2.2}$:}

\medskip\noindent $\underline{\mathsf{A}_{3.4}\oplus\mathsf{A}_{2.2}}$:

\begin{align*}
\mathsf{A}_{3.4}\oplus\mathsf{A}_{2.2}&=\langle \gen x, \gen t, t\gen t + [t+x]\gen x, \gen u, u\gen u\rangle\\
\mathsf{A}_{3.4}\oplus\mathsf{A}_{2.2}&=\langle \gen x, \gen t, t\gen t + [t+x]\gen x + u\gen u, u\gen x, -u\gen u\rangle.
\end{align*}

\medskip\noindent $\underline{\mathsf{A}_{3.5}\oplus\mathsf{A}_{2.2}}$:
\begin{align*}
\mathsf{A}_{3.5}\oplus\mathsf{A}_{2.2}&=\langle \gen t, \gen x, t\gen t + x\gen x, \gen u, u\gen u\rangle\\
\mathsf{A}_{3.5}\oplus\mathsf{A}_{2.2}&=\langle \gen t, \gen x, t\gen t + x\gen x + u\gen u, u\gen x, -u\gen u\rangle.
\end{align*}

\medskip\noindent $\underline{\mathsf{A}_{3.6}\oplus\mathsf{A}_{2.2}}$:

\begin{align*}
\mathsf{A}_{3.6}\oplus\mathsf{A}_{2.2}&=\langle \gen t, \gen x, t\gen t - x\gen x, \gen u, u\gen u\rangle\\
\mathsf{A}_{3.6}\oplus\mathsf{A}_{2.2}&=\langle \gen t, \gen x, t\gen t - x\gen x-u\gen u, u\gen x, -u\gen u\rangle.
\end{align*}

\medskip\noindent $\underline{\mathsf{A}_{3.7}\oplus\mathsf{A}_{2.2}}$:

\begin{align*}
\mathsf{A}_{3.7}\oplus\mathsf{A}_{2.2}&=\langle \gen u, x\gen u, (1-q)x\gen x + u\gen u, \gen t, t\gen t + \kappa x^{1/(1-q)}\gen u\rangle\\
\mathsf{A}_{3.7}\oplus\mathsf{A}_{2.2}&=\langle \gen t, \gen x, t\gen t + qx\gen x, \gen u, u\gen u\rangle\\
\mathsf{A}_{3.7}\oplus\mathsf{A}_{2.2}&=\langle \gen x, \gen t, qt\gen t + x\gen x, \gen u, u\gen u\rangle\\
\mathsf{A}_{3.7}\oplus\mathsf{A}_{2.2}&=\langle \gen t, \gen x, t\gen t + qx\gen x + u\gen u, u^q\gen x, -\frac{1}{q}u\gen u\rangle\\
\mathsf{A}_{3.7}\oplus\mathsf{A}_{2.2}&=\langle \gen x, \gen t, qt\gen t + x\gen x + u\gen u, u\gen x, -u\gen u\rangle.
\end{align*}

\medskip\noindent $\underline{\mathsf{A}_{3.8}\oplus\mathsf{A}_{2.2}}$:

\begin{align*}
\mathsf{A}_{3.8}\oplus\mathsf{A}_{2.2}&=\langle \gen u, x\gen u, -(1+x^2)\gen x - xu\gen u, \gen t, t\gen t + \kappa\sqrt{1+x^2}\gen u\rangle.
\end{align*}

\bigskip\noindent {\bf Non-linearizing realizations of $\mathsf{A}_{4.k}\oplus\mathsf{A}_1$:}

\begin{align*}
\mathsf{A}_{4.7}\oplus\mathsf{A}_1&= \langle \gen u, x\gen u, -\gen x, x\gen x + \left[2u-\frac{x^2}{2}\right]\gen u, \gen t \rangle\\
\mathsf{A}_{4.7}\oplus\mathsf{A}_1&= \langle \gen u, x\gen u, -\gen x, t\gen t + x\gen x + \left[2u-\frac{x^2}{2}\right]\gen u, t\gen t \rangle.
\end{align*}

\begin{align*}
\mathsf{A}_{4.8}\oplus\mathsf{A}_1&=\langle \gen u, x\gen u, -\gen x, t\gen t + qx\gen x + (1+q)u\gen u, t\gen t\rangle\\
\mathsf{A}_{4.8}\oplus\mathsf{A}_1&=\langle \gen u, x\gen u, -\gen x, qx\gen x + (1+q)u\gen u, \gen t\rangle\\
\mathsf{A}_{4.8}\oplus\mathsf{A}_1&=\langle \gen u, x\gen u, -(\gen t + \gen x), qt\gen t + qx\gen x + (1+q)u\gen u, t\gen t\rangle, \;\; q\neq 0\\
\mathsf{A}_{4.8}\oplus\mathsf{A}_1&=\langle \gen x, \gen u, u\gen x, t\gen t + (1+q)x\gen x + u\gen u, t\gen t\rangle\\
\mathsf{A}_{4.8}\oplus\mathsf{A}_1&=\langle \gen x, \gen u, u\gen x, (1+q)x\gen x + u\gen u, \gen t\rangle\\
\mathsf{A}_{4.8}\oplus\mathsf{A}_1&=\langle \gen x, \gen u, t\gen t + u\gen x, x\gen x + u\gen u, t\gen t\rangle,\;\; q=0.
\end{align*}

\[
\mathsf{A}_{4.10}\oplus\mathsf{A}_1=\langle \gen x, \gen u, x\gen x + u\gen u, u\gen x-x\gen u, \gen t\rangle.
\]

\subsection{Indecomposable solvable Lie algebras.}

\bigskip\noindent{\bf Inequivalent non-linearizing realizations of $\mathsf{A}_{5.19}$:}

\begin{align*}
\mathsf{A}_{5.19}&=\langle \gen u, x\gen u, -\gen x, \gen t, qt\gen t + px\gen x + (p+1)u\gen u\rangle,\\
\mathsf{A}_{5.19}&=\langle \gen x, \gen t, (t+u)\gen x, \gen t-\gen u, t\gen t + (1+p)x\gen x + u\gen u\rangle,\;\; q=1,\\
\mathsf{A}_{5.19}&=\langle \gen x, \gen u, u\gen x, \gen t, qt\gen t + (1+p)x\gen x + u\gen u\rangle\\
\mathsf{A}_{5.19}&=\langle \gen x, \gen t, t\gen x-\gen u, \gen t-u\gen x + \kappa \gen u, t\gen t + (1+p)x\gen x + pu\gen u\rangle,\;\;q=1\\
\mathsf{A}_{5.19}&=\langle \gen x, \gen u, t\gen t+u\gen x, t\gen t, qt\ln|t|\gen t + (1+q)x\gen x + u\gen u\rangle,\;\; p=q\neq 0.
\end{align*}
In the above $p, \kappa\in \mathbb{R}$ and $q\neq 0$ unless otherwise indicated.

\bigskip\noindent{\bf Inequivalent non-linearizing realizations of $\mathsf{A}_{5.20}$:}

\begin{align*}
\mathsf{A}_{5.20}&=\langle \gen u, x\gen u, -\gen x, \gen t, (1+p)t\gen {t} + px\gen {x} + [(1+p)u+t]\gen {u}\rangle,\quad p\in \mathbb{R}\\
\mathsf{A}_{5.20}&=\langle \gen x, \gen t, (t+u)\gen x, \gen t-\gen u, t\gen t + (x-u)\gen x + u\gen u\rangle,\quad p=0\\
\mathsf{A}_{5.20}&=\langle\gen x, \gen u, u\gen x, \gen t, (1+p)t\gen t + [(1+p)x+t]\gen x + u\gen u\rangle, \quad p\in \mathbb{R}\\
\mathsf{A}_{5.20}&=\langle \gen x, \gen t, t\gen x-\gen u, \gen t-u\gen x, t\gen t + x\gen x + \gen u\rangle,\quad p=0.
\end{align*}

\bigskip\noindent {\bf Non-linearizing realizations of $\mathsf{A}_{5.22}$:}

\[
\mathsf{A}_{5.22}=\langle \gen x, \gen u, u\gen x, \gen t, t\gen t + \frac{u^2}{2}\gen x\rangle.
\]

\bigskip\noindent{\bf Non-linearizing realizations of $\mathsf{A}_{5.23}$:}

\begin{align*}
\mathsf{A}_{5.23}&=\langle \gen x, \gen t, (t+u)\gen x, \gen t-\gen u, t\gen t + [2x+\frac{(t+u)^2}{2}]\gen x + u\gen u\rangle\\
\mathsf{A}_{5.23}&=\langle \gen x, \gen u, u\gen x, \gen t, pt\gen t+[2x+\frac{u^2}{2}]\gen x+ u\gen u\rangle.
\end{align*}

\bigskip\noindent{\bf Non-linearizing realizations of $\mathsf{A}_{5.24}$:}

\[
\mathsf{A}_{5.24}=\langle \gen x, \gen u, u\gen x, \gen t, 2t\gen t+[2x+\frac{u^2}{2}+\epsilon t]\gen x+ u\gen u\rangle,\;\; \epsilon=\pm 1.
\]

\bigskip\noindent {\bf Non-linearizing realizations of $\mathsf{A}_{5.30}$:}

\begin{align*}
\mathsf{A}_{5.30}&=\langle \gen x, \gen u , \gen t, u\gen x + t\gen u, pt\gen t+(p+2)x\gen x + (p+1)u\gen u\rangle\\
\mathsf{A}_{5.30}&=\langle \gen u, x\gen u , \gen t, -\gen x + tx\gen u, pt\gen t + x\gen x + (p+2)u\gen u\rangle,\; p\neq 0\\
\mathsf{A}_{5.30}&=\langle \gen u, x\gen u , \gen t, -\gen x + tx\gen u, x\gen x + 2u\gen u + \alpha(2\gen t -x^2\gen u)\rangle,\; \alpha\neq 0.
\end{align*}

\bigskip\noindent {\bf Non-linearizing realizations of $\mathsf{A}_{5.32}$:}

\[
\mathsf{A}_{5.32}=\langle \gen x, \gen u, \gen t, u\gen x + t\gen u, t\gen t + [x+pt]\gen x + u\gen u\rangle
\]

\bigskip\noindent{\bf Non-linearizing realizations of $\mathsf{A}_{5.33}$:}

\begin{align*}
\mathsf{A}_{5.33}&=\langle \gen t, \gen x, \gen u, t\gen t + pu\gen u, x\gen x + qu\gen u\rangle,\; p\neq 0\\
\mathsf{A}_{5.33}&=\langle \gen x, \gen t, \gen u, x\gen x + pu\gen u, t\gen t + qu\gen u\rangle,\; q\neq 0\\
\mathsf{A}_{5.33}&=\langle \gen x, \gen u, \gen t, pt\gen t + x\gen x, qt\gen t + u\gen u\rangle\\
\mathsf{A}_{5.33}&=\langle \gen t, \gen u, x\gen u, t\gen t - px\gen x, (1-q)x\gen x + u\gen u\rangle,\; q\neq 1,\; p\neq 0\\
\mathsf{A}_{5.33}&=\langle \gen t, \gen u, x\gen u, t\gen t + \kappa x^{1/(1-q)}\gen u, (1-q)x\gen x + u\gen u\rangle,\; \kappa\neq 0,\; q\neq 1\\
\mathsf{A}_{5.33}&=\langle \gen u, \gen t, x\gen u, (1-p)x\gen x + u\gen u, t\gen t - qx\gen x\rangle,\; p\neq 1,\; q\neq 0\\
\mathsf{A}_{5.33}&=\langle \gen u, \gen t, x\gen u, (1-p)x\gen x + u\gen u, t\gen t + \kappa x^{1/(1-p)}\gen u\rangle,\; \kappa\neq 0,\; p\neq 1\\
\mathsf{A}_{5.33}&=\langle \gen u, x\gen u, \gen t, pt\gen t + x\gen x + u\gen u, qt\gen t - x\gen x\rangle.
\end{align*}

\bigskip\noindent{\bf Non-linearizing realizations of $\mathsf{A}_{5.34}$:}

\begin{align*}
\mathsf{A}_{5.34}&=\langle \gen t, \gen x, \gen u, pt\gen t + x\gen x + u\gen u, t\gen t + u\gen x\rangle\\
\mathsf{A}_{5.34}&=\langle \gen x, \gen u, \gen t, t\gen t + px\gen x + u\gen u, x\gen x + t\gen u\rangle\\
\mathsf{A}_{5.34}&=\langle \gen t, \gen u, x\gen u, pt\gen t + u\gen u, t\gen t - \gen x\rangle\\
\mathsf{A}_{5.34}&=\langle \gen u, x\gen u, \gen t, t\gen t + (p-1)x\gen x + pu\gen u, x\gen x + [u+tx]\gen u\rangle.
\end{align*}

\bigskip\noindent{\bf Non-linearizing realizations of $\mathsf{A}_{5.35}$:}

\begin{align*}
\mathsf{A}_{5.35}&=\langle \gen t, \gen x, \gen u, pt\gen t + x\gen x + u\gen u, qt\gen t + u\gen x-x\gen u\rangle\\
\mathsf{A}_{5.35}&=\langle \gen t, \gen u, x\gen u, pt\gen t + u\gen u, qt\gen t - (1+x^2)\gen x-xu\gen u\rangle.
\end{align*}

\bigskip\noindent{\bf Non-linearizing realizations of $\mathsf{A}_{5.36}$:}

\begin{align*}
\mathsf{A}_{5.36}&=\langle\gen u, x\gen u, -\gen x, t\gen t + u\gen u, pt\gen t + x\gen x\rangle,\; p\in \mathbb{R}\\
\mathsf{A}_{5.36}&=\langle\gen u, x\gen u, -\gen x, t\gen t + u\gen u, x\gen x + t\gen u\rangle\\
\mathsf{A}_{5.36}&=\langle\gen x, \gen t, (t+u)\gen x, t\gen t + x\gen x + u\gen u, -t\gen t - u\gen u\rangle\\
\mathsf{A}_{5.36}&=\langle\gen x, \gen u, u\gen x, t\gen t + x\gen x +u\gen u, \alpha t\gen t - u\gen u\rangle,\; \alpha\in \mathbb{R}\\
\mathsf{A}_{5.36}&=\langle\gen x, \gen u, u\gen x, t\gen t + x\gen x +u\gen u, t\gen x - u\gen u\rangle.
\\\mathsf{A}_{5.36}&=\langle \gen x, \gen u, u\gen x+t\gen u, t\gen t + x\gen x +u\gen u, -2t\gen t - u\gen u\rangle\\
\mathsf{A}_{5.36}&=\langle \gen x, \gen t, t\gen x-\gen u, t\gen t + x\gen x, -t\gen t + u\gen u\rangle\\
\mathsf{A}_{5.36}&=\langle \gen x, \gen u, \gen t+u\gen x, x\gen x + u\gen u, t\gen t -u\gen u\rangle
\end{align*}

\bigskip\noindent{\bf Non-linearizing realizations of $\mathsf{A}_{5.38}$:}

\begin{align*}
\mathsf{A}_{5.38}&=\langle \gen u, \gen t, x\gen u, x\gen x+u\gen u, t\gen t + x\ln|x|\gen u\rangle\\
\mathsf{A}_{5.38}&=\langle \gen u, \gen t, x\gen u, t\gen t + \ln|x|\gen u, -x\gen x\gen u\rangle.
\end{align*}

\section{Admissible  invariant equations}
In this Section we list all admissible equations invariant under the obtained realizations of solvable algebras of dimension $\leq 5$. We mention that those equations invariant under algebras of dimension 5 depend on two arbitrary constants.
\subsection{Two-dimensional non-abelian solvable algebras}
$\mathsf{A}_{2.2}$:

\begin{enumerate}
\item
$$F=x^2f(u,\tau_1,\tau_2), \quad G=x^{-1}g(u,\tau_1,\tau_2), \quad \tau_1=u_2u_1^{-2}, \quad \tau_2=xu_1$$
\item
$$F=t^2f(u_1,\tau_1,\tau_2),  \quad G=t^{-1}f(u_1,\tau_1,\tau_2),  \quad \tau_1=u_2u_1^{-2}, \quad \tau_2=tu_1$$
\item
$$F=u^3f(t,u_1,\tau),  \quad G=ug(t,u_1,\tau),  \quad \tau=uu_2$$
\item
$$F=u_1^{-3}f(t,u,\tau), \quad G=g(t,u,\tau), \quad \tau=u_2u_1^{-2}$$

\end{enumerate}
\subsection{Three-dimensional solvable algebras}

$\mathsf{A}_{3.2}$:

\begin{enumerate}
\item
$$F=u_1^{-3}f(u,\tau),\quad G=g(u,\tau),\quad \tau=u_2u_1^{-2}$$

\item
$$F=u_1^{-3}f(t,\tau),\quad G=g(t,\tau),\quad \tau=u_2u_1^{-2}$$

\item
$$F=x^2 e^{u/\tau}\tau^{-3}f(u,\omega),\quad \tau=x u_1,\quad \sigma=x^2 u_2, \quad \omega=\tau^{-3}(\tau+\sigma)$$ $$G=xe^{u/\tau}[(-3\tau^2\omega^2+3\omega-2\tau^{-2})f(u,\omega)+\tau g(u,\omega)]$$

\item
$$F=t^2 \tau^{-3} f(u,\omega),\quad G=t^{-1}[-u-3(\tau \omega)^2 f(u,\omega)+\tau g(u,\omega)]$$
$$\tau=t u_1,\quad \sigma=t^2 u_2, \quad \omega=\sigma\tau^{-3}$$

\item
$$F=t^{-1}u_1^{-3}f(u,\tau),\quad G=t^{-1}f(u,\tau),\quad \tau=u_1^{-2}u_2$$

\item

$$F=u^3 u_1^{-3}f(t,\omega),\quad G=u u_1[-3\omega^2u_1f(t,\omega)+g(t,\omega)]$$
$$\omega=uu_1^{-3}u_2$$

\item
$$F=u_1^{-3}f(t,\omega),\quad G=\frac{1}{t}\ln|u_1|+g(t,\omega),\quad \omega=uu_1^{-2}u_2$$

\item
$$F=u_1^{-2}\tau^{1/2}f(u,\sigma),\quad G=u_1\tau^{1/2}g(u,\sigma)$$
$$\sigma=t^4u_2^2\tau^{-3},\quad \tau=tu_1$$

\item
$$F=u^3f(u_1,\tau),\quad G=ug(u_1,\tau),\quad \tau=uu_2$$
\end{enumerate}

$\mathsf{A}_{3.3}$:

\begin{enumerate}
\item
$$F=f(t,u_2),\quad G=g(t,u_2)$$

\item

$$F=f(\tau,u_2),\quad G=g(\tau,u_2),\quad \tau=x-t$$

\item

$$F=u_1^{-3}f(u,\tau),\quad \tau=u_2u_1^{-3}$$
$$G=-1-3u_1^2\tau^2 f(u,\tau)+u_1g(u,\tau)$$

\item

$$F=u_1^{-3}f(t,\tau),\quad \tau=u_2u_1^{-3}$$
$$G=u_1[-3u_1 \tau^2 f(t,\tau)+g(t,\tau)]$$

\item
$$F=u_1^{-3}f(t,\tau),\quad \tau=u_2u_1^{-3}$$
$$G=u_1[\frac{1}{2}u_1^{-2}-3u_1 \tau^2 f(t,\tau)]+u_1 g(t,\tau)$$

\item

$$F=f(u_1,u_2),  \quad G=uu_1+g(u_1,u_2)$$


\item

$$F=f(\tau,\omega),\quad G=u_1[-\ln|u_1|+g(\tau,\omega)]$$
$$\tau=u^{-1}u_2,\quad \omega=u^{-1}u_1$$
\end{enumerate}

$\mathsf{A}_{3.4}$:

\begin{enumerate}

\item
$$F=u^3u_1^{-3}f(t,\omega),\quad G=u_1[-3u_1u^3\tau^2 f(t,\omega)+u^{-1} g(t,\omega)]$$
$$\omega=u^3u_1^{-3}u_2$$

\item

$$F=u^4 u_1^{-3}f(\omega_1,\omega_2),\quad \omega_1=tu,\quad \omega_2=u^3u_1^{-3}u_2$$
$$G=ug(\omega_1,\omega_2)$$

\item

$$F=u_1^{-2}f(u,\tau),\quad G=u_1(\ln|u_1|+g(u,\tau)),\quad \tau=u_2u_1^{-2}$$

\item

$$F=u^2 f(u_1,\tau),\quad G=-u_1\ln|u_1|+g(u_1,\tau),\quad \tau=uu_2$$

\item

$$F=u_1^{-3}\exp{(3/u_1)}f(t,\tau),\quad \tau=\exp{(1/u_1)}u_1^{-3}u_2$$
$$G=u_1\exp{(1/u_1)}[-3u_1 \tau^3 f(t,\tau)+g(t,\tau)]$$

\item

$$F=u_1^{-3}\exp{(2/u_1)}f(\tau,\omega),\quad G=u_1[-3u_1 \tau^2 f(\tau,\omega)+g(\tau,\omega)]$$
$$\tau=\exp{(1/u_1)}u_1^{-3}u_2,\quad \omega=t\exp{(-1/u_1)}$$

\end{enumerate}

$\mathsf{A}_{3.5}$:

\begin{enumerate}

\item
$$F=f(t,x),\quad G=u_2g(t,x)$$

\item
$$F=t^{-1}f(x,\tau),\quad G=g(x,\tau),\quad \tau=u_2t^{-1}$$

\item

$$F=u_1^{-2}f(u,\tau),\quad G=u_1g(u,\tau),\quad \tau=u_2u_1^{-2}$$

\item

$$F=u^2f(u_1,\tau),\quad G=g(u_1,\tau),\quad \tau=uu_2$$

\item

$$F=t^2f(u_1,\tau),\quad G=f(u_1,\tau),\quad \tau=tu_2$$

\item

$$F=u_2^{-3}f(t,u_1),\quad G=u_2^{-1}g(t,u_1)$$

\end{enumerate}

$\mathsf{A}_{3.6}$:

\begin{enumerate}

\item
$$F=x^3 f(t,\tau),\quad \tau=x^{-2}u_2^{-3},$$
$$G=x^{1/2}g(t,\tau)$$
\item
$$F=x^{5/2}f(\tau,\omega),\quad \tau=x^{-2}u_2^{-3},\quad \omega=x t^{-2}$$
$$G=g(\tau,\omega)$$
\item
$$F=u_1^{-4}f(u,\tau),\quad \tau=u_1^{-2}u_2$$
$$G=u_1^{-1}g(u,\tau)$$

\item

$$F=u^{4}f(u_1,\tau),\quad \tau=uu_2$$
$$G=u^{2}g(u_1,\tau)$$

\item

$$F=t^2f(\tau,\omega),\quad G=t^{-2}g(\tau,\omega)$$
$$\tau=t^{2}u_1,\quad \omega=t^{3}u_2$$

\end{enumerate}

$\mathsf{A}_{3.7}$:

\begin{enumerate}
\item
$$F=x^3 f(t,\tau),\quad \tau=x^{q-1}u_2^{2q-1},$$
$$G=x^{1/(1-q)}g(t,\tau)$$
\item
$$F=x^{(3q-2)/(q-1)}f(\tau,\omega),\quad \tau=x^{q-1}u_2^{2q-1},\quad \omega=x t^{q-1}$$
$$G=g(\tau,\omega)$$
\item
$$F=u_1^{(1/q-3)}f(u,\tau),\quad \tau=u_1^{-2}u_2$$
$$G=u_1^{1/q}g(u,\tau)$$

\item
$$F=u^{3q-1}f(\tau,\omega),\quad \tau=u^{q-1}u_1,\quad \omega=u^{2q-1}u_2$$
$$G=g(\tau,\omega)$$

\item

$$F=u^{3-q}f(u_1,\tau),\quad \tau=uu_2$$
$$G=u^{1-q}g(u_1,\tau)$$

\item

$$F=u_1^{3/(q-1)}f(t,\tau),\quad \tau=u_1^{q-2}u_2^{1-q}$$
$$G=u_1^{q/(q-1)}g(t,\tau)$$

\item

$$F=t^2f(\tau,\omega),\quad G=t^{q-1}g(\tau,\omega)$$
$$\tau=t^{1-q}u_1,\quad \omega=t^{2-q}u_2$$

\end{enumerate}

$\mathsf{A}_{3.8}$:

\begin{enumerate}
\item
$$F=(1+x^2)^3f(t,\omega),\quad G=(1+x^2)^{1/2}[3x\omega f(t,\omega)+g(t,\omega)]$$

\item
$$F=(1+x^2)^3\exp(-\arctan x)f(\tau,\omega),\quad \tau=\arctan x-\ln|t|,\quad \omega=(1+x^2)^{3/2}u_2$$
$$G=(1+x^2)^{1/2}\exp(-\arctan x)[3x \omega f(\tau,\omega)+g(\tau,\omega)]$$

\item
$$F=u_2^{-1}f(t,\omega),\quad \omega=(1+u_1^2)^{3/2}u_2^{-1}$$
$$G=-3(1+u_1^2)u_1u_2f(t,\omega)+(1+u_1^2)^{1/2}g(t,\omega)$$

\item
$$F=(1+u_1^2)^{-3/2}\exp(\arctan u_1)f(\tau,\omega),\quad G=$$
$$G=\exp(\arctan u_1)[-3u_1u_2^2(1+u_1^2)^{-5/2}f(t,\omega)+(1+u_1^2)^{1/2}g(\tau,\omega)]$$
$$\tau=\arctan u_1+\ln|t|,\quad \omega=(1+u_1^2)^{3/2}u_2^{-1}$$
\end{enumerate}

$\mathsf{A}_{3.9}$:
\begin{enumerate}
\item
$$F=(1+x^2)^3 f(t,\tau),\quad \tau=\exp\{q \arctan x\} (1+x^2)^{3/2}u_2 $$
$$G=\exp\{-q \arctan x\}(1+x^2)^{1/2}[3x\tau f(t,\tau)+g(t,\tau) ] $$
\item
$$F=\exp\{-\arctan x\}(1+x^2)^3f(\tau,\omega),$$
$$\tau=\arctan x-\ln|t|,\quad \omega=\exp\{q \arctan x\} (1+x^2)^{3/2}u_2$$
$$G=(1+x^2)^{1/2}\exp\{(q+1) \arctan x\}[3x \omega f(\tau,\omega)+g(\tau,\omega)]$$
\item
$$F=\exp\{-3q\arctan u_1\}(1+u_1^2)^{-3/2}f(t,\tau),$$
$$G=\exp\{-q \arctan u_1\} (1+u_1^2)^{1/2}[-3\tau^2 u_1 f(t,\tau)+g(t,\tau)]$$
$$ \tau=\exp\{-q \arctan u_1\} (1+u_1^2)^{-3/2}u_2$$

\item
$$F=\exp\{-(3q-1)\arctan u_1\}(1+u_1^2)^{-3/2}f(\tau,\omega),$$
$$\tau=\exp\{-q \arctan u_1\} (1+x^2)^{-3/2}u_2,\quad \omega=\exp\{q \arctan u_1\}+\ln|t|$$
$$G=\exp\{-(q-1)\arctan u_1\}(1+u_1^2)^{1/2}[-3\tau^2 u_1f(\tau,\omega)+g(\tau,\omega)]$$
\end{enumerate}

\subsection{Four-dimensional solvable algebras}
\subsubsection{Decomposable four-dimensional solvable algebras}
$\mathsf{A}_{2.1}\oplus\mathsf{A}_{2.2}$:

\begin{enumerate}

\item
$$F=x^3f(\tau),\quad G=xg(\tau),\quad \tau=xu_2$$

\item
$$F=e^{-u_2/c''(x)}f(x),\quad G=e^{-u_2/c''(x)}[-\frac{c'''}{c''}u_2f(x)+g(x)]$$

\item
$$F=f(\tau),\quad G=u_1g(\tau),\quad \tau=u_2u_1^{-1}$$

\end{enumerate}

$\mathsf{A}_{2.2}\oplus\mathsf{A}_{2.2}$

\begin{enumerate}
\item
$$F=u^2u_1^{-3}f(\sigma),\quad G=g(\sigma),\quad \sigma=uu_1^{-2}u_2$$

\item
$$F=u_1^3u_2^{-3}f(t),\quad G=u_1^2u_2^{-1}g(t)$$

\item

$$F=t^2u_1^{-3}f(\sigma),\quad G=g(\sigma),\quad \sigma=tu_1^{-2}u_2$$

\item
$$F=x^2e^{1/\omega_1}\omega_1^{-3}f(\omega),$$
$$G=ux^{-1}e^{1/\omega_1}[(-3(\omega\omega_1)^2-2\omega_1^{-2}+3\omega) f(\omega)+\omega_1g(\omega)]$$
$$\omega_1=xu_1/u,\quad \omega_2=x^2 u_2/u,\quad \omega=(\omega_1+\omega_2)\omega_1^{-3}$$

\item
$$F=u^3(tu_1^3)^{-1}f(\omega),\quad G=u[-t^{-1}-3t^{-1}u^2u_1^{-4}u_2^2f(\omega)+u^{-1/(1+a)}u_1g(\omega)]$$
$$ \omega=tu^{(a+2)/(a+1)}u_1^{-3}u_2$$

\item
$$F=t^{-1}u^3  u_1^{-3}f(\rho),\quad G=t^{-1}u\omega[-3\omega f(\rho)+g(\rho) ]$$
$$\omega=uu_1^{-2}u_2,\quad \rho=t^{-1}uu_1^{-3}u_2$$

\item
$$F=x^2 \tau^{-q}f(\omega),\quad G=x^{-1} \tau^{1-q}g(\omega)$$
$$\tau=xu_1,\quad \sigma=x^2u_2,\quad \omega=\sigma \tau^{-1}$$

\item
$$F=t^2 f(\rho),\quad G=u_1[-q \ln|\tau|+g(\rho)]$$
$$\tau=tu_1,\quad \rho=tu_1^{-1}u_2$$

\item
$$F=t^2\tau^{-3a/(1+a)}f(\omega),\quad G=t^{-1}\tau^{1/(1+a)}g(\omega)$$
$$\tau=tu_1,\quad \sigma=t^2u_2,\quad \omega=\sigma\tau^{-(1+2a)/(1+a)}$$

\item
$$F=u^3 u_1^{-(1+3q)/q}f(\sigma),\quad G=uu_1^{-1/q}g(\sigma)$$
$$\sigma=uu_1^{-2}u_2$$

\item
$$F=u^3 u_1^{-3}e^{-1/u_1}f(\sigma),\quad G=u u_1e^{-1/u_1}[-3u_1 \sigma^2 f(\sigma)+g(\sigma)]$$
$$\sigma=uu_1^{-3}u_2$$
\end{enumerate}

$\mathsf{A}_{3.3}\oplus\mathsf{A_1}$

\begin{enumerate}
\item
$$F=f(u_2),\quad G=g(u_2)$$

\item
$$F=u_1^{-3}f(\tau),\quad G=-1-3\tau^2 u_1^2 f(\tau)+u_1 g(\tau),\quad \tau=u_2u_1^{-3}$$

\item
$$F= u_1^{-3}f(\tau),\quad G=u_1[-3\tau^2 u_1f(\tau)+g(\tau)],\quad \tau=u_2u_1^{-3}$$

\item
$$F=u_1^{-3}f(\tau),\quad G=uu_1+k-3\tau^2u_1^2f(\tau)+u_1g(\tau),\quad \tau=u_2u_1^{-3}$$

\item
$$F=(tu_1)^{-3}f(\tau),\quad G=t^{-1}uu_1[-3\tau^2 u_1 f(\tau)+g(\tau)],\quad \tau=u_2u_1^{-3}$$

\end{enumerate}

$\mathsf{A}_{3.4}\oplus\mathsf{A}_1$

\begin{enumerate}
\item
$$F=u_1^{-3}f(\sigma),\quad G=u_1[-3u_1\tau^2 f(\sigma)+ e^{-u}g(\sigma)],\quad \tau=u_2u_1^{-3}$$

\item
$$F=Ku_1^{-3}e^u,\quad G=u_1[-3Ku_1\tau^2+L\tau],\quad \tau=u_2u_1^{-3}$$

\item
$$F=u_1^{-2}f(\tau),\quad G=u_1[\ln|u_1|+g(\tau)],\quad \tau=u_2u_1^{-2}$$

\item
$$F=e^{2/u_1}u_1^{-3}f(\tau),\quad G=u_1[-3u_1\tau^2 f(\tau)+g(\tau)],\quad \tau=e^{2/u_1}u_1^{-3}u_2$$

\item
$$F=u_1^3u_2^{-2}f(\sigma),\quad G=u_1[-3u_1f(\sigma)-\ln|u|+g(\sigma)],\quad \sigma=uu_2u_1^{-3}$$

\item
$$F=\tau^{-2}f(\omega),\quad G=u\tau[\ln|\tau|+g(\omega)]$$
$$\tau=u_1u^{-1},\quad \sigma= u_2u^{-1},\quad  \omega=\sigma\tau^{-2}$$

\item
$$F=t^{-1}e^{3/u_1}u_1^{-3}f(\tau),\quad G=t^{-1}e^{1/u_1}u_1[-3\tau^2 u_1f(\tau)+g(\tau)]$$
$$\tau=e^{1/u_1}u_1^{-3}u_2$$

\end{enumerate}

$\mathsf{A}_{3.5}\oplus\mathsf{A}_1.$

\begin{enumerate}
\item
$$F=u_1^2f(\tau),\quad G=u_1g(\tau),\quad \tau=u_2u_1^{-2}$$

\item
$$F=u_2^{-3}f(u_1),\quad G=u_2^{-1}g(u_1)$$

\item
$$F=u_1^{-3}u^2f(\sigma),\quad G=u_1[-3u_1 \tau^2 u^2f(\sigma)+g(\sigma)]$$
$$\tau=u_2u_1^{-3},\quad \sigma=u \tau$$

\end{enumerate}

$\mathsf{A}_{3.6}\oplus\mathsf{A}_1.$

\begin{enumerate}
\item
$$F=x^3f(\tau),\quad G=x^{1/2}g(\tau),\quad \tau=x^{3/2}u_2$$

\item
$$F=u_1^{-4}f(\tau),\quad G=u_1^{-4}g(\tau),\quad \tau=u_2u^{-2}$$

\item
$$F=u_1^{-3/2}f(\tau),\quad G=u_1^{1/2}g(\tau),\quad \tau=u_2u_1^{-3/2}$$

\item
$$F=u_1^{-3}u^4f(\sigma),\quad G=u_1[-3\tau^2 u_1 u^4f(\sigma)+u^2g(\sigma)]$$
$$\tau=u_2u_1^{-3},\quad \sigma=u\tau$$
\end{enumerate}

$\mathsf{A}_{3.7}\oplus\mathsf{A}_1.$

\begin{enumerate}
\item
$$F=x^3f(\tau),\quad G=x^{1/(1-q)}g(\tau),\quad \tau=x^{(2q-1)/(q-1)}u_2$$

\item
$$F=u_1^{(1-3q)/q}f(\tau),\quad G=u_1^{1/q}g(\tau),\quad \tau=u_2u_1^{-2}$$

\item
$$F=x^{3-q}f(\tau),\quad G=x^{-q}g(\tau),\quad \tau=x^2 u_2$$

\item
$$F=u_1^{3/(q-1)}f(\tau),\quad G=u_1^{q/(q-1)}g(\tau),\quad \tau=u_2u_1^{(2-q)/(q-1)}$$

\item
$$F=u^{2}u_1^{-3}f(\omega),\quad G=-[(2q-1)(q-1)+3\tau^2\omega^2]f(\omega)+\tau g(\omega)$$
$$\tau=u_1 u^{q-1},\quad \sigma=u_2 u^{2q-1},\quad \omega=[\sigma+(q-1)\tau^2]\tau^{-3}$$

\item
$$F=u^{3q-1}\tau^{-3}f(\omega),\quad G=\tau[-3\omega^2 \tau f(\omega)+g(\omega)]$$
$$\tau=u^{q-1}u_1,\quad \sigma=u^{2q-1}u_2,\quad \omega=\sigma\tau^{-3}$$

\end{enumerate}

$\mathsf{A}_{3.8}\oplus\mathsf{A}_1.$

\begin{enumerate}
\item
$$F=(1+x^2)^3f(\omega),\quad G=(1+x^2)^{1/2}[3x\omega f(\omega)+g(\omega)]$$
$$\omega=(1+x^2)^{3/2}u_2$$

\item
$$F=u_2^{-1}f(\omega),\quad G=-3u_1u_2(1+u_1^2)^{-1}f(\omega)+(1+u_1^2)^{1/2}g(\omega)$$
$$\omega=(1+x^2)^{3/2}u_2$$

\end{enumerate}

$\mathsf{A}_{3.9}\oplus\mathsf{A}_1.$

\begin{enumerate}
\item
$$F=(1+x^2)^3f(\tau),\quad G=(1+x^2)^{1/2}e^{-q\arctan x}[3x\tau f(\tau)+g(\tau)]$$
$$\tau=(1+x^2)^{3/2}e^{q\arctan x}u_2$$

\item
$$F=(1+u_1^2)^{-3/2}e^{-3q\arctan u_1}f(\tau),\quad G=(1+u_1^2)^{1/2}e^{-q\arctan u_1}[-3\tau^2 u_1 f(\tau)+g(\tau)]$$
$$\tau=(1+u_1^2)^{-3/2}e^{-q\arctan u_1}u_2$$

\end{enumerate}

\section{Non-decomposable four-dimensional solvable algebras.}

$\mathsf{A}_{4.1}$:

\begin{enumerate}
\item
$$F=u_1^{-3}f(\tau),\quad G=\frac{1}{2}u_1^{-1}-3u_1^2\tau^2 f(\tau)+u_1g(\tau),\quad \tau=u_1^{-3}u_2$$

\item
$$F=f(u_2),\quad G=-\frac{x^2}{2}+g(u_2)$$

\end{enumerate}

$\mathsf{A}_{4.2}$:

\begin{enumerate}
\item
$$F=u_1^{-3}e^{(3-q)/u_1}f(\tau),\quad G=u_1e^{(1-q)/u_1}[-3\tau^2 u_1f(\tau)+G(\tau)]$$
$$\tau=u_1^{-3}e^{1/u_1}u_2$$

\item
$$F=u_1^{(3q-1)/(1-q)}f(\tau),\quad G=(1-q)\ln|u_1|+g(\tau),\quad \tau=u_1^{1-2q}u_2^{q-1}$$

\item
$$F=x^{q+3}f(\tau),\quad G=x^{q-1}g(\tau),\quad \tau=x^3u_2$$

\item
$$F=x^{(3q-4)/(q-1)}f(\tau),\quad G=x[\frac{1}{q-1}\ln|x|+g(\tau)],\quad \tau=x^{(q-2)/(q-1)}u_2$$

\end{enumerate}

$\mathsf{A}_{4.3}$:

\begin{enumerate}
\item
$$F=u_1^{-3}e^{-1/u_1}f(\tau),\quad G=u_1e^{-1/u_1}[3u_1\tau^2 f(\tau)+g(\tau)],\quad \tau=u_1^{-3}u_2$$

\item
$$F=f(u_2),\quad G=u_1+g(u_2)$$

\item
$$F=e^xf(u_2),\quad G=e^xg(u_2)$$

\item
$$F=x^3f(\tau),\quad G=x[\ln |x|+g(\tau)],\quad \tau=xu_2$$

\end{enumerate}

$\mathsf{A}_{4.4}$:

\begin{enumerate}
\item
$$F=u_1^{-3}e^{2/u_1}f(\tau),\quad G=u_1[\frac{1}{2}u_1^{-2}-3\tau^2 u_1f(\tau)+g(\tau)]$$
$$\tau=u_1^{-3}e^{1/u_1}u_2$$

\item
$$F=x^4f(\tau),\quad G=-\ln|x|+g(\tau),\quad \tau=x^3u_2$$
\end{enumerate}

$\mathsf{A}_{4.5}$:

\begin{enumerate}
\item
$$F=u_1^{1-3q}f(\tau),\quad G=u_1^{1-p}g(\tau),\quad \tau=u_1^{p-2q}u_2$$

\item
$$F=u_1^{(q-3)/(p-1)}f(\tau),\quad G=u_1^{(p-q)/(p-1)}g(\tau),\quad \tau=u_1^{(2-p)/(p-1)}u_2$$

\item
$$F=u_1^{(p-3)/(q-1)}f(\tau),\quad G=u_1^{(q-p)/(q-1)}g(\tau),\quad \tau=u_1^{(2-q)/(q-1)}u_2$$

\item
$$F=x^{(3(p-q)+1)/(p-q)}f(\tau),\quad G=x^{(1-q)/(p-q)}g(\tau),\quad \tau=x^{(2p-q)/(p-q)}u_2$$

\item
$$F=F=x^{(q+3p-3)/(p-1)}f(\tau),\quad G=x^{(q-1)/(p-1)}g(\tau),\quad \tau=x^{(2p-1)/(p-1)}u_2$$

\item
$$F=F=x^{(p+3q-3)/(q-1)}f(\tau),\quad G=x^{(p-1)/(q-1)}g(\tau),\quad \tau=x^{(2q-1)/(q-1)}u_2$$

\end{enumerate}

$\mathsf{A}_{4.6}$:

\begin{enumerate}
\item
$$F=(1+u_1^2)^{-3/2}e^{(q-3p)\arctan u_1}f(\tau),\quad G= (1+u_1^2)^{1/2}e^{(q-p)\arctan u_1}[3\tau^2 u_1 f(\tau)+g(\tau)]$$
$$\tau=(1+u_1^2)^{-3/2}e^{-p\arctan u_1}u_2$$

\item
$$F=(1+x^2)^3e^{q\arctan x}f(\tau),\quad G=(1+x^2)^{1/2}e^{(q-p)\arctan x}[3x\tau f(\tau)+g(\tau)]$$
$$\tau=(1+x^2)^{3/2}e^{p\arctan x}u_2$$

\end{enumerate}

$\mathsf{A}_{4.7}$:

\begin{enumerate}
\item
$$F=e^{-3u_2}f(t),\quad G=e^{-2u_2}g(t)$$

\item
$$F=t^2f(\tau),\quad G=tg(\tau),\quad \tau=u_2+\ln|t|$$

\item
$$F=e^{-2u_2}f(\tau),\quad G=e^{-u_2}g(\tau),\quad \tau=u_2+\ln|t|$$

\item
$$F=u_1^{-3}e^{-3t}f(\tau),\quad G=\frac{1}{2}u_1^{-1}-3u_1^2\tau^2e^{-3t}f(\tau)+u_1e^{-2t}g(\tau),\quad \tau=u_1^{-3}u_2$$

\end{enumerate}

$\mathsf{A}_{4.8}$:

\begin{enumerate}

\item
$$F=t^{3q-1}f(\tau),\quad G=t^{q}g(\tau),\quad \tau=t^{q-1}u_2$$

\item
$$F=u_2^{3q/(1-q)}f(t),\quad G=u_2^{(q+1)/(1-q)}g(t)$$

\item
$$F=u_2^{-3/2}f(t),\quad G=g(t),\quad \tau=t^{-2}u_2$$

\item
$$F=t^{-4}f(\tau),\quad G=t^{-1}g(\tau),\quad \tau=t^{-2}u_2$$

\item
$$F=(x-t)^2f(\sigma),\quad G=(x-t)^{1/q}g(\sigma),\quad \sigma=(x-t)^{(q-1)/q}u_2$$

\item
$$F=f(\sigma),\quad G=u_2g(\sigma),\quad \sigma=x-t+a\ln|u_2|$$

\item
$$F=u_1^{-3}u^2f(\sigma),\quad G=-1-3u_1^2\tau^2 u^2f(\sigma)+u_1u^qg(\sigma)$$
$$\sigma=u_1^{-3}u^{1-q}u_2$$

\item
$$F=u_1^{-3}u^2f(\sigma),\quad G=u_1[-3\tau^2 u_1 f(\sigma)+g(\sigma)],\quad \sigma=u_1^{-3}u_2$$

\item
$$F=u_1^{-3} t^2 f(\sigma),\quad G=u_1[-3\tau^2 u_1 t^2f(\sigma)+t^{q}g(\sigma)]$$
$$\sigma=t^{1-q}\tau,\quad \tau=u_1^{-3}u_2$$

\item
$$F=(u_1/u_2)^{3/2}f(t),\quad G=-3(u_1u_2)^{1/2}f(t)+u_1 g(t)$$

\item
$$F=u^3u_1^{-3}f(\sigma),\quad G=u_1[-3u_1\tau^2 u^3f(\sigma)-\ln|u|+g(\sigma) ]$$
$$\sigma=u^2u_1^{-3}u_2$$

\item
$$F=t^{(q+2)/(1-q)}u_1^{-3}f(\sigma),\quad G=\frac{1}{2}u_1^{-1}-3u_1^2\tau^2 t^{(q+2)/(1-q)}f(\sigma)+t^{2q/(1-q)}u_1g(\sigma)$$
$$\sigma=t\tau,\quad \tau=u_1^{-3}u_2$$

\item
$$F=u_1^{-(3q+2)}f(\tau),\quad G=uu_1+u_1^{1-q}g(\tau),\quad \tau=u_1^{-(q+2)}u_2$$

\item
$$F=u_1^{-2}f(\tau),\quad G=uu_1+u_1[c\ln|u_1|+g(\tau)],\quad \tau=u_1^{-2}u_2$$

\item
$$F=u_1^{-3}\exp[-1/u_1+\sigma]\sigma^{(3-q)/q}f(\omega),$$
$$G=u_1e^{-1/u_1}[-3\tau^2 u_1 \sigma^{(3-q)/q}e^{\sigma}f(\omega)+e^{\sigma}\sigma^{1/q}g(\omega)]$$
$$\sigma=u_1^{-1}-\ln|t|,\quad \omega=\sigma^{(1-q)/q}u_1^{-3}u_2$$

\item
$$F=u_1^{-3}\exp[-1/u_1]\tau^{\alpha-3}f(\omega),\quad \omega=\tau^{\alpha}e^{-\sigma},\quad \sigma=u_1^{-1}-\ln|t|,\quad \tau=u_1^{-3}u_2$$
$$G=u_1\exp[-1/u_1][-3\tau^2 u_1 \tau^{\alpha-3}f(\omega)+\tau^{\alpha-1}g(\omega)]$$

\end{enumerate}

$\mathsf{A}_{4.9}$:

$$F=u_1^{-3}(1+t^2)^{1/2}e^{-3q\arctan t}f(\sigma)$$
$$G=u_1[-3\tau^2u_1 (1+t^2)^{1/2}e^{-3q\arctan t}f(\sigma)+\frac{1}{2}u_1^{-2}+(1+t^2)^{-1}e^{-2q\arctan t}g(\sigma)]$$
$$\sigma=t+(1+t^2)\tau,\quad \tau=u_2 u_1^{-3}$$

$\mathsf{A}_{4.10}$:

\begin{enumerate}
\item
$$F=(1+x^2)^3f(\sigma),\quad G=(1+x^2)^{1/2}e^{a(\arctan x)}[3x\sigma f(\sigma)+g(\sigma)]$$
$$\sigma=t^{-1}(1+x^2)^{3/2}e^{-a(\arctan x)}u_2$$

\item
$$F=(1+u_1^2)^{3}u_2^{-3}f(t),\quad G=(1+u_1^2)^{2}u_2^{-1}[-3u_1f(t)+g(t)]$$

\item
$$F=(1+u_1^2)^{-3/2}t^{2}f(\sigma),\quad G=(1+u_1^2)^{1/2}[-3u_1\tau^2 f(\sigma)+g(\sigma)],$$
$$\sigma=t (1+u_1^2)^{-3/2}u_2$$

\item
$$F=(1+u_1^2)^{-3/2}e^{a \arctan u_1}\tau^{-2}f(\omega),\quad \omega=te^{a \arctan u_1}(1+u_1^2)^{-3/2}u_2$$
$$G=(1+u_1^2)^{1/2}e^{a \arctan u_1}\tau^{-2}[-3\tau^2 u_1 f(\omega)+g(\omega)]$$

\end{enumerate}

\section{Realizations of five-dimensional solvable Lie algebras.}

\subsection{Decomposable algebras:}

\bigskip\noindent $2\mathsf{A}_{2.2}\oplus\mathsf{A}_1$:

$$u_t=Ku_1^3u_2^{-3}u_3+Lu_1^2 u_2^{-1}$$

$\underline{\mathsf{A}_{3.4}\oplus\mathsf{A}_{2.2}}$:

\begin{enumerate}
\item
$$F=K e^{2/u_1}u_1^3u_2^{-2},\quad G=u_1[\ln|u_2u_1^{-2}|+L]$$

\item
$$F=K u_1^2 (uu_2)^{-1},\quad G=u_1[\ln|u_1u_2^{-1}|-3Kuu_1^{-2}u_2+L]$$
\end{enumerate}

$\underline{\mathsf{A}_{3.5}\oplus\mathsf{A}_{2.2}}$:

\begin{enumerate}
\item
$$F=Ku_1^2u_2^{-2},\quad G=Lu_1$$

\item
$$F=Kuu_2^{-1},\quad G=Lu_1$$

\end{enumerate}

$\underline{\mathsf{A}_{3.6}\oplus\mathsf{A}_{2.2}}$:

\begin{enumerate}

\item
$$u_t=Ku_1^{4}u_2^{-4}u_3+Lu_1^3u_2^{-2}$$

\item
$$u_t=Ku^{5}u_1^{-6}u_2u_3+u^4u_1^{-5}u_2^{2}[L-3Kuu_1^{-2}u_2]$$

\end{enumerate}

$\underline{\mathsf{A}_{3.7}\oplus\mathsf{A}_{2.2}}$:

\begin{enumerate}

\item
$$F=Kx^3\exp[-(q-1)^2/(q \kappa x^{\alpha} u_2)],\quad \tau=x^{(2q-1)/(q-1)}u_2$$
$$G=x^{1/(1-q)}\exp[-(q-1)^2/(q \kappa) x^{\alpha} \tau)][K\frac{2q-1}{q-1}\tau+L]$$

\item

$$F=K(u_1/u_2)^{3-(1/q)},\quad G=Lu_1^{2-(1/q)}u_2^{(1/q)-1}$$

\item
$$F=Lu_1^{2-q}u_2^{q-1},\quad G=Ku_1^{3-q}u_2^{q-3}$$

\item
$$F=Ku^{2}u_1^{-3}\omega^{-1/q},\quad G=-K[(2q-1)(q-1)+3\tau^2\omega^2]\omega^{-1/q}+L\tau \omega^{(q-1)/q}$$
$$\tau=u_1 u^{q-1},\quad \sigma=u_2 u^{2q-1},\quad \omega=[\sigma+(q-1)\tau^2]\tau^{-3}$$

\item
$$F=Ku^{3-q}\sigma^{-q},\quad G=u^{1-q}\sigma^{1-q}[-3K\sigma+L],\quad \sigma=uu_2u_1^{-2}$$

\end{enumerate}

\medskip\noindent $\underline{\mathsf{A}_{3.8}\oplus\mathsf{A}_{2.2}}$:

$$F=K(1+x^2)^3e^{-\omega/\kappa},\quad G=(1+x^2)^{1/2}e^{-\omega/\kappa}[3Kx+L],\quad \omega=(1+x^2)^{1/2}u_2$$

\bigskip\noindent$\underline{\mathsf{A}_{4.7}\oplus\mathsf{A}_1}$:

\begin{enumerate}
\item
$$u_t=Ke^{-3u_2}u_3+Le^{-2u_2}$$

\item
$$u_t=t^{-1}[Ke^{-3u_2}u_3+Le^{-2u_2}]$$

\end{enumerate}

\bigskip\noindent$\underline{\mathsf{A}_{4.8}\oplus\mathsf{A}_1}$:

\begin{enumerate}
\item
$$F=Kt^{3q-1}\tau^{3q/(1-q)},\quad G=Lt^q\tau^{(q+1)/(1-q)},\quad \tau=t^{q-1}u_2$$

\item
$$u_t=Ku_2^{3q/(1-q)}u_3+Lu_2^{(1+q)/(1-q)}$$

\item
$$\text{No invariant equation}$$

\item
$$F=Kt^{-1}u_1^{-3}\tau^{3/(q-1)},\quad G=t^{-1}u_1[-3Ku_1\tau^{(2q+1)/(q-1)}+L\tau^{(q+1)/(q-1)}],$$
$$\tau=u_2u_1^{-3}$$

\item
$$F=Ku_1^{-3}\tau^{3/(q-1)},\quad G=u_1[-3Ku_1\tau^{(2q+1)/(q-1)}+L\tau^{(q+1)/(q-1)}],$$
$$\tau=u_2u_1^{-3}$$

\item
$$u_t=Kt^{-1}u_1^{6}u_2^{-3}u_3+t^{-1}u_1^{4}u_2^{-1}[-3Ku_1+L]$$

\end{enumerate}

\bigskip\noindent$\underline{\mathsf{A}_{4.10}\oplus\mathsf{A}_1}$:

$$u_t=K(1+u_1^2)^3u_2^{-3/2}u_3+(1+u_1^2)^2u_2^{-1}[-3Ku_1+L]$$

\subsection{Indecomposable solvable Lie algebras.}

$\mathsf{A}_{5.19}$:

\begin{enumerate}
\item
$$u_t=Ku_2^{(3p-q)/(1-p)}u_3+Lu_2^{(1+p-q)/(1-p)}$$

\item
$$F=Ku_1^{-3}\tau^{2/(p-1)},\quad G=-1-3Ku_2^2u_2\tau^{(2p)/(1-p)}+L\tau^{p/(1-p)}$$
$$\tau=u_2u_1^{-3}$$

\item
$$F=Ku_1^{-3}\tau^{(q-3)/(1-p)},\quad G=u_1[-3u_1\tau^{(q-2p-1)/(1-p)}+L\tau^{p/(p-1)}]$$
$$\tau=u_2u_1^{-3}$$

\item
$$F=Ku_1^{-3}\tau^{2/(p-1)},\quad G=u_1[-3Ku_1\tau^{2p/(p-1)}+L\tau^{p/(p-1)}]$$
$$\tau=u_2u_1^{-3}$$

\item
$$F=Ku_1^{-3}t^{-1}\tau^{(q-3)/(1-q)},\quad G=u_1e^{-1/u_1}[-3Kt^{-1}u_1e^{1/u_1}\tau^{(q+1)/(q-1)}+L\omega^{-1/q}]$$
$$\tau=u_2u_1^{-3},\quad \sigma=u_1^{-1}-\ln|t|,\quad \omega=\sigma\tau^{q/(1-q)}$$

\end{enumerate}

$\mathsf{A}_{5.20}$:

\begin{enumerate}
\item
$$F=K u_2^{(1-2p)/(p-1)},\quad G=\frac{1}{1-p}\ln |u_2|+L,\quad p \ne 1$$

\item
$$\text{No invariant equation}$$

\item

Up to transformation $t'=t$, $x'=u$, $u'=x$
$$F=Ku_2^{(p-2)/(1-p)},\quad G=\frac{1}{p-1}\ln |u_2|+L,\quad p \ne 1$$

\item
$$\text{No invariant equation}$$

\end{enumerate}

$\mathsf{A}_{5.22}$:

$$F=Ku_1^{-3}e^{\tau},\quad G=u_1e^{\tau}(-3K\tau^2 u_1 +L),\quad \tau=u_2u_1^{-3}$$

$\mathsf{A}_{5.23}$:

\begin{enumerate}
\item
$$F=Ku_1^{-3}e^{-2\tau},\quad G=u_1[u_1^{-1}-3K\tau^2 u_1 e^{-2\tau}+ Le^{-\tau}],\quad \tau=u_2u_1^{-3}$$

\item
$$F=Ku_1^{-3}e^{(p-3)\tau},\quad G=u_1[-3K\tau^2 u_1 e^{(p-3)\tau}+Le^{(p-2)\tau}],\quad \tau=u_2u_1^{-3}$$

\end{enumerate}

$\mathsf{A}_{5.24}$:

$$F=Ku_1^{-3}e^{-\tau},\quad G=u_1[-3Ku_1 \tau^2 e^{-\tau}+\epsilon \tau],\quad \tau=u_2u_1^{-3}$$

$\mathsf{A}_{5.30}$:

\begin{enumerate}
\item
$$F=Ku_1^{-3}\tau^{-(3+2p)/p},\quad G=u_1[\frac{1}{2}u_1^{-2}-3Ku_1 \tau^{-3/p}+L\tau^{-2/p}]$$
$$\tau=u_2u_1^{-3}$$

\item
$$F=Ku_2^{-(p-3)/p},\quad G=-\frac{x^2}{2}+Lu_2^{2/p}$$

\item
$$F=Ke^{-3/(2\alpha)u_2},\quad G=-\frac{x^2}{2}+Le^{-u_2/\alpha}$$

\end{enumerate}

$\mathsf{A}_{5.32}$:

$$F=Ku_2^{-3}u_1^3,\quad G=u_1[\frac{1}{2}u_1^{-2}-3Ku_1+p\ln|\tau|+L]$$

$\mathsf{A}_{5.33}$:

\begin{enumerate}
\item
$$F=Ku_1^{3/(q-1)}\tau^{-(q+3p-1)/p},\quad G=Lu_1^{q/(q-1)}\tau^{-(q+p-1)/p}$$
$$\tau=u_2u_1^{-(q-2)/(q-1)}$$

\item
$$F=Ku_1^{3/(p-1)}\tau^{-(3+q(p-1))},\quad G=Lu_1^{p/(p-1)}\tau^{-(1+q(p-1))}$$
$$\tau=u_2u_1^{-(p-2)/(p-1)}$$

\item
$$F=Ku_1^{(p-3)}\tau^{p+q-3},\quad G=Lu_1^{p}\tau^{p+q-1}$$
$$\tau=u_2u_1^{-2}$$

\item
$$F=Kx^{1/p-(2q-3)},\quad G=Lx^{1/p-(q-1)}$$

\item
$$F=Kx^3\exp[-(q-1)^2/(\kappa q)\tau],\quad \tau=x^{(2q-1)/(q-1)}u_2$$
$$G=x^{1/(1-q)}\exp[-(q-1)^2/(\kappa q)\tau][K\frac{(2q-1)}{q-1}\tau+L]$$

\item
$$F=Kx^{1/q-(2p-3)},\quad G=Lx^{1/q-(p-1)}$$

\item
$$\text{same as case 5 with } p\leftrightarrow q$$

\item
$$F=Kx^{-q-2p+3}u_2^{-(p+q)},\quad G=Lx^{2-2p-q}u_2^{1-(p+q)}$$

\end{enumerate}

$\mathsf{A}_{5.34}$:

\begin{enumerate}
\item
$$F=Ke^{-1/u_1}u_1^{-3(p-2)}u_2^{p-3},\quad G=Le^{-1/u_1}u_1^{-(3p-4)}u_2^{p-1}$$

\item
$$F=Ku_1^{-7}u_2^2,\quad G=-\ln|u_1|+(p-1)\ln|\tau|+L,\quad \tau=u_2u_1^{-2}$$

\item
$$F=Ke^xu_2^{-p},\quad G=Le^xu_2^{1-p}$$

\item
$$F=Kx^2u_2^{-1},\quad G=x[\ln|x|+L(xu_2)^{1-p}]$$

\end{enumerate}

$\mathsf{A}_{5.35}$:

\begin{enumerate}
\item
$$F=K(1+u_1^2)^{-(3/2)(p-2)}e^{q\arctan u_1}u_2^{p-3},$$
$$G=(1+u_1^2)^{-(3/2)(p-4)}e^{q\arctan u_1}u_2^{p-1}(-3Ku_1+L)$$

\item
$$F=K(1+x^2)^{-(3/2)(p-2)}e^{q\arctan x}u_2^{-p},$$
$$G=(1+x^2)^{-(3/2)(p-4)}e^{q\arctan x}u_2^{1-p}(-3Ku_1+L)$$

\end{enumerate}

$\mathsf{A}_{5.36}$:

\begin{enumerate}
\item
$$F=Kt^{(1-p)/(2+p)}u_2^{-3/(p+2)},\quad G=Lt^{-p/(p+2)}u_2^{p/(p+2)}$$

\item
$$F=Kt^{1/2}u_2^{-3/2},\quad G=-\frac{1}{2}\ln|\tau|+L,\tau=u_2t^{-1}$$

\item
$$F=Ku_1^{-3},\quad G=-1-3Ku^2u_1^{-4}u_2^{2}+Luu_1^{-2}u_2$$

\item
$$F=Kt^2 u_1^{3\alpha}\sigma^{-3(\alpha+1)/(\alpha+2)},\quad \sigma=tu_1^{\alpha-1}u_2$$ $$G=u_1^{\alpha+1}\sigma^{-\alpha/(\alpha+2)}[-\frac{3K}{\alpha+3}\sigma^{(\alpha+3)/(\alpha+2)}+L]$$

\item
$$F=Kt^2\sigma^{-3/2},\quad G=u_1[\ln|u_1|-\frac{1}{2}\ln|\sigma|-3K\sigma^{1/2}+L]$$
$$\sigma=tu_1^{-1}u_2$$

\item
$$F=0$$
\item
$$F=Ku_1^5u_2^{-3},\quad G=uu_1+Lu_1^3u_2^{-1}$$

\item
$$F=Ku_1^{4}(\tau-1)^2,\quad G=u_1^3(\tau-1)u_2^{-1}[-3K\tau+L],\quad \tau=tu_1$$

\end{enumerate}

$\mathsf{A}_{5.37}$:

$$F=0: \quad \text{Inadmissible}$$

$\mathsf{A}_{5.38}$:

\begin{enumerate}
\item
$$F=Kx^3e^{-xu_2},\quad G=e^{-xu_2}[Kx^2u_2+Lx]$$

\item
$$F=Kx^3e^{x^2u_2},\quad G=e^{x^2u_2}[2Kx^2u_2+L]$$

\end{enumerate}

\section{Conclusions}
At this stage, all non-linearities have been specified (up to a multiplicative constant). In this sense, we have been able to give a complete specification of these equations by point symmetries. Not all the types admitting $\dim\mathsf{A}=4$ admit a solvable $\mathsf{A}$ with $\dim\mathsf{A}=5$ and these types of equation then require higher-order symmetries in order to specify the non-linearities (up to constants). Thus, we see that our classification touches upon the problem of finding complete symmetry groups for an evolution equation.


\end{document}